\documentclass[10pt,oneside, draft]{amsart}
\usepackage[a4paper, textwidth=16.5cm]{geometry}
\usepackage[T1]{fontenc}
\usepackage[english]{babel}
\usepackage{xcolor} % colors
\usepackage{tikz} % drawing
\usepackage[draft=false]{hyperref} % hyperlinks
\usepackage{mathtools}
\usepackage{amsmath,amsthm,amssymb}
\usepackage[all,cmtip]{xy}
\usepackage[capitalise]{cleveref}  % load after hyperref
\newcommand{\clevertheorem}[3]{%
	\newtheorem{#1}[thm]{#2}
	\crefname{#1}{#2}{#3}
}
\numberwithin{equation}{section} %% Comment out for sequentially-numbered
\numberwithin{figure}{section} %% Comment out for sequentially-numbered

\theoremstyle{plain} % bold environment name, italic text
\newtheorem{thm}{Theorem}[section]
\crefname{thm}{Theorem}{Theorems}
\newtheorem*{thm*}{Theorem}
\clevertheorem{proposition}{Proposition}{Propositions}
\clevertheorem{prop}{Proposition}{Propositions}
\newtheorem*{prop*}{Proposition}
\clevertheorem{theorem}{Theorem}{Theorems}
\clevertheorem{lemma}{Lemma}{Lemmas}

\clevertheorem{corollary}{Corollary}{Corollaries}
\clevertheorem{cor}{Corollary}{Corollaries}
\clevertheorem{conj}{Conjecture}{Conjectures}
\clevertheorem{questn}{Question}{Questions}

\theoremstyle{definition} % bold environment name, plain text
\clevertheorem{definition}{Definition}{Definitions}
\clevertheorem{assumption}{Assumption}{Assumptions}
\clevertheorem{defn}{Definition}{Definitions}
\clevertheorem{notn}{Notation}{Notations}
\clevertheorem{conv}{Convention}{Conventions}

\clevertheorem{remark}{Remark}{Remarks}
\clevertheorem{rmk}{Remark}{Remarks}
\clevertheorem{warn}{Warning}{Warnings}
\clevertheorem{example}{Example}{Examples}
\clevertheorem{summ}{Summary}{Summaries}

\usepackage{amsmath,amsthm,amssymb,amsfonts,extarrows}
\usepackage{enumerate,amscd,amsxtra,MnSymbol}
\usepackage{mathrsfs}
\usepackage{color}
\usepackage[cmtip,all]{xy}
\usepackage{eucal}
\usepackage{bbold}
\usepackage{upgreek}
\usepackage{paralist}
\usepackage{tikz-cd}
\usepackage{dsfont}
\usepackage{bbm}
\usepackage[bbgreekl]{mathbbol}

\DeclareMathSymbol\bbDelta \mathord{bbold}{"01}
\DeclareMathSymbol\bDelta \mathord{bbold}{"01}

\newtheorem{remark*}{Remark}

\newtheorem{construction}[thm]{Construction}

\newtheorem{notation}[thm]{Notation}

\newcommand{\bE}{{\mathbb E}}

\renewcommand{\P}{{\mathbb P}}

\newcommand{\mA}{{\mathcal A}}
\newcommand{\mB}{{\mathcal B}}
\newcommand{\mC}{{\mathcal C}}
\newcommand{\mD}{{\mathcal D}}
\newcommand{\mE}{{\mathcal E}}
\newcommand{\mF}{{\mathcal F}}
\newcommand{\mG}{{\mathcal G}}

\newcommand{\mJ}{{\mathcal J}}
\newcommand{\mK}{{\mathcal K}}

\newcommand{\mM}{{\mathcal M}}
\newcommand{\mN}{{\mathcal N}}
\newcommand{\mO}{{\mathcal O}}
\newcommand{\mP}{{\mathcal P}}

\newcommand{\mS}{{\mathcal S}}
\newcommand{\mT}{{\mathcal T}}

\newcommand{\mV}{{\mathcal V}}
\newcommand{\mW}{{\mathcal W}}
\newcommand{\mX}{{\mathcal X}}
\newcommand{\mY}{{\mathcal Y}}

\newcommand{\A}{{\mathrm A}}
\newcommand{\B}{{\mathrm B}}
\newcommand{\C}{{\mathrm C}}
\newcommand{\D}{{\mathrm D}}
\newcommand{\E}{{\mathrm E}}
\newcommand{\F}{{\mathrm F}}
\newcommand{\G}{{\mathrm G}}
\newcommand{\rH}{{\mathrm H}}

\newcommand{\K}{{\mathrm K}}
\renewcommand{\L}{{\mathrm L}}

\renewcommand{\P}{{\mathrm P}}

\newcommand{\R}{{\mathrm R}}
\newcommand{\rS}{{\mathrm S}}
\newcommand{\T}{{\mathrm T}}

\newcommand{\V}{{\mathrm V}}
\newcommand{\W}{{\mathrm W}}
\newcommand{\X}{{\mathrm X}}
\newcommand{\Y}{{\mathrm Y}}
\newcommand{\Z}{{\mathrm Z}}

\newcommand{\rc}{\mathrm{c}}

\newcommand{\bj}{\mathrm{j}}
\newcommand{\bi}{\mathrm{i}}
\newcommand{\m}{{\mathrm{m}}}
\newcommand{\bk}{\mathrm{k}}

\newcommand{\rt}{\mathrm{t}}

\newcommand{\g}{\mathrm{g}}
\newcommand{\n}{\mathrm{n}}
\newcommand{\br}{\mathrm{r}}

\newcommand{\op}{\mathrm{op}}
\newcommand{\s}{\mathrm{s}}
\newcommand{\colim}{\mathrm{colim}}
\newcommand{\Mod}{{\mathrm{Mod}}}
\newcommand{\LMod}{{\mathrm{LMod}}}
\newcommand{\RMod}{{\mathrm{RMod}}}

\newcommand{\rev}{{\mathrm{rev}}}

\newcommand{\Env}{{\mathrm{Env}}}  
\newcommand{\ot}{\otimes}

\newcommand{\co}{\mathrm{co}}

\newcommand{\coLMod}{{\mathrm{coLMod}}}  

\newcommand{\Hopf}{{\mathrm{Hopf}}}  
\newcommand{\sSeq}{{\mathrm{sSeq}}}  
\newcommand{\LaxFun}{{\mathrm{LaxFun}}}  
\newcommand{\loc}{{\mathrm{loc}}}  
\newcommand{\Nat}{{\mathrm{Nat}}}  
\newcommand{\Cmon}{{\mathrm{Cmon}}}

\newcommand{\id}{\mathrm{id}}
\newcommand{\act}{\mathrm{act}}
\newcommand{\Cat}{\mathrm{Cat}}

\newcommand{\Set}{\mathrm{Set}}

\newcommand{\Coalg}{\mathrm{Coalg}}
\newcommand{\Alg}{\mathrm{Alg}}
\newcommand{\Calg}{\mathrm{Calg}}
\newcommand{\Mon}{\mathrm{Mon}}
\newcommand{\Fun}{\mathrm{Fun}}
\newcommand{\Op}{{\mathrm{Op}}}
\newcommand{\cocart}{{\mathrm{cocart}}}

\newcommand{\cart}{{\mathrm{cart}}}

\newcommand{\lax}{{\mathrm{lax}}}
\newcommand{\oplax}{{\mathrm{oplax}}}

\newcommand{\tu}{{\mathbb 1}}

\newcommand{\FUN}{\mathrm{FUN}}

\newcommand{\Mul}{{\mathrm{Mul}}}
\newcommand{\ev}{{\mathrm{ev}}}

\newcommand{\End}{{\mathrm{End}}}
\newcommand{\Map}{{\mathrm{Map}}}

\newcommand{\Enr}{{\mathrm{Enr}}} 

\newcommand{\LinFun}{{\mathrm{LinFun}}}

\newcommand{\Mor}{{\mathrm{Mor}}} 

\newcommand{\Lax}{\mathrm{Lax}}

\newcommand{\Fin}{{\mathrm{Fin}}}

\newcommand{\mon}{{\mathrm{mon}}}

\newcommand{\f}{{\mathrm{f}}} 

\newcommand{\h}{{\mathrm{h}}} 
\newcommand{\Ho}{{\mathrm{Ho}}}

\begin{document}
	
%\title{A monadicity theorem for higher algebraic structures}

% \title{On the classification of structured monads}

% \title{On structured monads on $\infty$-categories}

%\title{The higher algebra of monads}

%\title{A monadic view on categorical algebra}

\title{A Duality Between Monads and monadic morphisms}

\author{Hadrian Heine, \\ Max Planck Institute of Mathematics, Bonn, \\ heine@mpim-bonn.mpg.de}
\maketitle

\begin{abstract} 

We establish a duality between monads and monadic morphisms in any $(\infty,2)$-category and characterize monadic morphisms in a wide class of examples.
%Specifically, we show that monads in $(\infty,2)$-categories are classified
%by their Eilenberg-Moore objects, which internalize the $\infty$-category of algebras of the monad.
This duality unifies several dualities between algebraic structures and their representations, and provides a general mechanism for transferring structure from a monad to its $\infty$-category of algebras.
This transfer of structure yields uniform constructions of tensor products for algebras over lax symmetric monoidal and oplax symmetric monoidal monads, extending classical tensor products for modules and operadic algebras.
Using this framework, we construct a relative tensor product for algebras
over lax monoidal monads, a tensor product for algebras over Hopf $\infty$-operads and equip the $\infty$-category of operadic algebras with canonical enrichment.
%and produce an equivalence between enriched 
%$\infty$-operads and marked presentably symmetric monoidal $\infty$-categories.

\end{abstract}

\tableofcontents

\section{Introduction}

%\subsubsection{Motivation}

A central principle in mathematics is that algebraic structures are most naturally understood through their representations.
Groups, for instance, are studied via their actions on vector spaces and manifolds, and rings are encoded by their categories of modules, giving rise to $K$-theory.
This perspective naturally extends to higher category theory.
An associative algebra $A$ in a presentably monoidal $\infty$-category $\mathcal V$ can be recovered from its $\infty$-category of right $A$-modules in $\mV$, provided one remembers its natural $\mathcal V$-linear structure and the forgetful functor from modules to $\mV$.
More precisely, sending an associative algebra $A$ in $\mV$ to the left $\mV$-tensored $\infty$-category of $A$-modules in $\mV$ promotes to an equivalence
\begin{equation}\label{eq1} \{ \text{associative} \ \text{algebras} \ \text{in} \ \mV\} \simeq \{ \text{presentably} \ \text{left} \ \mV-\text{tensored} \ \text{categories} \  \text{with} \ \text{a} \ \text{monadic} \ \text{functor} \ \text{to} \ \mV\} \end{equation}

This duality (\ref{eq1}) is a manifestation of a more fundamental duality that becomes only visible at the level of $(\infty,2)$-categories.
A monoidal $\infty$-category can be regarded as an $(\infty,2)$-category $B\mV$ with a single object and endomorphisms $\mV$, and associative algebras in $\mV$ are monads in $B\mV$ - associative algebras in the endomorphisms of $B\mV$. The $\infty$-category of modules over the algebra then coincides with the $\infty$-category of algebras - or representations - of the corresponding monad.
The equivalence (\ref{eq1}) is thus an instance of a more fundamental duality between monads and monadic morphisms in any suitable $(\infty,2)$-category.
%In this generality a morphism is monadic if it admits a left adjoint and is representably monadic - using the natural enrichment of $(\infty,2)$-categories in $(\infty,1)$-categories.
The purpose of this paper is to establish this equivalence between monads and monadic morphisms in full generality, thereby unifying previously known dualities and producing new ones.

A key conceptual feature of this duality is that it works as a mechanism to 
transport
%additional structure on monads corresponds to 
structure from a monad to its $\infty$-category of algebras.
For instance, when $\mathcal V$ is a symmetric monoidal $\infty$-category, equivalence (\ref{eq1}) is symmetric monoidal, and so transports the following structures:

\begin{center}\label{uj}
\begin{tabular}{ | l | l | l | p{5cm} |}
\hline
Algebras & Representations \\ \hline
cocommutative bialgebras & presentably symmetric monoidal $\mV$-linear $\infty$-categories over $\mV$ \\ \hline
$\bE_{\n+1}$-algebras & presentably $\bE_\n$-monoidal $\mV$-linear $\infty$-categories  \\ 
\hline
\end{tabular}
\end{center}

The $(\infty,2)$-category $B\mV$ is the full subcategory spanned by
$\mV$ of the $(\infty,2)$-category of presentably left $\mV$-tensored $\infty$-categories.
Under this identification, associative algebras $A$ in $\mV$ correspond to the $\mV$-linear monads $ A \ot (-): \mV \to \mV$. 
Passing from associative algebras in a symmetric monoidal $\infty$-category $\mV$ to monads in $B\mV$, cocommutative bialgebras in $\mV$ become monads in the $(\infty,2)$-category of symmetric monoidal $\mV$-linear $\infty$-categories and oplax symmetric monoidal $\mV$-linear functors. Similarly, $\bE_{\n+1}$-algebras in $\mV$ become monads in the $(\infty,2)$-category of $\mV$-linear $\bE_{\n}$-monoidal $\infty$-categories and lax $\bE_{\n}$-monoidal $\mV$-linear functors:
\begin{center}\label{uj}
\begin{tabular}{ | l | l | l | p{5cm} |}
\hline
$\mV$ & $B\mV$ \\ \hline
associative algebras & $\mV$-linear monads  \\ \hline
cocommutative bialgebras & oplax symmetric monoidal $\mV$-linear monads \\ \hline 
$\bE_{\n+1}$-algebras & lax $\bE_\n$-monoidal $\mV$-linear monads \\
\hline
\end{tabular}
\end{center}

To capture this transfer of structure systematically, it is essential to pass from associative algebras to the more flexible notion of monads. At the level of monads, the transfer of structure becomes intrinsic to the framework: for example, transferring symmetric monoidal structure to the $\infty$-category of representations of an oplax symmetric monoidal monad amounts to applying the duality between monads and monadic morphisms within the $(\infty,2)$-category of symmetric monoidal $\mV$-linear $\infty$-categories and oplax symmetric monoidal $\mV$-linear functors.

%\subsection{Main results}

A monad on some object $\X$ of an $(\infty,2)$-category $\mC$
is an associative algebra in the monoidal $\infty$-category $\Mor_\mC(\X,\X)$ of endomorphisms of $\X$. 
Similarly, the notion of Eilenberg-Moore category of representations of a monad makes sense in any $(\infty,2)$-category.
A morphism $\Y \to \X$ in an $(\infty,2)$-category $\mC$ is an Eilenberg-Moore object for a monad $\T$ on $\X $ if there is an equivalence 
\begin{equation}\label{pli}
\begin{xy}
\xymatrix{
\Mor_\mC(\Z, \Y) \ar[rd]^{} \ar[rr]^{\simeq}
&& \LMod_\T(\Mor_\mC(\Z,\X)) \ar[ld]^{} 
\\ & \Mor_\mC(\Z,\X)
}
\end{xy} 
\end{equation}
natural in $\Z \in \mC$ between the $\infty$-category of morphisms $\Z \to \Y$ in $\mC$
and the category of left modules taken with respect to the natural left action of the monoidal $\infty$-category $\Mor_\mC(\X,\X)$ on $\Mor_\mC(\Z,\X).$

We establish the following fundamental duality:

\begin{theorem}\label{uikp}(Theorem \ref{uhnggfdaa})
Let $\mC$ be an $(\infty,2)$-category that admits Eilenberg-Moore objects.
There is a duality 
$$ \Fun([1],\mC)^\mon \simeq \Mon(\mC)^\op $$ 
between the $\infty$-category of monadic morphisms in $\mC$ and the
$\infty$-category of monads in $\mC$.

This duality arises from a localization $$ \Fun([1],\mC)^\R\rightleftarrows \Mon(\mC)^\op$$ 
between the $\infty$-category of right adjoint morphisms in $\mC$ and the
$\infty$-category of monads in $\mC$ whose local objects are the monadic morphisms.
\end{theorem}

\cref{uikp} for $\mC= \infty\Cat$ and more generally
cosmifiable $\infty$-categories was independently proven by Haugseng \cite[Theorem 1.2.]{Haugseng2020OnLT} %for the case that $\mC= \infty\Cat$ (or more generally if $\mC$ admits a model by an $\infty$-cosmos
relying on fundamental results of Riehl-Verity \cite{RIEHL2016802} and Zaganidis \cite{zaganidis2017towards}.
%fundamental results 
%developed in the setting of $\infty$-cosmoi).

Specializing Theorem \ref{uikp} to the fiber over any chosen object of $\mC$
we obtain the following corollary:
\begin{corollary}\label{coro1}

Let $\mC$ be an $(\infty,2)$-category that admits Eilenberg-Moore objects and $\X \in \mC.$ There is a duality $$ \mC_{/\X}^\mon \simeq \Alg(\Mor_\mC(\X,\X))^\op$$ 
between the $\infty$-category of monadic morphisms in $\mC$ over $\X$ and the $\infty$-category of monads in $\mC$.

This duality arises from a localization $$ \mC_{/\X}^\R\rightleftarrows \Alg(\Mor_\mC(\X,\X))^\op$$ 
between the $\infty$-category of right adjoint morphisms in $\mC$ over $\X$ and the $\infty$-category of monads in $\mC$ whose local objects are the monadic morphisms over $\X.$

\end{corollary}
% \cref{coro1} for $\mC= \infty\Cat$, the $\infty$-category of small $\infty$-categories, was conjectured by Lurie \cite[Remark 4.7.3.8.]{lurie.higheralgebra}. 

Applying \cref{coro1} to the $(\infty,2)$-category of presentably left $\mV$-tensored $\infty$-categories we recover equivalence (\ref{eq1}).

Moreover we recover a many object version of equivalence (\ref{eq1}), which was used by \cite{reutter2025enriched} to identify a new model of enriched $\infty$-categories.
% several fundamental extensions, which apply to enriched
% $\infty$-categories and enriched $\infty$-operads, organized by the following table:
% \begin{center}\label{uj}
% \begin{tabular}{ | l | l | l | p{5cm} |}
% \hline
% One object & Many objects \\ \hline
% associative algebras in $\mV$ & $\mV$-enriched $\infty$-categories \\ \hline
% $\infty$-operads in $\mV$ & $\mV$-enriched $\infty$-operads  \\ 
% \hline
% \end{tabular}
% \end{center}
Gepner-Haugseng \cite{GEPNER2015575} define enriched $\infty$-categories as a many object version of associative algebras. They define $\infty$-categories enriched in any presentably monoidal $\infty$-category $\mV$ with fixed space of objects
$\X$ as algebras in $\mV$ over a many object version of the associative
$\infty$-operad. Hinich \cite{HINICH2020107129} and MacPherson \cite{MR4185309} identify Gepner-Haugseng's model of
$\mV$-enriched $\infty$-categories with space of objects $\X$ with associative algebras in the monoidal $\infty$-category of small colimits preserving $\mV$-linear endofunctors of $\mV^\X$, the $\infty$-category of functors
$\X \to \mV.$ In other words, $\mV$-enriched $\infty$-categories with space of objects $\X$ in the sense of Gepner-Haugseng are precisely $\mV$-linear monads on $\mV^\X$ preserving small colimits.

\cref{coro1} implies the following corollary:

\begin{corollary}\label{coro10}

Let $\mV$ be a presentably monoidal $\infty$-category and $\X$ a small space.
There is a duality $$ \{ \text{Presentably} \ \text{left} \ \mV-\text{tensored} \ \infty-\text{categories} \  \text{with} \ \text{a} \ \text{monadic} \ \text{functor} \ \text{to} \ \mV^\X\} \simeq $$$$ \{ \mV-\text{enriched} \ \infty-\text{categories} \ \text{with} \ \text{space} \ \text{of} \ \text{objects} \ \X \}$$ 

\end{corollary}

\cref{coro10} extends equivalence (\ref{eq1}) from the final space to an arbitrary space. Sending left to right adjoints, \cref{coro10}
provides an equivalence
\begin{equation}\label{eqpl}
\{ \text{Presentably} \ \text{left} \ \mV-\text{tensored} \ \infty-\text{categories} \  \text{with} \ \text{a} \ \text{free} \ \text{functor} \ \text{from} \ \mV^\X\} \simeq \end{equation}
$$
\{ \mV-\text{enriched} \ \infty-\text{categories} \ \text{with} \ \text{space} \ \text{of} \ \text{objects} \ \X \}$$

Equivalence (\ref{eqpl}) recovers an equivalence of \cite[Theorem 4.7.]{reutter2025enriched},
which was used by \cite{reutter2025enriched} to
identify a new model for enriched $\infty$-categories:
by the universal property of $\mV^\X$ the left hand side of (\ref{eqpl})
corresponds to certain pairs 
$(\mM,\alpha: \X \to \mM)$ consisting of a presentably left $\mV$-tensored $\infty$-category $\mM$ and a functor $\alpha: \X \to \mM$ such that
$\alpha^*: \mM \to \mV^\X$ is conservative, small colimits preserving and $\mV$-linear. 
Such pairs were introduced by \cite{reutter2025enriched} as a model for enriched $\infty$-categories, which only requires the technology of presentably left tensored $\infty$-categories.

In practice it is crucial to know \cref{coro10} not only for a fixed space of objects - but for varying spaces of objects - since morphisms between enriched $\infty$-categories do rarely preserve the space of objects.
\cref{uikp} provides such a global version of \cref{coro10},
where the space of objects does not need to be fixed.

As another important corollary of \cref{uikp} we obtain a multicategorical version of \cref{coro10}, which will be explored in detail in forthcoming work.
% is a previously unknown duality between enriched $\infty$-operads and suitable presentably symmetric monoidal $\infty$-categories equipped with a map of spaces.
Haugseng \cite[Corollary 4.2.9.]{Rune} and Brantner-Campos-Nuiten
\cite[Definition 3.4]{brantner2021pd} define $\infty$-operads enriched in a presentably symmetric monoidal $\infty$-category $\mV$ with space of objects $\X$
as algebras for the composition product on the $\infty$-category $\sSeq_\X(\mC) $ of $\X$-colored symmetric sequences in $\mV$.
These sequences are functors 
$\X \times \Sigma(\X) \to \mV$, where $\Sigma(\X):= \coprod_{\n \geq 0} \X^{\times \n} $ is the free symmetric monoidal $\infty$-category generated by $\X.$
By construction of \cite{brantner2021pd} the monoidal $\infty$-category 
of $\X$-colored symmetric sequences in $\mV$ is the opposite of the monoidal $\infty$-category of small colimits preserving symmetric monoidal endofunctors
of $\mV^{\Sigma(\X)}$ under $\mV$.
In other words, $\mV$-enriched $\infty$-operads with space of objects $\X$ are precisely symmetric monoidal $\mV$-linear monads on $\mV^{\Sigma(\X)}$ preserving small colimits.

\cref{coro1} implies the following corollary:

\begin{corollary}\label{coro11}

Let $\mV$ be a presentably symmetric monoidal $\infty$-category and $\X$ a small space.
There is a duality $$ \{ \text{Presentably} \ \text{symmetric} \  \text{monoidal} \ \mV-\text{linear} \ \infty-\text{categories} \  \text{with} \ \text{a} \ \text{monadic} \ \text{functor} \ \text{to} \ \mV^{\Sigma(\X)}\} \simeq $$$$\{ \mV-\text{enriched} \ \infty-\text{operads} \ \text{with} \ \text{space} \ \text{of} \ \text{objects} \ \X \}$$
\end{corollary}

Again sending left to right adjoints, \cref{coro11} provides an equivalence
\begin{equation}\label{eqpl2}
\{ \text{Presentably} \ \text{symmetric} \  \text{monoidal} \ \mV-\text{linear} \ \infty-\text{categories} \  \text{with} \ \text{a} \ \text{free} \ \text{functor} \ \text{from} \ \mV^{\Sigma(\X)}\} \simeq \end{equation}
$$ \{ \mV-\text{enriched} \ \infty-\text{operads} \ \text{with} \ \text{space} \ \text{of} \ \text{objects} \ \X \}$$

By the universal property of $\mV^{\Sigma(\X)}$ the left hand side of (\ref{eqpl2})
corresponds to certain pairs 
% Presentably left $\mV$-tensored $\infty$-categories $\mM$ equipped with a left adjoint $\mV$-linear functor $\mV^\X \to \mM$ precisely correspond to pairs
$(\mM,\alpha: \X \to \mM)$ consisting of a presentably symmetric monoidal
$\mV$-linear $\infty$-category $\mM$ and a functor $\alpha: \X \to \mM$ such that $\alpha^*: \mM \to \mV^\X$ is conservative, small colimits preserving, $\mV$-linear and symmetric monoidal.
We call such pairs marked presentably symmetric monoidal $\mV$-linear $\infty$-categories. We think of marked presentably symmetric monoidal $\mV$-linear $\infty$-categories as a convenient model for $\mV$-enriched $\infty$-operads, which only requires the technology of presentably symmetric monoidal $\infty$-categories.

Again \cref{uikp} provides a global version of \cref{coro11},
where the space of objects does not need to be fixed:

\begin{corollary}\label{coro12}

Let $\mV$ be a presentably symmetric monoidal $\infty$-category.
There is an equivalence $$ \{\text{Marked} \ \text{presentably} \ \text{symmetric} \  \text{monoidal} \ \mV-\text{linear} \ \infty-\text{categories} \} \simeq $$$$\{ \mV-\text{enriched} \ \infty-\text{operads}\}$$
\end{corollary}

To apply \cref{uikp} and \cref{coro1} it is necessary to know that Eilenberg-Moore objects exist. 
Classically, the category of algebras of a monad $\T$ on a category $\X$ is the lax limit of the lax functor classifying the diagram:
\begin{equation}\label{shape}
\begin{tikzcd}[column sep=huge]
\X
\arrow[r, bend left=40, "\id"]
\arrow[r, bend right=40, swap, "\T"{name=G}]
& \X
\arrow[l, "\id"{name=S, description}]
\arrow[Rightarrow, from=S, to=G, ""]
\end{tikzcd}\end{equation}

We extend this construction to any $(\infty,2)$-category and prove the following existence result:
\begin{theorem}\label{th001} (Theorem \ref{stru})
Let $\mC$ be an $(\infty,2)$-category and $\T$ a monad in $ \mC$.
The diagram (\ref{shape}) canonically extends to a lax functor
$$\theta_\T: \Delta^\op \to \mC$$
whose lax limit if it exists, is the Eilenberg-Moore object of $\T.$
\end{theorem}

\cref{stru} guarantees that \cref{uikp} applies to a wide class of $(\infty,2)$-categories, including all presentable $(\infty,2)$-categories,
but also the $(\infty,2)$-categories of (symmetric) monoidal $(\infty$-categories and lax (oplax) (symmetric) monoidal functors,
which is not presentable but admits the lax limit of \cref{stru}.

\vspace{1mm}

\cref{uikp} also raises the question how to identify monadic morphisms in $(\infty,2)$-categories of interest.
For the case of $\infty\Cat$, the $(\infty,2)$-category of small $\infty$-categories, this question is answered by the monadicity theorem
\cite[Theorem 4.7.3.5]{lurie.higheralgebra}, \cite[Theorem 7.2.7.]{RIEHL2016802}. It is of special importance in higher category theory since it is a key tool to prove rectification results \cite[Theorem 4.1.8.4., Theorem 4.5.4.7.]{lurie.higheralgebra}, \cite[Theorem 1.1]{MR3402334}, \cite[Theorem 4.1.1.]{hinich2015rectification} and prove equivalences between models \cite[Theorem 4.1.]{Haugseng_2017}.

We extend the monadicity theorem to a wide class of $(\infty,2)$-categories of higher algebraic structures. These include various types of $\infty$-operads, (symmetric) monoidal $\infty$-categories, enriched $\infty$-categories, tensored $\infty$-categories, double $\infty$-categories and more generally, $\bk$-fold $\infty$-categories.
To treat this generality, we use the language of algebraic pattern \cite{chu2021homotopy}.
This formalism gives a uniform way to talk about operad-like structures,
which we call $\mathfrak{P}$-operads if $\mathfrak{P}$ is an algebraic pattern. We prove a monadicity theorem for $\mathfrak{P}$-operads. To state this monadicity theorem, we introduce the concept of $\mathfrak{P}$-operadic colimit,
which generalizes Lurie's notion of operadic colimit but also the notion of conical colimit in enriched $\infty$-categories.
We prove the following monadicity theorem for higher algebraic structures:

\begin{theorem}\label{th9}(\cref{strmon}) Let $\mathfrak{P}$ be an algebraic pattern.

A morphism $\G: \mD \to \mC$ in the $(\infty,2)$-category of $\mathfrak{P}$-operads is monadic if and only if the following hold:

\begin{enumerate}
\item The functor $\G: \mD \to \mC$ admits a left adjoint in the $(\infty,2)$-category of $\mathfrak{P}$-operads.

\item The map $\G: \mD \to \mC$ is fiberwise a monadic functor of $\infty$-categories.
	
\item The colimit of every fiberwise $\G$-split simplicial object is a  $\mathfrak{P}$-operadic colimit.
	
\end{enumerate}

\end{theorem}

\cref{th9} specializes to the following characterization of monadic morphisms in the 
$(\infty,2)$-category of (symmetric) monoidal $\infty$-categories and
lax (symmetric) monoidal functors:
\begin{corollary}(\cref{relbart})Let $\mC,\mD$ be (symmetric) monoidal $\infty$-categories.

A lax (symmetric) monoidal functor $\G: \mD \to \mC$ is monadic if and only if the following conditions hold: 
\begin{enumerate}
\item The functor $\G$ admits a (symmetric) monoidal left adjoint.

\item The functor $\G$ is conservative.

\item Every $\G$-split simplicial object of $\mD$ admits a colimit that is preserved by $\G$.

\item For every $\n \geq 0$ the functor $\ot: \mD^{\times \n} \to \mD $
preserves geometric realizations of $\G^{\times \n}$-split simplicial objects.

\end{enumerate}
\end{corollary}

\cref{th9} also specializes to the following characterization of monadicity for enriched $\infty$-categories, which was independently proven by Stefanich \cite[Theorem 7.4.10.]{stefanich2021higher}.

\begin{corollary}\label{EnrBarr}(\cref{EnrBar}) 
Let $\mV$ be a monoidal $\infty$-category. A $\mV$-enriched functor $\G: \mN \to \mM$ is monadic if and only if the following conditions hold:

\begin{enumerate}
\item The $\mV$-enriched functor $\G: \mN \to \mM$ admits a $\mV$-enriched left adjoint.

\item The $\mV$-enriched functor $\G: \mN \to \mM$ is conservative.

\item Every $\G$-split simplicial object of $\mN$ admits a $\mV$-enriched colimit that is preserved by $\G$.

\end{enumerate}

\end{corollary}	

%\subsection{Applications}

In view of \cref{uikp} we can transport structure from a monad to its $\infty$-category of algebras by applying the duality to a chosen $(\infty,2)$-category of $\infty$-categories with structure - provided we know 
that Eilenberg-Moore objects exist.

As the first application we prove that the $(\infty,2)$-category of $\bE_\n$-monoidal $\infty$-categories compatible with geometric realizations
and lax $\bE_\n$-monoidal functors preserving geometric realizations, 
admits Eilenberg-Moore objects.
This way we construct a relative tensor product for 
algebras over lax monoidal monads that extends the relative tensor product of modules:
\begin{theorem}\label{plob}(\cref{Rell}, \cref{atono})
Let $0 \leq \n \leq \infty$ and $\mC$ an $\bE_\n$-monoidal $\infty$-category compatible with geometric realizations. Let $\T$ be a lax $\bE_\n$-monoidal monad on $\mC$ that preserves geometric realizations.

There is an $\bE_\n$-monoidal structure on the $\infty$-category
$\LMod_\T(\mC)$ of $\T$-algebras in $\mC$ and the free functor $$\mC \to \LMod_\T(\mC)$$ refines to an $\bE_\n$-monoidal functor.
%and the forgetful functor $\LMod_\T(\mC) \to \mC $ is an Eilenberg-Moore object for $\T$ in the $(\infty,2)$-category of $\bE_\n$-monoidal $\infty$-categories and lax $\bE_\n$-monoidal functors.

\end{theorem}

For every $\bE_{\n+1}$-algebra $\A$ in an
$\bE_{\n+1}$-monoidal $\infty$-category $\mV$ the associated monad
$\A \ot(-): \mV \to \mV$ is an $\bE_{\n}$-monoidal monad.
So Theorem \ref{plob} recovers the result that in every $\bE_{\n+1}$-monoidal $\infty$-category compatible with geometric realizations the $\infty$-category of left $\A$-modules in $\mV$ is $\bE_{\n}$-monoidal via the relative tensor product.

Lurie's construction of the relative tensor product for modules over $\bE_{\n+1}$-algebras \cite[4.4.]{lurie.higheralgebra} uses substantial machinery, including the theory of operadic Kan extensions and the theory of enriched strings. 
By contrast, our construction for monadic algebras avoids these techniques and yields an alternative conceptually simpler approach to the relative tensor product.
To construct the relative tensor product of monadic algebras 
we make the idea rigorous that the relative tensor product agrees on free algebras with the given tensor product, and extends to all algebras by
resolving an arbitrary algebra by free algebras.
	
As a second application we use \cref{th001} to show that the $(\infty,2)$-category of symmetric monoidal $\infty$-categories and oplax symmetric monoidal functors admits Eilenberg-Moore objects. This way we construct a tensor product for algebras over oplax monoidal monads that lifts the monoidal structure along the forgetful functor:

\begin{theorem}\label{opl}(\cref{patty}, \cref{exkq})
Let $0 \leq \n \leq \infty$ and $\mC$ an $\bE_\n$-monoidal $\infty$-category.	Let $\T$ be an oplax $\bE_\n$-monoidal monad on $\mC$. 
There is an $\bE_\n$-monoidal structure on the $\infty$-category
$\LMod_\T(\mC)$ of $\T$-algebras in $\mC$ and the forgetful functor $$\LMod_\T(\mC) \to \mC $$ refines to an $\bE_\n$-monoidal functor.

%and an Eilenberg-Moore object for $\T$ in the $(\infty,2)$-category of $\bE_\n$-monoidal $\infty$-categories and oplax $\bE_\n$-monoidal functors.
\end{theorem}

We use the tensor product of algebras over oplax monoidal monads of Theorem \ref{opl} to construct a tensor product for algebras over Hopf $\infty$-operads enriched in any symmetric monoidal $\infty$-category $\mC,$
which are monochromatic $\infty$-operads enriched in the symmetric monoidal $\infty$-category of $\bE_\infty$-coalgebras in $\mC.$

Via the composition product the $\infty$-category $\sSeq(\mC)$ of
(single colored) symmetric sequences in $\mC$ acts on $\mC$ from the left, and the $\infty$-category $\Alg_\mO(\mC)$ of algebras over a $\mC$-enriched $\infty$-operad $\mO$ is the $\infty$-category of left $\mO$-modules with respect to this action.
The left action of $\mathrm{sSeq}(\mC)$ on $\mC$ induces a monoidal functor
\begin{equation}\label{monil}
\mathrm{sSeq}(\mC) \to \Fun(\mC,\mC), \mO \mapsto (X \mapsto \mO \circ X \simeq \bigoplus_{\n \geq 0} \mO_\n \ot_{\Sigma_\n} X^{\ot\n})\end{equation}
that sends a $\mC$-enriched $\infty$-operad $\mO$ to its associated monad $\T_\mO$ on $\mC$ whose $\infty$-category of algebras is $\Alg_\mO(\mC).$
We lift the monoidal functor (\ref{monil}) to a monoidal functor
from the composition product on $\mathrm{sSeq}(\Coalg_{\bE_\infty}(\mC))$ to the endomorphism monoidal structure on the $\infty$-category $ \Fun^{\ot, \oplax}(\mC, \mC)$ of oplax symmetric monoidal endofunctors:

\begin{theorem}\label{sqat}(Theorem \ref{dfghjbbvl})
Let $\mC$ be a symmetric monoidal $\infty$-category compatible with small colimits. There is a commutative square of monoidal $\infty$-categories
\begin{equation*}\label{hhgbhhjvcff} 
\begin{xy}
\xymatrix{
\mathrm{sSeq}(\Coalg_{\bE_\infty}(\mC)) \ar[d] 
\ar[r]^{  } 
& \Fun^{\ot, \oplax}(\mC, \mC) \ar[d] 
\\
\mathrm{sSeq}(\mC) \ar[r]^{ }  &  \Fun(\mC, \mC).
}
\end{xy} 
\end{equation*} 
\end{theorem}
Theorem \ref{sqat} guarantees that the monad associated to a $\mC$-enriched Hopf $\infty$-operad underlies an oplax symmetric monoidal monad, whose
$\infty$-category of algebras carries a canonical symmetric monoidal structure
by Theorem \ref{opl}. We obtain the following theorem as an important application:

\begin{corollary}\label{Hoopfft}(\cref{Hoopf})
Let $\mC$ be a symmetric monoidal $\infty$-category compatible with small colimits
and $\mO$ a $\mC$-enriched Hopf $\infty$-operad.
The $\infty$-category $\Alg_\mO(\mC) $ carries a canonical symmetric monoidal structure such that the forgetful functor $\Alg_\mO(\mC) \to \mC$ is symmetric monoidal.
	
\end{corollary}

\cref{Hoopfft} was applied by \cite[Proposition 2.55.]{heine2024restricted} to endow the $\infty$-category of
algebras over the spectral $\bE_\infty$-operad with a symmetric monoidal structure in order to prove an equivalence between algebras over the spectral $\bE_\infty$-operad and commutative algebras in the sense of Lurie \cite[Proposition 2.56.]{heine2024restricted}.

As another application we endow the $\infty$-category of algebras in every
presentably symmetric monoidal $\mV$-enriched $\infty$-category with $\mV$-enrichment, where $\mV$ is a cartesian closed presentable $\infty$-category:

\begin{theorem}\label{sqat0}(\cref{enrlift})
Let $\mV$ be a cartesian closed presentable $\infty$-category and 
$\mC$ a presentably symmetric monoidal $\mV$-enriched $\infty$-category.
The monoidal functor $\sSeq(\mC) \to \Fun(\mC,\mC)$ lifts to a monoidal functor
$$\sSeq(\mC) \to \mV\mathrm{-}\Fun(\mC,\mC)$$
from the $\infty$-category of symmetric sequences in $\mC$ 
to the $\infty$-category of $\mV$-enriched endofunctors of $\mC.$

\end{theorem}

\cref{sqat0} guarantees that the monad associated to a $\mC$-enriched $\infty$-operad underlies a $\mV$-enriched monad, whose
$\infty$-category of algebras underlies a $\mV$-enriched $\infty$-category. We obtain the following result:

\begin{corollary}(\cref{enrist})
Let $\mV$ be a cartesian closed presentable $\infty$-category and $\mC$ a presentably symmetric monoidal $\mV$-enriched $\infty$-category.

For every $\mC$-enriched $\infty$-operad $\mO$ the $\infty$-category of $\mO$-algebras in $\mC$ underlies a $\mV$-enriched $\infty$-category that is tensored and cotensored over $\mV.$
	
\end{corollary}

\cref{enrist} was applied by \cite{mocchetti2024coproduct} to stable motivic homotopy theory, where $\mV$ is the $\infty$-category of motivic spaces,
$\mC$ is the $\infty$-category of motivic spectra and $\mO$ is the
$\bE_\infty$-operad in motivic spectra.
In this case \cref{enrist} gives that the $\infty$-category of motivic $\bE_\infty$-ring spectra is enriched in motivic spaces,
which was used in \cite{mocchetti2024coproduct} to study the tensor of any motivic sphere with rational motivic cohomology.

%\subsection{Strategy of proof}

The equivalence of \cref{uikp} at the level of objects 
is a rather formal consequence of the presentation of Eilenberg-Moore objects as a lax limit.
By contrast, extending the correspondence to the level of morphisms between 
monads on different objects is highly non-trivial.
To deal with morphisms, we send diagrams of monads to diagrams of Eilenberg-Moore objects. Here our main insight is a classification of diagrams of monads as a single monad in a suitable $(\infty,2)$-category of diagrams.
This is technically challenging since this $(\infty,2)$-category of diagrams
is not an $(\infty,2)$-category of functors and transformations but of oplax transformations. We use the identification of diagrams of monads with monads in diagrams to reduce \cref{uikp} to the analysis of Eilenberg-Moore objects in the $(\infty,2)$-category of diagrams.

To prove the classification of diagrams of monads by monads in a suitable $(\infty,2)$-category of diagrams,
we construct for every pair $(\mB,\mC)$ of $(\infty,2)$-categories
a specific enrichment of the $(\infty,2)$-category of
functors $\mB \to \mC$ and oplax transformations within the $\infty$-category of $(\infty,2)$-categories over $\mB$, which might be of independent interest.

\subsubsection{History of this work and relation to other work}

The first version \cite{heine2017equivalence} of this work, titled “About the equivalence between monads and monadic functors”, appeared on the arXiv in 2017. Subsequent developments in enriched $\infty$-category theory allowed us to simplify several proofs and constructions, making the framework more conceptual and broadly applicable.
After our work has been made public, Haugseng \cite[Theorem 1.2.]{Haugseng2020OnLT} gave an alternative proof of \cref{uikp}
for the case of $\mC= \infty\Cat$
or more generally if the $(\infty,2)$-category $\mC$ is modeled by an $\infty$-cosmos in the sense of \cite{riehl2016infinity}. 
Haugseng's proof relies on two strong results proven in the setting of $\infty$-cosmoi: the coherent refinement of adjunctions to coherent adjunctions of Riehl-Verity \cite{RIEHL2016802}
and an equivalence of Zaganidis \cite{zaganidis2017towards}
between monadic coherent adjunctions and monads in the sense of Riehl-Verity \cite{RIEHL2016802}, which Haugseng identifies with monads in our sense.
In contrast, our equivalence between monads and monadic morphisms is direct and applies under weaker assumptions.

In 2023 we added a monadicity theorem for higher algebraic structures 
to our first version and changed the title to “A monadicity theorem for higher algebraic structures”. This monadicity theorem characterized the essential image of the duality between monads and monadic morphisms in a wide class of examples. 

\subsubsection*{Notation and terminology}

We fix a hierarchy of Grothendieck universes whose objects we call small, large, very large, etc.
We call a space small, large, etc. if its set of path components and its homotopy groups are small, large, respectively, for any choice of base point. We call an $\infty$-category small, large, etc if its maximal subspace and all its mapping spaces are small, large, respectively.

\vspace{2mm}

We write 
\begin{itemize}
\item $\Set$ for the category of small sets.
\item $\Delta$ for the category of finite, non-empty, partially ordered sets and order preserving maps, whose objects we denote by $[\n] = \{0 < ... < \n\}$ for $\n \geq 0$.

\item $\Fin_*$ for the category of finite pointed sets and base point preserving maps, whose objects we denote by $\langle\n\rangle = \{*, 1, ..., \n\}$ for $\n \geq 0$.

\item $\mS$ for the $\infty$-category of small spaces.
\item $ \infty\Cat$ for the $\infty$-category of small $\infty$-categories.
\item $\infty\Cat^{\rc \rc} $ for the $\infty$-category of large $\infty$-categories with small colimits and small colimits preserving functors.
\item $\Fun(\mC,\mD)$ for the $\infty$-category of functors between two $\infty$-categories $\mC,\mD.$
\end{itemize}

We indicate $\infty$-categories of large objects by $\widehat{(-)}$, for example we write $\widehat{\mS}, \infty\widehat{\Cat}$ for the $\infty$-categories of large spaces, $\infty$-categories, respectively.

%Note that $\Ho(\infty\Cat)$ is cartesian closed and for small $\infty$-categories $\mC,\mD$ we write $\Fun(\mC,\mD)$ for the internal hom, the $\infty$-category of functors $\mC \to \mD.$ 

%We denote by $\infty\Cat^{\rc \rc} \subset \infty\widehat{\Cat}$ the subcategory of $\infty$-categories with small colimits and small colimits preserving functors.
%By \cite{lurie.higheralgebra} Proposition 4.8.1.3. the $\infty$-category $\infty\Cat^{\rc \rc}$ carries a closed symmetric monoidal structure such that the subcategory inclusion $\infty\Cat^{\rc \rc} \subset \infty\widehat{\Cat}$ is lax symmetric monoidal.

\vspace{1mm}

For any $\infty$-category $\mC$ containing objects $\A, \B$ we write
\begin{itemize}
\item $\mC(\A,\B)$ for the space of maps $\A \to \B$ in $\mC$,
\item $\mC_{/\A}$ for the $\infty$-category of objects over $\A$,
\item $\Ho(\mC)$ for the homotopy category,
\item $\mC^{\triangleleft}, \mC^{\triangleright}$ for the $\infty$-category arising from $\mC$ by adding an initial, final object, respectively,
\item $\mC^\simeq $ for the maximal subspace in $\mC$.
\end{itemize}

\subsubsection*{Inclusions and embeddings}
We call a fully faithful functor an embedding.
We call a functor $\phi: \mC \to \mD$ a  (subcategory) inclusion if one of the following equivalent conditions holds:
\begin{itemize}
\item For every $\infty$-category $\mB$ the induced map
$\infty\Cat(\mB,\mC) \to \infty\Cat(\mB,\mD)$ is an embedding.
\item $\phi: \mC \to \mD$ induces an embedding on maximal subspaces and on all mapping spaces.
\item The functor $\Ho(\phi):\Ho(\mC) \to \Ho(\mD)$ is an inclusion and the functor $\mC \to \Ho(\mC) \times_{\Ho(\mD)} \mD$ is an equivalence.
\end{itemize}

%In this case $\phi$ is uniquely determined by $\mD$ and $\Ho(\phi): \Ho(\mC) \to \Ho(\mD).$\vspace{2mm}

%We call a monoidal $\infty$-category compatible with small colimits if it admits small colimits, which are preserved by the tensor product in each component.
%We call a monoidal $\infty$-category presentable if it is compatible with small colimits and its underlying $\infty$-category is presentable.
%\vspace{2mm}

\subsubsection*{Relative cocartesian fibrations}

\begin{itemize}
\item A functor $\phi: \X \to \rS$ is a locally (co)cartesian fibration relative to some full subcategory $\mE \subset \Fun([1],\rS)$ if for every functor $\alpha: [1] \to \rS$ that belongs to $\mE$ the pullback $[1] \times_\rS \X \to [1]$ along $\alpha$ is a (co)cartesian fibration.

\item A (co)cartesian fibration $\phi: \X \to \rS$ relative to $\mE$ is a locally (co)cartesian fibration relative to $\mE$ such that for every functor $\alpha: [1] \to \rS$ that belongs to $\mE$ the pullback $[1] \times_\rS \X \to [1]$ along $\alpha$ is a (co)cartesian fibration whose (co)cartesian morphisms are $\phi$-(co)cartesian.

\item A functor $\phi: \X \to \rS$ is a (co)cartesian fibration relative to some full subcategories $\mE \subset \Fun([1],\rS)$ and $ \T \subset \Fun([2],\T)$ if it is a locally (co)cartesian fibration relative to $\mE $ and for every functor $\sigma: [2] \to \rS$ that belongs to $\T$ the pullback $[2] \times_\rS \X \to [2]$ along $\sigma$ is a (co)cartesian fibration relative to $\sigma^{-1}(\mE).$

\end{itemize}

\vspace{1mm}

For every cocartesian fibration $\mC \to \rS$ classifying a functor $\alpha: \rS \to \infty\Cat $ the opposite cocartesian fibration
$\mC^\rev \to \rS$ is the cocartesian fibration classifying the functor
$\rS \xrightarrow{\alpha} \infty\Cat \xrightarrow{(-)^\op} \infty\Cat.$

\subsubsection{Relative adjunctions}

A functor $\mC \to \mD$ over an $\infty$-category $\mB$ admits a right adjoint relative to $\mB$ if it admits a right adjoint and the counit (or equivalently by the triangle identities the unit) lies over an equivalence in $\mB$. 
There is the dual notion of left adjoint relative to $\mB.$

A map $\mC \to \mD$ of cocartesian fibrations over $\mB$ admits a right adjoint relative to $\mB$ if and only if it induces on the fiber over every object of $\mB$ a left adjoint functor \cite[Proposition 7.3.2.6.]{lurie.HTT}.

% \begin{definition}Let $\mB$ be an $\infty$-category.
% A functor $\mC \to \mD$ over $\mB$ admits a right adjoint relative to $\mB$ if it admits a right adjoint and the counit (or equivalently by the triangle identities the unit) lies over an equivalence in $\mB$. 
% There is the dual notion of left adjoint relative to $\mB.$

% \end{definition}

% \begin{remark}\label{cocartes}
% A functor $\mC \to \mD$ over $\mB$ that admits a right adjoint relative to $\mB$, preserves morphisms cocartesian over $\mB$. This follows immediately from the definition of cocartesian morphism.

% \end{remark}

% \begin{remark}\label{cocartesi}
% A map $\mC \to \mD$ of cocartesian fibrations over $\mB$ admits a right adjoint relative to $\mB$ if and only if it induces on the fiber over every object of $\mB$ a left adjoint functor \cite[Proposition 7.3.2.6.]{lurie.HTT}.

% \end{remark}

% \begin{remark}\label{twoca}
% A functor $\mC \to \mD$ over $\mB$
% that admits a right adjoint relative to $\mB $ uniquely refines to an adjunction in the $(\infty,2)$-category $\infty\Cat_{/\mB}$ \cite[Proposition 7.3.2.1.]{lurie.higheralgebra}.    
% \end{remark}

\subsubsection{Acknowledgements}

%Especially I thank my advisor Markus Spitzweck for all the time he spent answering my questions and whoshowed me the basic ideas on which this topic is based on. 

We thank Peter Arndt, Hongyi Chu, David Gepner, Thomas Nikolaus 
and Markus Spitzweck for helpful and inspiring discussions.

\label{rel}

\section{Enriched $\infty$-categories}\label{Enrici}

\subsection{Monoidal \texorpdfstring{$\infty$}{∞}-categories and their actions}

Next we define monoidal $\infty$-categories.

\begin{notation}\label{ooop}
A map $[\n] \to [\m]$ in $\Delta$
\begin{itemize}
\item is inert if it is of the form $[\m] \simeq \{\bi, \bi+1,..., \bi+\m \} \subset [\n]$ for some $\bi \geq 0.$

\item is active if it preserves the minimum and maximum.

\item preserves the minimum if it is a map $[\m] \to [\n]$ sending $0$ to $0$.

\item preserves the maximum if it is a map $[\m] \to [\n]$ sending $\m$ to $\n$.

\end{itemize}

\end{notation}

\begin{remark}

For every $\n \geq 0$ there are $\n$ inert maps $[1] \to [\n]$, where the $\bi$-th inert map $[1] \to [\n]$ for $1 \leq \bi \leq \n$ is the map $[1] \simeq \{\bi-1, \bi\} \subset [\n]$.

\end{remark}

\begin{definition}\label{ek}
A monoidal $\infty$-category is a cocartesian fibration $ \mV^\ot \to \Delta^\op$ such that for every $\n \geq 0$ the family of inert maps $[1]\to [\n]$ in $\Delta$ induces an equivalence $\mV^\ot_{[\n]} \to \mV^{\times \n}, $ where we set $\mV:=\mV^\ot_{[1]}.$
	
\end{definition}

\begin{notation}

For every monoidal $\infty$-category $\phi: \mV^\ot \to \Delta^\op $ we call $ \mV := \mV^\ot_{[1]} $ the underlying $\infty$-category of $\phi$.

\end{notation}

\begin{notation}

We usually denote a monoidal $\infty$-category $\mV^\ot \to \Delta^\op$ by $\mV$ leaving the monoidal structure implicite.

We write $\tu_{\mV}$ for the tensor unit, the object $ * \simeq \mV^\ot_{[0]} \to \mV^\ot_{[1]}= \mV $ induced by the map $[1] \to [0]$ in $\Delta.$

We write $\ot:\mV \times \mV \simeq \mV^\ot_{[2]} \to \mV^\ot_{[1]}= \mV$ for the functor induced by the map $[1]=\{0,2\} \subset [2]$ in $\Delta.$
    
\end{notation}

\begin{definition}\label{compost} Let $\mK \subset \infty\Cat$ be a full subcategory.

\begin{enumerate}
\item A monoidal $\infty$-category $\mV$ is compatible with $\mK$-indexed colimits if the underlying $\infty$-category admits $\mK$-indexed colimits and the tensor product $\ot: \mV \times \mV \to \mV$ preserves $\mK$-indexed colimits componentwise.

\item A monoidal $\infty$-category $\mV$ is closed if the tensor product $\ot: \mV \times \mV \to \mV$ is left adjoint componentwise.
\item A presentably monoidal $\infty$-category is a closed monoidal $\infty$-category whose underlying $\infty$-category is presentable.

%\item A $\kappa$-compactly generated monoidal $\infty$-category is a closed monoidal $\infty$-category whose underlying $\infty$-category is $\kappa$-presentable and such that the tensor unit is $\kappa$-compact and the tensor product of two $\kappa$-compact objects is $\kappa$-compact.
\end{enumerate}
\end{definition}
% \begin{remark}
% By \cite[Proposition 7.15.]{Rune} every presentably monoidal $\infty$-category is $\kappa$-compactly generated for some small regular cardinal $\kappa.$
% \end{remark}

\begin{definition}
Let $\mV, \mW$ be monoidal $\infty$-categories.

\begin{itemize}
\item A lax monoidal functor $\mV \to \mW $ is a functor $\mV^\ot \to \mW^\ot$ over $\Delta^\op$ that preserves cocartesian lifts of inert maps.

\item A monoidal functor $\mV \to \mW $ is a functor $\mV^\ot \to \mW^\ot$ over $\Delta^\op$ that preserves morphisms cocartesian over $\Delta^\op.$

\end{itemize}

\end{definition}

\begin{definition}
A lax monoidal functor $\mV \to \mW $ is an embedding 
if the functor $\mV^\ot \to \mW^\ot$ is fully faithful.
    
\end{definition}

\begin{definition}Let $\mV$ be a monoidal $\infty$-category.
An associative algebra in $\mV$ is a lax monoidal functor $* \to \mV$.
    
\end{definition}

\begin{notation} Let $\mV$ be a monoidal $\infty$-category.

Let $$\Alg(\mV) \subset \Fun_{\Delta^\op}(\Delta^\op, \mV^\ot) $$ be the full subcategory of associative algebras in $\mV.$

\end{notation}

\begin{notation}We fix the following notation:
\begin{enumerate}

\item Let $\Lax\Mon \subset \infty\Cat_{/\Delta^\op}$ be the subcategory of monoidal $\infty$-categories and lax monoidal functors.

\item Let $\Mon \subset \Lax\Mon $ be the subcategory of monoidal $\infty$-categories and monoidal functors.

%\item Let $\rc\rc\Mon \subset \widehat{\Mon}$ be the subcategory of monoidal $\infty$-categories compatible with small colimits and monoidal functors preserving small colimits.

%\item Let $\Pr\Mon \subset \rc\rc\Mon$ be the full subcategory of presentably monoidal $\infty$-categories.

\end{enumerate}

\end{notation}

% \begin{definition}
% Let $\mC$ be an $\infty$-category that admits finite products.
% A monoid in $\mC$ is a functor $\X: \Delta^\op \to \mC$
% such that for every
% $\n \geq 0$ the family of inert maps $[1] \to [\n]$ induces an equivalence $$\X([\n]) \to \X([1])^{\times \n}.$$

% \end{definition}

% \begin{notation}

% Let $$ \Mon(\mC) \subset \Fun(\Delta^\op, \mC) $$
% be the full subcategory of monoids in $\mC.$
 
% \end{notation}

% \begin{remark}

% The Grothendieck construction provides an equivalence
% $$ \Mon(\infty\Cat) \simeq \Mon.$$
% \end{remark}

\begin{notation}\label{invo}
The category $\Delta^\op$ carries a canonical involution sending $[\n] $ to $[\n]$ and a map $\f:[\n] \to [\m]$ to the map $[\n] \to [\m], \bi \mapsto \m-\f(\n-\bi).$ The involution on $\Delta^\op$ induces an involution on $\Mon(\mC)$ for every $\infty$-category $\mC$ that admits finite products, which forms the reversed monoid.

Moreover it induces an involution on $\infty\Cat_{/\Delta^\op}$ that restricts to involutions $(-)^\rev$ on $\Lax\Mon \subset \infty\Cat_{/\Delta^\op} $
and $\Mon \subset \Lax\Mon, $ which form the reversed monoidal $\infty$-category.

\end{notation}

Next we define left, right and bitensored $\infty$-categories.

\begin{definition} Let $\mV,\mW$ be monoidal $\infty$-categories.

\begin{enumerate}
\item An $\infty$-category left tensored over $\mV$ is a map $ \mM^\circledast \to \mV^\ot$ of cocartesian fibrations over $\Delta^\op$ such that for every $\n \geq 0$ the map $[0]\simeq \{\n\}\subset [\n]$ in $\Delta^\op$ induces an equivalence,
where we set $\mM:= \mM^\circledast_{[0]}:$
$$ \mM^\circledast_{[\n]} \to \mV^\ot_{[\n]} \times \mM. $$

\item An $\infty$-category right tensored over $\mV$ is a functor $ \mM^\circledast \to \mV^\ot$ such that the functor
$$ \mM^\circledast \xrightarrow{} \mV^\ot \simeq (\mV^\ot)^\rev $$
is an $\infty$-category left tensored over $\mV^\rev$.

\item An $\infty$-category bitensored over $\mV, \mW$ is a map $ \mM^\circledast \to \mV^\ot \times \mW^\ot$ of cocartesian fibrations over $\Delta^\op \times \Delta^\op$ such that for every $\n,\m \geq 0$ the maps $[0]\simeq \{\n\}\subset [\n]$ and $[0]\simeq \{0\}\subset [\m]$ induce an equivalence, where we set $\mM:= \mM^\circledast_{[0][0]}:$
$$ \mM^\circledast_{[\n][\m]} \to \mV^\ot_{[\n]} \times \mM \times \mW^\ot_{[\m]}.$$

\end{enumerate}
\end{definition} 

\begin{notation}

For every $\infty$-category $\phi: \mM^\circledast \to \mV^\ot $ left (right) tensored over $\mV$ we call $\mM:= \mM^\circledast_{[0]} $ the underlying $\infty$-category of $\phi$. We say that $\phi$ exhibits $\mM$ as left (right) tensored over $\mV$.

For every $\infty$-category $\phi: \mM^\circledast \to \mV^\ot \times \mW^\ot $ bitensored over $\mV, \mW$ we call $\mM:= \mM^\circledast_{[0] [0]} $ the underlying $\infty$-category of $\phi$. We say that $\phi$ exhibits $\mM$ as bitensored over $\mV, \mW$.

\end{notation}

\begin{definition}\label{tenkap} Let $\mV$ be a monoidal $\infty$-category.

\begin{enumerate}

\item An $\infty$-category $\mM^\circledast \to \mV^\ot $ left tensored over $\mV$ is compatible with small colimits if $\mM$ admits small colimits, $\mV$ is a monoidal $\infty$-category compatible with small colimits
and for every $\V \in \mV, \X \in \mC$ the functors $(-) \ot \X: \mV \to \mM, \V \ot (-): \mM \to \mM$ preserve small colimits.

\item An $\infty$-category presentably left tensored over $\mV$ is an $\infty$-category left tensored over $\mV$ compatible with small colimits whose underlying $\infty$-category is presentable and such that $\mV$ is a presentably monoidal $\infty$-category.

% \item An $\infty$-category $\mM$ left tensored over $\mV$ is $\kappa$-compactly generated if the underlying $\infty$-category is $\kappa$-presentable, $\mV$ is a $\kappa$-compactly generated monoidal $\infty$-category and the left $\mV$-action on $\mM$ restricts to a left $\mV^\kappa$-action on $\mM^\kappa$. 

\end{enumerate}	
\end{definition}

We make a similar definition for right and bitensored $\infty$-categories.

\begin{notation}
Let $$\L\Mod \subset \Fun([1],\infty\Cat_{/ \Delta^\op}), \ \R\Mod \subset \Fun([1],\infty\Cat_{/ \Delta^\op}) $$ be the subcategory of left (right) tensored $\infty$-categories and linear functors.

% \item Let $$ \B\Mod \subset \Fun([1],\infty\Cat_{/ \Delta^\op \times \Delta^\op}) $$ be the subcategory of bitensored $\infty$-categories and linear functors.

% \item Let $$\rc\rc\L\Mod \subset \widehat{\LMod}, \ \rc\rc\R\Mod \subset \widehat{\RMod}, \ \rc\rc\B\Mod \subset \widehat{\BMod} $$ be the respective subcategories of left, right, bitensored $\infty$-categories compatible with small colimits and linear functors preserving small colimits.

% \item Let $$\Pr\L\Mod \subset \rc\rc\L\Mod,\ \Pr\R\Mod \subset \rc\rc\R\Mod, \ \Pr\B\Mod \subset \rc\rc\B\Mod $$ be the respective full subcategories of presentably left, right, bitensored $\infty$-categories.

\end{notation}

\begin{remark}

The involution on $\Delta^\op$ taking the opposite category induces an involution on $\infty\Cat_{/\Delta^\op}$
that restricts to an equivalence $\LMod \simeq \RMod$.
\end{remark}

\subsection{Enriched \texorpdfstring{$\infty$}{∞}-categories}

In this subsection we recall the theory of enriched $\infty$-categories.

\begin{definition}\label{bla}
Let $\mV$ be a monoidal $\infty$-category.
An $\infty$-category weakly enriched in $\mV$ is a functor
$\phi: \mM^\circledast \to \mV^\ot $ 
such that the following conditions hold, where we set $ \mM:= \mM^\circledast_{[0]}:$
\begin{enumerate}

\item for every $\X \in \mM^\circledast$ lying over $[\n]$
every inert map $[\n] \to [\m]$ in $\Delta^\op$ preserving the maximum admits a cocartesian lift $\X \to \X'$ in $ \mM^\circledast.$

\item for every $\n \geq 0$ the unique map $[0]\to [\n]$ in $\Delta$ preserving the maximum induces an equivalence
$$ \theta: \mM^\circledast_{[\n]} \to \mV^\ot_{[\n]} \times \mM,$$
\vspace{1mm}
\item for every $\X,\Y \in \mM^\circledast$, where $[\n] \in \Delta$ is the image of $\X$, the cocartesian lift $\X \to \X'$ in $ \mM^\circledast $ of the unique inert map
$[0] \to [\n]$ of $\Delta$ preserving the maximum induces a pullback square:
\begin{equation*} 
\begin{xy}
\xymatrix{
\Map_{\mM^\circledast}(\Y,\X)  \ar[d]^{} \ar[r]^{ }
& \Map_{\mM^\circledast}(\Y,\X') \ar[d]^{} 
\\ \Map_{\mV^\ot}(\phi(\Y),\phi(\X))  
\ar[r]^{}  & \Map_{\mV^\ot}(\phi(\Y),\phi(\X')). 
}
\end{xy} 
\end{equation*}
\end{enumerate}

\end{definition}

\begin{notation}

For every weakly $\mV$-enriched $\infty$-category $\phi: \mM^\circledast \to \mV^\ot $ we call $\mM:= \mM^\circledast_{[0]} $ the underlying $\infty$-category of $\phi$. We say that $\phi$ exhibits $\mM$ as weakly enriched in $\mV$. We usually denote an $\infty$-category $\phi: \mM^\circledast \to \mV^\ot $ weakly enriched in $\mV$ by $\mM$ leaving the functor $\phi$ notationally implicite.

\end{notation}

\begin{notation}\label{mult}Let $\mV$ be a monoidal $\infty$-category,
$\mM$ a weakly $\mV$-enriched $\infty$-category and $\V_1,..., \V_{\n} \in \mV, \X, \Y \in \mM$ for $\n \geq 0 $. The space of multimorphisms
$\V_1,..., \V_{\n},  \X \to \Y $ in $\mM$ is the full subspace
$$\Mul_{\mM}(\V_1,..., \V_\n,\X, \Y) \subset \Map_{\mM^\circledast}(\Z,\Y)$$ spanned by the morphisms $\Z \to \Y$ in $\mM^\circledast$ lying over the unique map 
$[0] \to [\n]$ in $\Delta$ preserving the minimum, 
where $\Z \in \mM_{[\n]}^\circledast \simeq \mV^{\times \n} \times \mM$ corresponds to $(\V_1,..., \V_\n,\X) $.

\end{notation}

The next definition is \cite[Definition 2.120.]{heine2024bienriched}:

\begin{definition}Let $\mV$ be a monoidal $\infty$-category.
An $\infty$-category pseudo-enriched in $\mV$ is an $\infty$-category
$\mM^\circledast \to \mV^\ot $ weakly enriched in $\mV$ such that for every $\X,\Y \in \mM $ and every $\n \geq 0$ the universal multimorphism
$\V_1, ..., \V_\n \to \V_1 \ot ... \ot \V_\n $ in $\mV$ induces an equivalence
$$ \Mul_\mM(\V_1 \ot ... \ot \V_\n, \X, \Y) \to  \Mul_\mM(\V_1, ..., \V_\n, \X, \Y).$$

\end{definition}

\begin{example}

Every $\infty$-category left tensored over $\mV$ is an $\infty$-category pseudo-enriched in $\mV$.
    
\end{example}

Next we define enriched $\infty$-categories. The next definition is \cite[Definition 2.128.]{heine2024bienriched}:

\begin{definition}Let $\mV$ be a monoidal $\infty$-category
and $\mM $ an $\infty$-category pseudo-enriched in $\mV.$
A morphism object of $ \X, \Y \in \mM$ is an object $$\Mor_{\mC}(\X, \Y) \in \mV $$ such that there is a multimorphism $\beta \in \Mul_{\mM}(\Mor_{\mM}(\X, \Y), \X, \Y) $ that induces for every $\V \in \mV $ an equivalence
$$ \Map_{\mV}(\V, \Mor_{\mM}(\X, \Y)) \simeq \Mul_{\mM}(\V, \X, \Y).$$	
	
\end{definition}

The next definition is \cite[Definition 2.133.]{heine2024bienriched}:

\begin{definition}Let $\mV$ be a monoidal $\infty$-category.
A $\mV$-enriched $\infty$-category is an $\infty$-category $\mM$ pseudo-enriched in $\mV$ such that for every $\X, \Y \in \mM $ there is a morphism object $\Mor_{\mM}( \X, \Y) \in \mV $.

\end{definition}

\begin{notation}

For every $\infty$-category $\phi: \mM^\circledast \to \mV^\ot $ 
(pseudo-) enriched in $\mV$ we call $\mM:= \mM^\circledast_{[0]} $ the underlying $\infty$-category of $\phi$ and say that $\phi$ exhibits $\mM$ as (pseudo-) enriched in $\mV$. We usually denote an $\infty$-category $\phi: \mM^\circledast \to \mV^\ot $ (pseudo-) enriched in $\mV$ by $\mM$ leaving the functor $\phi$ notationally implicite.

\end{notation}

\begin{example}\label{Biolk}

Every $\infty$-category $\mM$ presentably left tensored over $\mV$ is a $\mV$-enriched $\infty$-category.

For every $\X \in \mM$ the functor $(-) \ot \X : \mV \to \mM$ admits a right adjoint, which sends $\Y \in \mM$ to the morphism object $\Mor_{\mM}(\X, \Y) \in \mV $.
In particular, for every presentably monoidal $\infty$-category $\mV$ the left action of $\mV$ on itself of \cref{selfbit} exhibits $\mV$ as $\mV$-enriched.

\end{example}

\begin{example}\label{enrsub}Let $\mV$ be a monoidal $\infty$-category,
$\mM $ a weakly $\mV$-enriched $\infty$-category and $\mN \subset \mM$ a full subcategory.
Let $\mN^\circledast \subset \mM^\circledast$ be the full subcategory of objects lying over some $\V \in \mV^\ot $ corresponding to some object of $\mN \subset \mM.$
The restriction $\mN^\circledast \subset \mM^\circledast \to \mV^\ot$
is a weakly $\mV$-enriched $\infty$-category whose underlying $\infty$-category is $\mN$.
If $\mM $ is a $\mV$-enriched $\infty$-category, then $\mN^\circledast \subset \mM^\circledast \to \mV^\ot$ is a $\mV$-enriched $\infty$-category.
We call $\mN^\circledast \subset \mM^\circledast \to \mV^\ot$ the full (weakly) $\mV$-enriched subcategory of $ \mM$ spanned by $\mN.$

\end{example}

% \begin{definition} Let $\mV$ be a monoidal $\infty$-category. A $\mV$-enriched $\infty$-category $ \mM^\circledast \to \mV^\ot $ is small if the collection of equivalence classes of objects of $\mM$ is small.

% \end{definition}

% \begin{remark}\label{smalo} Let $\mV$ be a monoidal $\infty$-category and $ \mM $ a $\mV$-enriched $\infty$-category.
% If $\mV$ is locally small, also $\mM$ is locally small since for every $\X, \Y \in \mM$ there is a canonical equivalence
% $$ \Map_{\mM}(\X,\Y) \simeq \Map_{\mV}(\tu_\mV,\Mor_\mC(\X,\Y)). $$

% Thus if $\mV$ is locally small, a $\mV$-enriched $\infty$-category $ \mM^\circledast \to \mV^\ot $ is small if and only if $\mM$ is small.

% \end{remark}

Next we generalize the notion of left tensored $\infty$-category by introducing the notion of tensor.

\begin{definition}

Let $\mV$ be a monoidal $\infty$-category, $\mM$ an $\infty$-category
pseudo-enriched in $\mV$ and $\V \in \mV, \X,\Y \in \mM$.

\begin{enumerate}
\item A multimorphism $\psi : \V, \X \to \Y $ in $\mM$ exhibits $\Y$ as the tensor of $\V, \X$, denoted by $\V \ot \X$, if for every $\Z \in \mM, \W \in \mV $ the following induced map is an equivalence:
\begin{equation*}
\Mul_\mM(\W, \Y, \Z) \to \Mul_\mM(\W, \V, \X, \Z).
\end{equation*}

\item A multimorphism $\psi : \V, \Y \to \X$ in $\mM$ exhibits $\Y$ as the cotensor of $\V, \X$, denoted by $^\V \X$, if for every $\Z \in \mM, \W \in \mV $ the following induced map is an equivalence:
\begin{equation*}
\Mul_\mM(\W, \Z, \Y) \to \Mul_\mM(\W, \V,\Z, \X).
\end{equation*}
\end{enumerate}
\end{definition}

We usually consider tensors and cotensors in $\infty$-categories enriched in a presentably monoidal $\infty$-category, which by \cref{Biolk} is enriched in itself.
In this case the definition of (co)tensors simplifies:

\begin{example}

Let $\mV$ be a presentably monoidal $\infty$-category, $\mM$ a $\mV$-enriched $\infty$-category and $\V \in \mV, \X,\Y \in \mM$.
	
\begin{enumerate}
\item A morphism $\psi : \V \to \Mor_\mM(\X, \Y)$ in $\mV$ exhibits $\Y$ as the tensor of $\V, \X$, denoted by $\V \ot \X$, if for every $\Z \in \mM$ the following induced morphism is an equivalence:
\begin{equation*}
\Mor_\mM(\Y, \Z) \to \Mor_\mV(\V,\Mor_\mM(\X, \Z)).
\end{equation*}

\item A morphism $\psi : \V \to \Mor_\mM(\Y, \X)$ in $\mV$ exhibits $\Y$ as the cotensor of $\V, \X$, denoted by $^\V \X$, if for every $\Z \in \mM$ the following induced morphism is an equivalence:
\begin{equation*}
\Mor_\mM(\Z, \Y) \to \Mor_\mV(\V,\Mor_\mM(\Z, \X)).
\end{equation*}
\end{enumerate}

\end{example}

The following is \cite[Lemma 2.122.]{heine2024bienriched}:

\begin{example}\label{exiq} Let $\mV$ be a monoidal $\infty$-category.
An $\infty$-category pseudo-enriched in $\mV$ is an $\infty$-category left tensored over $\mV$ if and only if it admits tensors.
    
\end{example}

%\subsection{Enriched functors}

Next we define enriched functors. 
The next definition is \cite[Definition 2.50.]{heine2024bienriched}:

\begin{definition}\label{linmapp}
Let $\mV, \mW$ be monoidal $\infty$-categories, $\alpha: \mV \to \mW$ a lax monoidal functor and $\mM, \mN$ weakly enriched $\infty$-categories.
An enriched functor $\mM \to \mN$ is a commutative square of $\infty$-categories
\begin{equation}\label{curv}
\begin{xy}
\xymatrix{
\mM^\circledast  \ar[d]^{} \ar[rr]^{\gamma}
&&\mN^\circledast \ar[d]^{} 
\\ \mV^\ot
\ar[rr]^{\alpha}  && \mW^\ot
}
\end{xy} 
\end{equation}
such that $\gamma$ preserves cocartesian lifts of inert morphisms of $\Delta^\op$ that preserves the maximum.

\begin{itemize}
 
\item An enriched functor is linear if it preserves tensors.

\item An enriched functor $\mM \to \mN$ is $\mV$-enriched if $\alpha$ is the identity.

\item An enriched functor is $\mV$-linear if it is $\mV$-enriched and linear.

\item A $\mV$-enriched functor admits a $\mV$-enriched left (right) adjoint if it admits a left (right) adjoint relative to $\mV^\ot.$

\item A $\mV$-enriched functor $\mM \to \mN $ is an embedding if $\gamma$ is fully faithful.

\item An enriched functor is essentially surjective if it induces an 
essentially surjective functor on underlying $\infty$-categories.

\end{itemize}
\end{definition}

\begin{notation}Let $\mV$ be a monoidal $\infty$-category and $\mM,\mN$ weakly $\mV$-enriched $\infty$-categories.
Let $$ \mV\mathrm{-}\Fun(\mM,\mN) \subset \Fun_{\mV^\ot}(\mM^\circledast,\mN^\circledast)$$ be the full subcategory of $\mV$-enriched functors.

\end{notation} 

\begin{notation}\label{linea} Let $\mV$ be a monoidal $\infty$-category and $\mM,\mN$ weakly $\mV$-enriched $\infty$-categories. Let $$ \mV-\LinFun(\mM,\mN) \subset \mV-\Fun(\mM,\mN) $$ be the full subcategory of $\mV$-linear functors.

Let $$ \mV-\LinFun^\L(\mM,\mN) \subset \mV-\LinFun(\mM,\mN) $$ be the full subcategory of $\mV$-linear functors that admit a right adjoint.

\end{notation}

% \begin{notation}Let $\mV$ be a monoidal $\infty$-category and $\mM,\mN$ weakly $\mV$-enriched $\infty$-categories. Let $$ \mV-\Fun^\L(\mM,\mN) \subset \mV-\Fun(\mM,\mN) $$ be the full subcategory spanned by the $\mV$-enriched functors that admit a $\mV$-enriched right adjoint. 

% Let $$ \mV-\Fun^\R(\mM,\mN) \subset \mV-\Fun(\mM,\mN) $$ be the full subcategory spanned by the $\mV$-enriched functors that admit a $\mV$-enriched left adjoint. 

% \end{notation}

% \begin{example}
% Let $\mV$ be a monoidal $\infty$-category and $\mM, \mN$ be $\infty$-categories left tensored over $\mV$.
% A $\mV$-enriched functor $\mM \to \mN $
% is $\mV$-linear if and only if it is a left $\mV$-linear functor of \cref{linea}.

%\end{example}

\begin{remark}\label{lemuta}Let $\mV$ be a monoidal $\infty$-category.
For every left $\mV$-tensored $\infty$-category $\mM$ the functor $\mV-\LinFun(\mV,\mM) \to \mM $ evaluating at the tensor unit is an equivalence by \cite[Corollary 4.2.4.7.]{lurie.higheralgebra}.

\end{remark}

% The following lemma is \cite[Lemma 10.1.]{HEINE2023108941}.

% \begin{lemma}\label{laxxl}Let $\mV$ be a monoidal $\infty$-category,
% $\mM$ a $\mV$-enriched $\infty$-category and $\X \in \mM$. The functor $$\Mor_\mM(\X,-):\mM \to \mV $$ refines to a $\mV$-enriched functor. 
% \end{lemma}

The following is \cite[example 2.30.]{heine2024bienriched}:

\begin{example}\label{euuz}Let $\mV$ be a monoidal $\infty$-category, 
$\mM $ a weakly $\mV$-enriched $\infty$-category and $\K$ an $\infty$-category.
\begin{itemize}
\item The functor $$ \mM^\circledast \times \K \to \mV^\ot $$	
is a weakly $\mV$-enriched $\infty$-category that exhibits $\mM \times \K$ as weakly enriched in $\mV.$
	
\item The pullback $$(\mM^\K)^\circledast:=\mV^\ot \times_{\Fun(\K,\mV^\ot)} \Fun(\K,\mM^\circledast) \to \mV^\ot $$ along the diagonal functor is a weakly $\mV$-enriched $\infty$-category that exhibits $\Fun(\K,\mM)$ as weakly enriched in $\mV.$
\end{itemize}
\end{example}

For the next remark we use the notation of \cref{euuz}:

\begin{remark}\label{2-cat2}Let $\mV$ be a monoidal $\infty$-category,
$\mM, \mN$ weakly $\mV$-enriched $\infty$-categories and $\K$ an $\infty$-category.
The canonical equivalences $$ \Fun(\K,\Fun_{\mV^\ot}(\mM^\circledast,\mN^\circledast)) \simeq \Fun_{\mV^\ot}(\mM^\circledast \times \K,\mN^\circledast)
\simeq \Fun_{\mV^\ot}(\mM^\circledast,(\mN^\K)^\circledast)$$
restrict to equivalences
$$ \Fun(\K,\mV\mathrm{-}\Fun(\mM,\mN)) \simeq \mV\mathrm{-}\Fun(\mM \times \K,\mN) \simeq\mV\mathrm{-}\Fun(\mM,\mN^\K).$$

\end{remark}

The following is \cite[Lemma 2.77.]{heine2024bienriched}:

\begin{lemma}\label{Adj}\label{remqa} Let $\mV$ be a monoidal $\infty$-category and let $\mM,\mN$ be $\mV$-enriched $\infty$-categories such that $\mM$ admits tensors.
A $\mV$-enriched functor $\mM \to \mN $ admits a $\mV$-enriched right adjoint if and only if it is $\mV$-linear and the underlying functor admits a right adjoint.

\end{lemma}

\begin{notation}We fix the following notation:

\begin{itemize}

\item Let $$\omega\Enr \subset \Fun([1],\infty\Cat_{/ \Delta^\op}) $$ be the subcategory of weakly enriched $\infty$-categories and enriched functors.

% \item Let $$ \P\Enr \subset \omega\Enr $$
% be the full subcategory of pseudo-enriched $\infty$-categories.

\item Let $$\Enr \subset \omega\Enr $$ be the full subcategory of enriched $\infty$-categories.

% \item Let $$\Enr\Cat \subset \widehat{\Enr} $$ be the subcategory of small $\infty$-categories enriched in a presentably monoidal $\infty$-category, and enriched functors lying over a left adjoint monoidal functor.

\end{itemize}

\end{notation}

\begin{remark}
Evaluation at the target restricts to forgetful functors $$ \omega\Enr \to \Lax\Mon, \ \Enr \to \Lax\Mon.$$

\end{remark}

\begin{notation}\label{prenri} Let $\mV$ be a monoidal $\infty$-category.

\begin{itemize}
\item We denote the fiber of the forgetful functor $ \omega\Enr \to \Lax\Mon $ over $\mV$ by $$ \mV\mathrm{-}\omega\Enr.$$

% \item We denote the fiber of the forgetful functor $ \P\Enr \to \Lax\Mon $
% over $\mV$ by $$ \mV\mathrm{-}\P\Enr.$$

\item We denote the fiber of the forgetful functor $ \Enr \to \Lax\Mon $ over $\mV$ by $$ \mV\mathrm{-}\Cat.$$

% \item We denote the fiber of the forgetful functor $ \Enr\Cat \to \Pr\Mon $
% over every presentably monoidal $\infty$-category
% $\mV$ by $$ \mV\mathrm{-}\Cat.$$

\end{itemize}

\end{notation}

\begin{notation}

Evaluation at the source followed by taking the fiber over 
$[0] \in \Delta^\op$ gives functors
$$\omega\Enr \to \infty\Cat, \ \Enr \to \infty\Cat.$$
    
\end{notation}

\begin{example}\label{unenr}

By \cite[Corollary 3.23.]{heine2024bienriched} the forgetful functor
$ \mS \mathrm{-} \Cat \to \infty\Cat $
is an equivalence.
\end{example}

The following is \cite[Corollary 3.62.]{heine2024bienriched}:

\begin{proposition}\label{present} Let $\mV$ be a presentably monoidal $\infty$-category.
The $\infty$-category $\mV\mathrm{-}\Cat$ is presentable.
    
\end{proposition}

\begin{remark}\label{enr2cat} Let $\mV$ be a monoidal $\infty$-category.
The finite products preserving functor
$$\infty\Cat \to \mV\mathrm{-}\omega\Enr, \ \K \mapsto \mV \times \K$$
makes $\mV\mathrm{-}\omega\Enr$ to a right $\infty\Cat$-tensored $\infty$-category. This action is closed by \cref{2-cat2} and so exhibits 
$\mV\mathrm{-}\omega\Enr$ as an $\infty\Cat$-enriched $\infty$-category.
In particular, the full subcategory
$\mV\mathrm{-}\Cat \subset \mV\mathrm{-}\omega\Enr$ is an $\infty\Cat$-enriched $\infty$-category.
Moreover this left action restricts to a closed right action of
$\infty\Cat$ on $\mV\mathrm{-}\L\Mod$, which exhibits $\mV\mathrm{-}\L\Mod$
as an $\infty\Cat$-enriched $\infty$-category.

\end{remark}

For every enriched $\infty$-category there is an opposite enriched $\infty$-category by the following proposition, which is \cite[Proposition 5.45.]{heine2024bienriched}:

\begin{proposition}\label{oppoen} Let $\mV$ be a monoidal $\infty$-category. There is an involution $$(-)^\op: {\mV\mathrm{-}\Cat} \simeq {\mV^\rev}\mathrm{-}\Cat$$
fitting into a commutative square 
\begin{equation*} 
\begin{xy}
\xymatrix{
{\mV\mathrm{-}\Cat} \ar[d] \ar[rr]^{(-)^\op}
&& {\mV^\rev}\mathrm{-}\Cat \ar[d] 
\\ \infty\Cat
\ar[rr]^{(-)^\op}  &&  \infty\Cat
}
\end{xy} 
\end{equation*}
such that for every $\mV$-enriched $\infty$-category $\mC $ and $\X,\Y \in \mC$ there is an equivalence in $\mV$: $$\Mor_{\mC^\op}(\X,\Y) \simeq \Mor_\mC(\Y,\X).$$	
\end{proposition}

Next we consider functoriality of change of enrichment.
The next proposition is \cite[Theorem 3.67., Proposition 2.65.]{heine2024bienriched}:

\begin{proposition}\label{adjures}\label{bica1}\label{bica}

Let $\alpha: \mV \to \mW$ be a lax monoidal functor.

\begin{enumerate}
\item There is a $\infty\Cat$-linear adjunction
$$ \alpha_!: \mV\mathrm{-}\omega\Enr \rightleftarrows \mW\mathrm{-}\omega\Enr : \alpha^*.$$

\vspace{1mm}
\item The left adjoint restricts to a functor
$ \alpha_!: \mV\mathrm{-}\Cat \to \mW\mathrm{-}\Cat. $
For every $\mV$-enriched $\infty$-category $\mC$ the unit
$\rho:\mC \to \alpha^*(\alpha_!(\mC)) $
is essentially surjective and induces for every $\X, \Y \in \mC$
an equivalence $$ \alpha(\Mor_\mC(\X,\Y)) \simeq \Mor_{\alpha_!(\mC)}(\rho(\X),\rho(\Y)).$$

\vspace{1mm}
\item If $\alpha$ is a monoidal functor, the right adjoint restricts to a functor
$ \alpha^*: \mW\mathrm{-}\LMod \to \mV\mathrm{-}\LMod, $
which identifies with the functor that restricts the action along $\alpha.$

\vspace{1mm}
\item The right adjoint restricts to a functor
$ \alpha^*: \mW\mathrm{-}\Cat \to \mV\mathrm{-}\Cat $
if $\alpha$ is a monoidal functor and admits a right adjoint.

\end{enumerate}
\end{proposition}

% \begin{lemma}

% Let $\bj: \mV \to \mW$ be an embedding of monoidal $\infty$-categories.
% The induced functor 
% $ \mV\mathrm{-}\Cat \to \mW\mathrm{-}\Cat $
% is fully faithful. The essential image precisely consists
% of the $\mW$-enriched $\infty$-categories whose morphism objects lie in $\mV.$

% \end{lemma}

\begin{definition}\label{ind}

Let $\alpha: \mV \to \mW$ be a lax monoidal functor.
The weak $\mV$-enrichment on $\mW$ induced by $\alpha$ is the restriction along $\alpha$ of the canonical weak $\mW$-enrichment on $\mW.$

\end{definition}

\begin{remark}\label{classi}

Let $\mV$ be a presentably monoidal $\infty$-category and 
$\mM$ a small $\mV$-enriched $\infty$-category.
By \cref{unenr} and \cref{bica1} an object $\X \in \mM$ corresponds to an enriched functor
$* \to \mM$ lying over the unique left adjoint monoidal functor
$\mS \to \mV.$
    
\end{remark}

\begin{notation}

Let $\n \geq 2$. For every $1 \leq \bi \leq \n$ let $\mV_\bi$ be a presentably monoidal $\infty$-category and $\mM_\bi$ a $\mV_\bi$-enriched $\infty$-category.
Let $$ \mM_1 \boxtimes ... \boxtimes \mM_\n  := \alpha_!(\mM_1 \times ... \times \mM_\n) \in \mV_1 \ot ... \ot \mV_\n \mathrm{-}\Cat.$$
where $$\alpha: \mV_1 \times ...  \times \mV_\n \to \mV_1 \ot ... \ot \mV_\n $$ is the universal monoidal functor preserving small colimits componentwise.

\end{notation}

% \begin{remark}

% Let $1 \leq \bi \leq \n.$
% By \cref{classi} an object $\X \in \mM_\bi$ corresponds to an enriched functor $* \to \mM_\bi$ lying over the left adjoint monoidal functor
% $\mS \to \mV_\bi.$
% The latter gives rise to an enriched
% functor $ \mM_1 \times ... \times \mM_{\bi-1} \times * \times \mM_{\bi+1} \times ... \times \mM_\n \to \mM_1 \times ... \times \mM_\n \to \mM_1 \boxtimes ... \boxtimes \mM_\n $
% lying over the induced left adjoint monoidal functor
% $ \mV_1 \times ... \times \mV_{\bi-1} \times \mS \times \mV_{\bi+1} \times ... \times \mV_\n \to \mV_1 \times ... \times \mV_\n \to \mV_1 \otimes ... \otimes \mV_\n, $
% which corresponds to an enriched functor 
% $$ \mM_1 \boxtimes ... \boxtimes \mM_{\bi-1} \boxtimes \mM_{\bi+1} \boxtimes ... \boxtimes \mM_\n \to \mM_1 \boxtimes ... \boxtimes \mM_\n $$
% which lies over the induced left adjoint monoidal functor
% $ \mV_1 \otimes ... \otimes \mV_{\bi-1} \otimes \mV_{\bi+1} \otimes ... \otimes \mV_\n \to \mV_1 \otimes ... \otimes \mV_\n. $

% \end{remark}

\begin{example}\label{enrtenso}

Let $\mV$ be a presentably $\bE_2$-monoidal $\infty$-category.
%Then $\mV\mathrm{-}\Cat$ carries a monoidal structure by \cref{weakten}.
The presentably $\bE_2$-monoidal $\infty$-category $\mV$ provides a 
left adjoint monoidal tensor product functor $$ \ot: \mV \ot \mV \to \mV. $$ 
The tensor product of $\mV$-enriched $\infty$-categories 
$\mM, \mN$ is 
the $\mV$-enriched $\infty$-category $$ \mM \otimes^\mV \mN:=\ot_!(\mM \boxtimes \mN).$$
    
\end{example}

\subsection{Enriched Yoneda-embedding}

In this subsection we introduce the enriched $\infty$-category of enriched functors between two enriched $\infty$-categories.

The next theorem is \cite[Proposition 4.11., Theorem 4.86.]{heine2024bienriched}:

\begin{theorem}\label{psinho} Let $\mV, \mW$ be presentably monoidal $\infty$-categories and $\mN$ a $\mV \ot \mW$-enriched $\infty$-category.

\begin{enumerate}
\item Let $\mM$ be a small $\mV$-enriched $\infty$-category. 
The $\infty$-category $\mV\mathrm{-}\Fun(\mM, \mN) $
is enriched in $\mW$ and there is a $\mV \ot \mW$-enriched functor
$$ \mM \boxtimes (\mV\mathrm{-}\Fun(\mM, \mN)) \to \mN $$ 
such that for every $\mW$-enriched $\infty$-category $\mO$ the following induced functor is an equivalence
$$ \mW \mathrm{-}\Fun(\mO,\mV\mathrm{-}\Fun(\mM, \mN)) \to \mV \ot \mW\mathrm{-}\Fun(\mM \boxtimes \mO,\mN).$$

\item Let $\mO$ be a small $\mW$-enriched $\infty$-category. 
The $\infty$-category $\mW\mathrm{-}\Fun(\mO, \mN) $
is enriched in $\mV$ and there is a $\mV \ot \mW$-enriched functor
$$ (\mW\mathrm{-}\Fun(\mO, \mN)) \boxtimes \mO \to \mN $$ 
such that for every $\mV$-enriched $\infty$-category $\mM$ the following induced functor is an equivalence
$$ \mV \mathrm{-}\Fun(\mM,\mW\mathrm{-}\Fun(\mO, \mN)) \to \mV \ot \mW\mathrm{-}\Fun(\mM \boxtimes \mO,\mN).$$
\end{enumerate}

\end{theorem}

The following is \cite[Lemma 2.58.]{heine2024bienriched}:

\begin{theorem}\label{psinho2} Let $\mV, \mW $ be presentably monoidal $\infty$-categories and $ \mM $ a small $\mV$-enriched $\infty$-category. 
Let $\mN$ be an $\infty$-category presentably left tensored over $\mV \ot \mW$.
The $\mW$-enriched $\infty$-category $$\mV\mathrm{-}\Fun(\mM, {\mN}) $$ is a presentably left $\mW$-tensored $\infty$-category.

\end{theorem}

For the next definition we use that every monoidal $\infty$-category $\mV$ compatible with small colimits gives rise to a canonical left action of $\mV^\rev \otimes \mV$ on $\mV$ compatible with small colimits 
\cite[Lemma 3.77.]{heine2024bienriched}.

\begin{definition}\label{unitol} Let $\mV$ be a presentably monoidal $\infty$-category and $\mM $ a small $\mV$-enriched $\infty$-category. The $\mV$-enriched $\infty$-category of $\mV$-enriched presheaves on $\mM$ is
$$\mP_\mV(\mM) := \Fun_{\mV^\rev}(\mM^\op,\mV^\rev). $$	
	
\end{definition}

By \cite[Notation 4.56.]{heine2024bienriched} there is a $\mV$-enriched Yoneda-embedding $$\iota: \mM \to \mP_\mV(\mM), $$
which is a $\mV$-enriched embedding, and sends every $\X \in \mM $ to a $\mV^\rev$-enriched functor $ \mM^\op \to \mV^\rev $ whose underlying functor is the functor $$ \Mor_{\mM}(-,\X): \mM^\op \to \mV.$$

\vspace{2mm}

The next is \cite[Theorem 3.41.]{heine2024bienriched}.
It is an enriched version of the universal property of presheaves \cite[Theorem 5.1.5.6.]{lurie.HTT}, \cite[Proposition 3.16.]{heine2024local}:

\begin{theorem}\label{envvcor} Let $\mV $ be a presentably monoidal $\infty$-category, $\mM $ a small $\mV$-enriched $\infty$-category and $\mN$ a left tensored $\infty$-category compatible with small colimits.
The induced functor $$\mV\mathrm{-}\Fun(\mP_{\mV}(\mM), \mN) \to \mV\mathrm{-}\Fun(\mM, \mN)$$
admits a fully faithful left adjoint that lands in the full subcategory $ \mV\mathrm{-}\LinFun^\L(\mP_{\mV}(\mM), \mN) $ of $\mV$-linear functors that admit a right adjoint. In particular, the following functor is an equivalence: $$\mV\mathrm{-}\LinFun^\L(\mP_{\mV}(\mM), \mN) \to \mV\mathrm{-}\Fun(\mM, \mN).$$

\end{theorem}

% The next follows from \cite[Corollary 4.71.]{heine2024bienriched}:

% \begin{proposition}\label{yoadj}

% Let $\mV $ be a presentably monoidal $\infty$-category and $\mM, \mN $ small $\mV$-enriched $\infty$-categories.
% The left adjoint $\mV$-enriched functor $$
% \F_!: \mP_\mV(\mM) \to \mP_\mV(\mN) $$
% extending $\mM \xrightarrow{\F} \mN \subset \mP_\mV(\mN)$ of \cref{envvcor} is $\mV$-enriched left adjoint to the $\mV$-enriched functor $ \F^*: \mP_\mV(\mN) \xrightarrow{} \mP_\mV(\mM),$
% which itself admits a $\mV$-enriched right adjoint.

% \end{proposition}

The next is \cite[Theorem 4.43.]{heine2024bienriched} and the enriched version of the Yoneda-lemma:

\begin{theorem}\label{enryol}

Let $\mV$ be a presentably monoidal $\infty$-category and $\mM $ a small $\mV$-enriched $\infty$-category.
The $\mV$-enriched $\infty$-category $\mP_\mV(\mM)$ is generated under small colimits and tensors by the essential image of the 
$\mV$-enriched Yoneda-embedding $\iota: \mM \to \mP_\mV(\mM).$
Moreover for every $\X \in \mM$
the $\mV$-enriched functor
$$ \Mor_{\mP_\mV(\mM)}(\Mor_\mM(-,\X),-):\mP_\mV(\mM) \to \mV $$
is equivalent to the $\mV$-enriched functor evaluating at $\X$
and preserves small colimits and tensors.

\end{theorem}

% The next is \cite[Corollary 3.54.]{heine2024bienriched}:

% \begin{proposition}\label{yoadj2}

% Let $\mV $ be a presentably monoidal $\infty$-category, $\mM $ a small $\mV$-enriched $\infty$-category, $\mN$ a left tensored $\infty$-category compatible with small colimits and $\F: \mM \to \mN$ a $\mV$-enriched embedding such that for every $\X \in \mM$
% the $\mV$-enriched functor $\Mor_\mN(\F(\X),-): \mN \to \mV$ preserves small colimits and tensors.
% The left adjoint $\mV$-linear functor $$
% \bar{\F}: \mP_\mV(\mM) \to \mN $$
% extending $\F:\mM \to \mN $ of \cref{envvcor} 
% induces a $\mV$-enriched equivalence to the full $\mV$-enriched subcategory of $\mN$ generated under small colimits and tensors by the essential image of $\F.$

% \end{proposition}

% The next follows from \cite[Corollary 4.71.]{heine2024bienriched}:

% \begin{proposition}
% Let $\mV$ be a presentably monoidal $\infty$-category, $\mM$ a small $\mV$-enriched $\infty$-category, $\mN$ a $\mV$-enriched and left tensored $\infty$-category compatible with small colimits and $\F: \mM \to \mN$ a $\mV$-enriched functor.
% The left adjoint $\mV$-enriched functor $ \bar{\F}: \mP_\mV(\mM) \to \mN $
% extending $\F: \mM \to \mN$ of \cref{envvcor} is $\mV$-enriched left adjoint to the $\mV$-enriched functor $ \mN \subset \mP_\mV(\mN) \xrightarrow{\F^*} \mP_\mV(\mM).$

% \end{proposition}

\begin{notation}\label{delopi}

Let $\mV$ be a monoidal $\infty$-category and $\A$ an associative algebra in $\mV$.
Let $\B \A \subset \RMod_\A(\mV) $ be the full $\mV$-enriched subcategory spanned by $\A.$
    
\end{notation}

\begin{remark}\label{unito}

Let $\mV$ be a monoidal $\infty$-category.
The $\mV$-enriched embedding $\B \A \subset \RMod_\A(\mV) $ induces a $\mV$-enriched functor
$\mP_\mV(\B\A) \to \RMod_\A(\mV)$
using \cref{envvcor}, which is an equivalence by \cite[Corollary 3.56.]{heine2024bienriched}.
     
\end{remark}

\begin{notation}\label{freene}
Let $\mV$ be a presentably monoidal $\infty$-category.
For every small $\infty$-category $\K$ we set
$$ \K_\mV:= \K \boxtimes \B\tu_\mV \in \mV\mathrm{-}\Cat. $$
    
\end{notation}

\begin{proposition}\label{freeenr}

Let $\mV$ be a presentably monoidal $\infty$-category, $\K$ an $\infty$-category and $\mN$ a small $\mV$-enriched $\infty$-category.
The canonical functor $\K \to \K_\mV$
induces an equivalence
$$ \mV\mathrm{-}\Fun(\K_\mV,\mN) \simeq \Fun(\K,\mN).$$

\end{proposition}

Next we define modules.

\begin{notation}

Let $\mV, \mW$ be presentably monoidal $\infty$-categories,
$\A$ an associative algebra in $\mV$ and $\mC$ a $\mV \ot \mW$-enriched
$\infty$-category.
The $\mW^\rev$-enriched $\infty$-category of left $\A$-modules in 
$\mM$ is $$ \L\Mod_\A(\mC) := \mV\mathrm{-}\Fun(\B\A,\mC).$$

\end{notation}

Evaluation at the unique object of $\B\A$, i.e. the
canonical $\mV$-enriched functor $\B\tu \to \B\A$, induces a $\mW^\rev$-enriched functor
$$ \L\Mod_\A(\mC) \to \mC.$$

The latter admits a $\mW^\rev$-enriched left adjoint if
$\mC$ is a left $\mV \ot \mW$-tensored $\infty$-category \cite[Proposition 4.2.4.2.]{lurie.higheralgebra}.

\subsection{Endomorphism actions}\label{endoo}

\begin{definition}Let $\mV$ be a monoidal $\infty$-category, $ \mM $ a weakly $\mV$-enriched $\infty$-category and $\X\in\mM$.
The endomorphism left action on $\X$ if it exists, is the final object of
$ \{\X \} \times_{\mM} \LMod(\mM).$
\end{definition}

Let $\mV$ be a monoidal $\infty$-category, $ \mM $ a weakly $\mV$-enriched $\infty$-category and $\X\in\mM$.
By \cite[Proposition 4.62.]{HEINE2023108941} the endomorphism left action on $\X$ exists if $\X$ admits an endomorphism object.

\vspace{1mm}
Endomorphism actions are functorial in the following way:
a map of weakly enriched $\infty$-categories like (\ref{linmapp})
%\begin{equation*} 
%\begin{xy}
%\xymatrix{
%\mM^\circledast
%\ar[d] 
%\ar[r]^\F  &  \mN^\circledast 
%\ar[d]
%\\
% \mV^\ot\ar[r] & \mW^\ot 
%}
%\end{xy}
%\end{equation*}
induces for every $\X \in \mM $ a commutative square 
\begin{equation*} 
\begin{xy}
\xymatrix{
\{\X \} \times_{\mM } \LMod(\mM)
\ar[d] 
\ar[r]  & \{\F(\X)\} \times_{\mN} \LMod(\mN) 
\ar[d]
\\
\Alg(\mV) \ar[r] & \Alg(\mW). 
}
\end{xy}
\end{equation*}

The endomorphism left actions are by definition the final objects of the top $\infty$-categories of the square.
Consequently, if $\X $ and $\F(\X) $ admit endomorphism objects, $\F$ sends the endomorphism left action on $\X$ to a left action on $\F(\X)$ that is the pullback of the endomorphism left action on $\F(\X) $ along a canonical map of associative algebras in $\mW: $$$\F(\Mor_\mC(\X,\X)) \to \Mor_\mD(\F(\X),\F(\X)).$$

%\subsubsection{Canonical action of endomorphisms on morphisms}
\begin{construction}
For every $\mV$-enriched $\infty$-category $\mM$
the $\mV$-enriched Yoneda-embedding $$\mM \to\mV\mathrm{-}\Fun(\mM^\op, \mV)$$ sends the endomorphism left $\Mor_\mM(\X,\X)$-action on $\X$ to a left $\Mor_\mM(\X,\X)$-action on $\Mor_\mM(-,\X): \mM^\op \to \mV $, which by
\cref{psinho} corresponds to a lift $$ \mM^\op \to \LMod_{\Mor_\mM(\X,\X)}(\mV) $$ of $\Mor_\mM(-,\X).$
So for every object $\Y$ of $\mM$ the morphism object $\Mor_\mM(\Y,\X) \in \mV$ carries a left $\Mor_\mM(\X,\X)$-action and for every morphism $\Y \to \Z $ in $\mM$ the morphism $\Mor_\mM(\Z,\X) \to \Mor_\mM(\Y,\X)$ is $\Mor_\mM(\X,\X)$-linear.
\end{construction}

There is the following functorality:
Let $\F: \mM \to \mN$ be a $\mV$-enriched functor
and $\X,\Y \in \mM.$
\vspace{2mm}
By the $\mV$-enriched Yoneda-lemma (\cref{enryol}) there is a $\mV$-enriched natural transformation $$ \Mor_\mM(\Y,-) \to \Mor_\mN(\F(\Y),-)\circ \F$$
of $\mV$-enriched functors $\mM \to \mV$ that sends the identity of $\Y$ to the identity of $\F(\Y).$
Consequently, the canonical morphism  
$$ \Mor_\mM(\Y,\X) \to \Mor_\mN(\F(\Y),\F(\X))$$
in $\mV$ is $\Mor_\mM(\X,\X)$-linear when the left hand side carries the canonical left $\Mor_\mM(\X,\X)$-action and the right hand side carries the left $\Mor_\mM(\X,\X)$-action that is the pullback of 
the canonical left $\Mor_\mN(\F(\X),\F(\X))$-action along a canonical map of associative algebras $\Mor_\mM(\X,\X)\to \Mor_\mN(\F(\X),\F(\X))$ in $\mV.$

\section{Monads and monadicity}

\subsection{$(\infty,2)$-categories}\label{2ca}

In this section we introduce $(\infty,2)$-categories, the natural context
to consider monads and monadic morphisms.

% \begin{definition}
% %We define the large $\infty$-category of small $(\infty,2)$-categories as the
% %$%\infty$-category $(\infty,2)\Cat$ of small $\infty$-categories enriched in the cartesian structure of $\infty\Cat.$	
% An $(\infty,2)$-category is an $\infty$-category enriched in the cartesian symmetric monoidal structure on $\infty\Cat.$ A 2-functor is a $\infty\Cat$-enriched functor.

% We set $(\infty,2)\Cat:=(\infty,2)\Cat.$
% \end{definition}

\begin{definition}
The $\infty$-category of small $(\infty,2)$-categories is $(\infty,2)\Cat:=\infty\Cat\mathrm{-}\Cat.$
\end{definition}

%an embedding is a $\infty\Cat$-enriched embedding, a 2-left adjoint
%is a $\infty\Cat$-enriched left adjoint. 
	
\begin{notation}
Let $\mC, \mD$ be $(\infty,2)$-categories.
Let $$\FUN(\mC, \mD) :=\infty\Cat\mathrm{-}\Fun(\mC,\mD).$$	 
	 
\end{notation}

\begin{example}

The canonical example of an $(\infty,2)$-category is $\infty\Cat$ with its natural enrichment in itself.
\end{example}

We call the $\infty\Cat$-enriched Yoneda-embedding the 2-Yoneda-embedding.

% \begin{notation}
% For every $(\infty,2)$-category $\mC$
% we call the $\infty\Cat$-enriched Yoneda-embedding $ \iota: \mC \hookrightarrow \FUN(\mC^\op, \infty\Cat)$
% of Definition \ref{nnoo} the 2-Yoneda-embedding.
% \end{notation}

\begin{notation}\label{oppost}

We define involutions
$$ (-)^\op, (-)^\co: (\infty,2)\Cat \to (\infty,2)\Cat,$$
where $(-)^\op$ is the opposite $\infty\Cat$-enriched $\infty$-category,
and $(-)^\co$ is the transfer of enrichment along the involution
$(-)^\op: \infty\Cat \to \infty\Cat.$

%:= ((-)^\op)_! : (\infty,2)\Cat=\Cat^{\infty\Cat} \to (\infty,2)\Cat=\Cat^{\infty\Cat}.$

\end{notation}

\begin{remark}
The opposite $\infty$-category involution
$$ (-)^\op: \infty\Cat \to \infty\Cat$$
is symmetric monoidal with respect to cartesian structures and so 
induces an equivalence of $(\infty,2)$-categories
$$ (-)^\op: \infty\Cat^\co \to \infty\Cat. $$

\end{remark}

% \begin{notation}\label{oppost}

% Let $\mC$ be an $(\infty,2)$-category.

% \begin{itemize}
% \item We write $(\mC^\mathrm{co})^\circledast \to \infty\Cat^\times$
% for the pullback of $\mC^\circledast \to \infty\Cat^\times$
% along the opposite $\infty$-category involution 
% $ (-)^\op : \infty\Cat^\times \to \infty\Cat^\times$.

% \item We write $(\mC^\mathrm{op})^\circledast \to \infty\Cat^\times$
% for the opposite $\infty\Cat$-enriched $\infty$-category.

% \end{itemize}
% \end{notation}

%\begin{definition}Let $\rS$   an $\infty$-category. A (locally) cocartesian $\rS$-family of$(\infty,2)$-categories is a (locally) cocartesian $\rS$-family of $\infty\Cat$-enriched $\infty$-categories.
%We call a (locally) cocartesian $\rS$-family of categories enriched in$\infty\Cat^\times $ a (locally) cocartesian $\rS$-family of 2-categories and a map of (locally) cocartesian $\rS$-families of $\infty\Cat^\times $-enriched categories a map of (locally) cocartesian $\rS$-families of 2-categories.\end{definition}

\begin{definition}

Let $\rS$ be an $\infty$-category. The canonical left action of $\infty\Cat$ on $\infty\Cat_{/ \rS}$ is the left action induced (Definition \ref{ind}) by the right adjoint functor $(-)\times\rS:\infty\Cat \to \infty\Cat_{/ \rS}$.

%This left action is provided by the left tensored $\infty$-category $ \infty\Cat^\times \times_{(\infty\Cat_{/ \rS})^\times} (\infty\Cat_{/ \rS})^\circledast \to \infty\Cat^\times.$

\end{definition}

\begin{example}Let $\rS$ be an $\infty$-category. The canonical left action of $\infty\Cat$ on $\infty\Cat_{/ \rS}$ is closed since for every $\infty$-category $\mB$ and functors $\mC \to \rS, \mD \to \rS$ there is a canonical equivalence
\begin{equation}\label{clokl}
\infty\Cat(\mB, \Fun_\rS(\mC,\mD)) \simeq \infty\Cat_{/ \rS}(\mB \times \mC,\mD).
\end{equation}
Thus $\infty\Cat_{/ \rS}$ is an $(\infty,2)$-category.
	
\end{example}

\begin{notation}\label{notabes}
Let $\infty\Cat_{/ \rS}^{\mE} \subset \infty\Cat_{/ \rS}$ be the subcategory of cocartesian fibrations relative to $\mE$ and maps of locally cocartesian fibrations relative to $\mE.$
Let $\infty\Cat_{/ \rS}^{\mE, \mT} \subset \infty\Cat_{/ \rS}$ be the subcategory of cocartesian fibrations relative to $\mE, \mT$ and maps of locally cocartesian fibrations relative to $\mE.$

\end{notation}

\begin{example}\label{Imporso}
	
Let $\rS$ be an $\infty$-category and $\mE \subset \Fun([1],\rS)$ a full subcategory. The canonical left $\infty\Cat$-action on $\infty\Cat_{/ \rS}$
restricts to a left $\infty\Cat$-action on $\infty\Cat^\mE_{/ \rS}$
(Definition \ref{notabes}).
For any cocartesian fibrations $\mC \to \rS, \mD \to \rS$ relative to $\mE$
let $\Fun^\mE_\rS(\mC,\mD) \subset \Fun_\rS(\mC,\mD)$ be the full subcategory of functors over $\rS$ preserving cocartesian lifts of morphisms of $\mE.$
For every $\infty$-category $\K$ equivalence (\ref{clokl}) restricts to an equivalence $$\Fun(\K, \Fun_\rS^\mE(\mC,\mD)) \simeq \Fun^\mE_\rS(\K \times \mC, \mD).$$
Hence the left $\infty\Cat$-action on $\infty\Cat^\mE_{/ \rS}$ is closed
so that $\infty\Cat^\mE_{/ \rS}$ is an $(\infty,2)$-category.
In particular, the functor $\Fun^\mE_\rS(\mC,-): \infty\Cat^\mE_{/ \rS} \to \infty\Cat$ refines to a 2-functor.

\end{example}

\subsection{Monads}\label{endomol}

In this section we define monads in every $(\infty,2)$-category
and associate a monad to any right adjoint morphism (Proposition \ref{rightadj}). Monads in $(\infty,2)$-categories were studied by \cite{RIEHL2016802}, \cite{Haugseng2020OnLT}, \cite[7.2.]{stefanich2021higher}.

\begin{definition}\label{moNN}
Let $\mC$ be an $(\infty,2)$-category,  $\X \in \mC.$
A monad on $\X$ in $\mC$ is an associative algebra in
the endomorphism monoidal structure on $\Mor_\mC(\X,\X).$
A comonad on $\X$ in $\mC$ is a monad on $\X$ in $\mC^\mathrm{co}.$
	
\end{definition}

\begin{example}\label{exck}
	
Let $\mC$ be an $(\infty,2)$-category,
$\mV$ a monoidal $\infty$-category and $\X$ a left $\mV$-module
in $\mC$ corresponding to a monoidal functor $\mV \to \Mor_\mC(\X,\X).$
Every associative algebra in $\mV$ gives rise to a monad on $\X$ in $\mC.$

\end{example}

\begin{remark}
Definition \ref{moNN} for $\mC=\infty\Cat$ is Lurie's definition of monad
\cite[Definition 4.7.0.1.]{lurie.higheralgebra}
and was generalized by \cite{Haugseng2020OnLT} to define monads in double
$\infty$-categories.
There are alternative definitions of a monad in an $(\infty,2)$-category:
Riehl-Verity \cite{RIEHL2016802} define monads in an $(\infty,2)$-category $\mC$ by 2-functors
from the walking monad, a 2-category encoding the shape of a monad, to $\mC$.
Zaganidis \cite{zaganidis2017towards} constructs a model for the $(\infty,2)$-category of monads for this definition.
Haugseng \cite[Theorem 1.4.]{Haugseng2020OnLT} identifies all these definitions.
\end{remark}

Next we assign a monad to any right adjoint morphism in any $(\infty,2)$-category.

%There is the notion of adjunction in any $(\infty,2)$-category:

\begin{definition}Let $\mC $ be an $(\infty,2)$-category and $\F: \X \to \Y, \G: \Y \to \X$ morphisms of $\mC.$ 
We say that $\F $ is left adjoint to $\G$ (or $\G$ is right adjoint to $\F$) 
if there are morphisms $ \eta: \id_\X \to \G \circ \F $ in $\Mor_\mC(\X,\X)$
and $ \varepsilon : \F \circ \G \to \id_\Y $ in $\Mor_\mC(\Y,\Y)$
such that the following triangle identities hold:
$$ (\varepsilon \circ \F) \circ (\F \circ \eta) = \id_\F, \ (\G \circ \varepsilon) \circ (\eta \circ \G ) = \id_\G.$$
\end{definition}

\begin{definition}
Let $\mC$ be an $(\infty,2)$-category, $\T$ a monad on some object $\X$ of $\mC$ and $\G: \Y \to \X$ a right adjoint morphism
of $\mC$ equipped with a left $\T$-action. The left action map $\mu: \T \circ \G \to \G$ in $\Mor_\mC(\Y,\X)$ exhibits $\T$ as the monad associated to $\G$ if $\mu$ exhibits $\T$ as the endomorphism object of $\G$ 
with respect to the canonical left $\Mor_\mC(\X,\X)$-action on $\Mor_\mC(\Y,\X).$

\end{definition}

\begin{lemma}\label{trwgqaxy}

Let $\mC $ be an $(\infty,2)$-category and $\G: \Y \to \X $ a morphism of $\mC$ that admits a left adjoint $ \F: \X \to \Y $. 
Let $\eta: \id_\X \to \G \circ \F $ be the unit and $\varepsilon: \F \circ \G \to \id_\Y $ the counit.
% of the adjunction. 

\begin{enumerate}

\item For every morphism $\h: \X \to \X $ of $\mC $ the following two maps are inverse to each other:
$$ \hspace{7mm}\alpha: \Mor_\mC(\X,\X)(\h, \G \circ \F) \to \Mor_\mC(\Y,\X) (\h \circ \G, \G \circ \F \circ \G ) \xrightarrow{\Mor_\mC(\Y,\X) (\h \circ \G, \G \circ \varepsilon ) } \Mor_\mC(\Y,\X) (\h \circ \G, \G) $$
$$\hspace{7mm} \beta: \Mor_\mC(\Y,\X) (\h \circ \G, \G) \to \Mor_\mC(\X,\X) (\h \circ \G \circ \F, \G \circ \F) \xrightarrow{\Mor_\mC(\X,\X) (\h \circ \eta, \G \circ \F) } \Mor_\mC(\X,\X)(\h, \G \circ \F).$$

\item Let $\T: \X \to \X $ be a morphism of $\mC$ and $\varphi : \T \circ \G \to \G $ a morphism in $  \Mor_\mC(\Y,\X). $

Let $ \psi $ be the composition $  \T \xrightarrow{ \T \circ \eta} \T \circ \G \circ \F \xrightarrow{ \varphi \circ \F } \G \circ \F $
in $\Mor_\mC(\X,\X) $.

Then $ \varphi$ factors as $ \T \circ \G \xrightarrow{\psi \circ \G } \G \circ \F \circ \G \xrightarrow{\G \circ \epsilon }  \G. $

\vspace{1mm}

So for every morphism $\h: \X \to \X  $ of $\mC$
the map $ \Gamma: $$$ \hspace{8mm} \Mor_\mC(\X,\X) (\h, \T) \to \Mor_\mC(\Y,\X) (\h \circ \G, \T \circ \G ) \xrightarrow{\Mor_\mC(\Y,\X)(\h \circ \G, \varphi ) } \Mor_\mC(\Y,\X) (\h \circ \G, \G)  $$ factors as 
$$ \hspace{15mm} \Mor_\mC(\X,\X) (\h, \T) \xrightarrow{ \Mor_\mC(\X,\X)(\h, \psi)} \Mor_\mC(\X,\X) (\h, \G \circ \F ) \xrightarrow{\alpha} \Mor_\mC(\Y,\X) (\h \circ \G, \G).  $$
So $\psi$ is an equivalence if and only if for any $\h: \X \to \X  $ of $\mC$ the map $ \Gamma $ is an equivalence.

\vspace{2mm}

\item Let $\G: \Y \to \X, \h: \Z \to \X $ be morphisms of $\mC$ that admit left adjoints $\F:\X \to \Y$, $\bk: \X \to \Z $, respectively, and let $\phi: \Y \to \Z$ be a morphism in $\mC$ over $\X.$ 
Let $ \omega$ be the morphism $$\h \circ \bk \to \h \circ \bk \circ \G \circ \F \simeq \h \circ \bk \circ \h \circ \phi \circ \F \to \h \circ \phi \circ \F \simeq \G \circ \F   $$ in $\Mor_\mC(\X,\X). $ 
Then $\h \circ \bk \circ \G \xrightarrow{\omega \circ \G} \G \circ \F \circ \G \to \G $ factors as $$\h \circ \bk \circ \G \simeq \h \circ \bk \circ \h \circ \phi \to \h \circ \phi \simeq \G. $$

\end{enumerate}

\end{lemma}

\begin{proof}

(1): The composition $\beta \circ \alpha $ factors as
$$ \Mor_\mC(\X,\X) (\h, \G \circ \F ) \xrightarrow{\Mor_\mC(\X,\X) (\h, \G \circ \F  \circ \eta )  } \Mor_\mC(\X,\X) (\h, \G \circ \F \circ \G \circ \F)$$$$ \xrightarrow{\Mor_\mC(\X,\X) (\h, \G \circ \varepsilon \circ \F)}
\Mor_\mC(\X,\X) (\h, \G \circ \F) $$ and $\alpha \circ \beta$ factors as
$$ \Mor_\mC(\Y,\X) (\h \circ \G, \G) \xrightarrow{\Mor_\mC(\Y,\X) (\h \circ \G \circ \varepsilon, \G )  } \Mor_\mC(\Y,\X) (\h \circ \G \circ \F \circ \G, \G )$$$$ \xrightarrow{\Mor_\mC(\Y,\X) (\h \circ \eta \circ \G, \G)} \Mor_\mC(\Y,\X) (\h \circ \G, \G).$$ 

Therefore statement (1) follows from the triangle identities.

%The compositions $  \F \xrightarrow{ \F \circ \eta} \F \circ \G \circ \F \xrightarrow{ \varepsilon \circ \F } \F  $ and$  \G \xrightarrow{ \eta \circ \G } \G \circ \F \circ \G \xrightarrow{ \G \circ \varepsilon } \G  $of morphisms in $\Mor_\mC(\X, \Y) , \Mor_\mC(\Y,\X) $ are the identities,respectively.

(2): The composition $ \psi: \T \circ \G \xrightarrow{ \T \circ \eta \circ \G} \T \circ \G \circ \F \circ \G \xrightarrow{ \varphi \circ \F \circ \G } \G \circ \F \circ \G  \xrightarrow{\G \circ \epsilon }  \G $ factors as
$$ \T \circ \G \xrightarrow{ \T \circ \eta \circ \G } \T \circ \G \circ \F \circ \G \xrightarrow{ \T \circ \G \circ \epsilon } \T \circ \G  \xrightarrow{\varphi } \G $$ and is thus equivalent to $\varphi$ by the triangle identities.

(3): The composition $$\h \circ \bk  \circ \G \to \h \circ \bk \circ \G \circ \F  \circ \G \simeq \h \circ \bk \circ \h \circ \phi \circ \F  \circ \G \to \h \circ \phi \circ \F  \circ \G \simeq \G \circ \F   \circ \G \to \G  $$ 
factors as
$\h \circ \bk  \circ \G \to \h \circ \bk \circ \G \circ \F  \circ \G \simeq \h \circ \bk \circ \h \circ \phi \circ \F  \circ \G \to \h \circ \phi \circ \F  \circ \G \to \h \circ \phi\simeq \G  $ and so factors as
$$\h \circ \bk  \circ \G \to \h \circ \bk \circ \G \circ \F  \circ \G \to
\h \circ \bk  \circ \G \simeq \h \circ \bk \circ \h \circ \phi \to \h \circ \phi\simeq \G,  $$ which is equivalent to 
$\h \circ \bk \circ \G \simeq \h \circ \bk \circ \h \circ \phi \to \h \circ \phi \simeq \G $ by the triangle identities. 	

\end{proof}

Lemma \ref{trwgqaxy} implies the following proposition:

\begin{proposition}\label{rightadj}

Let $\mC$ be an $(\infty,2)$-category and $\G: \Y \to \X $ a morphism of $\mC$ that admits a left adjoint $ \F: \X \to \Y $.
Let $\eta: \id_\X \to \G \circ \F $ be the unit and $\varepsilon: \F \circ \G \to \id_\Y $ the counit.

% of the adjunction. \vspace{2mm}

\begin{enumerate}
\item For every morphism $\h: \X \to \X $ of $\mC $ the following map $\alpha$ is an equivalence:
$$\hspace{10mm} \Mor_\mC(\X,\X)(\h, \G \circ \F) \to \Mor_\mC(\Y,\X) (\h \circ \G, \G \circ \F \circ \G ) \xrightarrow{\Mor_\mC(\Y,\X) (\h \circ \G, \G \circ \varepsilon)} \Mor_\mC(\Y,\X) (\h \circ \G, \G).$$

Hence $\G \circ \varepsilon : \G \circ \F \circ \G \to \G $ exhibits $\G \circ\F $ as the endomorphism object of $\G: \Y \to \X $ with respect to the canonical
left $\Mor_\mC(\X,\X)$-action on $\Mor_\mC(\Y,\X).$

\vspace{1mm}

\item Let $\T: \X \to \X $ be a morphism of $\mC$ and $\varphi : \T \circ \G \to \G $ a morphism in $  \Mor_\mC(\Y,\X). $

Let $ \psi $ be the composition $ \T \xrightarrow{ \T \circ \eta} \T \circ \G \circ \F \xrightarrow{ \varphi \circ \F} \G \circ \F $ in $\Mor_\mC(\X,\X) $ and $\Gamma$ the composition 	
$$\hspace{17mm} \Mor_\mC(\X,\X) (\h, \T) \to \Mor_\mC(\Y,\X) (\h \circ \G, \T \circ \G ) \xrightarrow{\Mor_\mC(\Y,\X)(\h \circ \G, \varphi )} \Mor_\mC(\Y,\X) (\h \circ \G, \G).  $$ 

Then $\psi$ is an equivalence if and only if for every morphism $\h: \X \to \X  $ of $\mC$ the map $ \Gamma $ is an equivalence. 
Thus $\varphi : \T \circ \G \to \G $ exhibits $\T$ as the endomorphism object of $\G: \Y \to \X $ with respect to the left $\Mor_\mC(\X,\X)$-action on $\Mor_\mC(\Y,\X)$ if and only if $\psi$ is an equivalence. 

\end{enumerate}	

\end{proposition}

\begin{corollary}\label{rightadj2}
Let $\mC $ be an $(\infty,2)$-category and $\G: \Y \to \X $ a morphism of $\mC$ that admits a left adjoint $ \F: \X \to \Y $.
The counit $\G \circ \F\circ \G \to \G \circ \id = \G$ exhibits $\G \circ \F$ as an endomorphism object of $\G$ with respect to the canonical left $\Mor_\mC(\X,\X)$-action on $\Mor_\mC(\Y,\X)$.
	
\end{corollary}

\begin{remark}

After the first version of this paper appeared, \cref{rightadj} was also proven by \cite[Proposition 7.2.4.]{stefanich2021higher}.
    
\end{remark}

We obtain the following corollary:

\begin{corollary}\label{sdfghjxcvbnggggg} 

Let $\F: \mC \to \mD$ be a 2-functor, $\T$ a monad in $\mC$
on an object $\X$ of $\mC$ and $\phi: \Y \to \X $ a morphism of $\mC$ that carries a left $\T$-action, which is sent by $\F$ to a left $\F(\T)$-action on $\F(\phi): \F(\Y ) \to \F(\X) $.
If $\phi $ carries the endomorphism left $\T$-action, then $\F(\phi)$ carries the endomorphism left $ \F(\T) $-action.

\end{corollary}

\begin{remark}\label{funcmon}

Let $\mC$ be an $(\infty,2)$-category and
\begin{equation}\label{pli}
\begin{xy}
\xymatrix{
\Y \ar[rd]_{\g'} \ar[rr]^{\h}
&& \Y' \ar[ld]^{\g} 
\\ & \X
}
\end{xy} 
\end{equation}
a triangle in $\mC$ such that $\g$ admits a left adjoint $\f$ and
$\g'$ admits a left adjoint $\f'$.
Let $\T$ be the monad associated to $\g$ and $\T'$ the monad associated to $\g'.$

The functor $\Mor_\mC(\h,\X): \Mor_\mC(\Y',\X) \to \Mor_\mC(\Y,\X)$
is $\Mor_\mC(\X,\X)$-linear and so sends the endomorphism left $\T$-action 
on $\g$ to a left $\T$-action on $\g'$, which corresponds to a map $\T \to \T'$
of monads on $\X $ by the universal property of the endomorphism algebra.

\end{remark}

\subsection{Eilenberg-Moore objects}\label{mondl}

In this section we develop the theory of Eilenberg-Moore objects in every $(\infty,2)$-category, which for $\infty\Cat$ specialize to the
$\infty$-category of algebras of a monad. 
%We prove a monadicity theorem in any $(\infty,2)$-category that admits Eilenberg-Moore objects (\cref{monaT}).
Eilenberg-Moore objects were also studied by 
\cite[7.3.]{stefanich2021higher}.

\begin{construction}\label{endoi}
	
Let $\mC$ be an $(\infty, 2)$-category. The 2-Yoneda-embedding
$ \iota: \mC \hookrightarrow \FUN(\mC^\op, \infty\Cat) $
sends the endomorphism left $\Mor_{\mC}(\X,\X) $-action on any object $\X \in \mC$ to a left $\Mor_{\mC}(\X,\X) $-action on $\iota(\X)$,
which by \cref{psinho} corresponds to a 2-functor $$\bar{\iota}(\X): \mC^\op \to \Mor_{\mC}(\X,\X)\mathrm{-}\LMod $$ lifting $\iota(\X): \mC^\op \to \infty\Cat$ along the forgetful 2-functor $\nu: \Mor_{\mC}(\X,\X)\mathrm{-}\LMod \to \infty\Cat$.
\end{construction}
Let $\T$ be a monad on $\X \in \mC.$ There is a 2-functor $\LMod_\T: \Mor_{\mC}(\X,\X)\mathrm{-}\LMod \to \infty\Cat$
and a 2-transformation $\LMod_\T \to \nu$.

\begin{definition}\label{monadi} 
Let $\mC$ be an $(\infty, 2)$-category and $\T$ a monad on some $\X \in \mC.$
A morphism $\phi: \Y \to \X$ of $\mC$ is an Eilenberg-Moore object of $\T $ if there is a commutative triangle in $\FUN(\mC^\op, \infty\Cat):$
\begin{equation*}
\begin{xy}
\xymatrix{
\iota(\Y) \ar[rd]^{\iota(\phi)}
\ar[rr]^{\simeq} 
&&\LMod_\T \circ \bar{\iota}(\X) \ar[ld]
\\
& \iota(\X).}
\end{xy} 
\end{equation*}  

\end{definition}

\begin{remark}\label{charos}
Let $\mC$ be an $(\infty,2)$-category, $\phi: \Y \to \X$ a morphism in $\mC$ and $\T$ a monad on $\X$ in $\mC.$
By the 2-Yoneda-lemma a left $\T$-action on $\phi$ corresponds to a map 
$\theta: \iota(\Y) \to \LMod_\T \circ \bar{\iota}(\X)$ in $\FUN(\mC^\op, \infty\Cat)$ over $\iota(\X).$ 
The morphism $\phi$ is an an Eilenberg-Moore object of $\T$
if and only if for every $\Z \in \mC$ the functor $$\theta_\Z: \Mor_\mC(\Z,\Y) \to \LMod_\T(\Mor_\mC(\Z,\X)) $$ is an equivalence.

\end{remark}

\begin{remark}\label{funcem}

Let $\mC$ be an $(\infty,2)$-category, $\X \in \mC$ and $\theta: \T \to \T'$ a map
of monads on $\X $ in $\mC.$
Let $\Y \to \X$ an Eilenberg-Moore object of $\T$ in $\mC$ and
$\Y' \to \X$ an Eilenberg-Moore object of $\T'$ in $\mC$.
The restriction of the left $\T'$-action on $\Y' \to \X$ along $\theta$
is a left $\T$-action on $\Y' \to \X$, which by universal property of the
Eilenberg-Moore object $\Y \to \X$ corresponds to a morphism $\Y' \to \Y$ over $\X.$

\end{remark}

\begin{lemma}\label{rembras0}
Let $\G: \mD \to \mC$ be a fully faithful 2-functor that admits a 2-left adjoint. Every Eilenberg-Moore object of $\mC$
whose target belongs to $\mD$, is an Eilenberg-Moore object of $\mD$.
\end{lemma}

\begin{proof}

Let $\F$ be a left adjoint of $\G$ and $\G$ is a fully faithful 2-functor.
Let $\X \in \mD$ and $\psi: \Y \to \G(\X)$ a right adjoint morphism of $\mC.$ 
We like to see that for every morphism $\A \to \B$ of $\mC$ that is inverted by
$\F$ the induced functor
$\Mor_\mC(\B, \Y) \to \Mor_\mC(\A, \Y)$
is an equivalence.
This morphism identifies with the functor
$$\LMod_\T(\Mor_\mC(\B, \G(\X))) \to \LMod_\T(\Mor_\mC(\A, \G(\X)))$$
induced by the functor $\Mor_\mC(\B, \G(\X)) \to \Mor_\mC(\A, \G(\X)),$
which identifies with the equivalence 
$$\Mor_\mD(\F(\B), \X) \to \Mor_\mD(\F(\A), \X).$$
\end{proof}

\begin{definition}\label{monadi} 
Let $\mC$ be an $(\infty, 2)$-category.
A morphism $ \Y \to \X$ of $\mC$ is monadic if it is an
Eilenberg-Moore object of some monad $\T$ on $\X.$
	
\end{definition}

There are three kind of opposite $(\infty,2)$-categories (Notation \ref{oppost}) giving the following three notions:

\begin{notation}
Let $\mC$ be an $(\infty, 2)$-category.

\begin{itemize}
\item A coEilenberg-Moore object in $\mC$ for a comonad $\R$ on $\X$ is an Eilenberg-Moore object $\Y \to \X$ for the monad $\R$ in $\mC^\mathrm{co}.$

\item A (co)Kleisli object in $\mC$ for a (co)monad $\T$ on $\X$ is a (co)Eilenberg-Moore object $\X \to \Y$ for the (co)monad $\T$ in $\mC^\mathrm{op}.$

\end{itemize}	
	
\end{notation}

%\begin{lemma}Let $\G: \mD^\circledast \to \mC^\circledast$ be an embedding that admits a 2-left adjoint and $\psi: \Y \to \X$ a  monadic morphism.Then $\Y$ belongs to the essential image of $\G$ if $\X$ does.\end{lemma}

%\begin{proof}An object $\Y \in \mC$ belongs to the essential image of$\G$ if and only if for any morphism $\Z \to \Z'$ of $\mC$ inverted bythe left adjoint of $\G$ the induced functor $\Mor_\mC(\Z',\Y) \to \Mor_\mC(\Z,\Y)$ is an equivalence.

%Let $\T$ be the monad associated to $\psi: \Y \to \X$.Since $\X $ belongs to the essential image of $\G$,the induced functor $\Mor_\mC(\Z',\X) \to \Mor_\mC(\Z,\X)$,which is $\Mor_\mC(\X,\X)$-linear, is an equivalence and so induces an equivalence$\LMod_\T(\Mor_\mC(\Z',\X)) \to \LMod_\T(\Mor_\mC(\Z,\X))$.
%This equivalence identifies with the functor$\Mor_\mC(\Z',\Y) \to \Mor_\mC(\Z,\Y).$\end{proof}

We will prove the following:

\begin{proposition}\label{monaT}
Let $\mC$ be an $(\infty,2)$-category.
A morphism $\phi: \Y \to \X$ in $\mC$ is monadic if and only if the following conditions hold:
\begin{enumerate}
\item The morphism $\phi: \Y \to \X$ admits a left adjoint.
	
\item For every $\Z \in \mC$ the induced functor $\phi_*: \Mor_\mC(\Z,\Y) \to \Mor_\mC(\Z,\X)$ is monadic.
		
\end{enumerate}
	
\end{proposition}

We prepare the proof of \cref{monaT}.
Let $\mB,\mC$ be $\infty$-categories and $\T$ a monad on $\mC$. By Corollary \cref{2-cat2} there is a canonical equivalence over $\Fun(\mB,\mC)$: $$ \Fun(\mB, \LMod_\T(\mC)) \simeq \LMod_{\T}(\Fun(\mB,\mC)), $$
where the right hand side are modules with respect to the left action of $\Fun(\mC,\mC)$ on $\Fun(\mB,\mC).$

Consequently, for every functor $\G:\mD \to \mC$ a left $\T$-action on $\G$ corresponds to a lift $\bar{\G}: \mD \to \LMod_\T(\mC)$ of $\G$.
Proposition \ref{rightadj} gives the following corollary:
\begin{corollary}\label{Hgb} Let $\G:\mD \to \mC$ be functor that carries a left $\T$-action corresponding to a lift $\bar{\G}: \mD \to \LMod_\T(\mC)$ of $\G$.
If $\G$ has a left adjoint $\F: \mC \to \mD, $ 
the left $\T$-action on $\G$ is the endomorphism left $\T$-action if and only if $\bar{\G}$ preserves the left adjoints, i.e. the canonical map $ \T \circ \id \to \T \circ \G \circ \F \to \bar{\G} \circ \F $ is an equivalence.
	
\end{corollary}

% \begin{remark}\label{inpos}

% Let $\mC$ be an $(\infty,2)$-category and $\G:\Y \to \X$ a morphism in $\mC$ that admits a left adjoint $\F.$ Let $\T $ be a monad on $\X$ that acts on $\G$. Let $\mu:\T \circ \G \to \G$ be the action morphism in $ \Mor_\mC(\Y,\X).$

% By \cref{endoi} the functor $\Mor_\mC(\F,\X): \Mor_\mC(\Y,\X) \to \Mor_\mC(\X,\X)$ is left $\Mor_\mC(\X,\X)$-linear and so sends the left
% $\T$-action on $\G$ to a left $\T$-action on $\G\circ \F$ whose action morphism is $\mu \circ \F: \T \circ \G\circ \F \to \G\circ \F.$
% By \cref{rightadj} (2) the left $\T$-action on $\G$ is the endomorphism $\T$-action if and only if $\T \to \T \circ \G\circ \F \xrightarrow{\mu\circ\F} \G\circ \F$  is an equivalence, where the first morphism is induced by the unit $\id \to \G \circ \F.$
% This is equivalent to say that the unit $\id \to \G \circ \F$ exhibits
% $\G\circ\F$ as the free left $\T$-module generated by $\id: \X \to \X,$
% which is $\T.$
% In other words, a left $\T$-action on $\G:\Y \to \X$ is the endomorphism action if and only if it is sent by $\Mor_\mC(\F,\X): \Mor_\mC(\Y,\X) \to \Mor_\mC(\X,\X)$ to $\T$ viewed as module over itself.
    
% \end{remark}

% \begin{remark}

% After the first version of this paper appeared, \cref{monaT} was also proven by \cite[Theorem 2.5.2.]{stefanich2021higher}

% \end{remark}

Let $\G:\mD \to \mC$ be functor that admits a left adjoint $\F: \mC \to \mD.$ 
By Corollary \ref{rightadj2} there is an endomorphism left $\T$-action on $\G$
for some monad $\T$ on $\mC$ such that $\T \simeq \G \circ \F.$
So by Corollary \ref{Hgb} the corresponding lift $\bar{\G}: \mD \to \LMod_\T(\mC)$ of $\G$ preserves the left adjoints.

\begin{definition}
A functor $\G:\mD \to \mC$ is monadic if $\G$ admits a left adjoint $\F: \mC \to \mD$ and $\bar{\G}: \mD \to \LMod_\T(\mC)$ is an equivalence, where $\T$ is the endomorphism algebra of $\G.$
	
\end{definition}

\begin{remark}\label{monchar}
Let $\G:\mD \to \mC$ be functor that admits a left adjoint, $\T$ the endomorphism algebra of $\G$ acting on $\G$ via the endomorphism left action
and $\rH: \mD \to \LMod_\T(\mC)$ a lift of $\G$ that preserves the left adjoints.	
Then $\G:\mD \to \mC$ is monadic if and only if $\rH$ is an equivalence. 

\end{remark}

\begin{proof}
	
%Let $\G:\mD \to \mC$ be functor that admits a left adjoint $\F$ and$\T$ the endomorphism algebra of $\G$ (Corollary \ref{rightadj2}) acting on $\G$ via the endomorphism left action. 
%Let $\bar{\G}: \mD \to \LMod_\T(\mC)$ be the lift of $\G$ corresponding to the endomorphism left $\T$-action.

The lift $\rH: \mD \to \LMod_\T(\mC)$ of $\G$ corresponds by Corollary \ref{Hgb} to the endomorphism left $\T$-action on $\G.$
By the uniqueness of endomorphism left $\T$-actions on $\G$ there is an autoequivalence $\alpha$ of $\T$ in $\Alg(\Fun(\mC,\mC))$ such that the one endomorphism left $\T$-action on $\G$ is the pullback along $\alpha$ of the other. %endomorphism left $\T$-action on $\G$. %that identifies both endomorphism left $\T$-actions on $\G.$
The equivalence $\alpha$ gives rise to an autoequivalence $\alpha^*$ of $\LMod_\T(\mC)$ such that $\rH \simeq \alpha^* \circ \bar{\G}.$

\end{proof}

\begin{lemma}\label{Examp}
Let $\T$ be a monad on some $\infty$-category $\mC.$ The functor $\LMod_\T(\mC) \to \mC$ is monadic.	
	
\end{lemma}

\begin{proof}
The identity of $\LMod_\T(\mC)$ corresponds to a left $\T$-action on the functor $\LMod_\T(\mC) \to \mC$, which is the endomorphism left action since the identity preserves the left adjoints (Corollary \ref{Hgb}). 

%We apply Lemma \ref{monchar}.	
	
\end{proof}

\begin{remark}\label{sub}

Let $\mC$ be an $(\infty,2)$-category and $\mD \subset \mC$ a locally full subcategory, $\X \in \mD$ and $\T$ a monad on $\X$ in $\mD$. Let $\G: \Y \to \X$ be a morphism in $\mD$ that is an Eilenberg-Moore object for $\T$ in $\mC$.
Then $\G: \Y \to \X$ is an Eilenberg-Moore object for $\T$ in $\mD$
if a morphism $\Z \to \Y$ in $\mC$ lies in $\mD$ if the composition
$\Z \to \Y \xrightarrow{\G} \X$ lies in $\mD$.

\end{remark}

\begin{notation}
Let $\mK \subset \infty\Cat$ be a full subcategory. Let
$$\infty\Cat(\mK) \subset \infty\Cat$$ be the locally full subcategory 
of the $(\infty,2)$-category $\infty\Cat$ of small $\infty$-categories that admit $\mK$-indexed colimits and functors preserving $\mK$-indexed colimits.

\end{notation}

\begin{example}Let $\mK \subset \infty\Cat$ be a full subcategory,
$\mC$ an $\infty$-category that admits $\mK$-indexed colimits and 
$\T$ a monad on $\mC$ that preserves $\mK$-indexed colimits.
By \cite[Corollary 4.2.3.5.]{lurie.higheralgebra} the $\infty$-category $\LMod_\T(\mC)$ admits $\mK$-indexed colimits 
and the functor $\LMod_\T(\mC) \to \mC$, which is monadic by \cref{Examp},
preserves $\mK$-indexed colimits.
Since the functor $\LMod_\T(\mC) \to \mC$ is conservative, it detects $\mK$-indexed colimits. This implies by \cref{sub} that 
$\LMod_\T(\mC) \to \mC$ is a monadic morphism in
the $(\infty,2)$-category $\infty\Cat(\mK).$

\end{example}

\begin{lemma}\label{eilmu}
Let $\mC$ be an $(\infty,2)$-category and $\G: \Y \to \X$ a
morphism in $\mC$ that admits a left adjoint $\F$ and $\T$ a monad on $\X$ that acts on $\G.$ The following conditions are equivalent:

\begin{enumerate}

\item For every $\Z \in \mC$ the induced functor $$ \Mor_\mC(\Z,\Y) \to \LMod_\T(\Mor_\mC(\Z,\X)) $$
over $\Mor_\mC(\Z,\X)$ preserves the free functors.
    
\item The induced functor $$ \Mor_\mC(\X,\Y) \to \LMod_\T(\Mor_\mC(\X,\X)) $$
over $\Mor_\mC(\X,\X)$ preserves the free object on the identity.
    
\item The left $\T$-action on $\G$ is the endomorphism action.

\end{enumerate}

\end{lemma}

\begin{proof}

Since $\F: \X \to \Y $ is right adjoint to $\G: \Y \to \X,$ the functor
$\Mor_\mC(\Z,\F)$ is left adjoint to $\Mor_\mC(\Z,\G)$. 
Condition (3) is equivalent to say that for every morphism $\h: \Z \to \X$
the $\Mor_\mC(\X,\X)$-linear functor $$\Mor_\mC(\F \circ \h, \X): \Mor_\mC(\Y,\X) \to \Mor_\mC(\Z,\X) $$ sends the left $\T$-action on 
$\G : \Y \to \X$ to the free left $\T$-action on $\h: \Z \to \X.$

The $\Mor_\mC(\X,\X)$-linear functor $\Mor_\mC(\h,\X): \Mor_\mC(\X,\X) \to \Mor_\mC(\Z,\X) $ preserves free left $\T$-modules.
Therefore the $\Mor_\mC(\X,\X)$-linear functor $\Mor_\mC(\F \circ \h, \X) $ sends the left $\T$-action on 
$\G : \Y \to \X$ to the free left $\T$-action on $\h: \Z \to \X$ for every $\h: \Z \to \X$ if and only if the $\Mor_\mC(\X,\X)$-linear functor $\Mor_\mC(\F, \X) $ sends the left $\T$-action on $\G : \Y \to \X$ to the free left $\T$-action on $\id: \X \to \X$. This shows that (1) and (2) are equivalent.

We prove that (2) and (3) are equivalent.
Let $\mu:\T \circ \G \to \G$ be the action morphism in $ \Mor_\mC(\Y,\X).$
The $\Mor_\mC(\X,\X)$-linear functor $\Mor_\mC(\F,\X): \Mor_\mC(\Y,\X) \to \Mor_\mC(\X,\X)$ sends the left
$\T$-action on $\G$ to a left $\T$-action on $\G\circ \F$ whose action morphism is $\mu \circ \F: \T \circ \G\circ \F \to \G\circ \F.$
By \cref{rightadj} (2) the left $\T$-action on $\G$ is the endomorphism $\T$-action if and only if $\T \to \T \circ \G\circ \F \xrightarrow{\mu\circ\F} \G\circ \F$  is an equivalence, where the first morphism is induced by the unit $\id \to \G \circ \F.$
This is equivalent to say that the unit $\id \to \G \circ \F$ exhibits
$\G\circ\F$ as the free left $\T$-module generated by $\id: \X \to \X,$
which is $\T.$
So a left $\T$-action on $\G:\Y \to \X$ is the endomorphism action if and only if it is sent by $\Mor_\mC(\F,\X): \Mor_\mC(\Y,\X) \to \Mor_\mC(\X,\X)$ to $\T$ viewed as module over itself.

\end{proof}

\begin{proof}[Proof of \cref{monaT}]
	
If $\phi$ is monadic, $\phi$ is an Eilenberg-Moore object of a monad $\T$ on some $\X \in \mC$. So for every $\Z \in \mC$
the induced functor $\Mor_\mC(\Z,\Y) \to \Mor_\mC(\Z,\X)$ identifies with the functor $$\LMod_\T(\Mor_\mC(\Z,\X)) \to \Mor_\mC(\Z,\X), $$ which is monadic by \cref{monchar}. So (2) holds.
For every morphism $\Z \to \Z' $ in $\mC$ the induced $\Mor_\mC(\X,\X)$-linear functor $\Mor_\mC(\Z',\X) \to \Mor_\mC(\Z,\X)$ 
gives a functor $$\alpha: \LMod_\T(\Mor_\mC(\Z',\X)) \to \LMod_\T(\Mor_\mC(\Z,\X))$$
covering the functor $\Mor_\mC(\Z',\X) \to \Mor_\mC(\Z,\X)$ and preserving the free left $\T$-modules.
The functor $\alpha$ is equivalent over $\Mor_\mC(\Z',\X) \to \Mor_\mC(\Z,\X)$ to the functor $\beta: \Mor_\mC(\Z',\Y) \to \Mor_\mC(\Z,\Y)$. So $\beta$ preserves the free functors. Proposition \ref{prosta} implies that $\phi$ has a left adjoint. So (1) holds.

\vspace{1mm}

We prove that $\phi$ is monadic if (1) and (2) hold.
If $\phi$ admits a left adjoint $\F$, by Corollary \ref{rightadj2} we find that $\phi$ admits an endomorphism algebra $\T \simeq \phi \circ \F$
that universally acts on $\phi$.
The endomorphism left $\T$-action on $\phi$ corresponds to a map 
$\theta: \iota(\Y) \to \LMod_\T \circ \bar{\iota}(\X)$ in $\FUN(\mC^\op, \infty\Cat)$ over $\iota(\X).$  
We like to see that for every $\Z \in \mC$ the functor $$\theta_\Z: \Mor_\mC(\Z,\Y) \to \LMod_\T(\Mor_\mC(\Z,\X)) $$ is an equivalence.
The functor $\theta_\Z $ preserves the left adjoints by \cref{eilmu}.
% since the canonical functor $$\Mor_\mC(\Z,\T) \to \theta_\Z \circ \Mor_\mC(\Z,\F)$$ is an equivalence: it identifies with the canonical map $$\Mor_\mC(\Z,\T) \circ \id \to \Mor_\mC(\Z,\T) \circ \Mor_\mC(\Z,\T) \to \Mor_\mC(\Z,\T),$$ where the latter map is the multiplication of the monad $\Mor_\mC(\Z,\T).$
By \cref{monchar} the functor $\theta_\Z$ is an equivalence since
%corresponds to a left $\Mor_\mC(\Z,\T)$-action on $\Mor_\mC(\Z,\phi): \Mor_\mC(\Z,\Y) \to \Mor_\mC(\Z,\X),$ which is the endomorphism left action since $\theta_\Z $ preserves the left adjoints.
$\theta_\Z$ preserves the left adjoints and $\Mor_\mC(\Z,\phi)$ is monadic.
% and corresponds to the endomorphism left $\Mor_\mC(\Z,\T)$-action, $\theta_\Z$ is an equivalence.So $\phi$ is the Eilenberg-Moore object of $\T$ and monadic.

\end{proof}

The $\infty$-categorical monadicity theorem \cite[Theorem 4.7.3.5]{lurie.higheralgebra} and \cref{monaT} imply the following:

\begin{corollary}\label{monat0}
Let $\mC$ be an $(\infty,2)$-category.
A morphism $\phi: \Y \to \X$ in $\mC$ is monadic if and only if the following conditions hold:
\begin{enumerate}
\item The morphism $\phi: \Y \to \X$ admits a left adjoint.
		
\item For every $\Z \in \mC$ the induced functor $\phi_*: \Mor_\mC(\Z,\Y) \to \Mor_\mC(\Z,\X)$ is conservative.
		
\item For every $\Z \in \mC$ every $\phi_*$-split simplicial object of
$\Mor_\mC(\Z,\Y)$ admits a colimit preserved by $\phi_*.$
		
\end{enumerate}
	
\end{corollary}

\cref{monaT} implies the following corollary:

\begin{corollary}\label{rembras}
Let $\G: \mD \to \mC$ be a 2-functor that admits a 2-left adjoint.
Then $\G$ preserves monadic morphisms and Eilenberg-Moore objects.

\end{corollary}

\begin{proof}

Let $\F$ be a left adjoint of $\G.$ 
% For every $\X \in \mC, \Y \in \mD $ the induced functor 
% $$ \Mor_\mD(\F(\X), \Y) \to \Mor_\mC(\G(\F(\X)), \G(\Y)) \to \Mor_\mC(\X, \G(\Y)) $$ is an equivalence.
Let $\psi: \Y \to \X$ be a right adjoint morphism of $\mD.$ Then $\G(\psi)$ is a right adjoint morphism of $\mC.$ For every $\Z \in \mC$ the functor $$\Mor_\mC(\Z, \G(\Y)) \to \Mor_\mC(\Z, \G(\X))$$ is equivalent to the functor $\Mor_\mD(\F(\Z), \Y) \to \Mor_\mD(\F(\Z), \X).$ So $\G(\psi)$ is monadic if $\psi$ is monadic by \cref{monaT}.
\end{proof}

% \begin{corollary}

% Let $\mC$ be an $(\infty,2)$-category, $\X \in \mC$ and $\T$ a monad on $\X $ in $\mC$ and $\G:\Y \to \X$ an Eilenberg-Moore object of $\T$ in $\mC$. 
% The left $\T$-action on $\G$ is the endomorphism action.

% \end{corollary}

% \begin{proof}

% By definition for every $\Z \in \mC$ the induced functor
% $$\theta: \Mor_\mC(\Z,\Y) \to \LMod_\T(\Mor_\mC(\Z,\X))$$ is an equivalence
% over $\Mor_\mC(\Z,\X).$
% By definition $\G: \Y \to \X $ is monadic and so by \cref{monaT} admits a left adjoint $\F: \X \to \Y.$
% Hence $\Mor_\mC(\Z,\F)$ is left adjoint to $\Mor_\mC(\Z,\G)$. 
% Since $\theta$ is an equivalence, it preserves the left adjoints.
% In particular, taking $\Z=\X$ and the image of the identity of $\X$
% under the functor $\theta \circ \Mor_\mC(\Z,\F)$, which is the image of the 
% left $\T$-action on $\G$ under the functor $\Mor_\mC(\F,\X): \Mor_\mC(\Y,\X) \to \Mor_\mC(\X,\X)$ and a left $\T$-action on $\G \circ \F$, canonically identifies with the free $\T$-module generated by $\id_\X$, which is $\T$.
% Hence by \cref{inpos} the left $\T$-action on $\G$ is the endomorphism action.

% \end{proof}

\begin{proposition}\label{preloc}

Let $\mC$ be an $(\infty,2)$-category, $\X \in \mC$ and $\T$ a monad on $\X $ in $\mC$ and $\G:\Z \to \X$ a right adjoint morphism whose associated monad we denote by $\L$ and $\Y \to \X$ a right adjoint morphism that is an Eilenberg-Moore object for $\T$.
There is a canonical equivalence
$$ \Mor_{\mC_{/\X}}(\Z,\Y) \simeq \Map_{\Alg(\Mor_\mC(\X,\X))^\op}(\L,\T).$$

\end{proposition}

\begin{proof}

There is a canonical equivalence
$$ \Mor_{\mC_{/\X}}(\Z,\Y) \simeq \{ \h\}\times_{\Mor_\mC(\Z,\X)} \Mor_\mC(\Z,\Y) \simeq \{ \h\}\times_{\Mor_\mC(\Z,\X)} \LMod_\T(\Mor_\mC(\Z,\X)) $$$$ \simeq \Map_{\Alg(\Mor_\mC(\X,\X))}(\T,\L)\simeq \Map_{\Alg(\Mor_\mC(\X,\X))^\op}(\L,\T).$$

\end{proof}

\begin{remark}

\cref{preloc} suggests that if Eilenberg-Moore objects exist, 
there is a localization
$$ \mC_{/\X}^\R \rightleftarrows \Alg(\Mor_\mC(\X,\X))^\op, $$
where the left adjoint sends a right adjoint morphism to $\X$ to its associated monad on $\X$ and where the right adjoint sends a monad on $\X$ to its Eilenberg-Moore object.
To make this precise one has to make the process of associating a monad to a right adjoint morphism functorial or likewise one has to make the process of assigning the Eilenberg-Moore object functorial.
This is highly non trivial although we know how to map morphisms to morphisms (\cref{funcmon}, \cref{funcem}).
We will make functoriality rigorous in \cref{Classa} to prove this aforementioned localization result (\cref{uhnggfdaa}).

\end{remark}

\subsection{Kleisli objects}

For the next proposition we use that every Eilenberg-Moore object admits a left adjoint by \cref{monaT}.

\begin{proposition}\label{endoii}

Let $\mC$ be an $(\infty,2)$-category, $\X \in \mC$ and $\T$ a monad in $\mC$
on $\X$. Let $\G: \Y \to \X$ be an Eilenberg-Moore object for $\T$ whose left adjoint we denote by $\F: \X \to \Y.$
Then $\F$ carries a canonical right $\T$-action that is the endomorphism action
with respect to the right action of $\Mor_\mC(\X,\X)$ on $\Mor_\mC(\X,\Y).$ 

\end{proposition}

\begin{proof}

\label{charos}

By the 2-Yoneda-lemma the left $\T$-action on $\phi$ corresponds to a map 
$\theta: \iota(\Y) \to \LMod_\T \circ \bar{\iota}(\X)$ in $\FUN(\mC^\op, \infty\Cat)$ over $\iota(\X).$ 
Since $\G: \Y \to \X$ is an an Eilenberg-Moore object of $\T$,
for every $\Z \in \mC$ the functor $$\theta_\Z: \Mor_\mC(\Z,\Y) \to \LMod_\T(\Mor_\mC(\Z,\X)) $$ is an equivalence.
Thus $\theta$ is an equivalence.
The equivalence $\theta$ in $\FUN(\mC^\op, \infty\Cat)$ over $\iota(\X)$
sends the endomorphim left action of $\Mor_\mC(\X,\X)$ on $\X$
to a right $\Mor_\mC(\X,\X)$-linear functor
$$ \Mor_\mC(\X,\Y) \to \LMod_\T(\Mor_\mC(\X,\X))$$
over $\Mor_\mC(\X,\X)$ refining $\theta_\X,$ which therefore is a right $\Mor_\mC(\X,\X)$-linear equivalence. 
Since $\theta_\X$ is an equivalence, it preserves free objects.
Thus the free left $\T$-action on the identity of $\X$, which is $\T$ viewed as left module over itself, is the image of the identity of $\X$ under $\theta_\X \circ \Mor_\mC(\X,\F)$. The latter is the image of the left $\T$-action on $\G$ under the left $\Mor_\mC(\X,\X)$-linear functor $\Mor_\mC(\F,\X).$
The associative algebra structure on $\T$ in $\Mor_\mC(\X,\X)$ gives rise to a
biaction of $\T$ on itself whose underlying left and right actions are the free ones on the tensor unit $\id_\X.$
A $\T$-biaction corresponds to a right $\T$-action in $\LMod_\T(\Mor_\mC(\X,\X))$ with respect to the right action of
$\Mor_\mC(\X,\X)$ on $\LMod_\T(\Mor_\mC(\X,\X))$.
Since $\theta_\X$ is a right $\Mor_\mC(\X,\X)$-linear equivalence, 
the right $\T$-action on the free left $\T$-module generated by $\id_\X$ corresponds via $\theta_\X$ to a right $\T$-action on $\F:\X \to \Y.$

We prove that this right $\T$-action on $\F$ is the endomorphism action.
By \cref{rightadj} it suffices to prove that the canonical morphism
\begin{equation}\label{epolt}
\T \xrightarrow{\eta \circ \T} \G \circ \F \circ \T \xrightarrow{\G \circ \mu} \G \circ \F \end{equation} is an equivalence, where $\mu: \F \circ \T \to \F$ is the action morphism and $\eta: \id \to \G \circ \F$ is the unit.

By definition the right $\T$-action on $\F$ is sent by $\theta_\X$ to
the canonical right $\T$-action on the free left $\T$-module generated by $\id_\X.$
The image in $\Mor_\mC(\X,\X)$ of the right $\T$-action on $\F$ under $\theta_\X$ is the image of the right $\T$-action on $\F$ under $\Mor_\mC(\X,\G),$ which is by definition is the free right $\T$-action on the tensor unit of $\id_\X.$
In particular, the induced action morphism $\G \circ \mu:\G \circ \F \circ \T \to \G \circ \F$ of $\Mor_\mC(\X,\G)(\F) = \G \circ \F$
is the action morphism of the free right $\T$-action on the tensor unit of $\id_\X$, which is the multiplication morphism of the associative algebra $\T.$
Thus the morphism (\ref{epolt}) is an equivalence.

\end{proof}

Replacing an $(\infty,2)$-category $\mC$ by $\mC^\co$ we obtain the following:

\begin{corollary}\label{Kleicom}
Let $\mC$ be an $(\infty,2)$-category, $\X \in \mC$ and $\T$ a monad in $\mC$
on $\X$. Let $\F: \X \to \Y$ be a Kleisli object for $\T$ whose right adjoint we denote by $\G: \Y \to \X.$
Then $\G$ carries a canonical left $\T$-action that is the endomorphism action
with respect to the left action of $\Mor_\mC(\X,\X)$ on $\Mor_\mC(\Y,\X).$ 

\end{corollary}

\begin{remark}

Let $\mC$ be an $(\infty,2)$-category, $\X \in \mC$ and $\T$ a monad in $\mC$
on $\X$. Let $\F: \X \to \Y$ be a Kleisli object for $\T$ whose right adjoint we denote by $\G: \Y \to \X,$ and let $\G': \Y \to \X$ be an Eilenberg-Moore object for $\T$.
By \cref{Kleicom} the morphism $\G: \Y \to \X $ carries a canonical left $\T$-action that is the endomorphism action with respect to the left action of $\Mor_\mC(\X,\X)$ on $\Mor_\mC(\Y,\X).$ 
By the universal property of the Eilenberg-Moore object there is a canonical morphism $\Y \to \Y'$ over $\X$.

We will prove later if Eilenberg-Moore objects exist, every right adjoint morphism to $\X$ can universally be turned into a monadic morphism to $\X$ and the local equivalences are precisely the morphisms over $\X$ that induce an equivalence on endomorphism monads (\cref{uhnggfdaa}).
Consequently, the Eilenberg-Moore object of the monad $\T$ on $\X$
is the monadic approximation of the Kleisli-object of $\T.$

\end{remark}

\begin{proposition}
Let $\mC$ be an $(\infty,2)$-category, $\X \in \mC$ and $\T$ a monad in $\mC$
on $\X$. Let $\G: \Y \to \X$ be an Eilenberg-Moore object for $\T$ whose left adjoint we denote by $\F: \X \to \Y.$ 

For every $\Z \in \mC$ the endomorphism right $\T$-action on $\F: \X \to \Y$ of \cref{endoii}
induces a functor
$$ \Mor_\mC(\Y,\Z) \to \RMod_\T(\Mor_\mC(\X,\Z)) $$
over $ \Mor_\mC(\X,\Z)$ that preserves the left adjoints.
This functor is an equivalence if and only if the induced functor
$\Mor_\mC(\Y,\Z) \to \Mor_\mC(\X,\Z)$ is monadic.

\end{proposition}

\begin{proof}

By the 2-Yoneda-lemma and \cref{endoi} the right $\T$-action on $\F: \X \to \Y$ of \cref{endoii} corresponds to a morphism
$$\theta: \Mor_\mC(\Y,-) \to \RMod_\T(\Mor_\mC(\X,-)) $$
over $ \Mor_\mC(\X,-)$ in $\Fun(\mC, \infty\Cat).$
We prove that for every $\Z \in \mC$ the induced functor
$$\theta_\Z: \Mor_\mC(\Y,\Z) \to \RMod_\T(\Mor_\mC(\X,\Z)) $$
over $ \Mor_\mC(\X,\Z)$ is an equivalence.
The functor $\theta_\Z$ preserves the left adjoints. This holds because
for every morphism $\h: \X \to \Z$ the canonical morphism
$ \T \circ \h \to \G \circ \F \circ \h$
identifies with the morphism
$$\rho: \T \circ \h \xrightarrow{\T \circ \eta \circ \h}  \T \circ \G \circ \F \circ \h \xrightarrow{\mu \circ  \F \circ \h} \G \circ \F \circ \h,$$
where $\eta: \id \to \G \circ \F$ is the unit and $\mu: \T \circ \G \to \G$ is the action morphism on $\G.$
The morphism $\rho $ is an equivalence by \cref{rightadj} because the right $\T$-action on $\F$ is the endomorphism action by \cref{endoii}.
Since the target of the functor $\theta_\Z$ is monadic over $ \Mor_\mC(\X,\Z)$,
and $\theta_\Z$ preserves the left adjoints, the functor $\theta_\Z$ is an equivalence by \cref{monchar} since the source of $\theta_\Z$ is monadic over $ \Mor_\mC(\X,\Z)$.

% This holds since the functor
% $\Mor_\mC(\Y,\Z) \to \Mor_\mC(\X,\Z)$ is conservative and preserves geometric realizations 

\end{proof}

\subsection{Construction of Eilenberg-Moore objects}\label{Exiis}

Next we study existence of Eilenberg-Moore objects
(Corollary \ref{EMex}) that we obtain from a description of Eilenberg-Moore objects as a lax limit (Theorem \ref{stru}).
We start with defining lax functors. 
Lax functors in different models were studied by \cite{GAGNA2021107986}, \cite{abellan2023comparing}. Lax limits were studied by \cite{abellan2024straightening}, \cite{gagna2020fibrations}, \cite{berman2020lax}.

\begin{notation}Let $\rS$ be an $\infty$-category and $\mE \subset \Fun([1],\rS)$ a full subcategory.
	
\begin{enumerate}
\item Let $\infty\Cat_{/ \rS}^{\mathrm{loc}} \subset \infty\Cat_{/ \rS}$ be the subcategory of locally cocartesian fibrations over $\rS$ and functors over $\rS$ preserving locally cocartesian lifts.

\item Let $\infty\Cat^\mE_{/ \rS} \subset \infty\Cat_{/ \rS} $ be the subcategory of cocartesian fibrations relative to $\mE$ and functors over $\rS$ preserving cocartesian lifts of morphisms of $\mE.$

\item Let $\infty\Cat_{/ \rS}^{\mathrm{loc}, \mE} \subset \infty\Cat_{/ \rS} $ be the intersection of the subcategories 
$ \infty\Cat_{/ \rS}^{\mathrm{loc}}$ and $\infty\Cat^\mE_{/ \rS}.$
\end{enumerate}	

\end{notation}

%The canonical left $\infty\Cat$-action on $ \infty\Cat_{/ \rS}$ restricts to a left$\infty\Cat$-action on $\infty\Cat_{/ \rS}^{\mathrm{locart}}$.Moreover the $\infty\Cat$-linear inclusion $\infty\Cat_{/ \rS}^{\mathrm{locart}} \subset \infty\Cat_{/ \rS}$ induces a $\infty\Cat$-linear inclusion $\infty\Cat_{/ \rS}^{\cart} \subset \infty\Cat_{/ \rS}^{\cart,\oplax}$.

\begin{notation}
	
Let $\rS$ be an $\infty$-category, $\mE \subset \Fun([1],\rS)$ a full subcategory and $\mC$ an $(\infty,2)$-category.
Let $$\mE\mathrm{-}\LaxFun(\rS,\mC) \subset \FUN(\mC^\op, \infty\Cat_{/ \rS}^{\loc, \mE}) $$ be the full subcategory of 2-functors
$\mC^\op \to \infty\Cat_{/ \rS}^{\mathrm{loc}, \mE}$ such that for every $\s \in \rS$ the 2-functor $\mC^\op \to \infty\Cat_{/ \rS}^{\mathrm{loc}, \mE} \to \infty\Cat $ is representable, where the latter 2-functor takes the fiber over $\s \in \rS.$

For $\mE$ the full subcategory spanned by the equivalences we remove $\mE$ from the notation and terminology.

\end{notation}

\begin{definition}Let $\rS$ be an $\infty$-category, $\mE \subset \Fun([1],\rS)$ a full subcategory and $\mC$ an $(\infty,2)$-category.

\begin{itemize}
\item An $\mE$-lax functor $\rS \to \mC$ is an object of $\mE\mathrm{-}\LaxFun(\rS,\mC)$.

\vspace{1mm}

\item An $\mE$-oplax functor $\F: \rS \to \mC$ is an $\mE$-lax functor $\rS \to \mC^\co,$ which we denote by $\F^\co.$
\end{itemize}

\end{definition}

\begin{definition}
Let $\F: \rS \to \mC$ be an $\mE$-lax functor.
For every $\s \in \rS$ the image $\F(\s) \in \mC$ represents the presheaf $\mC^\op \xrightarrow{\F} \infty\Cat_{/ \rS}^{\loc, \mE} \xrightarrow{ } \infty\Cat,$ where the latter 2-functor takes the fiber over $\s \in \rS.$
For every morphism $\alpha: \s \to \rt$ in $\mE$ the image $\F(\alpha): \F(\s) \to \F(\rt) $ in $\mC$ represents the presheaf $$\mC^\op \xrightarrow{\F} \infty\Cat_{/ \rS}^{\loc, \mE} \xrightarrow{ } \infty\Cat_{/ [1]}^{\cocart} \simeq \Fun([1],\infty\Cat),$$ where the latter 2-functor takes the pullback along the functor $[1] \to \rS$ taking $\alpha.$

\end{definition}

\begin{remark}
Let $\mC$ be an $(\infty, 2)$-category
and $\rS \to \T$ a functor.
The functor $$\FUN(\mC^\op, \infty\Cat_{/ \T}^{\loc, \mE}) \to \FUN(\mC^\op, \infty\Cat_{/ \rS}^{\loc, \mE})$$
restricts to a functor
$$\mE\mathrm{-}\LaxFun(\T,\mC) \to \mE\mathrm{-}\LaxFun(\rS,\mC).$$
\end{remark}

\begin{remark}

Let $\mC $ be an $(\infty,2)$-category.
The 2-Yoneda-embedding $$\mC \hookrightarrow \FUN(\mC^\op,\infty\Cat)$$ gives rise to an embedding of $(\infty,2)$-categories $$\FUN(\rS,\mC) \subset \FUN(\mC^\op,\infty\Cat)^\rS \simeq 
\FUN(\mC^\op, \FUN(\rS,\infty\Cat))
\simeq $$$$ \FUN(\mC^\op, \infty\Cat_{/ \rS}^{\cocart}) \subset \FUN(\mC^\op,\infty\Cat_{/ \rS}^{\mathrm{loc}, \mE}). $$
The latter induces an embedding of $(\infty,2)$-categories $$\FUN(\rS,\mC) \hookrightarrow \mE\mathrm{-}\LaxFun(\rS,\infty\Cat).$$
The latter induces an equivalence of $(\infty,2)$-categories $$\FUN(\rS,\mC) \simeq \Fun([1],\rS)\mathrm{-}\LaxFun(\rS,\infty\Cat)$$
identifying 2-functors $\rS \to \mC$ with $\Fun([1],\rS)$-lax functors.

\end{remark}

\begin{lemma}
Let $\rS$ be an $\infty$-category. There is an equivalence of $(\infty, 2)$-categories:
	
$$ \LaxFun(\rS,\infty\Cat) \simeq \infty\Cat_{/ \rS}^{\loc}.$$
	
\end{lemma}
\begin{proof}
The full subcategory
$$\LaxFun(\rS,\infty\Cat) \subset \FUN(\infty\Cat^\op, \infty\Cat_{/ \rS^\op}^{\loc})$$
agrees with the full subcategory $\FUN^{ \mathrm{coten}}(\infty\Cat^\op,\infty\Cat_{/ \rS}^{\mathrm{loc}})$
of 2-functors preserving cotensors.
So there is an equivalence of $(\infty,2)$-categories:
$$\infty\Cat_{/ \rS}^{\loc} \simeq (\LinFun_{\infty\Cat}(\infty\Cat, \infty\Cat_{/ \rS}^{\loc})^\op)^\op \simeq \FUN^{ \mathrm{coten}}(\infty\Cat^\op,\infty\Cat_{/ \rS}^{\mathrm{loc}})= \LaxFun(\rS,\infty\Cat). $$
	
\end{proof}

Let $\mE \subset \Fun([1],\rS)$ be a full subcategory. Recall the right adjoint 2-functor $$\Fun^\mE_\rS(\rS,-): \infty\Cat^\mE_{/ \rS} \to \infty\Cat$$ from Example \ref{Imporso} for the next definition:

\begin{definition}
Let $\rS$ be an $\infty$-category, $\mE \subset \Fun([1],\rS)$ a full subcategory, $\mC$ an $(\infty, 2)$-category and $\F: \rS \to \mC$ an $\mE$-lax functor.
The $\mE$-lax limit of $\F$ if it exists, is the object of $\mC$ representing the image of $\F$ under the 2-functor
$$\mE\mathrm{-}\LaxFun(\rS,\mC) \subset \FUN(\mC^\op, \infty\Cat_{/ \rS}^{\loc, \mE})
\xrightarrow{\FUN(\mC^\op, \Fun^\mE_\rS(\rS,-))} \FUN(\mC^\op, \infty\Cat).$$		
%is the object of $\mC$universal with respect to the diagonal 2-functor $\mC \simeq \LaxFun(*,\mC) \to \LaxFun(\rS,\mC).$
\end{definition}

\begin{definition}
Let $\rS$ be an $\infty$-category, $\mE \subset \Fun([1],\rS)$ a full subcategory, $\mC$ an $(\infty, 2)$-category and $\F: \rS \to \mC$ an $\mE$-oplax functor.
The $\mE$-oplax limit of $\F$ if it exists, is the $\mE$-lax limit of $\F^\co$.

\end{definition}

We have the following existence result for lax limits:

\begin{proposition}\label{Exist}
Let $\rS$ be an $\infty$-category, $\mE \subset \Fun([1],\rS)$ a full subcategory and $\mC$ an $(\infty,2)$-category.
If $\mC$ or $\mC^\op$ are presentable $(\infty,2)$-categories, every $\mE$-lax functor $\F: \rS \to \mC$ admits an $\mE$-lax limit.	
	
\end{proposition}

\begin{proof}
	
The $\mE$-lax functor $\F: \rS \to \mC$ is given by a 2-functor
$\alpha: \mC^\op \to \infty\Cat_{/ \rS}^{\loc, \mE}$
such that for every $\s \in \rS$ the 2-functor $\mC^\op \to \infty\Cat_{/ \rS}^{\mathrm{loc}, \mE} \to \infty\Cat$ is representable, where the latter 2-functor takes the fiber over $\s \in \rS.$
If $\mC$ is a presentable $(\infty,2)$-category, a 2-functor $\mC^\op\to \infty\Cat$ is representable if and only if it preserves small limits and cotensors
by \cref{remqa} and \cite[Corollary 5.5.2.9]{lurie.HTT}.
If $\mC^\op$ is a presentable $(\infty,2)$-category, a 2-functor $\mC^\op \to \infty\Cat$ is representable if and only if it preserves small limits, cotensors and is accessible by \cref{remqa}, \cref{remqa} and \cite[Corollary 5.5.2.9]{lurie.HTT}.

Since for every $\s \in \rS$ the 2-functor $\infty\Cat_{/ \rS}^{\mathrm{loc}, \mE} \to \infty\Cat $ taking the fiber over $\s$ preserves small limits and cotensors and is accessible, $\alpha: \mC^\op \to \infty\Cat_{/ \rS}^{\mathrm{loc}, \mE}$ preserves small limits and cotensors if $\mC$ is a presentable $(\infty,2)$-category, and $\alpha$ preserves small limits and cotensors and is accessible if $\mC^\op$ is a presentable $(\infty,2)$-category.
The inclusion $\infty\Cat_{/ \rS}^{\mathrm{loc}, \mE} \subset \infty\Cat^\mE_{/ \rS}$ preserves small limits and cotensors and is accessible.
The right adjoint 2-functor $\Fun^\mE_\rS(\rS,-): \infty\Cat^\mE_{/ \rS} \to \infty\Cat$ preserves small limits and cotensors and is accessible by \cref{remqa} and \cite[Corollary 5.5.2.9]{lurie.HTT} since source and target are presentable.
Hence the composition $$\rho: \mC^\op \xrightarrow{\alpha} \infty\Cat_{/ \rS}^{\loc, \mE} \subset \infty\Cat^\mE_{/ \rS} \xrightarrow{\Fun^\mE_\rS(\rS,-)} \infty\Cat$$ preserves small limits and cotensors and so is representable if $\mC$ is a presentable $(\infty,2)$-category. Similarly, $\rho$ preserves small limits and cotensors and is accessible, and so is representable if $\mC^\op$ is a presentable $(\infty,2)$-category.

\end{proof}

\begin{proposition}

Let $\mC$ be an $(\infty,2)$-category and $\mB$ a full subcategory such that the embedding $\mB \subset \mC$ of $(\infty,2)$-categories admits a left adjoint.
Let $\rS$ be an $\infty$-category and $\mE \subset \Fun([1],\rS)$ a full subcategory. The full subcategory $\mB$ is closed under $\mE$-lax limits.

\end{proposition}

\begin{proof}

Let $\F: \rS \to \mB$ be an $\mE$-lax functor such that the composition 
$\rS \to \mB \subset \mC$ admits an $\mE$-lax limit
$\Z \in \mC$. It suffices to see that for every local equivalence $\A \to \B$ in $\mC$ the induced functor $ \Mor_\mC(\B,\Z) \to \Mor_\mC(\A,\Z)$ is an equivalence.
By definition this functor is the image of the morphism $\A \to \mB$ in $\mC$ under the 2-functor $$\mC \xrightarrow{\alpha} \infty\Cat^{\loc, \mE} \xrightarrow{\Fun^\mE_\rS(\rS,-)} \infty\Cat,$$ where $\alpha$ classifies $\F.$
The 2-functor $\alpha$ sends the morphism $\A \to \B$ to an equivalence since for every $\s \in \rS$ the induced functor
on the fiber over $\s$ identifies with the functor 
$$ \Mor_\mC(\B,\F(\s)) \to \Mor_\mC(\A,\F(\s)) $$ 
which is an equivalence since $\F(\s) \in \mB$ by assumption.
\end{proof}

Next we apply the theory of lax limits to construct Eilenberg-Moore objects of monads.

\begin{notation}

Let $\max \subset \Fun([1],\Delta^\op)$ be the full subcategory of inert maps of $\Delta^\op$ preserving the maximum.

\end{notation}

\begin{construction}

Let $\mC$ be an $(\infty, 2)$-category and $\T$ a monad on some $\X\in \mC.$

Let $ \iota: \mC \hookrightarrow \FUN(\mC^\op, \infty\Cat) $ be the 2-Yoneda-embedding and $ \bar{\iota}: \mC \hookrightarrow \FUN(\mC^\op, \Mor_{\mC}(\X,\X)\mathrm{-}\LMod) $ the canonical lift of the 2-Yoneda-embedding.
Let $$\gamma: \infty\Cat^{\loc, \max}_{/\Mor_{\mC}(\X,\X)^\ot} \to \infty\Cat^{\loc, \max}_{/\Delta^\op}$$
be the 2-functor taking pullback along the lax monoidal functor $\Delta^\op \to \Fun_\rS(\X,\X)^\ot$ classifying $\T.$

Let $\theta_\T$ be the image of $\bar{\iota}(\X)$ under the 2-functor
$$\FUN(\mC^\op, \Mor_{\mC}(\X,\X)\mathrm{-}\LMod)
\subset \FUN(\mC^\op, \infty\Cat^{\loc, \max}_{/\Mor_{\mC}(\X,\X)^\ot})
\xrightarrow{\FUN(\mC^\op,\gamma)} \FUN(\mC^\op, \infty\Cat^{\loc, \max}_{/\Delta^\op}).$$

Then $\theta_\T$ is a $\max$-lax functor $\Delta^\op \to \mC.$

\end{construction}

\begin{remark}
    
The $\max$-lax functor $\theta_\T: \Delta^\op \to \mC$ sends every $[\n] \in \Delta$ to $\X$ and the morphism $[1] \simeq \{0,\n\} \subset [\n] $ to $\T^{\circ \n} : \X \to \X.$

This is because for every $[\n]\in \Delta^\op$ the composition $ \mC^\op \xrightarrow{\theta_\T}
\infty\Cat^{\loc, \max}_{/\Delta^\op} \to \infty\Cat$, where the latter functor takes the fiber over $[\n] \in \Delta$, is $\iota(\X).$
For every map $\alpha: [1] \simeq \{0,\n\} \subset [\n] $ in $\Delta$ the composition $$ \mC^\op \xrightarrow{\theta}
\infty\Cat^{\loc, \max}_{/\Delta^\op} \to \infty\Cat^{\cocart}_{/[1]} \simeq \Fun([1], \infty\Cat),$$ where the latter functor takes the pullback along the functor $[1]\to \Delta^\op$ taking $\alpha$, classifies the map $\iota(\T)^{\circ\n}: \iota(\X) \to \iota(\X).$
    
\end{remark}

\begin{proposition}\label{stru}
	
Let $\mC$ be an $(\infty, 2)$-category and $\T$ a monad on some $\X\in \mC.$
The $\max$-lax limit of $\theta_\T$ if it exists, is an Eilenberg-Moore object for $\T.$

\end{proposition}

\begin{proof}

% Let $ \iota: \mC \hookrightarrow \FUN(\mC^\op, \infty\Cat) $ be the 2-Yoneda-embedding and $$\gamma: \infty\Cat^{\loc, \max}_{/\Mor_{\mC}(\X,\X)^\ot} \to \infty\Cat^{\loc, \max}_{/\Delta^\op}$$
% the 2-functor taking pullback along $\T: \Delta^\op \to \Mor_{\mC}(\X,\X)^\ot.$

% Let $\theta$ be the image of $\bar{\iota}(\X)$ under the 2-functor
% $$\FUN(\mC^\op, \Mor_{\mC}(\X,\X)\mathrm{-}\LMod)
% \subset \FUN(\mC^\op, \infty\Cat^{\loc, \max}_{/\Mor_{\mC}(\X,\X)^\ot})
% \xrightarrow{\FUN(\mC^\op,\gamma)} \FUN(\mC^\op, \infty\Cat^{\loc, \max}_{/\Delta^\op}).$$

% For every $[\n]\in \Delta^\op$ the composition $ \mC^\op \xrightarrow{\theta_\T}
% \infty\Cat^{\loc, \max}_{/\Delta^\op} \to \infty\Cat$, where the latter functor takes the fiber over $[\n] \in \Delta$, is $\iota(\X).$
% For every map $\alpha: [1] \simeq \{0,\n\} \subset [\n] $ in $\Delta$ the composition $$ \mC^\op \xrightarrow{\theta}
% \infty\Cat^{\loc, \max}_{/\Delta^\op} \to \infty\Cat^{\cocart}_{/[1]} \simeq \Fun([1], \infty\Cat),$$ where the latter functor takes the pullback along the functor $[1]\to \Delta^\op$ taking $\alpha$, classifies the map $\iota(\T)^{\circ\n}: \iota(\X) \to \iota(\X).$

The $\max$-lax limit of $\theta_\T$ if it exists, represents the image of $\theta_\T$ under the functor 
$$\FUN(\mC^\op,\Fun^{\max}_{\Delta^\op}(\Delta^\op,-)): \FUN(\mC^\op, \infty\Cat^{\loc, \max}_{/\Delta^\op}) \to  \FUN(\mC^\op, \infty\Cat).$$

The 2-functor $\LMod_\T: \Mor_{\mC}(\X,\X)\mathrm{-}\LMod \to \infty\Cat$ 
factors as 
$$\Mor_{\mC}(\X,\X)\mathrm{-}\LMod \subset \infty\Cat^{\loc, \max}_{/\Mor_{\mC}(\X,\X)^\ot} \xrightarrow{\gamma} \infty\Cat^{\loc, \max}_{/\Delta^\op} \xrightarrow{\Fun^{\max}_{\Delta^\op}(\Delta^\op,-)} \infty\Cat. $$

Hence the $\max$-lax limit of $\theta_\T$ if it exists, represents the image of $\bar{\iota}(\X)$ under the 2-functor
$$ \FUN(\mC^\op,\LMod_\T): \FUN(\mC^\op, \Mor_{\mC}(\X,\X)\mathrm{-}\LMod) \to \FUN(\mC^\op, \infty\Cat)$$
and so by definition is an Eilenberg-Moore object for $\T.$
	
\end{proof}

\begin{corollary}\label{EMex}
	
Every presentable $(\infty,2)$-category $\mC$ admits Eilenberg-Moore objects, coEilenberg-Moore objects, Kleisli objects and coKleisli objects.
	
\end{corollary}

\begin{proof}
We apply \cref{Exist} and \cref{stru} to deduce that
$\mC$ admits Eilenberg-Moore objects and Kleisli objects.
To deduce that $\mC$ has coEilenberg-Moore objects and coKleisli objects
we use that trivially $\mC^\mathrm{co}$ is a presentable $(\infty,2)$-category if $\mC$ is a presentable $(\infty,2)$-category.

\end{proof}

In the following we study lax limits in the $(\infty,2)$-category $\infty\Cat_{/\T}$
for every $\infty$-category $\T.$ We apply this to construct and analyze Eilenberg-Moore objects in $\infty\Cat_{/\T}$.

\vspace{1mm}

Our strategy is the following: 

\begin{enumerate}
\item We first study how to produce lax functors from an $\infty$-category 
$\rS $ to the $(\infty,2)$-category $\infty\Cat_{/\T}.$
Via the Grothendieck construction maps $\mC \to \rS \times \T$ of cocartesian fibrations over $\rS$ give rise to functors 
$\rS \to \infty\Cat_{/\T}$.
We show that maps $\mC \to \rS \times \T$ of locally cocartesian fibrations over $\rS$ give rise to lax functors
$\rS \to \infty\Cat_{/\T}.$

\item As a second step we prove that the lax limit of a lax functor $\rS \to \infty\Cat_{/\T}$
arising from a map $\phi: \mC \to \rS \times \T$ of locally cocartesian fibrations over $\rS$ can be computed in terms of $\phi$
(\cref{simpli}).
    
\end{enumerate}

\vspace{1mm}

The following are \cite[Definition 3.21., Remark 3.23.]{heine2024local}:
\begin{definition}
A functor $\phi: \mC \to \rS$ is flat if the functor 
$ \mC \times_\rS (-): {{\infty\Cat}}_{/\rS} \to {{\infty\Cat}}_{/\mC} $
admits a right adjoint.

\end{definition}
\begin{remark}
The pullback and opposite of a flat functor are flat.
Cocartesian fibrations are flat  \cite[B.3.11.]{lurie.higheralgebra}.
\end{remark}
\begin{notation}\label{Parfun}
Let $\phi: \mC \to \rS, \beta: \mC \to \T$ be functors.
If $\phi$ is flat, the functor 
$$ {{\infty\Cat}}_{/\rS} \xrightarrow{(-)\times_\rS\mC}{{\infty\Cat}}_{/\mC} \xrightarrow{\beta_*} {{\infty\Cat}}_{/\T}$$
admits a right adjoint, which we denote by $\Fun^{\rS}_\T(\mC, -).$
If $\beta=\phi$, we write $\Fun^{\rS}(\mC, -)$ for $\Fun^{\rS}_\T(\mC, -).$

\end{notation}

%So for any cocartesian or cartesian fibration $\mC \to \rS$, functors $\mC \to \T,\mD \to \T$ and $\mB \to \rS $
%there is a canonical equivalence $$ \Fun_{\rS}(\mB, \Fun^{\rS}_\T(\mC, \mD)) \simeq \Fun_{\T}(\mB \times_\rS \mC,  \mD). $$

The next remark is \cite[Remark 3.71]{HEINE2023108941}:

\begin{remark}\label{zyaqcfd}
Let $\T \to \rS, \mB \to \rS, \mC \to \T, \mD \to \T$ be functors such that the composition $\mC \to \T \to \rS$ is flat.

\begin{enumerate} 

\item Let $\rS' \to \rS$ be a functor. There is a canonical equivalence $$ \rS' \times_\rS \Fun^{\rS}_\T(\mC, \mD) \simeq \Fun^{\rS'}_{\rS' \times_\rS \T}(\rS' \times_\rS \mC,\rS' \times_\rS \mD) $$ over $\rS'$.
In particular, for every object $\s $ of $\rS$ there is a canonical equivalence $$  \Fun^{\rS}_\T(\mC, \mD)_\s \simeq \Fun_{\T_\s}(  \mC_\s,  \mD_\s). $$

% \item There is a canonical equivalence $$ \Fun^{\rS}_\T(\mC, \mD) \simeq \rS \times_{\Fun^{\rS}(\mC, \T) } \Fun^{\rS}(\mC, \mD) $$ over $\rS$.

% More generally, for any functor $\mB \to \Fun^{/\rS}_\T(\mC, \mE) $ over $\rS$ there is a canonical equivalence$$\mB \times_{\Fun^{/\rS}_\T(\mC, \mE)} \Fun^{/\rS}_\T(\mC, \mD)\simeq \Fun^{/\mB}_{\mB \times_\rS\mE}(\mB \times_\rS\mC, \mB \times_\rS\mD) $$ over $\mB$ given by the composition $$ \mB \times_{\Fun^{/\rS}_\T(\mC, \mE)} \Fun^{/\rS}_\T(\mC, \mD) \simeq \mB \times_{(\mB \times_\rS \Fun^{/\rS}_\T(\mC, \mE)) } (\mB \times_\rS \Fun^{/\rS}_\T(\mC, \mD)) \simeq $$$$ \mB \times_{\Fun^{/\mB}_{\mB \times_\rS \T }(\mB \times_\rS\mC, \mB \times_\rS\mE)} \Fun^{/\mB}_{\mB \times_\rS \T }(\mB \times_\rS\mC, \mB \times_\rS\mD) \simeq \Fun^{/\mB}_{\mB \times_\rS\mE}(\mB \times_\rS\mC, \mB \times_\rS\mD) $$ of canonical equivalences over $\mB.$

\item For every functors $ \mB \to \T', \alpha: \T' \to \T $ there is a canonical equivalence over $\rS:$
$$ \Fun^{\rS}_{\T}(\alpha_\ast(\mB),  \mD) \simeq  \Fun^{\rS}_{\T'}(\mB,  \T' \times_{\T} \mD).$$

\end{enumerate}

\end{remark}

\begin{construction}
    
Let $\rS, \T$ be $\infty$-categories. Let
$$ \nu: (\infty\Cat^{\loc,\mE}_{/\rS})_{/\rS \times \T} \to \mE\mathrm{-}\LaxFun(\rS,\infty\Cat_{/\T})\subset\FUN((\infty\Cat_{/\T})^\op,\infty\Cat_{/ \rS}^{\mathrm{loc},\mE}) $$
be the 2-functor sending a functor
$\mC \to \rS \times \T$ to the functor $$(\mB \to \T) \mapsto \Fun_{\rS \times \T}^\rS(\rS \times \mB,\mC).$$

\end{construction}

% \begin{remark}
    
% Let $\rS, \T$ be $\infty$-categories. 
% The functor $$ \nu: (\infty\Cat_{/\rS})_{/\rS \times \T} \to \FUN^{ \mathrm{coten}}((\infty\Cat_{/\T})^\op,\infty\Cat_{/ \rS}) $$
% restricts to a functor
% $$ (\infty\Cat^{\loc,\mE}_{/\rS})_{/\rS \times \T} \to \FUN^{ \mathrm{coten}}((\infty\Cat_{/\T})^\op,\infty\Cat_{/ \rS}^{\mathrm{loc},\mE}) $$
% \end{remark}

\begin{remark}

For every functor $\rS' \to \rS, \T' \to \T$ there is a commutative square
\begin{equation*} 
\begin{xy}
\xymatrix{
(\infty\Cat^{\loc,\mE}_{/\rS})_{/\rS \times \T}
\ar[d] 
\ar[r]  & \mE\mathrm{-}\LaxFun(\rS,\infty\Cat_{/\T})\subset\FUN((\infty\Cat_{/\T})^\op,\infty\Cat_{/ \rS}^{\mathrm{loc},\mE})
\ar[d]
\\
(\infty\Cat^{\loc,\mE}_{/\rS'})_{/\rS' \times \T'} \ar[r] & \mE\mathrm{-}\LaxFun(\rS',\infty\Cat_{/\T'})\subset\FUN((\infty\Cat_{/\T'})^\op,\infty\Cat_{/ \rS'}^{\mathrm{loc},\mE}),
}
\end{xy}
\end{equation*}
where the left vertical functor is induced by pullback along the functor $\rS' \to \rS, \T' \to \T$ and the right functor is induced by pullback along the functor $\rS' \to \rS$ and the canonical forgetful functor $\infty\Cat_{/\T'} \to \infty\Cat_{/\T}.$

In particular, for every $\mC \in (\infty\Cat^{\loc,\mE}_{/\rS})_{/\rS \times \T}$ and $\s \in \rS$ the lax functor $\nu(\mC): \rS \to \infty\Cat_{/\T}$ sends any $\s \in \rS$ to the fiber $\mC_\s \to \T.$
    
\end{remark}

\begin{remark}\label{functola}

For $\mE=\Fun([1],\rS)$ the 2-functor $\nu$
identifies with the canonical equivalence
$$ (\infty\Cat^{\cocart}_{/\rS})_{/\rS \times \T} \simeq \Fun(\rS,\infty\Cat_{/\T})\subset\FUN((\infty\Cat_{/\T})^\op,\infty\Cat_{/ \rS}^{\cocart}) \simeq \Fun(\rS,\FUN(\infty\Cat_{/\T})^\op,\infty\Cat))$$
induced by the 2-Yoneda-embedding 
$$\iota: \infty\Cat_{/\T} \to \FUN((\infty\Cat_{/\T})^\op,\infty\Cat).$$
    
\end{remark}

\begin{notation}

Let $\rS, \T$ be $\infty$-categories, $\mE \subset \Fun([1],\rS)$ a full subcategory and
$\mC \to \rS \times \T$ a map of locally cocartesian fibrations over $\rS$.
Let $$\Fun_{\rS \times \T}^{\T,\mE}(\rS \times \T,\mC) \subset \Fun_{\rS \times \T}^\T(\rS \times \T,\mC)$$
be the full subcategory whose fiber over $\rt \in \T$ is the full subcategory of
$$\Fun^\mE_{\rS}(\rS,\mC_\rt) \subset \Fun_{\rS \times \T}^\T(\rS \times \T,\mC)_\rt \simeq \Fun_{\rS}(\rS,\mC_\rt) $$
spanned by the sections of $\mC_\rt \to \rS$ sending morphisms of 
$\mE$ to morphisms cocartesian over $\rS.$

\end{notation}

\begin{proposition}\label{simpli} Let $\rS, \T$ be $\infty$-categories, $\mE \subset \Fun([1],\rS)$ a full subcategory and
$\mC \to \rS \times \T$ a map of locally cocartesian fibrations over $\rS$.
The $\mE$-lax limit of the functor $\nu(\mC)$ is $$\Fun_{\rS \times \T}^{\T,\mE}(\rS \times \T,\mC) \to \T.$$

\end{proposition}

\begin{proof}
By definition the $\mE$-lax limit of the functor $\nu(\mC)$ is reprented by the functor $$(\mB \to \T) \mapsto \Fun^\mE_\rS(\rS,\Fun_{\rS \times \T}^\rS(\rS \times \mB,\mC)) \simeq \Fun^{\mE\times \T}_{\rS \times \T}(\rS \times \mB,\mC) \simeq \Fun_\T(\mB,\Fun_{\rS \times \T}^{\T,\mE}(\rS \times \T,\mC)).$$

\end{proof}

\begin{proposition}\label{koy}

Let $\rS,\T$ be $\infty$-categories and $\mF \subset \Fun([1],\T)$ a full subcategory and $\phi: \mC \to \rS \times \T$ a map of locally cocartesian fibrations over $\rS$.

\begin{enumerate}
\item The functor $\phi: \mC \to \rS \times \T$ is a map of (locally) cocartesian fibrations relative to $\mF$ if the lax functor $\nu(\mC): \rS \to \infty\Cat_{/\T}$ lands in the subcategory of (locally) cocartesian fibrations relative to $\mF$ and maps of such.

\vspace{1mm}

\item The functor $\phi: \mC \to \rS \times \T$ is a map of (locally) cartesian fibrations relative to $\mF$ if the lax functor $\nu(\mC): \rS \to \infty\Cat_{/\T}$ lands in the full subcategory of (locally) cartesian fibrations relative to $\mF$.

\end{enumerate}	

\end{proposition}

\begin{proof}
This is \cref{swi}.

\end{proof}

\begin{proposition}

Let $\rS,\T$ be $\infty$-categories and $\mE \subset \Fun([1],\rS), \mF \subset \Fun([1],\T)$ full subcategories and $\phi : \mC \to \rS \times \T$ a map of locally cocartesian fibrations over $\rS$.

\begin{enumerate}
\item If the lax functor $\nu(\mC): \rS \to \infty\Cat_{/\T}$ lands in the subcategory of (locally) cocartesian fibrations relative to $\mF$ and maps of such, the lax functor $\nu(\phi): \rS \to \infty\Cat_{/\T}$ admits an $\mE$-lax limit in the subcategory of (locally) cocartesian fibrations relative to $\mF$ and maps of such.

\vspace{1mm}

\item If the lax functor $\nu(\mC): \rS \to \infty\Cat_{/\T}$ lands in the full subcategory of (locally) cartesian fibrations relative to $\mF$ and sends morphisms of $\mE$ to maps of (locally) cartesian fibrations relative to $\mF$,
the lax functor $\nu(\phi): \rS \to \infty\Cat_{/\T}$ admits an $\mE$-lax limit in the full subcategory of (locally) cartesian fibrations relative to $\mF$.
\end{enumerate}	

\end{proposition}

\begin{proof}

(1): By \cref{koy} the functor $\phi: \mC \to \rS \times \T$ is a map of (locally) cocartesian fibrations relative to $\mF$. By 
\cref{simpli} the $\mE$-lax limit of $\nu(\phi)$ is the functor
$$\Fun_{\rS \times \T}^{\T,\mE}(\rS \times \T,\mC) \to \T.$$

By \cite[Theorem B.4.2.]{lurie.higheralgebra} the functor $$\Fun_{\rS \times \T}^{\T}(\rS \times \T,\mC) \to \T$$
is a (locally) cocartesian fibrations relative to $\mF$.
We will prove that the restriction $$\Fun_{\rS \times \T}^{\T,\mE}(\rS \times \T,\mC) \to \T.$$
is also a (locally) cocartesian fibrations relative to $\mF$. This holds if the fiber transport of $\Fun_{\rS \times \T}^{\T}(\rS \times \T,\mC) \to \T$ along any morphism $\rt \to \rt'$ in $\mF$ preserves the full subcategory $\Fun_{\rS \times \T}^{\T,\mE}(\rS \times \T,\mC) \to \T.$

This fiber transport identifies with the functor
$$ \Fun_{\rS}(\rS,\mC_\rt) \to \Fun_{\rS}(\rS,\mC_{\rt'})$$
induced by the functor $\mC_\rt \to \mC_{\rt'}$. The latter restricts to a functor $$ \Fun^\mE_{\rS}(\rS,\mC_\rt) \to \Fun^\mE_{\rS}(\rS,\mC_{\rt'})$$
if the functor $\mC_\rt \to \mC_{\rt'}$ is a map of cocartesian fibrations relative to $\mE.$
This is equivalent to say that for every morphism of $\mE$ the fiber transports of $\phi: \mC \to \rS \times \T$ preserve cocartesian lifts of morphisms of $\mF$, which holds by assumption.

(2): By \cref{koy} the functor $\phi: \mC \to \rS \times \T$ is a map of (locally) cartesian fibrations relative to $\mF$. By 
\cref{simpli} the $\mE$-lax limit of $\nu(\phi)$ is the functor
$$\Fun_{\rS \times \T}^{\T,\mE}(\rS \times \T,\mC) \to \T.$$

By \cite[Theorem B.4.2.]{lurie.higheralgebra} the functor $$\Fun_{\rS \times \T}^{\T}(\rS \times \T,\mC) \to \T$$
is a (locally) cartesian fibrations relative to $\mF$.
We will prove that the restriction $$\Fun_{\rS \times \T}^{\T,\mE}(\rS \times \T,\mC) \to \T.$$
is also a (locally) cartesian fibrations relative to $\mF$. This holds if the fiber transport of $\Fun_{\rS \times \T}^{\T}(\rS \times \T,\mC) \to \T$ along any morphism $\rt \to \rt'$ in $\mF$ preserves the full subcategory $\Fun_{\rS \times \T}^{\T,\mE}(\rS \times \T,\mC) \to \T.$

This fiber transport identifies with the functor
$$ \Fun_{\rS}(\rS,\mC_{\rt'}) \to \Fun_{\rS}(\rS,\mC_{\rt})$$
induced by the functor $\mC_{\rt'} \to \mC_{\rt}$. The latter restricts to a functor $$ \Fun^\mE_{\rS}(\rS,\mC_{\rt'}) \to \Fun^\mE_{\rS}(\rS,\mC_{\rt})$$
if the functor $\mC_{\rt'} \to \mC_{\rt}$ is a map of cocartesian fibrations relative to $\mE.$
This is equivalent to say that for every morphism of $\mE$ the fiber transports of $\phi: \mC \to \rS \times \T$ preserve cartesian lifts of morphisms of $\mF$, which holds by assumption.

\end{proof}

% \begin{notation}

% Let $\rS$ be an $\infty$-category. Let $$ \infty\Cat^{\cart, \lax}_{/\rS} \subset \infty\Cat_{/\rS}$$ the full subcategory of cartesian fibrations over $\rS.$

% \end{notation}

% % \begin{theorem}

% % Let $\rS, \T$ be $\infty$-categories and $\mE \subset \Fun([1],\rS)$ a full subcategory.
% % Every $\mE$-lax functor $ \rS \to \infty\Cat^{\cart, \lax}_{/\T}$ admits an $\mE$-lax limit.

% % \end{theorem}

% \begin{corollary}

% Let $\rS$ be $\infty$-categories and $\mE \subset \Fun([1],\rS)$ a full subcategory.
% The $(\infty,2)$-category $$ \infty\Cat^{\cart, \lax}_{/\rS}$$ admits $\mE$-lax limits.

% \end{corollary}

Next we apply the theory of lax limits in $\infty\Cat_{/\rS}$ to study Eilenberg-Moore objects in $\infty\Cat_{/\rS}$.

\begin{notation}

Let $\X \to \rS$ be a functor.
The endomorphism left action of $\Fun_\rS(\X,\X)$ on $\X \to \rS$
classifies a $\rS$-family $$ \mC^\circledast \to  \Fun_\rS(\X,\X)^\ot \times \rS$$ of $\infty$-categories left tensored over $\Fun_\rS(\X,\X)$.

%We apply Notation \ref{mmod4} to this situation to obtain for every monad $\T$ on $\mC \to \rS$ a functor $\LMod^\rS_\T(\mC) \to \rS.$

\end{notation}

\begin{notation}

Let $\rS$ be an $\infty$-category, $\X \to \rS$ a functor and $\T$ a monad on $\X \to \rS$ in $\infty\Cat_{/\rS}$. 

Let 
$$ \LMod_\T^\rS(\X) \subset \Fun^\rS_{\Delta^\op \times \rS}(\Delta^\op \times \rS,\Delta^\op \times_{\Fun_\rS(\X,\X)^\ot}\X^\circledast) $$
be the full subcategory whose fiber over $\s \in \rS$ is
the full subcategory
$$ \LMod_\T(\X_\s) \subset \Fun^\rS_{\Delta^\op \times \rS}(\Delta^\op \times \rS,\Delta^\op \times_{\Fun_\rS(\X,\X)^\ot}\X^\circledast)_\s \simeq \Fun_{\Delta^\op}(\Delta^\op,\Delta^\op \times_{\Fun_\rS(\X,\X)^\ot}\X_\s^\circledast) $$
spanned by the $\T$-modules in $\X_\s.$

\end{notation}

\begin{notation}
Let $\X \to \rS$ be a functor and $\R$ a comonad on $\X \to \rS$ corresponding to a monad on $\X^\op \to \rS^\op$.
Let $$\coLMod_\R^\rS(\X):=\LMod_{\R^\op}^{\rS^\op}(\X^\op)^\op.$$

\end{notation}

%\begin{remark}A bitensored $\infty$-category $ \mM^\circledast \to \mV^\ot \times \mW^\ot $ is in particular a $\mW^\ot$-family of $\infty$-categories left tensored over $\mV$.So by Remark \ref{Fama} there is a canonical monoidal functor$ \mV \to \Fun_{\mW^\ot}(\mM_{[0]}^\circledast,\mM_{[0]}^\circledast)$ that induces a functor $\mV \to \mW\mathrm{-}\LinFun(\mM,\mM)$ that sends an associative algebra $\A$ to the $\mW$-linear monad $\A \ot (-)$.Consequently, for every associative algebra $\A$ in $\mV$ there is a canonical equivalence $\LMod_\A^{\mW^\ot}(\mM)\simeq \LMod_{\A\ot(-)}^{\mW^\ot}(\mM) $ of$\infty$-categories left tensored over $\mW$. 

\begin{remark}
Let $\X \to \rS$ be a functor.
The 2-functor $$\Cat^\cocart_{ \infty / \rS} \simeq \Cat^\cart_{ \infty / \rS^\op}, \X \to \rS \mapsto (\X^\rev)^\op \to \rS^\op$$ sends the endomorphism left action of $\Fun^\cocart_\rS(\X,\X)$ on $\X \to \rS$
to a left action of $\Fun^\cocart_\rS(\X,\X)$ on $(\X^\rev)^\op \to \rS^\op$
that is the pullback of the endomorphism left action of $\Fun^\cocart_\rS(\X,\X)$ on $(\X^\rev)^\op \to \rS^\op$ along a canonical monoidal equivalence
$$ \Fun^\cocart_\rS(\X,\X) \simeq \Fun^\cart_{\rS^\op}((\X^\rev)^\op,(\X^\rev)^\op).$$
Thus for every monad $\T$ on $\X \to \rS$ preserving cocartesian morphisms over $\rS$ there is a canonical equivalence 
$$\LMod^\rS_\T(\X) \simeq \coLMod^\rS_{\T^\rev}(\X^\rev).$$
	
\end{remark}

\begin{proposition}\label{classiso}
Let $\rS$ be an $\infty$-category, $\X \to \rS$ a functor and $\T$ a monad on $\X \to \rS$ in $\infty\Cat_{/\rS}$. 
Let $$ \phi: \X^\circledast \to \Fun_\rS(\X,\X)^\ot \times \rS$$ be the $\infty$-category left tensored over $\Fun_\rS(\X,\X)$ classifying the endomorphism left action of $ \Fun_\rS(\X,\X)$ on $\X \to \rS$ in $\infty\Cat_{/\rS}$. 
Let $$\phi': \Delta^\op \times_{\Fun_\rS(\X,\X)^\ot}\X^\circledast \to \Delta^\op \times \rS$$ be the pullback of $\phi$ along the lax monoidal functor $\Delta^\op \to \Fun_\rS(\X,\X)^\ot$ classifying $\T.$
The image of $\phi'$ under $$ \nu: (\infty\Cat^{\loc}_{/\Delta^\op})_{/\Delta^\op \times \rS} \to \FUN^{ \mathrm{coten}}((\infty\Cat_{/\rS})^\op,\infty\Cat_{/ \Delta^\op}^{\mathrm{loc}}) $$
is classified by the lax functor
$\theta_\T: \Delta^\op \to \infty\Cat_{/\rS}.$

\end{proposition}

\begin{proof}
The lax functor
$\theta_\T: \Delta^\op \to \infty\Cat_{/\rS} $
is classified by the image of $\bar{\iota}(\X)$ under the 2-functor 
$$ \FUN((\infty\Cat_{/\rS})^\op,\Fun_\rS(\X,\X)\mathrm{-}\LMod) \xrightarrow{ } \FUN((\infty\Cat_{/\rS})^\op,\infty\Cat_{/ \Delta^\op}^{\mathrm{loc},\max}),$$
where the last 2-functor restricts along the lax monoidal functor $\Delta^\op \to \Fun_\rS(\X,\X)^\ot$ classifying $\T.$

By \cref{functola} the 2-functor $$ \nu: (\infty\Cat^{\cocart}_{/\Delta^\op})_{/\Delta^\op \times \rS} \to \FUN((\infty\Cat_{/\rS})^\op,\infty\Cat_{/ \Delta^\op}^{\mathrm{cocart}}) $$
sends the morphism $\phi$ to a morphism
$$\nu(\X^\circledast) \to \nu(\Fun_\rS(\X,\X)^\ot \times \rS) \simeq \nu'(\Fun_\rS(\X,\X)^\ot),$$
where $\nu'$ is the 2-functor
$$ \infty\Cat_{/ \Delta^\op}^{\mathrm{cocart}}\simeq \Fun^{\mathrm{coten}}(\infty\Cat^\op,\infty\Cat_{/ \Delta^\op}^{\mathrm{cocart}}) \to \FUN((\infty\Cat_{/\rS})^\op,\infty\Cat_{/ \Delta^\op}^{\mathrm{cocart}}).$$

The 2-functor $\nu'$ identifies with the 2-functor
$$ \Fun(\Delta^\op,\infty\Cat) \simeq \Fun(\Delta^\op,\FUN^{\mathrm{coten}} (\infty\Cat^\op,\infty\Cat)) \to \Fun(\Delta^\op,\FUN^{ }((\infty\Cat_{/\rS})^\op,\infty\Cat))$$
induced by the 2-functor $$\infty\Cat \simeq \FUN^{\mathrm{coten}}(\infty\Cat^\op,\infty\Cat) \to \FUN^{ }((\infty\Cat_{/\rS})^\op, \infty\Cat),$$
which factors as 
$$\infty\Cat \xrightarrow{(-)\times \rS} \infty\Cat_{/\rS} \xrightarrow{\iota} \FUN((\infty\Cat_{/\rS})^\op, \infty\Cat).$$

% Let $$\delta: \infty\Cat^\cocart_{/\Delta^\op}\to \FUN((\infty\Cat_{/\rS})^\op,\infty\Cat^\cocart_{/\Delta^\op})$$ be the diagonal functor.

Let $$\delta: \infty\Cat\to \FUN((\infty\Cat_{/\rS})^\op,\infty\Cat)$$ be the diagonal 2-functor and $\iota: (\infty\Cat_{/\rS}\to \FUN((\infty\Cat_{/\rS})^\op,\infty\Cat)$
the 2-Yoneda-embedding.

There is a canonical natural 2-transformation
$$\delta \to \iota\circ ((-)\times \rS) $$
whose component as any $\mA \in \infty\Cat$ and $\mB \in \infty\Cat_{/\rS}$ is the diagonal functor
$$ \mA \to \Fun_\rS(\mB,\mA \times \rS) \simeq \Fun(\mB,\mA). $$

We apply \cref{heopo} to the latter natural transformation.
We find that the functor 
\begin{equation}\label{llkkmo}
\LMod_{\Fun_\rS(\X,\X)}(\infty\Cat_{/\rS}) \to \LMod_{\Fun_\rS(\X,\X)}(\FUN((\infty\Cat_{/\rS})^\op,\infty\Cat))\end{equation}
induced by the 2-Yoneda-embedding factors as
$$\LMod_{\Fun_\rS(\X,\X)}(\infty\Cat_{/\rS}) \to \LMod_{ \iota(\rS \times\Fun_\rS(\X,\X))}(\FUN((\infty\Cat_{/\rS})^\op,\infty\Cat)) \to$$$$ \LMod_{\delta(\Fun_\rS(\X,\X))}(\FUN((\infty\Cat_{/\rS})^\op,\infty\Cat))\simeq$$$$ \LMod_{\Fun_\rS(\X,\X)}(\FUN((\infty\Cat_{/\rS})^\op,\infty\Cat)).$$

Therefore the image of
the endomorphism left action of $\Fun_\rS(\X,\X) $ on $\X \to \rS$ under the functor (\ref{llkkmo})
% $$\LMod_{\Fun_\rS(\X,\X)}(\infty\Cat_{/\rS}) \to \LMod_{\Fun_\rS(\X,\X)}(\FUN((\infty\Cat_{/\rS})^\op,\infty\Cat)) \simeq$$$$
% \FUN((\infty\Cat_{/\rS})^\op,\Fun_\rS(\X,\X)\mathrm{-}\Mod) \subset $$$$
% \FUN((\infty\Cat_{/\rS})^\op,\infty\Cat_{/\Fun_\rS(\X,\X)^\ot}) \simeq$$$$ \FUN((\infty\Cat_{/\rS})^\op,\infty\Cat_{/\Delta^\op})_{/\delta(\Fun_\rS(\X,\X)^\ot)}$$
classifies the pullback of $$\nu(\phi): \nu(\X^\circledast) \to \nu(\Fun_\rS(\X,\X)^\ot \times \rS) \simeq \nu'(\Fun_\rS(\X,\X)^\ot) $$
along the canonical map
$$ \delta(\Fun_\rS(\X,\X)^\ot) \to \iota_*(\Fun_\rS(\X,\X)^\ot \times \rS)\simeq \nu'(\Fun_\rS(\X,\X)^\ot).$$ 

Hence by naturality the image of the endomorphism left action of $\Fun_\rS(\X,\X) $ on $\X \to \rS$ under the 2-functor
$$\Fun_\rS(\X,\X)\mathrm{-}\LMod \to \LMod_{\Fun_\rS(\X,\X)}(\FUN((\infty\Cat_{/\rS})^\op,\infty\Cat)) \simeq$$$$
\FUN((\infty\Cat_{/\rS})^\op,\Fun_\rS(\X,\X)\mathrm{-}\Mod) \subset $$$$
\FUN((\infty\Cat_{/\rS})^\op,\infty\Cat_{/\Fun_\rS(\X,\X)^\ot}) \to $$$$
\FUN((\infty\Cat_{/\rS})^\op,\infty\Cat_{/\Delta^\op}), $$
where the last functor restricts along the lax monoidal functor $\Delta^\op \to \Fun_\rS(\X,\X)^\ot$ classifying $\T,$
classifies the pullback of $$\nu(\phi'): \nu(\X^\circledast) \to \nu(\Delta^\op \times \rS) \simeq * $$
along the canonical map
$ \delta(\Delta^\op)\simeq * \to *,$ and so identifies with $\nu(\phi').$

\end{proof}

\begin{theorem}

Let $\rS$ be an $\infty$-category, $\X \to \rS$ a functor and $\T$ a monad on $\X \to \rS$ in $\infty\Cat_{/\rS}$. 
The Eilenberg-Moore object of $\T$ in $\infty\Cat_{/\rS}$ is the functor
$$ \LMod_\T^\rS(\X) \to \rS.$$

\end{theorem}

\begin{proof}

Let $ \phi: \X^\circledast \to \Fun_\rS(\X,\X)^\ot \times \rS$ be the $\infty$-category left tensored over $\Fun_\rS(\X,\X)$ classifying the endomorphism left action of $ \Fun_\rS(\X,\X)$ on $\X \to \rS$ in $\infty\Cat_{/\rS}$.
Let $$\phi': \Delta^\op \times_{\Fun_\rS(\X,\X)^\ot}\X^\circledast \to \Delta^\op \times \rS$$ be the pullback of $\phi$ along the lax monoidal functor $\Delta^\op \to \Fun_\rS(\X,\X)^\ot$ classifying $\T.$
By \cref{classiso} the lax functor $\theta_\T$ is classified by the
image of $\phi'$ under $$ \nu: (\infty\Cat^{\loc}_{/\Delta^\op})_{/\Delta^\op \times \rS} \to \FUN^{ \mathrm{coten}}((\infty\Cat_{/\rS})^\op,\infty\Cat_{/ \Delta^\op}^{\mathrm{loc}}).$$

Therefore by \cref{simpli} the $\mathrm{max}$-lax limit of $\theta_\T \simeq \nu(\phi')$ is $$\LMod_\T^\rS(\X) =\Fun^{\rS,{\mathrm{max}}}_{\Delta^\op \times \rS}(\Delta^\op \times \rS,\Delta^\op \times_{\Fun_\rS(\X,\X)^\ot}\X^\circledast).$$

By \cref{stru} the $\mathrm{max}$-lax limit of $\theta_\T$ is the Eilenberg-Moore object of $\T$ in $\infty\Cat_{/\rS}$.

\end{proof}

\begin{lemma}\label{lemru}

Let $\rS$ be an $\infty$-category, $\X \to \rS$ a functor and $\T$ a monad on $\X \to \rS$ in $\infty\Cat_{/\rS}.$
 The pullback $$ {\Delta^\op} \times_{\Fun_\rS(\X,\X)^\ot} \X^\circledast \to \rS \times {\Delta^\op} $$ of $\phi$ along $\T$ is a map of locally cocartesian fibrations over $\Delta^\op$.

\end{lemma}

\begin{proof}

We consider the $\infty$-category $\phi: \X^\circledast \to \Fun_\rS(\X,\X)^\ot \times \rS$ left tensored over $\Fun_\rS(\X,\X)$ classifying the endomorphism $\Fun_\rS(\X,\X)$-action on $\X \to \rS$ in $\infty\Cat_{/\rS}.$
The functor $\phi: \X^\circledast \to \Fun_\rS(\X,\X)^\ot \times \rS$ is a map of cocartesian fibrations over $\Delta^\op$.
The latter induces on the fiber over every $[\n]\in \Delta^\op$ the canonical functor $\X \times \Fun_\rS(\X,\X)^{\times \n} \to \Fun_\rS(\X,\X)^{\times \n} \times \rS, $ which is a map of locally cocartesian fibrations over
$\Fun_\rS(\X,\X)^{\times \n}.$
By \cite[Proposition 2.4.2.11.]{lurie.HTT} this implies that $\phi$ is a map of locally cocartesian fibrations over $\Fun_\rS(\X,\X)^\ot.$
Thus the pullback $ {\Delta^\op} \times_{\Fun_\rS(\X,\X)^\ot} \X^\circledast \to \rS \times {\Delta^\op} $ of $\phi$ along $\T$ is a map of locally cocartesian fibrations over $\Delta^\op$.

\end{proof}

\begin{proposition}\label{patter3}
Let $\rS$ be an $\infty$-category, $\mE \subset \Fun([1],\rS)$ a full subcategory, $\X \to \rS$ a functor and $\T$ a monad on $\X \to \rS$ in $\infty\Cat_{/\rS}.$
 	
\begin{enumerate}
		
\item If the functor $\X \to \rS$ is a cocartesian fibration relative to $\mE$ and for every $\n \geq 0$ the functor $\T^{\circ \n} : \X \to \X$ is a map of locally cocartesian fibrations relative to $\mE$, the
pullback $$ {\Delta^\op} \times_{\Fun_\rS(\X,\X)^\ot} \X^\circledast \to \rS \times {\Delta^\op} $$ along $\T$ is a map of cocartesian fibrations relative to $\mE$ whose fiber transports preserve cocartesian lifts of inert morphisms.

\vspace{1mm}
		
\item If the functor $\X \to \rS$ is a cartesian fibration relative to $\mE$, the
pullback $$ {\Delta^\op} \times_{\Fun_\rS(\X,\X)^\ot} \X^\circledast \to \rS \times {\Delta^\op} $$ along $\T$ is a map of cartesian fibrations relative to $\mE$ whose fiber transports preserve cocartesian lifts of inert morphisms.

\end{enumerate}
	
\end{proposition}

\begin{proof}

This follows from \cref{lemru}, \cref{swi}.

\end{proof}

\begin{corollary}\label{patter2}\label{corso}
Let $\rS$ be an $\infty$-category, $\mE \subset \Fun([1],\rS)$ a full subcategory and $\X \to \rS$ a functor and $\T$ a monad on $\X \to \rS$ in $\infty\Cat_{/\rS}.$
 	
\begin{enumerate}
		
\item If the functor $\X \to \rS$ is a cocartesian fibration relative to $\mE$ and for every $\n \geq 0$ the functor $\T^{\circ \n} : \X \to \X$ is a map of cocartesian fibrations relative to $\mE$, then $\LMod_\T^\rS(\X)\to \rS$ is a map of cocartesian fibration relative to $\mE$ and the functor $\LMod_\T^\rS(\X)\to \X$ is a map of cocartesian fibrations relative to $\mE$.

\vspace{1mm}
		
\item If the functor $\X \to \rS$ is a cartesian fibration relative to $\mE$, then $\LMod_\T^\rS(\X) \to \rS$ is a map of cartesian fibration relative to $\mE$ and the functor $\LMod_\T^\rS(\X)\to \X$ is a map of cartesian fibrations relative to $\mE$.

\end{enumerate}
	
\end{corollary}

\begin{remark}

Applying the opposite $\infty$-category involution
$$(-)^\op: \infty\Cat^\co \simeq \infty\Cat$$
statement (2) of \cref{patter3} is equivalent to say that for every $\infty$-category $\rS$ and full subcategory $\mE \subset \Fun([1],\rS)$,
cocartesian fibration $\X \to \rS$ relative to $\mE$
and monad $\T$ on $\X \to \rS$ in $\infty\Cat_{/\rS}$
the functor $\co\LMod_\T^\rS(\X) \to \rS$ is a map of cocartesian fibration relative to $\mE$ and the functor $\co\LMod_\T^\rS(\X)\to \X$ is a map of cocartesian fibrations relative to $\mE$.

\end{remark}

\begin{theorem}\label{tgvwxlkm}\label{patty}

Let $\rS$ be an $\infty$-category and $\mE \subset \Fun([1], \rS)$ a full subcategory.

\begin{enumerate}
\item For every monad $\T$ on $\X \to \rS$ in $\infty\Cat_{/\rS}^\mE$ the forgetful functor $$ \LMod^{\rS}_\T(\X) \to \X $$
is an Eilenberg-Moore object for $\T $ in $\infty\Cat_{/\rS}^\mE$.

\vspace{1mm}

\item For every comonad $\T$ on $\X \to \rS$ in $\infty\Cat_{/\rS}^\mE$ the forgetful functor $$ \co\LMod^{\rS}_\T(\X) \to \X $$
is a coEilenberg-Moore object for $\T $ in $\infty\Cat_{/\rS}^\mE$.

\vspace{1mm}

\item Let $\X \to \rS$ be a cocartesian fibration relative to $\mE$.
For every comonad $\T$ on $\X \to \rS$ in $\infty\Cat_{/\rS}$ the forgetful functor $$ \co\LMod^{\rS}_\T(\X) \to \X $$
is an Eilenberg-Moore object for $\T $ in the full subcategory of
$\infty\Cat_{/\rS}$ of cocartesian fibrations relative to $\mE.$
Moreover the functor $$ \co\LMod^{\rS}_\T(\X) \to \X $$ is a map
of cocartesian fibrations relative to $\mE.$

\end{enumerate}

\end{theorem}

\begin{corollary}\label{zhgbfdc}
	
Let $\rS$ be an $\infty$-category and $\mE \subset \Fun([1], \rS)$ a full subcategory.
\begin{enumerate}
\item The full subcategory of ${\infty\Cat}_{/ \rS} $ of cocartesian fibrations relative to $\mE$ is closed under coEilenberg-Moore objects.
Moreover coEilenberg-Moore objects are maps of cocartesian fibrations relative to $\mE$. 
	
\item Applying the opposite involution, the full subcategory of ${\infty\Cat}_{/ \rS} $ of cartesian fibrations relative to $\mE$ is closed under Eilenberg-Moore objects.
Moreover Eilenberg-Moore objects are maps of cartesian fibrations relative to $\mE$. 

\end{enumerate}
\end{corollary}

% \begin{corollary}\label{zhgbfdc}
	
% Let $\rS$ be an $\infty$-category and $\mE \subset \Fun([1], \rS), \mT \subset \Fun([2],\rS)$ full subcategories.
% The full subcategory of ${\infty\Cat}_{/ \rS} $ of cartesian fibrations relative to $\mE, \mT$ admits Eilenberg-Moore objects, which are preserved by the embedding to ${\infty\Cat}_{/ \rS}. $
	
% \end{corollary}

\subsection{Construction of Kleisli objects}

% \begin{notation}Let $\rS$ be an $\infty$-category and $\mE \subset \Fun([1], \rS), \mT \subset \Fun([2],\rS)$ full subcategories.
% Let $$\infty\Cat_{/ \rS}^{\mE,\mT} \subset \infty\Cat_{/ \rS}$$ be the subcategory of cocartesian fibrations relative to
% $\mE,\mT$ and maps of such.

% \end{notation}

Next we study Kleisli objects in $\infty\Cat^\mE_{/\rS}$.

\begin{theorem}\label{coKlei}

Let $\rS$ be an $\infty$-category, $\mE \subset \Fun([1],\rS)$ a full subcategory and $ \mC \to \rS$ a cocartesian fibration relative to $\mE.$

\begin{enumerate}
\item For every monad $\T$ on $\mC \to \rS$ in $\infty\Cat^{\mE}_{/\rS}$ let $ \LMod^{\rS}_\T(\mC)' \subset \LMod^{\rS}_\T(\mC)$ be the essential image of the free functor.
The functor $\mC \to \LMod^{\rS}_\T(\mC)'$ is a Kleisli-object in $\infty\Cat^{\mE}_{/\rS}.$

\vspace{1mm}

\item For every comonad $\R $ on $\mC \to \rS$ in $\infty\Cat^{\mE}_{/\rS}$ let $ \coLMod^{\rS}_\R(\mC)' \subset \coLMod^{\rS}_\R(\mC)$ be the essential image of the cofree functor. The functor $\mC \to \coLMod^{\rS}_\R(\mC)'$ is a coKleisli-object in $\infty\Cat^{\mE}_{/\rS}.$
\end{enumerate}

\end{theorem}

\begin{proof}
% The essential image of a map of $\mathfrak{P}$-operads is again
% $\mathfrak{P}$-operadic for every discrete algebraic pattern $\mathfrak{P}$.
% Thus we can reduce to proving the statement for algebraic pattern $ \mathfrak{P}$ on $\rS$ such that $\rS= \rS^\circ.$

By \cref{tgvwxlkm} the conservative functors $\LMod^{\rS}_\T(\mC) \to \mC $ and $ \coLMod^{\rS}_\R(\mC) \to \mC$ are maps of cocartesian fibrations over $\mE$.
The left adjoint relative to $\rS$ of the functor $\LMod^{\rS}_\T(\mC) \to \mC $
is a map of cocartesian fibrations over $\rS$.
The right adjoint relative to $\rS$ of the functor $\coLMod^{\rS}_\R(\mC) \to \mC $ is a map of cocartesian fibrations relative to $\mE$ since the composition of the right adjoint with the functor $\coLMod^{\rS}_\R(\mC) \to \mC $ is 
$\R$, which is a map of cocartesian fibrations relative to $\mE$.

So the free functor $\mC \to \LMod^{\rS}_\T(\mC)$ and the cofree functor
$\mC \to \coLMod^{\rS}_\R(\mC)$ are maps of cocartesian fibrations relative to $\mE$ and so induce maps $\mC \to \LMod^{\rS}_\T(\mC)'$ and
$\mC \to \coLMod^{\rS}_\R(\mC)'$ of cocartesian fibrations relative to $\mE$.
It remains to see that the functor $\mC \to \LMod^{\rS}_\T(\mC)'$ is a Kleisli-object in $\infty\Cat^{\mE}_{/\rS}$ and the functor $\mC \to \coLMod^{\rS}_\R(\mC)'$ is a coKleisli-object in $\infty\Cat^{\mE}_{/\rS}.$

We prove the first statement. The second statement is dual to the first one. 
The canonical equivalence
$$ \Fun^\mE_\rS(-,\LMod^\rS_\T(\mC)) \to \LMod_\T(\Fun^\mE_\rS(-,\mC))$$
of 2-functors $(\infty\Cat^\mE_{/ \rS})^\op \to \infty\Cat$
sends the endomorphism left $\Fun^\mE_\rS(\mC,\mC)$-action on $\mC \to \rS$ to a $\Fun^\mE_\rS(\mC,\mC)$-linear equivalence
$\Fun^\mE_\rS(\mC,\LMod^\rS_\T(\mC)) \to \LMod_\T(\Fun^\mE_\rS(\mC,\mC)).$
Under the latter equivalence the canonical $\T,\T$-biaction on $\T$ with respect to the canonical $\Fun^\mE_\rS(\mC,\mC), \Fun^\mE_\rS(\mC,\mC)$-biaction on $\Fun^\mE_\rS(\mC,\mC)$
corresponds to a right $\T$-action on the free functor $\mC \to \LMod^\rS_\T(\mC)$ over $\rS.$
The embedding $\LMod^\rS_\T(\mC)' \subset \LMod^\rS_\T(\mC)$ over $\rS$ induces an embedding $$\RMod_\T(\Fun^\mE_\rS(\mC,\LMod^\rS_\T(\mC)')) \subset \RMod_\T(\Fun^\mE_\rS(\mC,\LMod^\rS_\T(\mC)))$$
so that the right $\T$-action on the free functor $\mC \to \LMod^\rS_\T(\mC)$
induces a right $\T$-action on the functor $\mC \to \LMod^\rS_\T(\mC)'$ over $\rS.$

It remains to see that for every functor $\mD \to \rS$ the induced functor
$$ \Fun^\mE_\rS(\LMod^\rS_\T(\mC)',\mD) \to \RMod_\T(\Fun^\mE_{\rS}(\mC,\mD)) $$
is an equivalence.
The latter functor is a functor over $\Fun^\mE_\rS(\mC,\mD).$
So by \cite[Lemma 4.7.3.13.]{lurie.higheralgebra} it is enough to see that the functor $$\theta: \Fun^\mE_\rS(\LMod^\rS_\T(\mC)',\mD) \to \Fun^\mE_{\rS}(\mC,\mD)$$ is
conservative, admits a left adjoint $\mF$ such that the induced map
$(-)\circ \T \to \theta \circ \mF $ is an equivalence, and that every $\theta$-split simplicial object of $\Fun^\mE_\rS(\LMod^\rS_\T(\mC)',\mD)$ admits a colimit that is preserved by $\theta.$ The functor $\theta$ is conservative because the functor $\mC \to \LMod^\rS_\T(\mC)'$ is essentially surjective.
The functor $\theta$ is right adjoint to the functor
$$\Fun^\mE_{\rS}(\mC,\mD) \to \Fun^\mE_\rS(\LMod^\rS_\T(\mC)',\mD)$$ induced by the functor $\LMod^\rS_\T(\mC)' \to \mC$ so that the induced map
$(-)\circ \T \to \theta \circ \mF $ is an equivalence.

We prove that every $\theta$-split simplicial object of $\Fun^\mE_\rS(\LMod^\rS_\T(\mC)',\mD)$ admits a colimit that is preserved by $\theta$. We first consider the case that $\mE=\rS^\simeq.$
Let $\Env^\rS(\mD) \to \rS$ be the functor $$ \Fun([1], \rS) \times_{\Fun(\{1\}, \rS)} \mD \to \Fun([1], \rS) \to \Fun(\{0\}, \rS),$$ which is a cartesian fibration.
The diagonal embedding $\rS \to \Fun([1], \rS)$ yields an embedding $\mD \to \Env^\rS(\mD)$ over $\rS$.
The functor $\theta$ is the pullback of the similarly defined functor
$$\theta': \Fun_\rS(\LMod^\rS_\T(\mC)',\Env^\rS(\mD)) \to \Fun_{\rS}(\mC,\Env^\rS(\mD)).$$
So it is enough to see that every $\theta'$-split simplicial object of $\Fun_\rS(\LMod^\rS_\T(\mC)',\Env^\rS(\mD))$ admits a colimit preserved by $\theta'.$
Since $\Env^\rS(\mD) \to \rS$ is a cartesian fibration,
a simplicial object $\X$ of $\Fun_\rS(\LMod^\rS_\T(\mC)',\Env^\rS(\mD))$ 
admits a geometric realization if for every $\Z \in \mC$ lying over $\s \in \rS$
the induced functor $$\Delta^\op \to \Fun_\rS(\LMod^\rS_\T(\mC)',\Env^\rS(\mD)) \xrightarrow{\ev_{\T(\Z)}}\Env^\rS(\mD)_\s$$ admits a geometric realization.
But the latter functor splits and so admits a geometric realization
if $\X$ is $\theta'$-split. In this case $\theta'$ preserves this colimit.

Next we consider the case of arbitrary $\mE.$
Let $\X$ be a $\theta$-split simplicial object of $ \Fun^\mE_\rS(\LMod^\rS_\T(\mC)',\mD)$. By what we have proven,
$\X$ admits a geometric realization in $\Fun_\rS(\LMod^\rS_\T(\mC)',\mD)$ and we have to see that the geometric realization belongs to $\Fun^\mE_\rS(\LMod^\rS_\T(\mC)',\mD)$,
i.e. that for every $\A, \B \in \mC$ lying over objects $\s \in \rS$ and $\rt \in \rS$ and cocartesian lift $\f: \T(\A) \to \T(\B) $ in $\LMod^\rS_\T(\mC)'$ of a morphism $\g: \s \to \rt$ of $\mE$ the induced morphism $\X(\T(\A)) \to \X(\T(\B))$ is cocartesian over $\rS.$
We have also proven that geometric realizations in $\Fun_\rS(\LMod^\rS_\T(\mC)',\mD)$ are formed objectwise.
So we have to see that the fiber transport $\g_!: \mD_\s \to \mD_\rt$
preserves the geometric realization of the simplicial object 
$\X(\T(\A))$ in $\mD_\s $.
But $ \X(\T(\A)) \simeq \theta(\X)(\A)$ and $\theta(\X)$ is split so that also
$\theta(\X)(\A)$ is split in $ \mD_\s$ and so preserved by any functor.

\end{proof}

\begin{corollary}\label{zhgbfdc1}
	
Let $\rS$ be an $\infty$-category and $\mE \subset \Fun([1], \rS)$ a full subcategory.
The full subcategory of ${\infty\Cat}_{/ \rS} $ of cartesian fibrations relative to $\mE$ is closed under 
Kleisli objects.
\end{corollary}

\begin{proof}
	
Let $ \mC \to \rS$ be a functor and $\T$ a monad on $\mC \to \rS$ in $\infty\Cat_{/ \rS}$. By \cref{tgvwxlkm} the functor $ \LMod^{\rS}_\T(\mC) \to \mC $ is an Eilenberg-Moore object for $\T $ in ${\infty\Cat}_{/\rS}$
and by \cref{coKlei} the free functor $\mC \to \LMod^{\rS}_\T(\mC)'$ 
to its essential image is a Kleisli-object for $\T$ in ${\infty\Cat}_{/\rS}.$
If $ \mC \to \rS$ is a cartesian fibration relative to $(\mE, \mT),$ by \cref{corso} the functor $ \LMod^{\rS}_\T(\mC) \to \rS $ is a cartesian fibration relative to $(\mE,\mT)$ and $ \LMod^{\rS}_\T(\mC) \to \mC $ and its left adjoint are maps of locally cartesian fibrations relative to $\mE$.
So the restriction $ \LMod^{\rS}_\T(\mC)' \to \rS $ is a cartesian fibrations relative to $(\mE,\mT)$ and the forgetful functor $\LMod^{\rS}_\T(\mC)' \to \mC$ and free functor $\mC \to \LMod^{\rS}_\T(\mC)' $ are maps of locally cartesian fibrations relative to $\mE.$
 
\end{proof}

% \begin{proposition}
% Let $\T$ be an $\infty$-category, $\mV$ a monoidal $\infty$-category, $\A$ an associative algebra in $\mV$ and $\phi: \mM^\circledast \to \mV^\ot \times \T$ be a $\T$-family of left $\mV$-tensored $\infty$-categories.
% Let $$\phi': \Delta^\op \times_{\mV^\ot}\mM^\circledast \to \Delta^\op \times \T$$ be the pullback of $\phi$ along $\A.$
% The image of $\phi'$ under $$ \nu: (\infty\Cat^{\loc}_{/\Delta^\op})_{/\Delta^\op \times \T} \to \FUN^{ \mathrm{coten}}((\infty\Cat_{/\T})^\op,\infty\Cat_{/ \Delta^\op}^{\mathrm{loc}}) $$
% is the functor
% $$ \infty\Cat_{/\T} \to \FUN((\infty\Cat_{/\T})^\op,\mV\mathrm{-}\LMod) \xrightarrow{\LMod_\A)} \FUN((\infty\Cat_{/\T})^\op,\infty\Cat).$$

% \end{proposition}

%\section{The higher algebra of structured monads}

%\section{Structured monads}

\section{A duality between monads and monadic morphisms}\label{Classa}

\subsection{A classification of monads by Eilenberg-Moore objects}\label{Classa}

In this section we classify monadic morphisms by monads in any $(\infty,2)$-category that admits Eilenberg-Moore objects (\cref{uhnggfdaa}) and prove that monadic morphisms are the local objects for a localization
%on the $(\infty,2)$-category of right adjoint morphisms 
(\cref{monan}).

%\subsection{A classification of families of monads}

We fix the following notation:
\begin{notation}Let $\mC$ be an $(\infty,2)$-category.
Let $$ \Fun([1],\mC)^\R \subset \Fun([1],\mC) $$ be the full subcategory of right adjoint morphisms.
Let $$ \Fun([1],\mC)^\mon \subset \Fun([1],\mC) $$ be the full subcategory of monadic morphisms.
For every $\X \in \mC$ let $$(\mC_{/\X})^\mon \subset (\mC_{/\X})^\R \subset \mC_{/\X}$$ the full subcategories of objects over $\X$ whose morphism to $\X$ is monadic, right adjoint, respectively.

\end{notation}

%\subsubsection{Eilenberg-Moore objects among oplax natural transformations}

To classify monads by their Eilenberg-Moore objects we use the following strategy:

\begin{enumerate}
\item We construct for every $\infty$-category $\rS$
an $\infty$-category $\underline{\FUN}^\oplax(\rS,\mC)$ of functors $\rS \to \mC$ and oplax natural transformations (Notation \ref{oplll}) that is enriched in $\infty\Cat_{/ \rS^\op}$ (Notation \ref{opplas}).
\vspace{1mm}
\item We use the enrichment of $\underline{\FUN}^\oplax(\rS,\mC)$ in $\infty\Cat_{/ \rS^\op}$ to construct an $\infty$-category $\Mon(\mC) $ of monads in $\mC$ (\cref{monasol}).
\vspace{1mm}
\item We identify monads and Eilenberg-Moore objects in the underlying $(\infty,2)$-category of $\underline{\FUN}^\oplax(\rS,\mC)$
with functors $\rS \to \Mon(\mC)$ and $\rS \to \Fun([1],\mC)^\mon,$ respectively (\cref{equis}, \cref{propp}).

\vspace{1mm}
\item We finally reduce the coherent classification of monads by Eilenberg-Moore objects in $\mC$ to the non-coherent classification of monads by Eilenberg-Moore objects in $\underline{\FUN}^\oplax(\rS,\mC)$.
    
\end{enumerate}

%identify $\Fun([1], \mC)^\mon$
% with an $\infty$-category $\Mon(\mC)$ of monads on $\mC$ (Theorem \ref{uhnggfdaa}).
% To define $\Mon(\mC)$ we 

%In the following we construct $\underline{\FUN}^\oplax(\rS,\mC)$.

% where $\rS$ is an $\infty$-category.Next we organize monads on varying objects of an ambient$(\infty,2)$-category to an $\infty$-category.
%\subsubsection{Oplax natural transformations}

\vspace{1mm}

We start with defining oplax natural transformations.

\begin{notation}

For every $\infty$-category $\rS$ let $\infty\Cat_{/ \rS}^{\cart,\oplax} \subset \infty\Cat_{/ \rS}$ be the full subcategory of cartesian fibrations over $\rS$.

\end{notation}

The canonical left $\infty\Cat$-action on $ \infty\Cat_{/ \rS}$ restricts to a left
$\infty\Cat$-action on $\infty\Cat_{/ \rS}^{\cart,\oplax}$.
Moreover the $\infty\Cat$-linear inclusion $\infty\Cat_{/ \rS}^{\cart} \subset \infty\Cat_{/ \rS}$ induces a $\infty\Cat$-linear inclusion $\infty\Cat_{/ \rS}^{\cart} \subset \infty\Cat_{/ \rS}^{\cart,\oplax}$.

\begin{remark}

Let $\mC$ be an $(\infty,2)$-category.
The 2-Yoneda-embedding $$\mC \hookrightarrow \FUN(\mC^\op,\infty\Cat)$$ gives rise to an inclusion of $(\infty,2)$-categories: $$\theta: \FUN(\rS,\mC) \subset \FUN(\rS, \FUN(\mC^\op,\infty\Cat)) \simeq $$$$
\FUN(\mC^\op, \FUN(\rS,\infty\Cat))
\simeq \FUN(\mC^\op, \infty\Cat_{/ \rS^\op}^{\cart}) \subset \FUN(\mC^\op, \infty\Cat_{/ \rS^\op}^{\cart,\oplax}) .$$

\end{remark}

\begin{notation}\label{oplll}

Let $\rS$ be an $\infty$-category and $\mC $ an $(\infty,2)$-category.
Let $\FUN^\oplax(\rS,\mC) $ be the essential image of the inclusion $\theta.$

\end{notation}

So there is an essentially surjective inclusion
$\FUN(\rS,\mC) \hookrightarrow\FUN^\oplax(\rS,\mC)$
of $(\infty,2)$-categories
that induces for every functors $\F,\G:\rS \to\mC$
an embedding $$\Nat(\F,\G):=\Mor_{\Fun(\rS,\mC)}(\F,\G) \hookrightarrow \Nat_\oplax(\F,\G):=\Mor_{\FUN^\oplax(\rS,\mC)}(\F,\G).$$

Next we refine the $\infty\Cat$-enrichment of $\FUN^\oplax(\rS,\mC)$
to an enrichment in $\infty\Cat_{/ \rS^\op}.$

\begin{notation}
Let $\rS$ be an $\infty$-category and $$\delta_\rS:=(-)\times \rS:\infty\Cat \to \infty\Cat_{/ \rS}: \Fun_{\rS}(\rS,-)$$ the canonical adjunction.
\end{notation}

\begin{lemma}\label{enros}
Let $\rS$ be an $\infty$-category. The canonical left action of $\infty\Cat_{/ \rS}$
on $\infty\Cat_{/ \rS}$ exhibits the full subcategory $\infty\Cat_{/ \rS}^{\cart,\oplax}$ as enriched in $\infty\Cat_{/ \rS}$,
where $\Mor_{ \infty\Cat_{/ \rS}^{\cart,\oplax}}(-,-)=\Fun^\rS(-,-)$ (Notation \ref{Parfun}).
%Pulling back this enrichment via the left adjoint functor$\delta_\rS$ we obtain the canonical enrichment in $\infty\Cat.$

\end{lemma}
\begin{proof}
%The $\infty$-category $\infty\Cat_{/ \rS}$ is left tensored over itselfvia the product and in particular left tensored over $\infty\Cat$by pulling back along $\delta_\rS.$In particular, $\infty\Cat_{/ \rS}$ is pseudo-enriched in itself,where the pullback of the pseudo-enrichment along $\delta_\rS$gives the canonical enrichment in $\infty\Cat$.The full subcategory $\infty\Cat_{/ \rS}^{\cart,\oplax}$is pseudo-enriched in $\infty\Cat_{/ \rS}$, where the pullback of the pseudo-enrichment along $\delta$gives the enrichment in $\infty\Cat$.
For every cartesian fibrations $\mA \to \rS,\mB \to \rS, \mD \to \rS$ there is a canonical equivalence
$$\Map_{\infty\Cat_{/ \rS}}(\mD, \Fun^\rS(\mA,\mB)) \simeq
\Map_{\infty\Cat_{/ \rS}}(\mD \times_\rS \mA,\mB) \simeq \Mul_{\infty\Cat_{/ \rS}}(\mD,\mA;\mB) \simeq  \Mul_{\infty\Cat_{/ \rS}^{\cart,\oplax}}(\mD,\mA;\mB).$$
\end{proof}

\begin{notation}Let $\rS$ be an $\infty$-category and $\mC, \mD$ be
$\infty$-categories enriched in $ \infty\Cat_{/ \rS}.$
We set $$\FUN^\rS(\mC,\mD) :={\infty\Cat_{/ \rS}} \mathrm{-}\Fun(\mC,\mD).$$ 
\end{notation}

\begin{notation}\label{opplas}

Let $\rS$ be an $\infty$-category and $\mC$ an $(\infty,2)$-category.
Let $$\underline{\FUN}^\oplax(\rS,\mC)$$ be the full $\infty\Cat_{/\rS^\op}$-enriched subcategory 
of $$\FUN^{\rS^\op}((\delta_{\rS^\op})_!(\mC)^\op, \infty\Cat_{/ \rS^\op}^{\cart,\oplax})$$ spanned by $\Fun^\oplax(\rS,\mC)$.

\end{notation}

\begin{remark}
Let $$ (-) \times \rS: \infty\Cat \to \infty\Cat_{/\rS}: \sigma:= \Fun_\rS(\rS,-)$$ be the canonical adjunction.
There is a canonical equivalence of $(\infty,2)$-categories
$$\FUN(\mC^\op, \infty\Cat_{/ \rS^\op}^{\cart,\oplax}) \simeq \sigma_!(\FUN^{\rS^\op}((\delta_{\rS^\op})_!(\mC)^\op, \infty\Cat_{/ \rS^\op}^{\cart,\oplax})) $$
that restricts to an equivalence $$\FUN^\oplax(\rS,\mC) \simeq  \sigma_!(\underline{\FUN}^\oplax(\rS,\mC)).$$
\end{remark}

\begin{notation}

Let $\rS$ be an $\infty$-category, $\mC$ an $(\infty,2)$-category
and $\F,\G:\rS \to\mC$ functors.
Let $$\Nat_\oplax^{\rS^\op}(\F,\G)
:=\Mor_{\underline{\FUN}^\oplax(\rS,\mC)}(\F,\G).$$

\end{notation}

Thus $$\Nat_\oplax(\F,\G)
\simeq \Fun_{\rS^\op}(\rS^\op, \Nat_\oplax^{\rS^\op}(\F,\G)).$$

\vspace{1mm}
There is the following functoriality of $\Nat_\oplax(\F,\G)$:

\begin{remark}\label{FUNK}

Let $\varphi: \T \to \rS$ be a functor.
The induced finite products preserving functor 
$$\Phi:=\varphi^*: \infty\Cat_{/ \rS^\op}^{\cart,\oplax} \to \infty\Cat_{/\T^\op}^{\cart,\oplax}$$ gives rise to a $\infty\Cat_{/\T^\op}$-enriched functor
$$\Phi_!(\FUN^{\rS^\op}(\delta_!(\mC)^\op, \infty\Cat_{/ \rS^\op}^{\cart,\oplax})) \to \FUN^{\T^\op}(\delta_!(\mC)^\op, \infty\Cat_{/\T^\op}^{\cart,\oplax}).$$
The latter restricts to a $\infty\Cat_{/\T^\op}$-enriched functor
$$\Phi_!(\underline{\FUN}^\oplax(\rS,\mC)) \to \underline{\FUN}^\oplax(\T,\mC).$$

Thus for every functors $\F,\G:\rS \to\mC$ we obtain a functor over $\T^\op:$
$$\T^\op \times_{\rS^\op} \Nat_\oplax^{\rS^\op}(\F,\G) \to \Nat_\oplax^{\T^\op}(\F \circ \varphi,\G \circ \varphi).$$
\end{remark}

\begin{proposition}\label{poqay}
Let $\varphi: \T \to \rS$ be a functor, $\mC$ an $(\infty,2)$-category and $\F,\G:\rS \to\mC$ functors.
The induced functor 
$$\T^\op \times_{\rS^\op} \Nat_\oplax^{\rS^\op}(\F,\G) \to \Nat_\oplax^{\T^\op}(\F \circ \varphi,\G \circ \varphi) $$
over $\T^\op$ is an equivalence.
\end{proposition}

\begin{proof}

Let $\R \to \T$ be a functor.
It is enough to prove that the induced functor
\begin{equation}\label{zuhl}
\Fun_{\rS^\op}(\R^\op, \Nat_\oplax^{\rS^\op}(\F,\G)) \simeq \Fun_{\T^\op}(\R^\op,\T^\op \times_{\rS^\op} \Nat_\oplax^{\rS^\op}(\F,\G)) \to \end{equation}
$$ \Fun_{\T^\op}(\R^\op,\Nat_\oplax^{\T^\op}(\F \circ \varphi,\G \circ \varphi)) $$
is an equivalence.
Let $$\kappa_\T:=(-)\times \R^\op: \infty\Cat \rightleftarrows \Cat_{\infty / \T^\op}: \Fun_{\T^\op}(\R^\op,-) $$ and let $\kappa_\rS$ be defined similarly.
Let $$\Phi:=\varphi^*: \infty\Cat_{/ \rS^\op} \to \infty\Cat_{/\T^\op}.$$

To see that the functor (\ref{zuhl}) is an equivalence, we prove that for any $(\infty,2)$-category $ \mC$ the following 2-functor induces equivalences on morphism $\infty$-categories:
\begin{equation}\label{eqax}\kappa_\rS^*(\FUN^{\rS^\op}((\delta_\rS)_!(\mC), \infty\Cat_{/ \rS^\op}^{\cart,\oplax})) \end{equation}
$$\simeq  \kappa_\T^*(\Phi_!(\FUN^{\rS^\op}((\delta_\rS)_!(\mC), \infty\Cat_{/ \rS^\op}^{\cart,\oplax}))) \to $$
$$ \kappa_\T^*(\FUN^{\T^\op}((\delta_\T)_!(\mC), \infty\Cat_{/\T^\op}^{\cart,\oplax})).
$$

In the following we prove that the 2-functor
(\ref{eqax}) is natural in $\mC$: we construct an $ (\infty,2)\Cat$-enriched natural transformation
$\zeta$ of functors $(\infty,2)\Cat^\op \to (\infty,2)\Cat$
whose component at any $(\infty,2)$-category $ \mC$ is the 2-functor (\ref{eqax}).
The adjunction $$\delta= (-) \times \rS^\op: \infty\Cat \rightleftarrows \infty\Cat_{/ \rS^\op}: \gamma:=\Fun_{\rS^\op}(\rS^\op,-)$$
gives rise to an adjunction $$\delta_!: (\infty,2)\Cat \to {\infty\Cat_{/\rS^\op}}\mathrm{-}\Cat : \gamma_!.$$
Since the functor $\delta$ preserves finite products,
the functor $\delta_!$ preserves finite products, too, and so makes ${\infty\Cat_{/\rS^\op}}\mathrm{-}\Cat$ to an $\infty$-category left tensored over $ (\infty,2)\Cat.$
Tautologically, the left adjoint $\delta_!$ is $(\infty,2)\Cat$-linear and thus admits a 
$(\infty,2)\Cat$-enriched right adjoint.

The functor
$\kappa_\rS= (-) \times \R^\op: \infty\Cat \to \infty\Cat_{/\rS^\op}$ is $ \infty\Cat$-linear. Hence the induced functor
$$(\kappa_\rS)_!: (\infty,2)\Cat \to {\infty\Cat_{/\rS^\op}}\mathrm{-}\Cat$$
is ${\infty\Cat_{/\rS^\op}}\mathrm{-}\Cat$-linear and thus in particular
$(\infty,2)\Cat$-linear by restricting along $\delta_!$
and so admits a $(\infty,2)\Cat$-enriched right adjoint $\kappa_\rS^*$.

For every $ \infty\Cat_{/\rS^\op}$-enriched $\infty$-category $\mD$ the functor $$ \FUN^{\rS^\op}(-,\mD):{\infty\Cat_{/\rS^\op}}\mathrm{-}\Cat^\op \to {\infty\Cat_{/\rS^\op}}\mathrm{-}\Cat $$
is canonically $ (\infty,2)\Cat$-enriched and admits a $ (\infty,2)\Cat$-enriched left adjoint.
Moreover the natural transformation
$$\lambda: \Phi_! \circ \FUN^{\rS^\op}(-,\mD) \to \FUN^{\T^\op}(-,\Phi_!(\mD)) \circ \Phi_! $$
of functors ${\infty\Cat_{/\rS^\op}}\mathrm{-}\Cat^\op \to  {\infty\Cat_{/\T^\op}}\mathrm{-}\Cat $ is $ (\infty,2)\Cat$-enriched.

Composing with $\kappa_\rS^*$ and $\delta_!$ we find that the functor 
$$\kappa_\rS^* \circ \FUN^{\rS^\op}(-,\mD) \circ \delta_!: (\infty,2)\Cat^\op \to (\infty,2)\Cat$$
is canonically $ (\infty,2)\Cat$-enriched and has a $ (\infty,2)\Cat$-linear left adjoint.
Moreover we obtain a $ (\infty,2)\Cat$-enriched transformation
$\zeta:=\kappa_\T^* \circ \lambda \circ \delta_!:$$$ \kappa_\rS^* \circ \FUN^{\rS^\op}(-,\mD)\circ \delta_! \simeq \kappa_\T^* \circ \Phi_! \circ \FUN^{\rS^\op}(-,\mD) \circ \delta_! \to $$$$ \kappa_\T^* \circ \FUN^{\T^\op}(-,\Phi_!(\mD)) \circ\Phi_! \circ \delta_! \simeq \kappa_\T^* \circ \FUN^{\T^\op}(-,\Phi_!(\mD)) \circ \delta_!$$
of $(\infty,2)\Cat$-enriched functors $(\infty,2)\Cat^\op \to (\infty,2)\Cat$
whose component at any $(\infty,2)$-category $ \mC$ is the 2-functor (\ref{eqax}).

We like to see that for any $(\infty,2)$-category $\mC$ the 2-functor $\zeta(\mC)$ is an embedding.
Since source and target of $\zeta$ admit $ (\infty,2)\Cat$-linear left adjoints (and so preserve $ (\infty,2)\Cat$-enriched cotensors),
$\zeta(\mC)$ identifies with the 2-functor $\FUN(\mC, \zeta(*))$,
which is an embedding if $\zeta(*)$ is an embedding. It remains to prove that $\zeta(*)$ is an embedding.
For that we need to show that for any cartesian fibrations $\X \to \rS^\op,\Y \to \rS^\op$ the functor 
$$\T^\op \times_{\rS^\op} \Fun^{\rS^\op}(\X,\Y) \to \Fun^{\T^\op}(\T^\op \times_{\rS^\op} \X, \T^\op \times_{\rS^\op} \Y) $$
over $\T^\op$ is an equivalence. This follows from Remark \ref{zyaqcfd}.

\end{proof}

%\begin{proposition}\label{spa}Let $\rS$ be a space and $\mC^\circledast \to \infty\Cat^\times$ an $(\infty,2)$-category.There is an equivalence $$\underline{\FUN}^\oplax(\rS,\mC)^\circledast \simeq (\mC^\circledast)^{\rS^\op}$$ of $\infty$-categories enriched in $\Fun(\rS^\op,\infty\Cat) \simeq \infty\Cat_{/\rS^\op}.$\end{proposition}

%\begin{proof}If $\rS$ is a space, the embedding $ \infty\Cat_{/ \rS^\op}^{\cart,\oplax} \subset \infty\Cat_{/ \rS^\op}$ is an equivalence.Thus the embedding of $ \infty\Cat_{/ \rS^\op} \simeq \Fun(\rS^\op,\infty\Cat)$-enriched $\infty$-categories$$\FUN^{\rS^\op}(\delta_!(\mC)^\op, \infty\Cat_{/\rS^\op}^{\cart,\oplax})^\circledast \hookrightarrow \FUN^{\rS^\op}(\delta_!(\mC)^\op, \infty\Cat_{/\rS^\op})^\circledast\simeq$$$$\Enr\Fun_{\Fun(\rS^\op,\infty\Cat)}(\delta_!(\mC)^\op, \infty\Cat^{\rS^\op})^\circledast \simeq \FUN(\mC^\op, \infty\Cat^{\rS^\op})^\circledast \simeq (\FUN(\mC^\op, \infty\Cat)^\circledast)^{\rS^\op}$$is an equivalence. This equivalence restricts to an equivalence $\underline{\FUN}^\oplax(\rS,\mC)^\circledast \simeq (\mC^\circledast)^{\rS^\op}$.\end{proof}

\begin{proposition}

Let $\rS$ be an $\infty$-category.
There is an equivalence of $\infty\Cat_{/ \rS^\op}$-enriched $\infty$-categories:

$$ \underline{\FUN}^\oplax(\rS,\infty\Cat) \simeq \infty\Cat_{/\rS^\op}^{\cart,\oplax}$$
\end{proposition}

\begin{proof}

By definition the full subcategory
$$\FUN^\oplax(\rS,\infty\Cat) \subset \FUN(\infty\Cat^\op, \Cat_{ \infty /\rS^\op}^{\cart,\oplax})$$
agrees with the full subcategory $\FUN^ \mathrm{coten}(\infty\Cat^\op, \infty\Cat_{/ \rS^\op}^{\cart,\oplax})$
of 2-functors preserving cotensors.
The equivalence
$$\infty\Cat_{/\rS^\op}^{\cart,\oplax} \simeq \infty\Cat\mathrm{-}\LinFun(\infty\Cat, (\infty\Cat_{/ \rS^\op}^{\cart,\oplax})^\op)^\op \simeq $$$$ \FUN^{ \mathrm{coten}}(\infty\Cat^\op, \infty\Cat_{/ \rS^\op}^{\cart,\oplax}) = \FUN^\oplax(\rS,\infty\Cat)$$
underlies a $\infty\Cat_{/\rS^\op}$-enriched functor
$$\theta: \infty\Cat_{/\rS^\op}^{\cart,\oplax} \to \underline{\FUN}^\oplax(\rS,\infty\Cat), (\mX \to \rS^\op) \mapsto (\K \mapsto \mX^\K \to \rS^\op).$$
Let $\mX \to \rS^\op,\mY \to \rS^\op$ be cartesian fibrations classifying functors $\F,\G: \rS \to \infty\Cat.$
The functor $\theta$ induces on morphism
$\infty$-categories over $\rS^\op$ a functor
$$\Fun^{\rS^\op}(\mX,\mY) \to \Nat_\oplax^{\rS^\op}(\F,\G)$$
that yields for every functor $\varphi: \T \to \rS$ the equivalence induced by $\rho$ for $\rS$ replaced by $\T:$
$$ \Fun_{\rS^\op}(\T^\op, \Fun^{\rS^\op}(\mX,\mY))^\simeq \simeq$$$$ \Fun_{\T^\op}(\T^\op, \Fun^{\T^\op}(\T^\op \times_{\rS^\op}\mX, \T^\op \times_{\rS^\op}\mY))^\simeq \to \Fun_{\rS^\op}(\T^\op,\Nat_\oplax^{\rS^\op}(\F,\G))^\simeq \simeq$$$$ \Fun_{\T^\op}(\T^\op,\Nat_\oplax^{\T^\op}(\F \circ \varphi,\G \circ \varphi))^\simeq.$$

\end{proof}

%\subsubsection{Eilenberg-Moore objects in functor $(\infty,2)$-categories}

Now we are ready to study Eilenberg-Moore objects in the $(\infty,2)$-category $\FUN^\oplax(\rS,\mC)$ for every $\infty$-category $\rS$ (Proposition \ref{oplax}). We start with studying right adjoint morphisms in
the $(\infty,2)$-category $\FUN^\oplax(\rS,\mC)$, which is
Proposition \ref{rightfun}, and continue with studying Eilenberg-Moore objects in 
the $(\infty,2)$-category of 2-presheaves (Proposition \ref{Eilpre})
and general functor $(\infty,2)$-categories (Proposition \ref{fun}).

%For that we first study right adjoint morphisms in general functor$(\infty,2)$-categories.

\begin{proposition}\label{rightfun}
	
Let $\rS$ be an $\infty$-category and $\mC$ an $(\infty,2)$-category and $\alpha: \F \to \G$ a morphism in $\FUN^\oplax(\rS,\mC)	$.
The following conditions are equivalent:
	
\begin{enumerate}
		
\item The morphism $\alpha$ admits a left adjoint.

\vspace{1mm}

\item The morphism $\alpha$ belongs to $\Fun(\rS,\mC)$
and for every object $\Z \in\rS$ the map $\alpha_\Z: \F(\Z) \to \G(\Z)$
in $\mC$ admits a left adjoint.
		
\end{enumerate}
	
\end{proposition}

\begin{proof}
Let $\kappa$ be the embedding
$\FUN^\oplax(\rS,\mC) \subset \FUN(\mC^\op,\infty\Cat^{\cart,\oplax}_{/ \rS})$.
	
Assume that (1) holds and let $\Z \in \rS$.
Then $\alpha_\Z: \F(\Z) \to \G(\Z)$ admits a left adjoint.
Moreover for every $\Y \in \mC$ the 
functor $\kappa(\alpha)(\Y): \kappa(\F)(\Y) \to \kappa(\G)(\Y)$ over $\rS$ admits a left adjoint $\beta^\Y$ relative to $\rS$ and
so is a map of cartesian fibrations over $\rS$ by \cite[Proposition 7.3.2.6.]{lurie.higheralgebra}. Hence $\alpha$ belongs to $\Fun(\rS,\mC).$
Assume that 2. holds and for every $\Z \in \mC$ let $\Gamma^\Z$ be the left adjoint of $\alpha_\Z.$ We show that $\kappa(\alpha)$ admits a left adjoint.
By Corollary \ref{propos} it is enough to see that for every $\Y \in \mC$
the functor $\kappa(\alpha)(\Y): \kappa(\F)(\Y) \to \kappa(\G)(\Y)$ over $\rS$ admits a left adjoint $\beta^\Y$ relative to $\rS$ and
for every morphism $\rho: \Y \to \Y'$ in $\mC$
the following canonical natural transformation over $\rS$ is an equivalence:
$$ \beta^{\Y'} \circ \kappa(\G)(\rho) \to \kappa(\F)(\rho) \circ \beta^\Y.$$

By the first condition of (2) the functor $\kappa(\alpha)(\Y)$ is a map of cartesian fibrations over $\rS$.
Using \cite[Proposition 7.3.2.6.]{lurie.higheralgebra} it is enough to see that 
for every $\Z \in \rS$ the functor $\kappa(\alpha)(\Y)_\Z: \kappa(\F)(\Y)_\Z \to \kappa(\G)(\Y)_\Z$ on the fiber over $\Z$ admits a left adjoint $\beta^{\Y,\Z}$ and 
for every morphism $\rho: \Y \to \Y'$ in $\mC$
the following canonical natural transformation is an equivalence:
$$\lambda: \beta^{\Y',\Z} \circ \kappa(\G)(\rho)_\Z \to \kappa(\F)(\rho)_\Z \circ \beta^{\Y,\Z}.$$

This follows from the fact that 
$\kappa(\alpha)(\Y)_\Z: \kappa(\F)(\Y)_\Z \to \kappa(\G)(\Y)_\Z$ 
identifies with the functor
$$\Mor_\mC(\Y, \alpha(\Z)) : \Mor_\mC(\Y, \F(\Z)) \to \Mor_\mC(\Y, \G(\Z))$$
right adjoint to $\Mor_\mC(\Y, \Gamma^\Z) $ and $ \lambda $ identifies with the equivalence
$$ \Mor_\mC(\Y', \Gamma^\Z) \circ \Mor_\mC(\rho, \G(\Z)) \simeq \Mor_\mC(\rho,  \Gamma^\Z) \simeq \Mor_\mC(\rho, \F(\Z)) \circ \Mor_\mC(\Y, \Gamma^\Z).$$
	
\end{proof}

\begin{remark}
Proposition \ref{rightfun}	for a differrent model of $(\infty,2)$-category of functors and oplax natural transformations was proven by \cite[Theorem 4.6.]{Haugseng2020OnLT} 
and generalized by \cite[Proposition 7.16.]{heine2025categorification}.
	
\end{remark}

%\subsubsection{Eilenberg-Moore objects in functor $(\infty,2)$-categories}

\begin{construction}

Let $\mC$ be an $(\infty, 2)$-category
and $\F$ a 2-functor $\mC^\op \to \infty\Cat.$
The endomorphism left $\Mor_{\FUN(\mC^\op, \infty\Cat)}(\F,\F)$-action on $\F $ corresponds by \cref{psinho} to a 2-functor $$\bar{\F}: \mC^\op \to {\Mor_{\FUN(\mC^\op, \infty\Cat)}(\F,\F)}\mathrm{-}\LMod $$ lifting $\F.$
For every monad $\T$ on $\F \in \FUN(\mC^\op, \infty\Cat)$ there is a 2-functor $$\LMod_\T: \Mor_{\FUN(\mC^\op, \infty\Cat)}(\F,\F) \mathrm{-}\LMod \to \infty\Cat$$
and a 2-transformation $\LMod_\T \to \nu$,  where $\nu: \Mor_{\FUN(\mC^\op, \infty\Cat)}(\F,\F) \mathrm{-}\LMod \to \infty\Cat$ is the forgetful 2-functor.
\end{construction}

\begin{proposition}\label{Eilpre}

Let $\mC$ be an $(\infty,2)$-category
and $\T$ a monad on some presheaf $\F$ in $\FUN(\mC^\op, \infty\Cat).$
The morphism $\LMod_\T \circ \bar{\F} \to \F$ is an Eilenberg-Moore object for $\F.$

\end{proposition}

\begin{proof}

We need to give an equivalence $$\theta: \LMod_\T \circ \Mor_{\FUN(\mC^\op, \infty\Cat)}(-,\F) \simeq \Mor_{\FUN(\mC^\op, \infty\Cat)}(-, \LMod_\T \circ \bar{\F}) $$ over $\Mor_{\FUN(\mC^\op, \infty\Cat)}(-,\F)$ in
$\FUN(\FUN(\mC^\op, \infty\Cat)^\op, \infty\Cat).$
By \cite[Theorem 5.1.]{HEINE2023108941} the functor $$ \FUN(\FUN(\mC^\op, \infty\Cat)^\op, \infty\Cat) \to  \FUN(\mC^\op, \infty\Cat)$$ that restricts along the 
2-Yoneda-embedding $\iota: \mC \to \FUN(\mC^\op, \infty\Cat)$ is an equivalence, where the left hand side is the full subcategory of 2-functors preserving small limits and cotensors.
Since source and target of $\theta$ preserve small limits and cotensors,
$\theta$ is given by its restriction along $\iota$. So it is enough to find an 
equivalence $$\LMod_\T \circ \Mor_{\FUN(\mC^\op, \infty\Cat)}(-,\F) \circ \iota \simeq \Mor_{\FUN(\mC^\op, \infty\Cat)}(-, \LMod_\T \circ \bar{\F}) \circ \iota $$ over $\Mor_{\FUN(\mC^\op, \infty\Cat)}(-,\F) \circ \iota $ in
$\FUN(\mC^\op, \infty\Cat).$
By the 2-Yoneda-lemma (\cref{enryol}) there are canonical equivalences 
$$\Mor_{\FUN(\mC^\op, \infty\Cat)}(-, \LMod_\T \circ \bar{\F}) \circ \iota \simeq \LMod_\T \circ \bar{\F}, \ \Mor_{\FUN(\mC^\op, \infty\Cat)}(-,\F) \circ \iota \simeq \F.$$

\end{proof}

\begin{corollary}\label{charo}

Let $\mC$ be an $(\infty,2)$-category.
The $(\infty,2)$-category $\FUN(\mC^\op, \infty\Cat)$ admits (co)Eilenberg-Moore objects, which are formed objectwise, and the 2-Yoneda-embedding $\iota: \mC \hookrightarrow \FUN(\mC^\op, \infty\Cat) $
preserves (co)Eilenberg-Moore objects.
\end{corollary}

\begin{proposition}\label{fun}

Let $\mC, \mD$ be $(\infty,2)$-categories.
If $\mD$ admits Eilenberg-Moore objects, coEilenberg Moore objects, Kleisli objects,
coKleisli objects, respectively, the $(\infty,2)$-category $\FUN(\mC, \mD)$ admits the same, which are formed objectwise.

\end{proposition}

\begin{proof}
It is enough to show that if $\mD$ admits Eilenberg-Moore objects, $\FUN(\mC, \mD)$ admits Eilenberg-Moore objects, which are formed objectwise, since there are equivalences of $(\infty,2)$-categories $$\FUN(\mC, \mD)^\mathrm{co} \simeq \FUN(\mC^\mathrm{co}, \mD^\mathrm{co}), \ \FUN(\mC, \mD)^\op \simeq \FUN(\mC^\op, \mD^\op).$$

The 2-Yoneda-embedding $\iota: \mD \hookrightarrow \FUN(\mD^\op, \infty\Cat)$ induces an embedding
$$\iota_!: \FUN(\mC, \mD) \hookrightarrow \FUN(\mC, \FUN(\mD^\op, \infty\Cat)) \simeq \FUN(\mC \times \mD^\op, \infty\Cat),$$
where the last equivalence is by \cite[Proposition 4.40.]{heine2024higher}. 
By Corollary \ref{charo} the $(\infty,2)$-category $$ \FUN(\mC \times \mD^\op, \infty\Cat)$$ admits Eilenberg-Moore objects,
which are formed objectwise.
So if $\F \in \FUN(\mC, \mD)$ and $\T$ is a monad on $\F$,
there is an Eilenberg-Moore object $\alpha \to \iota_!(\F)$ of $\iota_!(\T)$ such that for every $\X \in \mC, \Y \in \mD$
the induced functor
$\alpha(\Y,\X) \to \Mor_\mD(\Y,\F(\X))$ is an Eilenberg-Moore object
for $\Mor_\mD(\Y,\T(\X)).$
We prove that $\alpha$ belongs to the essential image of $\iota_!.$
Let $\mW \to \iota(\F(\X))$ be an Eilenberg-Moore object
in $\FUN(\mD^\op, \infty\Cat) $ for $\iota(\T(\X)).$
The induced left $\iota(\T(\X))$-action on the map
$\alpha(-,\X) \to \iota(\F(\X))$ in $\FUN(\mD^\op, \infty\Cat) $ 
yields a map $\lambda: \alpha(-,\X) \to \mW$ over $\iota(\F(\X))$
in $\FUN(\mD^\op, \infty\Cat) $.
The morphism $\lambda$ is objectwise an equivalence by the uniqueness of Eilenberg-Moore objects.
Thus the map
$\alpha(-,\X) \to \iota(\F(\X))$ in $\FUN( \mD^\op, \infty\Cat) $ is an Eilenberg-Moore object for $\iota(\T(\X)).$
By Corollary \ref{charo} the 2-Yoneda-embedding $\iota: \mD \hookrightarrow \FUN(\mD^\op, \infty\Cat)$ preserves Eilenberg-Moore objects
since $\mD$ has Eilenberg-Moore objects.
By the uniqueness of Eilenberg-Moore objects $\alpha(-,\X)$ belongs to the essential image of $\iota.$
Thus $\alpha$ belongs to the essential image of $\iota_!.$

\end{proof}
 
\begin{proposition}\label{oplax}

Let $\rS$ be an $\infty$-category and $\mC$ an $(\infty,2)$-category that admits Eilenberg-Moore objects.

\begin{enumerate}
\item Then $\Fun^{\oplax}(\rS, \mC)$ admits Eilenberg-Moore objects, which are formed objectwise.

\item Let $\alpha: \Y \to \X$ be a monadic morphism and $\beta: \Z \to \Y$ any morphism
in $\Fun^{\oplax}(\rS, \mC)$.
Then $\alpha$ belongs to $\Fun(\rS,\mC).$
Moreover $\beta$ belongs to $\Fun(\rS,\mC)$ if and only if
$\alpha \circ \beta$ belongs to $\Fun(\rS,\mC)$.
\end{enumerate}
\end{proposition}

\begin{proof}
(1): By \cref{zhgbfdc} the $(\infty,2)$-category $\infty\Cat^{\cart,\oplax}_{/ \rS^\op}$ admits Eilenberg-Moore objects and for every functors $\T \to \rS$ the induced 2-functor
$\infty\Cat^{\cart,\oplax}_{/ \rS^\op} \to \infty\Cat^{\cart,\oplax}_{/ \T^\op}$ preserves Eilenberg-Moore objects.
Proposition \ref{fun} implies that $\FUN(\mC^\op, \infty\Cat^{\cart,\oplax}_{/ \rS^\op})$
admits Eilenberg-Moore objects, which are formed objectwise. Thus the induced 2-functor 
$$\FUN(\mC^\op, \infty\Cat^{\cart,\oplax}_{/\rS^\op}) \to \FUN(\mC^\op, \infty\Cat^{\cart,\oplax}_{/\T^\op})$$
preserves Eilenberg-Moore objects.

We prove that $\Fun^{\oplax}(\rS, \mC)$ is closed in $\FUN(\mC^\op, \infty\Cat^{\cart,\oplax}_{/ \rS^\op})$
under Eilenberg-Moore objects.
First note that for every monadic morphism $\Y \to \X$ in $\FUN(\mC^\op, \infty\Cat^{\cart,\oplax}_{/ \rS^\op}),$ where $\X \in \FUN(\mC^\op, \infty\Cat^{\cart}_{/ \rS^\op})$,
also $\Y \in \FUN(\mC^\op, \infty\Cat^{\cart}_{/ \rS^\op})$.
This follows via Proposition \ref{fun} from the fact that every monadic morphism in
$\infty\Cat^{\cart,\oplax}_{/ \rS^\op}$ is a map of cartesian fibrations over
$\rS^\op$ that is fiberwise conservative (\cref{zhgbfdc}).

Let $\Y \to \X$ be an Eilenberg-Moore object in $\FUN(\mC^\op, \infty\Cat^{\cart,\oplax}_{/ \rS^\op})$ for some monad $\T$ on $\X$, where $\X \in \Fun^{\oplax}(\rS, \mC).$
Then $\Y \in \FUN(\mC^\op, \infty\Cat^{\cart}_{/ \rS^\op})$
and it is enough to see that for every $\Z \in \rS$ 
the image of $\Y$ under the 2-functor 
$$\Gamma_\Z: \FUN(\mC^\op, \infty\Cat^{\cart,\oplax}_{/ \rS^\op}) \to \FUN(\mC^\op, \infty\Cat)$$ taking the fiber
over $\Z$ belongs to the essential image of the 2-Yoneda-embedding.
Since $\Gamma_\Z$ preserves Eilenberg-Moore objects, $\Gamma_\Z(\Y) \to \Gamma_\Z(\X)$ is an Eilenberg-Moore object for $\Gamma_\Z(\T)$ in $\FUN(\mC^\op, \infty\Cat)$.
As $\X \in \Fun^{\oplax}(\rS, \mC),$ the object $\Gamma_\Z(\X)$
belongs to the essential image of the 2-Yoneda-embedding.
%and $\Gamma_\Z(\T)$ is a monad in $\mC.$
So also $\Gamma_\Z(\Y) $ belongs to the essential image of the 2-Yoneda-embedding as the 2-Yoneda-embedding preserves 
Eilenberg-Moore objects by Corollary \ref{charo}.

(2): Every monadic morphism has a left adjoint and so belongs to
$\Fun(\rS,\mC)$ by Proposition \ref{rightfun}.
The second part of 2. follows via Proposition \ref{fun} 
from the fact that every monadic morphism of $\infty\Cat^{\cart,\oplax}_{/ \rS^\op}$ is a fiberwise conservative map of cartesian fibrations over
$\rS^\op$ (\cref{zhgbfdc}).

\end{proof}

\begin{corollary}\label{monadis}

Let $\rS$ be an $\infty$-category, $\mC$ an $(\infty,2)$-category
and $\alpha: \F \to \G$ a morphism in $\FUN^\oplax(\rS,\mC)	$.
The following conditions are equivalent:

\begin{enumerate}

\item The morphism $\alpha$ is monadic.

\item The morphism $\alpha$ belongs to $\Fun(\rS,\mC)$
and for every $\Z \in\rS$ the map $\alpha_\Z: \F(\Z) \to \G(\Z)$
in $\mC$ is monadic.

\end{enumerate}

\end{corollary} 

\begin{proof}
	
If (1) holds, Proposition \ref{oplax} implies that for every object $\Z \in\rS$ the map $\alpha_\Z: \F(\Z) \to \G(\Z)$ in $\mC$ is monadic.
%, where we use that every monadic morphism is an Eilenberg-Moore object of some monad.
Moreover by Proposition \ref{rightfun} the morphism $\alpha$ belongs to $\Fun(\rS,\mC)$.

If (2) holds, by Proposition \ref{rightfun} the morphism $\alpha$ admits a left adjoint and so has an associated monad $\T$ on $\G.$
By Proposition \ref{oplax} there is an Eilenberg-Moore object
$\beta: \rH \to \G$ for $\T,$ which induces a morphism
$\gamma: \F \to \rH$ over $\G.$
By Proposition \ref{oplax} the morphisms $\beta,\gamma$ belong to
$\Fun(\rS,\mC)$ since $\alpha$ does.
Since Eilenberg-Moore objects in $\FUN^\oplax(\rS,\mC)$ are formed objectwise, by (2) and the uniqueness of Eilenberg-Moore objects,
the natural transformation $\gamma$ of functors $\rS \to \mC$ is objectwise an equivalence and so an equivalence.	
	
\end{proof}

Proposition \ref{rightfun} and Corollary \ref{monadis} imply the following corollary:

\begin{corollary}\label{propp}
Let $\rS$ be an $\infty$-category, $\mC$ an $(\infty,2)$-category and $\X: \rS \to \mC$ a functor.

\begin{enumerate}
\item The canonical equivalence
$$\Fun_\mC(\rS,\Fun([1],\mC)) \simeq \Fun(\rS,\mC)_{/\X}$$
restricts to an equivalence
$$\Fun_\mC(\rS,\Fun([1],\mC)^\R) \simeq (\FUN^\oplax(\rS,\mC)_{/\X})^\R.$$
\item If $\mC$ admits Eilenberg-Moore objects, the canonical equivalence
$$\Fun_\mC(\rS,\Fun([1],\mC)) \simeq \Fun(\rS,\mC)_{/\X}$$
restricts to an equivalence
$\Fun_\mC(\rS,\Fun([1],\mC)^\mon) \simeq (\FUN^\oplax(\rS,\mC)_{/\X})^\mon.$

\end{enumerate}
\end{corollary}

%Proposition \ref{propp} 1. follows from Proposition \ref{rightfun} and Proposition \ref{propp} 2. follows from Corollary \ref{monadis}.\vspace{1mm}

%\subsection{A classification of monads}

In the following we use the enrichment of $\underline{\FUN}^\oplax(\mC,\mC)$ in $\infty\Cat_{/ \mC^\op}$ to define an $\infty$-category of monads in $\mC.$
For the next definition we use the following notation

\begin{definition}\label{Exn}Let $\rS$ be an $\infty$-category.
A $\rS$-family of monoidal $\infty$-categories
is a map $\mV^\ot \to \rS \times \Delta^\op $
of cocartesian fibrations over $\Delta^\op$ such that for every $\s \in \rS$
the induced functor $\mV_\s^\ot \to \Delta^\op $ is a monoidal $\infty$-category.

\end{definition}

\begin{remark}
A functor $\mV^\ot \to \rS \times \Delta^\op $ is a $\rS$-family of monoidal $\infty$-categories if and only if the classifying functor
$\Delta^\op \to \infty\Cat_{/\rS}$ is a monoid object \cite[Definition 4.1.2.5.]{lurie.higheralgebra}, which corresponds to
an associative algebra for the cartesian symmetric monoidal structure on 
$\infty\Cat_{/\rS}$ by \cite[Proposition 2.4.1.7.]{lurie.higheralgebra}.
\end{remark}

\begin{notation}\label{Nojk}Let $\rS$ be an $\infty$-category and $\mV^\ot \to \rS \times \Delta^\op $ is a $\rS$-family of monoidal $\infty$-categories.
Let $$\Alg^\rS(\mV^\ot) \subset \Fun^{\rS}_{\rS\times \Delta^\op}(\rS \times \Delta^\op,\mV^\ot) \to \rS$$
whose fiber over $\s \in \rS$ is the full subcategory of functors $\Delta^\op \to \mV^\ot_\s$ 
over $\Delta^\op$ that preserve cocartesian lifts of inert maps, i.e. the associative algebras in $\mV_\s.$

\end{notation}

%and that associative algebras in $\infty\Cat_{/ \mC^\op}$ are classified by $\mC^\op$-families of monoidal $\infty$-categories (Example \ref{Exn}).

\begin{definition}\label{monasol}
Let $\mC$ be an $(\infty,2)$-category.
The $\infty$-category of monads in $\mC$ is
$$ \Mon(\mC):= \Alg^{\mC^\op}(\Nat^{\mC^\op}_\oplax(\id,\id)). $$

\end{definition}

\begin{notation}

For every $\X \in \mC$ we set $$\Mon_\X(\mC):= \{\X\}\times_{\mC^\op}\Mon(\mC) \simeq \Alg(\Mor_\mC(\X,\X)).$$

\end{notation}

The $\infty$-category $\Mon(\mC)$ over $\mC^\op$ is characterized by the following universal property:

\begin{proposition}\label{equis}
Let $\rS$ be an $\infty$-category and $\mC$ an $(\infty,2)$-category.
For every functor $\X: \rS \to \mC$
there is a canonical equivalence
$$\Fun_{\mC^\op}(\rS^\op,\Mon(\mC)) \simeq \Mon_\X(\FUN^\oplax(\rS,\mC)).$$

\end{proposition}

\begin{proof}Using the equivalence of Proposition \ref{poqay} there is an equivalence
$$\Fun_{\mC^\op}(\rS^\op,\Mon(\mC)) = \Fun_{\mC^\op}(\rS^\op, \Alg^{\mC^\op}(\Nat_\oplax^{\mC^\op}(\id,\id))) \simeq$$$$ \Fun_{\rS^\op}(\rS^\op, \rS^\op \times_{\mC^\op} \Alg^{\mC^\op}(\Nat^{\mC^\op}_\oplax(\id,\id))) \simeq $$$$\Fun_{\rS^\op}(\rS^\op, \Alg^{\rS^\op}(\rS^\op \times_{\mC^\op}\Nat^{\mC^\op}_\oplax(\id,\id))) \simeq $$$$\Fun_{\rS^\op}(\rS^\op, \Alg^{\rS^\op}(\Nat^{\rS^\op}_\oplax(\X,\X))) \simeq$$$$\Alg(\Fun_{\rS^\op}(\rS^\op, \Nat^{\rS^\op}_\oplax(\X,\X))) \simeq  \Alg(\Nat_\oplax(\X,\X)).$$

\end{proof}

%\begin{remark}By ... for any $(\infty,2)$-categories $\mB, \mC$there is an $(\infty,2)$-category $\FUN^\oplax(\mB,\mC)$and ... define the  $(\infty,2)$-category of monads in $\mC$ by$\Mon(\mC):= \FUN^\oplax(B\Delta_+, \mC)$,where $\Delta_+=\Delta^{\triangleleft}.$So by ... there is a canonical equivalence$\FUN^\oplax(\mB,\Mon(\mC)) \simeq \Mon(\FUN^\oplax(\mB,\mC)).$

%By ... the fiber over any $\Y \in \mC$ of the forgetful functor$\Mon(\mC):= \FUN^\oplax(B\Delta_+, \mC) \to \mC$ is $\Alg(\Mor_\mC(\Y,\Y)).$So for any functor $\X : \mC \to \rS$ there is a canonical equivalence$$\Fun_{\mC}(\rS,\Mon(\mC))^\simeq \simeq \{\alpha\}\times_{\FUN^\oplax(\rS, \mC)^\simeq} \FUN^\oplax(\rS,\Mon(\mC))^\simeq$$$$\simeq \{\alpha\}\times_{\FUN^\oplax(\rS, \mC)^\simeq} \Mon(\FUN^\oplax(\rS,\mC))^\simeq \simeq \Mor_{\FUN^\oplax(\rS,\mC)}(\X,\X)^\simeq.$$\end{remark}

%By Proposition \ref{equis} there is a canonical equivalence$$\Fun_{\mC^\op}(\rS^\op,\Mon(\mC)) \simeq \Alg(\Nat_\oplax(\X,\X)) = \Alg(\Mor_{\FUN^\oplax(\rS,\mC)}(\X,\X)).$$

If $\mC$ admits Eilenberg-Moore objects, by Proposition \ref{oplax} the $(\infty,2)$-category $\FUN^\oplax(\rS,\mC)$ admits Eilenberg-Moore objects, which are formed objectwise.
So we can make the following: % construction:

\begin{construction}Let $\mC $ be an $(\infty,2)$-category that admits Eilenberg-Moore objects and $\beta: \rS \to \Mon(\mC)^\op$ a functor.
Let $\X$ be the composition of $\beta$ and the functor $\Mon(\mC)^\op \to \mC$.
The functor $\beta^\op$ corresponds to a monad on $\X$ in $\FUN^\oplax(\rS,\mC)$,
whose Eilenberg-Moore object $\Y \to \X$ in $\FUN^\oplax(\rS,\mC)$
corresponds to a functor $$\Alg(\beta): \rS \to \Fun([1],\mC)^\mon$$ over $\mC.$
\end{construction}

\begin{notation}
Let $\mC$ be an $(\infty,2)$-category that admits Eilenberg-Moore objects.
For $\beta$ the identity of $\rS=  \Mon(\mC)^\op$ we write
$\Alg: \Mon(\mC)^\op \to  \Fun([1],\mC)^\mon$ for $\Alg(\beta).$	
\end{notation}

\begin{construction}Let $\mC $ be an $(\infty,2)$-category and $\gamma: \rS \to \Fun([1],\mC)^\R$ a functor.
Let $\X$ be the composition of $\gamma$ and evaluation at the target $\Fun([1],\mC)^\R \to \mC$.
The functor $\gamma$ over $\mC$ corresponds to a right adjoint morphism
$\Z \to \X$ in $\FUN^\oplax(\rS,\mC)$,
whose associated monad $\T$ on $\X$ corresponds to a functor $\End(\gamma): \rS \to \Mon(\mC)^\op$ over $\mC.$

\end{construction}

\begin{notation}
Let $\mC$ be an $(\infty,2)$-category.
For $\gamma$ the identity of $\rS= \Fun([1],\mC)^\R$ we write
$$\End: \Fun([1],\mC)^\R \to \Mon(\mC)^\op $$ for $\End(\gamma).$	
\end{notation}

Now we are ready to prove the main theorem of this section:

\begin{theorem}\label{uhnggfdaa}
Let $\mC$ be an $(\infty,2)$-category that admits Eilenberg-Moore objects. There is an adjunction $$ \End: \Fun([1],\mC)^\R\rightleftarrows \Mon(\mC)^\op: \Alg,$$ where $\Alg$ sends a monad to its Eilenberg-Moore object, and $\End$ sends a right adjoint morphism to its monad. The functor $\Alg$ is fully faithful and the essential image is $\Fun([1],\mC)^\mon $. For any right adjoint morphism $\g: \Z \to \X$ the unit $\Z \to \Alg(\End(\g))$ corresponds to the endomorphism left $\End(\g)$-action on $\g.$
\end{theorem}

\begin{proof}

Let $ \X: \rS \to \mC$ be a functor.
For every functor $ \alpha: \rS' \to \rS$ the induced 2-functor
$$\alpha^*: \FUN^\oplax(\rS,\mC) \to \FUN^\oplax(\rS',\mC)$$ preserves Eilenberg-Moore objects.
Consequently, for every functor $\beta: \rS \to \Mon(\mC)^\op$ over $\mC$
there is a canonical equivalence
$\Alg(\beta) \circ \alpha \simeq \Alg(\beta \circ \alpha)$
of functors $\rS' \to \Fun([1],\mC)^\mon.$
In particular, for $\beta$ the identity there is a canonical equivalence
$ \Alg \circ \alpha \simeq \Alg(\alpha): \rS \to \Fun([1],\mC)^\mon.$
By Remark \ref{FUNK} for every 
functor $\gamma: \rS \to \Fun([1], \mC)^\R$ over $\mC$ there is a canonical equivalence
$\End(\gamma) \circ \alpha \simeq \End(\gamma \circ \alpha)$ of functors
$\rS' \to \Mon(\mC)^\op$ and so an equivalence $\End \circ \alpha \simeq \End(\alpha): \rS \to \Mon(\mC)^\op$.

A functor $\beta: \rS \to \Mon(\mC)^\op$ over $\mC$ 
corresponds to a monad $\T$ on $\X$ in $\FUN^\oplax(\rS,\mC)$
that admits an Eilenberg-Moore object $\Y \to \X$ in $\FUN^\oplax(\rS,\mC)$
(corresponding to $\Alg(\beta): \rS \to \Fun([1], \mC)^\mon$), 
whose associated monad (corresponding to $\End(\Alg(\beta))$ is $\T.$
Thus $\End(\Alg(\beta)) \simeq \beta$. For $\beta$ the identity
we obtain an equivalence $\epsilon: \End \circ \Alg \simeq \id$
of functors $\Mon(\mC)^\op \to \Mon(\mC)^\op.$

A functor $\gamma: \rS \to \Fun([1], \mC)^\R$ over $\mC$ 
corresponds to a right adjoint morphism $\Z \to \X$ in
$\FUN^\oplax(\rS,\mC)$, whose monad admits an Eilenberg-Moore object
$\Y \to \X$ in $\FUN^\oplax(\rS,\mC)$
(corresponding to $\Alg(\End(\gamma)): \rS \to \Fun([1], \mC)^\mon$).
The induced morphism $\rho: \Z \to \Y$ over $\X$ in $\FUN^\oplax(\rS,\mC)$
corresponds to a natural transformation $\gamma \to \Alg(\End(\gamma)) \simeq \Alg \circ \End \circ \gamma$ of functors $\rS \to \Fun([1],\mC)^\R$ over $\mC.$
For $\gamma$ the identity we obtain a natural transformation $\eta: \id \to \Alg \circ \End $ of functors $\Fun([1],\mC)^\R \to \Fun([1],\mC)^\R$ over $\mC.$
If $\gamma$ lands in $\Fun([1], \mC)^\mon,$
the corresponding right adjoint morphism $\Z \to \X$ is monadic
so that $\rho: \Z \to \Y$ is an equivalence.
This implies that the restriction $\eta'$ of $\eta$ to
$\Fun([1],\mC)^\mon$ is an equivalence.
Consequently, $\eta'$ and $\epsilon$ exhibit $$\Fun([1], \mC)^\mon \subset \Fun([1], \mC)^\R \xrightarrow{\End} \Mon(\mC)^\op$$ as inverse to $$\Alg: \Mon(\mC)^\op \to \Fun([1], \mC)^\mon.$$
Moreover $\eta \circ \Alg$ is an equivalence.

The natural transformation $\End \circ \eta$ is also an equivalence
since for every right adjoint morphism $\G: \Z \to \X$ 
the canonical morphism $\rho: \Z \to \Y$ over $\X$
induces an equivalence on associated monads,
which follows from Lemma \ref{trwgqaxy} 3.
This implies by \cite[Proposition 5.2.7.4.]{lurie.HTT} that $\eta: \id \to \Alg \circ \End $ exhibits $\End$ as left adjoint to $\Alg$, where we use that $\Alg$ is fully faithful.

\end{proof}

\begin{remark}
Let $\G: \Y \to \X, \h: \Z \to \X $ be morphisms in $\mC$ that admit left adjoints $\F: \X \to \Y, \bk: \X \to \Z, $ 
respectively.
By Lemma \ref{trwgqaxy} 3. a morphism $\phi: \Y \to \Z$ in $\mC$ over $\X$ is a local equivalence if and only if the morphism 
$$\h \circ \bk \to  \h \circ\bk \circ \G \circ \F \simeq \h \circ \bk \circ \h \circ \phi \circ \F \to  \h \circ \phi \circ \F \simeq  \G \circ \F $$ in $\Mor_\mC(\X,\X) $ is an equivalence.

In particular, a morphism $\phi: \Y \to \Z$ in $ \mC$ over $\X$, where $\Z$ is a local object, is a local equivalence 
if and only if the morphism $\bk \to  \bk \circ \G \circ \F \simeq \bk \circ \h \circ \phi \circ \F \to  \phi \circ \F $ in $\Mor_\mC(\X, \Z)$ is an equivalence.	

\end{remark}

\begin{remark}

With a bit extra work one can enhance $\Mon(\mC)$ to an  $(\infty,2)$-category
and the equivalence of \cref{uhnggfdaa}
to an equivalence of $(\infty,2)$-categories.
The proof of \cref{uhnggfdaa} remains unchanged but one needs to introduce
cartesian fibrations of $(\infty,2)$-categories to extend the construction of
$\underline{\FUN}^\oplax(\rS,\mC)$ to the situation that $\rS$ is an $(\infty,2)$-category. A similar construction gives that 
$\underline{\FUN}^\oplax(\rS,\mC)$ is enriched in $(\infty,2)\Cat_{/\rS}$
and one defines $\Mon(\mC)$ to be the $\infty$-category of associative algebras relative to $\mC$ of the endomorphism algebra $$\Mor_{\underline{\FUN}^\oplax(\mC,\mC)}(\id_\mC,\id_\mC) \to \mC $$ of the identity in
$\underline{\FUN}^\oplax(\mC,\mC)$.
A similar argument gives that $\underline{\FUN}^\oplax(\rS,\mC)$ admits Eilenberg-Moore objects for any $(\infty,2)$-category $\rS$
and that functors $\rS \to \Mon(\mC)$ correspond to monads in
$\underline{\FUN}^\oplax(\rS,\mC)$.

Having this the same proof as \cref{uhnggfdaa} gives the  $(\infty,2)$-categorical refinement.
The only reasons why we did not present the enhanced version of \cref{uhnggfdaa} are that we did not want to introduce fibrations of $(\infty,2)$-categories,
which most readers are not familiar with, and that we do not have any application that makes use of the $(\infty,2)$-categorical refinement.

\end{remark}

\begin{corollary}\label{Corollis}
Let $\mV$ be a monoidal $\infty$-category.
There 
 a canonical localization
$$ \Alg(\mV)^\op \hookrightarrow (\mV\mathrm{-}\RMod_{/\mV})^\R$$
sending $\A$ to $\LMod_\A(\mV) \to \mV.$
\end{corollary}

\begin{proof}
The closed left $\infty\Cat$-action on $\mV\mathrm{-}\RMod$ of Remark \ref{enr2cat} makes $\mV\mathrm{-}\RMod$ to an $(\infty,2)$-category.
%By \cite[Corollary 4.2.4.7.]{lurie.higheralgebra} 
By Remark \ref{Fama} (3) there is a monoidal equivalence $ \mV \simeq \mV\mathrm{-}\LinFun(\mV,\mV), \A \mapsto \A \ot(-):\mV \to \mV, $ where $\mV$ carries the right action over itself, that sends an associative algebra $\A$ in $\mV$ 
to a monad $\A \ot(-):\mV \to \mV$ in $\RMod_{\mV}.$
By Remark \ref{Fama} (2) there is a $\mV$-linear equivalence  
$\LMod_\A(\mV) \simeq \LMod_{\A \ot (-)}(\mV).$ 
By Corollary \ref{exk} (2) the $(\infty,2)$-category $\mV\mathrm{-}\RMod$ admits Eilenberg-Moore objects and the $\mV$-linear functor $\LMod_\A(\mV)
\simeq \LMod_{\A \ot (-)}(\mV)$ is an Eilenberg-Moore object for $\A \ot (-): \mV \to \mV $ in $\mV\mathrm{-}\RMod.$ 
We apply Theorem \ref{uhnggfdaa}.

\end{proof}

% \begin{remark}

% By our construction the $\infty$-category of monads is 

% \cref{uhnggfdaa}
% \end{remark}

\begin{corollary}\label{monan}
Let $\mC$ be an $(\infty,2)$-category that admits Eilenberg-Moore objects. The embedding
$$ \Fun([1], \mC)^\mon \subset \Fun([1], \mC)^\R $$ 
of $(\infty,2)$-categories admits a left adjoint. 
\end{corollary}

\begin{proof}%[Proof of Theorem \ref{monan}]
Since the 2-Yoneda-embedding $\iota: \mC \hookrightarrow \FUN(\mC^\op, \infty\Cat)$ preserves Eilenberg-Moore objects by Corollary \ref{charo}, it is enough to prove the statement for $\FUN(\mC^\op, \infty\Cat)$, which is a presentable $(\infty,2)$-category. Thus we may assume that $\mC$ is tensored and cotensored over $\infty\Cat.$
In this case also $\mC^{[1]} $ exhibits $\Fun([1], \mC)$ as tensored and cotensored over $\infty\Cat$ with objectwise cotensor.
By Theorem \ref{uhnggfdaa} the embedding $$ \Fun([1], \mC)^\mon  \subset  \Fun([1], \mC)^\R $$ admits a left adjoint.
Therefore by \cref{remqa} (2) it is enough to see that the full subcategories $ \Fun([1], \mC)^\mon$ and $ \Fun([1], \mC)^\R $ of $ \Fun([1], \mC)$ are closed under cotensors. 

Because $\mC$ is tensored over $\infty\Cat$, for every $\K \in \infty\Cat$ there is a $\infty\Cat$-linear functor $\K \ot (-):\mC \to \mC$ left adjoint to the functor $(-)^\K: \mC \to \mC$ taking the cotensor,
which inherits a unique structure of a 2-functor from its left adjoint.
So for every right adjoint morphism $\h: \W \to \V$ in $\mC$
the induced morphism on cotensors $\h^\K: \W^\K \to \V^\K$ in $\mC$ is right adjoint.
For every $\Z \in \mC$ the functor $\Mor_\mC(\Z, \W^\K) \to \Mor_\mC(\Z, \V^\K)$ identifies with the functor
$\Fun(\K,\Mor_\mC(\Z, \W)) \to \Fun(\K,\Mor_\mC(\Z, \V)),$
which is monadic if $\h:\W \to \V$ is monadic.

\end{proof}

\subsection{A classification of enriched $\infty$-categories}

In this section we apply \cref{uhnggfdaa} to obtain a classification of
$\mV$-enriched $\infty$-categories with space of objects $\X$ by presentably left $\mV$-tensored $\infty$-categories equipped with a functor starting at $\X$ recovering a result of Reuter-Zetto \cite[Theorem 4.7.]{reutter2025enriched}.

\begin{notation}

Let $\mV$ be a presentably monoidal $\infty$-category.
Let $\mV \mathrm{-}\Mod^\L$ be the $(\infty,2)$-category of presentably
left $\mV$-tensored $\infty$-categories and left adjoint $\mV$-linear functors.

\end{notation}

\begin{remark}

For every $\mM, \mN \in \mV \mathrm{-}\Mod^\L$ the $\infty$-category of morphisms $\mM \to \mN$ in $\mV \mathrm{-}\Mod^\L$
is the $\infty$-category $\mV\mathrm{-}\LinFun^\L(\mM,\mN)$ of \cref{linea}.

\end{remark}

\begin{notation}
Let $\mV$ be a presentably monoidal $\infty$-category and $\X$ a small space.
We write $\mV^\X$ for the cotensor in $\mV \mathrm{-}\Mod^\L$,
which is the $\infty$-category of functors $\X \to \mV$ with diagonal left $\mV$-action. 
\end{notation}

\begin{remark}

By \cite[Proposition 4.8.1.17.]{lurie.higheralgebra} the presentably left $\mV$-tensored $\infty$-category $\mV^\X$ is the free presentably left $\mV$-tensored $\infty$-category 
generated by the $\infty$-category $\X.$
Thus for every presentably left $\mV$-tensored $\infty$-category $\mM$
the Yoneda-embedding $$ \X \to \mS^\X \xrightarrow{(-)\ot \tu_\mV} \mV^\X$$ induces an equivalence
\begin{equation}\label{lineq}
\mV\mathrm{-}\LinFun^\L(\mV^\X,\mM) \to \Fun(\X,\mM).
\end{equation}
The inverse sends a functor $\theta: \X \to \mM$ to a left adjoint $\mV$-linear functor
$\bar{\theta}: \mV^\X \to \mM$ whose right adjoint is the functor
$$\mM \to \mV^\X, \Z \mapsto \Mor_\mM(-,\Z) \circ \theta^\op,$$ which is canonically a $\mV$-enriched functor.

\end{remark}

\begin{notation}
Let $\mV$ be a presentably monoidal $\infty$-category and $\X$ a small space.

Let $$ \mathrm{Quiv}_\X(\mV):=\mV\mathrm{-}\LinFun^\L(\mV^\X,\mV^\X) $$ the endomorphism algebra of $\mV^\X$ in the $(\infty,2)$-category $\mV \mathrm{-}\Mod^\L$,
whose underlying $\infty$-category is equivalent to
$ \Fun(\X \times \X,\mV)$ by equivalence (\ref{lineq}).
\end{notation}

\begin{definition}
Let $\mV$ be a presentably monoidal $\infty$-category and $\X$ a small space.

A non-univalent $\mV$-enriched $\infty$-category with space of objects $\X$
(in the sense of Gepner-Haugseng) is an associative algebra in
$ \mathrm{Quiv}_\X(\mV)$.

We set $$ \mV \mathrm{-}\Cat^{\mathrm{nu}}_\X:= \Alg(\mathrm{Quiv}_\X(\mV)).$$

\end{definition}

\begin{remark}\label{GHAU}
By \cite[Proposition 4.5.3.]{HINICH2020107129} there is a canonical equivalence between $\mV \mathrm{-}\Cat_\X$ and the $\infty$-category of non-univalent $\mV$-enriched $\infty$-categories with space of objects $\X$ in the sense of
Hinich, which by \cite{MR4185309} is further identified with the $\infty$-category of non-univalent $\mV$-enriched $\infty$-categories with space of objects $\X$ in the sense of Gepner-Haugseng.

\end{remark}
\vspace{2mm}

\begin{lemma}\label{markmo}

Let $\mV$ be a presentably monoidal $\infty$-category, $\X$ a small space,
$\mM$ a presentably left $\mV$-tensored $\infty$-category and $\theta: \X \to \mM$ a functor. 
%Let $\bj: \theta(\X) \subset \mM$ be the $\mV$-enriched essential image of $\theta$.
The following conditions are equivalent:

\begin{enumerate}

\item The $\infty$-category $\mM$ is generated under small colimits and the left $\mV$-action by the essential image of $\theta.$
For every $\Y \in \X$ the $\mV$-enriched functor
$$\Mor_\mM(\theta(\Y),-): \mM \to \mV$$ is $\mV$-linear and preserves small colimits.

\vspace{1mm}

\item The functor
$$ \gamma: \mM \to \mV^\X, \ \Z \mapsto \Mor_\mM(-,\Z) \circ \theta^\op$$
right adjoint to $\bar{\theta}: \mV^\X \to \mM$, is conservative, preserves small colimits and is $\mV$-linear.
In other words, $ \gamma: \mM \to \mV^\X$ is a monadic morphism in 
$\mV \mathrm{-}\Mod^\L$.

\end{enumerate}

\end{lemma}

\begin{proof}

By construction of $\gamma$ condition (1) evidently implies (2).
We prove that (2) implies (1).
% By definition the $\mV$-enriched functor $\gamma$ preserves small colimits and is $\mV$-linear if and only if for every $\Y \in \X$ the $\mV$-enriched functor
% $$\Mor_\mM(\theta(\Y),-): \mM \to \mV$$ is $\mV$-linear and preserves small colimits.
% Moreover $\gamma$ is conservative if $\mM$ is generated under small colimits and the left $\mV$-action by the essential image of $\theta.$

% It remains to prove that $\mM$ is generated under small colimits and the left $\mV$-action by the essential image of $\theta$ if $\gamma$ is conservative.
% Let $\theta': \X \to \theta(\X)$ be the induced functor to the essential image.
% The $\mV$-enriched embedding $\bj: \theta(\X) \subset \mM$ uniquely extends to a
% left adjoint $\mV$-linear functor 
% $\bar{\bj}: \mP_\mV(\theta(\X)) \to \mM$.
% By \cref{enryol} the $\infty$-category $\mP_\mV(\theta(\X))$ is generated under small colimits and the left $\mV$-action by the essential image of the $\mV$-enriched Yoneda-embedding $\theta(\X) \to \mP_\mV(\theta(\X)). $
% So to see that $\mM$ is generated under small colimits and the left $\mV$-action by the essential image of $\theta$, it suffices to see that $\bar{\bj}: \mP_\mV(\theta(\X)) \to \mM$ is an equivalence.
We assume that (2) holds. By construction of $\gamma$ it suffices to prove that $\mM$ is generated under small colimits and the left $\mV$-action by the essential image of $\theta.$
Let $\mN \subset \mM$ be a full subcategory containing the essential image of $\theta$ and closed under small colimits and the left $\mV$-action.
We have to see that $\mN =\mM.$
Since $\gamma: \mM \to \mV^\X$ is monadic by assumption, $\mM$ is generated under small colimits by the free objects. So it suffices to see that
$\mN$ contains the free objects or equivalently that $\mV^\X$ is the full subcategory $\mN' \subset \mV^\X$ spanned by the objects whose image under the free functor belongs to $\mN$. Since $\mN$ contains the essential image of
$\theta: \X \to \mM$ and $\theta$ factors as $\X \to \mS^\X \xrightarrow{\tu\ot(-)} \mV^\X$ followed by the free functor $\mV^\X \to \mM$,
the full subcategory $\mN'$ contains the essential image of the functor
$\X \to \mS^\X \xrightarrow{\tu\ot(-)} \mV^\X$.
The full subcategory $\mN'$ is closed under small colimits and the left $\mV$-action since the free functor is $\mV$-linear and preserves small colimits. Thus $\mN'$ contains all objects of the form
$\Map_\X(-,\Y) \ot \V.$
So the result follows from the fact that $\mV^\X$ is generated under small colimits by the full subcategory of objects of the form $\Map_\X(-,\Y) \ot \V$
for $\Y \in \X$ and $\V \in \mV,$ which is \cite[Lemma 4.42.]{HEINE2023108941}.

\end{proof}

The following definition is due to Reuter-Zetto \cite[Definition 5.6.]{reutter2025enriched}:

\begin{definition}
Let $\mV$ be a presentably monoidal $\infty$-category.
A marked $\mV$-module is a pair $(\mM,\theta: \X \to \mM)$, where $\mM$ is a presentably left $\mV$-tensored $\infty$-category and $\theta: \X \to \mM$ is a functor starting at a small space $\X$, that satisfies the equivalent conditions of \cref{markmo}.
We call $\X$ the space of objects of $(\mM,\theta: \X \to \mM)$.

\end{definition}

\begin{notation}Let $\mV$ be a presentably monoidal $\infty$-category.
The $\infty$-category of marked $\mV$-modules $$(\mV \mathrm{-}\Mod^\L)^{\#}$$ is the full subcategory of the pullback $$\mS \times_{\Fun(\{0\},\infty\widehat{\Cat})} \Fun([1],\mV \mathrm{-}\Mod^\L)$$
spanned by the marked $\mV$-modules.
\end{notation}

\begin{remark}
Let $\mV$ be a presentably monoidal $\infty$-category and $(\mM,\theta: \X \to \mM)$ a marked $\mV$-module. The right adjoint functor
$$ \gamma: \mM \to \mV^\X, \ \Z \mapsto \Mor_\mM(-,\Z) \circ \theta^\op$$
is a monadic morphism in 
$\mV \mathrm{-}\Mod^\L$ and so $\mM$ is the $\infty$-category of algebras for a small colimits preserving $\mV$-linear monad on $\mV^\X.$
This monad is $\gamma \circ \bar{\theta}: \mV^\X \to \mV^\X$
and by (\ref{lineq}) corresponds to a functor $\X \times \X \to \mV,$ 
which is $\gamma \circ \theta: (\A,\B) \mapsto \Mor_\mM(\theta(\A),\theta(\B)).$
By \cref{GHAU} this functor $$\X \times \X \to \mV: (\A,\B) \mapsto \Mor_\mM(\theta(\A),\theta(\B))$$ is the graph of a $\mV$-enriched $\infty$-category with space of objects $\X$ in the sense of Gepner-Haugseng.

The $\infty$-category of algebras of the small colimits preserving $\mV$-linear monad $\gamma \circ \bar{\theta}$ on $\mV^\X$
is $\mP_\mV(\theta(\X)).$ 
%the presentably left $\mV$-tensored $\infty$-category 
%of $\mV$-enriched presheaves on the full $\mV$-enriched subcategory of $\mM$ spanned by the essential image of $\theta: \X \to \mM.$

\end{remark}

\begin{notation}

Let $\mV$ be a presentably monoidal $\infty$-category.
Let $ \xi: \mS \to \mV \mathrm{-}\Mod^\L$ be the functor
$\X \mapsto \mV^\X$.
For every functor $\mC \to \mV \mathrm{-}\Mod^\L$ we write 
$\xi^*\mC \to \mS$ for the pullback of $\mC \to \mV \mathrm{-}\Mod^\L$ along
$\xi.$

\end{notation}

\begin{remark}

Let $\mV$ be a presentably monoidal $\infty$-category.
The $\infty$-category $\xi^*\Fun^\mon(\mV \mathrm{-}\Mod^\L) \to \mS$
is the $\infty$-category of pairs $(\X, \gamma: \mM \to \mV^\X),$
where $\gamma$ is a monadic morphisms in $\mV \mathrm{-}\Mod^\L.$
There is a canonical equivalence
$$ \xi^*\Fun^\mon(\mV \mathrm{-}\Mod^\L)^\op \simeq (\mV \mathrm{-}\Mod^\L)^{\#}, $$$$ (\X, \gamma: \mM \to \mV^\X) \mapsto (\mM, \theta: \X \to \mM)$$
that sends a right adjoint to the restriction of its left adjoint.

\end{remark}

\begin{notation}

Let $\mV$ be a presentably monoidal $\infty$-category.
Let $$ \mV \mathrm{-}\Cat^{\mathrm{nu}}:= \xi^*\Mon(\mV \mathrm{-}\Mod^\L) \to \mS.$$

\end{notation}

\begin{remark}

The $\infty$-category $ \mV \mathrm{-}\Cat^{\mathrm{nu}} \to \mS$
is an $\infty$-category over $\mS$ whose fiber over any small space
$\X$ is $\mV \mathrm{-}\Cat^{\mathrm{nu}}_\X.$
By Gepner-Haugseng's fiberwise definition of the $\infty$-category of $\mV$-enriched $\infty$-categories one can show that 
$ \mV \mathrm{-}\Cat^{\mathrm{nu}}$ is the $\infty$-category of $\mV$-enriched $\infty$-categories in the sense of Gepner-Haugseng.

\end{remark}

Applying \cref{uhnggfdaa} to the $(\infty,2)$-category $\mV \mathrm{-}\Mod^\L$
and pulling back along $\xi: \mS \to \mV \mathrm{-}\Mod^\L$ we obtain the following corollary, which was proven by Reuter-Zetto \cite[Theorem 4.7.]{reutter2025enriched}:

\begin{corollary}\label{compmod} Let $\mV$ be a presentably monoidal $\infty$-category.
There is a canonical equivalence 
$$ \mV \mathrm{-}\Cat^{\mathrm{nu}} \simeq (\mV \mathrm{-}\Mod^\L)^{\#} $$
over $\mS$ sending a $\mV$-enriched $\infty$-category in the sense of Gepner-Haugseng to the presentably left $\mV$-tensored $\infty$-category of
$\mV$-enriched presheaves on it.

\end{corollary}

Restricting the equivalence of \cref{compmod} to the fiber over the point
we obtain the following important equivalence, which was proven by Lurie \cite[Theorem 4.8.5.11.]{lurie.higheralgebra},  and characterizes associative algebras by their $\infty$-category of modules pointed by the algebra:

\begin{corollary}\label{compmod2} Let $\mV$ be a presentably monoidal $\infty$-category.
There is a canonical equivalence 
$$\Alg(\mV) \simeq \mV \mathrm{-}\Mod^\L_{\mV/}$$
sending $\A $ to $(\RMod_\A(\mV),\A).$

\end{corollary}

\subsection{A classification of enriched $\infty$-operads}

In this section we apply \cref{uhnggfdaa} to obtain a classification of
$\mV$-enriched $\infty$-operads with space of colors $\X$ by presentably symmetric monoidal $\infty$-categories under $\mV$ equipped with a functor starting at $\X$.

\begin{notation}

Let $\mV$ be a presentably monoidal $\infty$-category.
Let $$\Cmon^\L$$ be the $(\infty,2)$-category of presentably symmetric monoidal
$\infty$-categories and left adjoint symmetric monoidal functors.

Let $$\mV \mathrm{-}\Cmon^\L$$ be the $(\infty,2)$-category 
$\Cmon_{\mV/}^\L.$
\end{notation}

We refer to objects of $\mV \mathrm{-}\Cmon^\L$ as presentably symmetric monoidal $\mV$-linear $\infty$-categories.
This terminology is also reflected in the following remark, which is \cite[Corollary 3.4.1.7.]{lurie.higheralgebra}:

\begin{remark}

Let $\mV$ be a presentably symmetric monoidal $\infty$-category.
There is a canonical equivalence
$$ \mV \mathrm{-}\Cmon^\L \simeq \Calg(\mV \mathrm{-}\Mod^\L).$$ 

\end{remark}

\begin{remark}

For every $\mM, \mN \in \Cmon^\L$ the $\infty$-category of morphisms $\mM \to \mN$ in $\Cmon^\L$
is the $\infty$-category $$\Fun^{\ot, \L}(\mM,\mN)$$ of left adjoint symmetric monoidal functors $\mM \to \mN.$

For every $\mM, \mN \in \mV \mathrm{-}\Cmon^\L$ the $\infty$-category of morphisms $\mM \to \mN$ in $\mV \mathrm{-}\Cmon^\L$
is the $\infty$-category $$ \mV \mathrm{-}\Fun^{\ot, \L}(\mM,\mN):= * \times_{\Fun^{\ot, \L}(\mV,\mN)} \Fun^{\ot, \L}(\mM,\mN) $$ of left adjoint symmetric monoidal functors $\mM \to \mN$ under $\mV.$

\end{remark}

\begin{notation}Let $\mV$ be a presentably symmetric monoidal $\infty$-category and $\mM \in \mV \mathrm{-}\Cmon^\L.$
For every $\Y_1,..., \Y_\n,\Z \in \mM$ for $\n \geq 0$ let
$$ \Mul\Mor_\mM(\Y_1,..., \Y_\n, \Z):=\Mor_\mM(\Y_1 \ot ...\ot \Y_\n, \Z)  \in \mV,$$
where on the right hand side we use the tensor product of $\mM$ and the enrichment of $\mM$ in $\mV.$

\end{notation}

\begin{notation}
Let $\X$ be a small space. We write $\Sigma(\X)$ for the free symmetric monoidal $\infty$-category on $\X$ and have $$\Sigma(\X)\simeq \coprod_{\n \geq 0}\X^{\times \n}.$$

\end{notation}

\begin{notation}Let $\mV$ be a presentably symmetric monoidal $\infty$-category.
We write $\mV^{\Sigma(\X)}$ for the free presentably smmetric monoidal
$\infty$-category under $\mV$ generated by the $\infty$-category $\X,$
which is the Day-convolution monoidal structure on the
$\infty$-category of functors $\Sigma(\X) \to \mV.$

\end{notation}

Thus for every presentably smmetric monoidal
$\infty$-category $\mM$ under $\mV$
the Yoneda-embedding $$ \X \to \mS^\X \xrightarrow{(-)\ot \tu_\mV} \mV^\X \to \mV^{\Sigma(\X)}$$ induces an equivalence
\begin{equation}\label{lineq2}
\mV \mathrm{-}\Fun^{\ot,\L}(\mV^{\Sigma(\X)},\mM) \to \Fun(\X,\mM).
\end{equation}
The inverse sends a functor $\theta: \X \to \mM$ to a left adjoint symmetric monoidal functor
$\bar{\theta}: \mV^{\Sigma(\X)} \to \mM$ under $\mV$ whose right adjoint is the functor
$$\mM \to \mV^{\Sigma(\X)} \simeq \prod_{\n \geq 0}\Fun(\X^{\times\n},\mV), \ \Z \mapsto \{\Mul\Mor_\mM(-,\Z) \circ \theta^{\times \n}\}_{\n \geq 0},$$ which is canonically a lax symmetric monoidal and $\mV$-enriched functor.

\begin{notation}\label{operat}Let $\mV$ be a presentably symmetric monoidal $\infty$-category and $\X$ a small space.

Let $$ \sSeq_\X(\mV):=\mV \mathrm{-}\Fun^{\ot,\L}(\mV^{\Sigma(\X)},\mV^{\Sigma(\X)})$$
be the endomorphism algebra of $\mV^{\Sigma(\X)}$ in the $(\infty,2)$-category $\mV\mathrm{-}\Calg^\L$, whose underlying $\infty$-category identifies with the $\infty$-category $ \Fun(\X \times \Sigma(\X),\mV)$
by equivalence (\ref{lineq2}).
\end{notation}

We make the following definition following \cite{brantner2021pd} and
\cite{Rune}:

\begin{definition}\label{GHAU2}Let $\mV$ be a presentably symmetric monoidal $\infty$-category and $\X$ a small space.

A (non-univalent) $\mV$-enriched $\infty$-operad with space of objects $\X$ is an 
associative algebra in $ \sSeq_\X(\mV).$

We set $$\mV\mathrm{-}\Op_\X^{\mathrm{nu}}:= \Alg(\sSeq_\X(\mV)).$$

\end{definition}

\vspace{2mm}

\begin{proposition}\label{markmo}

Let $\mV$ be a presentably symmetric monoidal $\infty$-category, $\X$ a small space,
$\mM$ a presentably symmetric monoidal $\mV$-linear $\infty$-category and $\theta: \X \to \mM$ a functor.
The following conditions are equivalent:

\begin{enumerate}

\item The $\infty$-category $\mM$ is generated under small colimits, tensor products and the left $\mV$-action by the essential image of $\theta.$
For every $\Y_1,..., \Y_\n \in \X$ for $\n \geq 0$ the $\mV$-enriched functor
$$\Mul\Mor_\mM(\theta(\Y_1),...,\theta(\Y_\n),-): \mM \to \mV$$ is $\mV$-linear
and preserves small colimits.

For every $\A, \B \in \mM$ and $\Y_1,..., \Y_\n \in \X$ for $\n \geq 0$ the canonical morphism 
\begin{equation}\label{erth}
\colim_{0 \leq \bi \leq \n} \Mul\Mor_\mM(\theta(\Y_1),...,\theta(\Y_\bi),\A) \ot \Mul\Mor_\mM(\theta(\Y_{\bi+1}),...,\theta(\Y_\n),\B) \to\end{equation}
$$\Mul\Mor_\mM(\theta(\Y_1),...,\theta(\Y_\n),\A \ot \B)$$
is an equivalence.

\vspace{1mm}

\item The functor
$$ \gamma: \mM \to \mV^{\Sigma(\X)} \simeq \prod_{\n \geq 0}\Fun(\X^{\times\n},\mV), \ \Z \mapsto \{\Mul\Mor_\mM(-,\Z) \circ \theta^{\times \n}\}_{\n \geq 0}$$
right adjoint to $\bar{\theta}: \mV^{\Sigma(\X)} \to \mM$, is conservative, preserves small colimits, is symmetric monoidal and $\mV$-linear.
In other words, $ \gamma: \mM \to \mV^{\Sigma(\X)}$ is a monadic morphism in 
$\mV \mathrm{-}\Calg^\L$.

\end{enumerate}

\end{proposition}

\begin{proof}

We first observe that by construction of the Day-convolution tensor product
\cite[Remark 2.2.6.15.]{lurie.higheralgebra} for every $\A, \B \in \mM$ and $\Y_1,..., \Y_\n \in \X$ for $\n \geq 0$ the canonical morphism 
$$ \gamma(\A) \ot \gamma(\B) \to \gamma(\A \ot \B)$$
evaluated at $(\Y_1,..., \Y_\n)$
identifies with the morphism (\ref{erth}).

Thus by construction of $\gamma$ condition (1) evidently implies (2).
We prove that (2) implies (1).
% By definition the $\mV$-enriched functor $\gamma$ preserves small colimits and is $\mV$-linear if and only if for every $\Y \in \X$ the $\mV$-enriched functor
% $$\Mor_\mM(\theta(\Y),-): \mM \to \mV$$ is $\mV$-linear and preserves small colimits.
% Moreover $\gamma$ is conservative if $\mM$ is generated under small colimits and the left $\mV$-action by the essential image of $\theta.$

% It remains to prove that $\mM$ is generated under small colimits and the left $\mV$-action by the essential image of $\theta$ if $\gamma$ is conservative.
% Let $\theta': \X \to \theta(\X)$ be the induced functor to the essential image.
% The $\mV$-enriched embedding $\bj: \theta(\X) \subset \mM$ uniquely extends to a
% left adjoint $\mV$-linear functor 
% $\bar{\bj}: \mP_\mV(\theta(\X)) \to \mM$.
% By \cref{enryol} the $\infty$-category $\mP_\mV(\theta(\X))$ is generated under small colimits and the left $\mV$-action by the essential image of the $\mV$-enriched Yoneda-embedding $\theta(\X) \to \mP_\mV(\theta(\X)). $
% So to see that $\mM$ is generated under small colimits and the left $\mV$-action by the essential image of $\theta$, it suffices to see that $\bar{\bj}: \mP_\mV(\theta(\X)) \to \mM$ is an equivalence.
We assume that (2) holds. By construction of $\gamma$ it suffices to prove that $\mM$ is generated under small colimits, tensor products and the left $\mV$-action by the essential image of $\theta.$
Let $\mN \subset \mM$ be a full subcategory containing the essential image of $\theta$ and closed under small colimits, tensor products and the left $\mV$-action.
We have to see that $\mN =\mM.$
Since $\gamma: \mM \to \mV^{\Sigma(\X)}$ is monadic by assumption, $\mM$ is generated under small colimits by the free objects. So it suffices to see that
$\mN$ contains the free objects or equivalently that $\mV^{\Sigma(\X)}$ is the full subcategory $\mN' \subset \mV^{\Sigma(\X)}$ spanned by the objects whose image under the free functor belongs to $\mN$. Since $\mN$ contains the essential image of
$\theta: \X \to \mM$ and $\theta$ factors as \begin{equation}\label{funzyl}
\X \to \Sigma(\X) \to \mS^{\Sigma(\X)} \xrightarrow{\tu\ot(-)} \mV^{\Sigma(\X)}\end{equation} followed by the free functor $\mV^{\Sigma(\X)} \to \mM$,
the full subcategory $\mN'$ contains the essential image of the functor
(\ref{funzyl}).
The full subcategory $\mN'$ is closed under small colimits, tensor products and the left $\mV$-action since the free functor is $\mV$-linear, symmetric monoidal and preserves small colimits. Thus $\mN'$ contains the full subcategory of $\mV^{\Sigma(\X)}$ generated by the essential image of 
the functor (\ref{funzyl}) under tensor products, which is the essential image of
the symmetric monoidal functor $\Sigma(\X) \to \mS^{\Sigma(\X)} \xrightarrow{\tu\ot(-)} \mV^{\Sigma(\X)}$.
Thus $\mN'$ contains all objects of the form
$\Map_{\Sigma(\X)}(-,\Y) \ot \V.$
So the result follows from the fact that $\mV^{\Sigma(\X)}$ is generated under small colimits by the full subcategory of objects of the form $\Map_{\Sigma(\X)}(-,\Y) \ot \V$
for $\Y \in \Sigma(\X)$ and $\V \in \mV,$ which is \cite[Lemma 4.42.]{HEINE2023108941}.

\end{proof}

\begin{definition}
Let $\mV$ be a presentably symmetric monoidal $\infty$-category.
A marked presentably symmetric monoidal $\mV$-linear $\infty$-category is a pair $(\mM,\theta: \X \to \mM)$, where $\mM$ is a presentably symmetric monoidal $\mV$-linear $\infty$-category and $\theta: \X \to \mM$ is a functor starting at a small space $\X$, that satisfies the equivalent conditions of \cref{markmo}.
We call $\X$ the space of colors of $(\mM,\theta: \X \to \mM)$.

\end{definition}

\begin{notation}Let $\mV$ be a presentably symmetric monoidal $\infty$-category.
The $\infty$-category of marked presentably symmetric monoidal $\mV$-linear
$\infty$-categories
$$(\mV \mathrm{-}\Calg^\L)^{\#}$$ is the full subcategory of the pullback $$\mS \times_{\Fun(\{0\},\infty\widehat{\Cat})} \Fun([1],\mV \mathrm{-}\Calg^\L)$$
spanned by the marked presentably symmetric monoidal $\mV$-linear
$\infty$-categories.
\end{notation}

\begin{remark}
Let $\mV$ be a presentably symmetric monoidal $\infty$-category and $(\mM,\theta: \X \to \mM)$ a marked presentably symmetric monoidal $\mV$-linear
$\infty$-category. The right adjoint functor
$$ \gamma: \mM \to \mV^{\Sigma(\X)} \simeq \prod_{\n \geq 0}\Fun(\X^{\times\n},\mV), \ \Z \mapsto \{\Mul\Mor_\mM(-,\Z) \circ \theta^{\times \n}\}_{\n \geq 0}$$
is a monadic morphism in 
$\mV \mathrm{-}\Calg^\L$ and so $\mM$ is the $\infty$-category of algebras for a small colimits preserving symmetric monoidal $\mV$-linear monad on $\mV^{\Sigma(\X)}.$
This monad is $\gamma \circ \bar{\theta}:\mV^{\Sigma(\X)}\to \mV^{\Sigma(\X)}$
and by (\ref{lineq2}) corresponds to a functor $\X \times \Sigma(\X) \to \mV,$ 
which is $$\gamma \circ \theta: (\A,\Y_1,...,\Y_\n) \mapsto \Mul\Mor_\mM(\theta(\A),\theta(\Y_1),..., \theta(\Y_\n)).$$
By \cref{GHAU} this functor is the multi-graph of a $\mV$-enriched $\infty$-operad with space of objects $\X$.

%The $\infty$-category of algebras of the small colimits preserving symmetric monoidal $\mV$-linear monad $\gamma \circ \bar{\theta}$ on $\mV^{\Sigma(\X)}$
%is $\Rep_\mO.$ 

\end{remark}

\begin{notation}

Let $\mV$ be a presentably symmetric monoidal $\infty$-category.
Let $ \xi: \mS \to \mV \mathrm{-}\Calg^\L$ be the functor
$\X \mapsto \mV^{\Sigma(\X)}$.
For every functor $\mC \to \mV \mathrm{-}\Calg^\L$ we write 
$\xi^*\mC \to \mS$ for the pullback of $\mC \to \mV \mathrm{-}\Calg^\L$ along
$\xi.$

\end{notation}

\begin{remark}

Let $\mV$ be a presentably symmetric monoidal $\infty$-category.
The $\infty$-category $$\xi^*\Fun^\mon(\mV \mathrm{-}\Calg^\L) \to \mS$$
is the $\infty$-category of pairs $(\X, \gamma: \mM \to \mV^{\Sigma(\X)}),$
where $\gamma$ is a monadic morphisms in $\mV \mathrm{-}\Calg^\L.$
There is a canonical equivalence
$$ \xi^*\Fun^\mon(\mV \mathrm{-}\Calg^\L)^\op \simeq (\mV \mathrm{-}\Calg^\L)^{\#}, $$$$ (\X, \gamma: \mM \to \mV^{\Sigma(\X)}) \mapsto (\mM, \theta: \X \to \mM)$$
that sends a right adjoint to the restriction of its left adjoint.

\end{remark}

\begin{notation}

Let $\mV$ be a presentably symmetric monoidal $\infty$-category.
Let $$ \mV\mathrm{-}\Op^{\mathrm{nu}}:= \xi^*\Mon(\mV \mathrm{-}\Calg^\L) \to \mS.$$

\end{notation}

\begin{remark}

The $\infty$-category $\mV\mathrm{-}\Op^{\mathrm{nu}} \to \mS$
is an $\infty$-category over $\mS$ whose fiber over any small space
$\X$ is $\mV\mathrm{-}\Op_\X^{\mathrm{nu}}.$

\end{remark}

Applying \cref{uhnggfdaa} to the $(\infty,2)$-category $\mV \mathrm{-}\Calg^\L$
and pulling back along $\xi: \mS \to \mV \mathrm{-}\Calg^\L$ we obtain the following result:

\begin{corollary}\label{compmod} Let $\mV$ be a presentably symmetric monoidal $\infty$-category.
There is a canonical equivalence 
$$\mV\mathrm{-}\Op^{\mathrm{nu}} \simeq (\mV \mathrm{-}\Calg^\L)^{\#} $$
over $\mS$.

%sending a $\mV$-enriched $\infty$-operad $\mO$ to the presentably symmetric monoidal $\mV$-linear $\infty$-category of representations of $\mO.$

\end{corollary}

% Restricting the equivalence of \cref{compmod} to the fiber over the point
% we obtain the following important equivalence, which was proven by Lurie \cite[Theorem 4.8.5.11.]{lurie.higheralgebra},  and characterizes associative algebras by their $\infty$-category of modules pointed by the algebra:

% \begin{corollary}\label{compmod2} Let $\mV$ be a presentably monoidal $\infty$-category.
% There is a canonical equivalence 
% $$\Alg(\mV) \simeq \mV \mathrm{-}\Calg^\L$$
% sending $\A $ to $(\RMod_\A(\mV),\A).$

% \end{corollary}

\section{A monadicity theorem for higher algebraic structures}

In this section we prove a monadicity theorem for higher algebraic structures (Theorem \ref{strmon}).

\subsection{Higher algebraic structures}\label{catpat}

%\subsection{Monads on higher algebraic structures}\label{catpat}

% For the next definition we use the notion of relative cocartesian fibration and relative limit \cite[Definition 4.3.1.1.]{lurie.HTT}:

%\begin{remark}A functor $\rH: \K^{\triangleleft} \to \mC$ corresponds toan object $\X \in \mC$ equipped with a map $\delta(\X)\to \T$ in $\Fun(\K,\mC)$,where $-\infty \in\K^{\triangleleft}$ is the cone point, $\X \simeq \rH(-\infty), \T \simeq \rH_{\mid\K}$ and $\delta$ is the diagonal functor. Since $\delta$ is left adjoint to evaluation at the cone point, we can rephrase Definition \ref{gggbv}:a map $\delta(\X)\to \T$ in $\Fun(\K,\mC)$ exhibits an object $\X \in \mC$ as the $\phi$-limit of a functor $\T:\K \to \mC$ if for every$\Y \in \mC$ the induced square
%\begin{equation*}\label{zhkq}
%\begin{xy}
%\xymatrix{
%\mC(\Y,\X) \ar[d]
%\ar[r]^{}& \Fun(\K,\mC)(\delta(\Y),\T) \ar[d]^{} 
%\\
%\rS(\phi(\Y),\phi(\X)) \ar[r]^{}  &  \Fun(\K,\rS)(\phi \circ \delta(\Y),\phi \circ \T).}\end{xy} 
%\end{equation*}is a pullback square.\end{remark}

%\begin{corollary}\label{uhnbv}Let $\phi: \mC \to \rS$ be a cocartesian fibration and $\K$ an $\infty$-category.The functor $$\Fun^\cocart_\rS(\K^{\vartriangleleft},\mC) \to \Fun^\cocart_\rS(\K,\mC) $$is fully faithful if and only if every functor $\K^{\vartriangleleft} \to \mC $ over $\rS$sending all morphisms to $\phi$-cocartesian ones is a $\phi$-limit diagram.\end{corollary}

We start by introducing operad-like structures using the framework of algebraic pattern \cite{chu2021homotopy}.
The following definition is \cite[Definition 2.1.]{chu2021homotopy}:

\begin{definition}

An algebraic pattern is a triple $(\rS,\mE, \rS^\circ)$
consisting of the following:

\begin{enumerate}

\item an $\infty$-category $\rS,$

\item a wide subcategory $\mE \subset \rS$ that is the left class of a factorization system on $\rS$.

\item a full subcategory $\rS^\circ \subset \rS$.

\end{enumerate}

\begin{notation}
Let $\mathfrak{P}=(\rS,\mE, \rS^\circ)$ be an algebraic pattern.
We call $\mE$ the class of inert morphisms of $\rS$.
We call the right class of the factorization system the class of active morphisms of $\rS$.
We call $\rS^\circ$ the full subcategory of elementary objects.

\end{notation}

\end{definition}

\begin{definition}\label{fibrous}
Let $\mathfrak{P}=(\rS,\mE, \rS^\circ)$ be an algebraic pattern.
A $\mathfrak{P}$-operad is a functor $\phi: \mC \to \rS$ such that the following conditions hold: 
\begin{enumerate}
\item The functor $\phi$ is a cocartesian fibration relative to $\mE$.

\item For every $\s \in \rS$ the canonical functor
$$ \mC_\s \to \lim_{\rt \in \rS^\circ \times_{\rS}\mE_{\s /}} \mC_\rt$$
is an equivalence.

\item For every $\s \in \rS$ and $\X \in \mC$ lying over $\s$ the canonical functor
$$(\rS^\circ \times_{\rS}\mE_{\s /})^{\triangleleft} \to \mC, \ (\alpha: \s \to \rt) \mapsto \alpha_!(\X) $$ 
induces for every $\Y \in \mC$ a pullback square
\begin{equation*} 
\begin{xy}
\xymatrix{
\Map_{\mC}(\Y,\X)
\ar[d] 
\ar[r]  & \lim_{\alpha \in \rS^\circ \times_{\rS}\mE_{\s /}}\Map_{\mC}(\Y,\alpha_!(\X))
\ar[d]
\\
\Map_{\mC}(\phi(\Y),\phi(\X)) \ar[r] & \lim_{\alpha \in \rS^\circ \times_{\rS}\mE_{\s /}}\Map_{\rS}(\phi(\Y),\phi(\alpha_!(\X))).
}
\end{xy}
\end{equation*}

%is a $\phi$-limit diagram.

\end{enumerate}

\end{definition}

\begin{notation}Let $\mathfrak{P}=(\rS,\mE, \rS^\circ)$ be an algebraic pattern
and $\phi: \mC \to \rS$ a $\mathfrak{P}$-operad.
We set $$\mC^\circ := \rS^\circ \times_\rS \mC.$$

\end{notation}

\begin{notation}
Let $$\infty\Cat_{/ \rS}^{ \mathfrak{P}  } \subset \infty\Cat_{/ \rS}^{\mE} $$ be the full subcategory of $\mathfrak{P}$-operads.
	
\end{notation}

\begin{example}

Let $\rS$ be an $\infty$-category.
Then $(\rS, \rS, \rS)$ is an algebraic pattern, which we denote by $\rS$.
Every functor $\mC \to \rS$ is fibrous for this pattern.

\end{example}

\begin{definition}
    
An algebraic pattern $\mathfrak{P}=(\rS,\mE, \rS^\circ)$ is discrete if for every $\s \in \rS$ the $\infty$-category $\rS^\circ \times_{\rS}\mE_{\s /}$
is a set.
    
\end{definition}

% \begin{notation}Let $\rS$ be an $\infty$-category, $\mathfrak{P} $ an algebraic pattern on $\rS$ and $\mC \to \rS, \mD \to \rS$ be $\mathfrak{P}$-operads.
% Let $$\Alg_{\mC/\mathfrak{P}}(\mD) \subset \Alg_{\mC/\rS}(\mD)$$ be the full subcategoy of functors over $\rS$ preserving cocartesian lifts of morphisms of $\mE$.
% Moreover we set $$\Alg_{\mathfrak{P}}(\mD):= \Alg_{\rS/\mathfrak{P}}(\mD).$$
	
% \end{notation}

We are mainly interested in the following generalization of monoidal $\infty$-categories:

\begin{definition}
Let $\rS$ be an $\infty$-category and $\mathfrak{P} $ an algebraic pattern on $\rS$.
A $\mathfrak{P}$-monoidal $\infty$-category is a cocartesian fibration $\mC \to \rS$ that is a $\mathfrak{P}$-operad.
	
\end{definition}

\begin{notation}
Let $$\infty\Cat_{/ \rS}^{\mathfrak{P}, \mon} \subset \infty\Cat_{/ \rS}^{\cocart} $$ be the full subcategory of $\mathfrak{P}$-monoidal $\infty$-categories.
	
\end{notation}

\begin{definition}Let $\rS$ be an $\infty$-category and $\mathfrak{P} $ an algebraic pattern on $\rS$.
The opposite of a $\mathfrak{P}$-monoidal $\infty$-category $\mC \to \rS$ is the fiberwise opposite cocartesian fibration of $\mC \to \rS$.
    
\end{definition}

We have the following basic operations to produce new algebraic pattern:
\begin{notation}
The product of two algebraic pattern $\mathfrak{P}= (\rS, \mE, \rS^\circ), \ \mathfrak{P}'= (\rS', \mE', \rS'^\circ) $
is the algebraic pattern $$ \mathfrak{P} \times \mathfrak{P}' := (\rS \times \rS',\mE \times \mE',\rS^\circ \times \rS' \cup \rS \times \rS'^\circ).$$
\end{notation}

We can pullback algebraic pattern along any relative cocartesian fibration:

\begin{definition}\label{nooor}

Let $\mathfrak{P}= (\rS, \mE, \rS^\circ) $ be an algebraic pattern and $\psi: \rS' \to \rS $ a cocartesian fibration relative to $\mE.$
Then $\phi^{-1}\mathfrak{P}:= (\rS', \mE', \rS'^\circ) $ is an algebraic pattern, where $\mE'$ is the wide subcategory of $\rS'$ of $\phi$-cocartesian lifts of morhisms of $\mE$ and $\rS'^\circ$ is the full subcategory of
$\rS'$ spanned by the objects of $\rS$ that lie over objects of $\rS^\circ.$ 
\end{definition}

\begin{remark}\label{recogn0}
Let $\mathfrak{P}= (\rS, \mE, \rS^\circ) $ be an algebraic pattern and $\psi: \rS' \to \rS $ a cocartesian fibration relative to $\mE.$
A functor $\mC \to \rS'$ is $\psi^{-1}\mathfrak{P}$-fibrous if and only if
it is a map in $\infty\Cat_{/ \rS}^{\mathfrak{P}}$ \cite[Proposition B.2.7.]{lurie.higheralgebra}. So there is a canonical equivalence 
$$(\infty\Cat_{/ \rS}^{\mathfrak{P}})_{/\rS'} \simeq \Cat_{\infty / \rS'}^{\mathfrak{P}'}.$$

%Hence a functor $\mC \to \mO'$ is a generalized $\mO'$-monoidal $\infty$-category if and only if it is a map of generalized $\mO$-operads and cocartesian fibration.

\end{remark}

Now we are ready to consider examples of algebraic pattern:

\subsubsection{Non-symmetric $\infty$-operads}\label{opa} 

%Next we use algebraic pattern to define non-symmetric $\infty$-operads.

% \begin{notation}
% Let $\Delta$ be the category of finite, non-empty, totally ordered sets and order preserving maps and let $\Delta^\op:= \Delta^\op.$
% \end{notation}

% \begin{definition}\label{ooop}
% A map $[\m] \to [\n]$ in $\Delta$ is inert if it is of the form $[\m] \simeq \{\bi, \bi+1,..., \bi+\m \} \subset [\n]$ for $\bi \geq 0.$
% \end{definition}

% For every $\n \geq 0$ there are $\n$ inert morphisms $[\n] \to [1]$ in $\Delta^\op$.

\begin{example}
There are two algebraic pattern on $\Delta^\op$, where $\mE \subset \Delta^\op$ is the wide subcategory of inert morphisms and $(\Delta^\op)^\circ$ is one of the following two full subcategories of $\Delta^\op:$

\begin{itemize}
\item The full subcategory spanned by $[1].$

\item The full subcategory spanned by $[1],[0].$

\end{itemize}

We call these algebraic pattern the algebraic pattern for non-symmetric $\infty$-operads and generalized non-symmetric $\infty$-operads.

\begin{example}\emph{}
	
\begin{itemize}
\item A (generalized) non-symmetric $\infty$-operad is a $\mathfrak{P}$-operad for the algebraic pattern $\mathfrak{P}$ for (generalized) non-symmetric $\infty$-operads.
		
\item A monoidal $\infty$-category is a $\mathfrak{P}$-monoidal $\infty$-category for the algebraic pattern $\mathfrak{P}$ for non-symmetric $\infty$-operads.

\item A double $\infty$-category is a $\mathfrak{P}$-monoidal $\infty$-category for the algebraic pattern $\mathfrak{P}$ for generalized non-symmetric $\infty$-operads.
		
\end{itemize}

\end{example}

\end{example}
\begin{example}Let $\bk \geq 1$.
The algebraic pattern for (generalized) $\bE_\bk$-operads is the product pattern $\bE_\bk:= (\Delta^\op)^{\times\bk}$ of the algebraic pattern for (generalized) non-symmetric $\infty$-operads.
%By convention the categorical pattern for (generalized) $\bE_\infty$-operads is the categorical pattern for (generalized) symmetric $\infty$-operads.

\end{example}

\begin{example}Let $\bk \geq 1.$
	
\begin{itemize}
\item A (generalized) $\bE_\bk$-operad is a $\mathfrak{P}$-operad for the algebraic pattern $\mathfrak{P}$ for (generalized) $\bE_\bk$-operads.
		
\item An $\bE_\bk$-monoidal $\infty$-category is a $\mathfrak{P}$-monoidal $\infty$-category for the algebraic pattern $\mathfrak{P}$ for $\bE_\bk$-operads.

\item A $\bk+1$-fold $\infty$-category is a $\mathfrak{P}$-monoidal $\infty$-category for the algebraic pattern $\mathfrak{P}$ for generalized $\bE_\bk$-operads.
		
\end{itemize}

\end{example}

\subsubsection{Symmetric $\infty$-operads}\label{opa} 

%Next we use categorical pattern to define symmetric $\infty$-operads.

\begin{notation}
Let $\Fin_*$ be the category of finite pointed sets.
\end{notation}

We write pointed finite sets as $\langle \n \rangle:= \{ \ast, 1, ..., \n\}$ for $\n \geq 0$, where $\ast$ is the base point.

\begin{definition}
A map of pointed finite sets $\langle \n \rangle \to \langle \m \rangle$ is inert if for every $1 \leq \bi \leq \m$ the fiber over $\bi$ consists precisely of one element.
%\begin{itemize}\item inert if for every $1 \leq \bi \leq \m$ the fiber of $\theta$ over $\bi$ consists precisely of one element.
%\item standard inert if $\m=1$.\item active if it sends only the base point to the base point.\end{itemize}
\end{definition}
For every $\n \geq 0$ there are $\n$ inert maps $\langle \n \rangle \to \langle 1 \rangle$, where the $\bi$-th map $ \langle \n \rangle \to \langle 1 \rangle$ for $1 \leq \bi \leq \n$ sends $\bi$ to 1.

% \begin{definition}\label{ooop}
% A map $[\m] \to [\n]$ in $\Delta$ is inert if it is of the form $[\m] \simeq \{\bi, \bi+1,..., \bi+\m \} \subset [\n]$ for $\bi \geq 0.$
% \end{definition}

%For every $\n \geq 0$ there are $\n$ inert morphisms $[\n] \to [1]$ in $\Delta^\op$.

\begin{example}
There are two algebraic pattern on $\Fin_*$, where $\mE \subset \Fin_*$ is the wide subcategory of inert morphisms and $\Fin_*^\circ$ is one of the following two full subcategories of $\Fin_*:$

\begin{itemize}
\item The full subcategory spanned by $\langle 1 \rangle.$

\item The full subcategory spanned by $\langle 1 \rangle, \langle 0 \rangle.$

\end{itemize}

We call these algebraic pattern the algebraic pattern for symmetric $\infty$-operads and generalized symmetric $\infty$-operads.

\end{example}

\begin{example}\emph{}

\begin{itemize}
\item A (generalized) symmetric $\infty$-operad is a $\mathfrak{P}$-operad for the algebraic pattern $\mathfrak{P}$ for (generalized) symmetric $\infty$-operads.
	
\item A symmetric monoidal $\infty$-category is a $\mathfrak{P}$-monoidal $\infty$-category for the algebraic pattern $\mathfrak{P}$ for symmetric $\infty$-operads.
		
\end{itemize}
	
\end{example}

%\subsubsection{$\mO$-operads}

\begin{example}Let $\psi: \mO \to \Fin_*$ a cocartesian fibration relative to the collection of inert morphisms. Let $\mathfrak{P}$ be the algebraic pattern for (generalized) symmetric $\infty$-operads.
The algebraic pattern for (generalized) $\mO$-operads is $\psi^{-1}\mathfrak{P}.$ 

\end{example}

\begin{example}Let $\psi: \mO \to \Fin_*$ a cocartesian fibration relative to the collection of inert morphisms.
	
\begin{itemize}
\item A (generalized) $\mO$-operad is a $\mathfrak{P}$-operad for the algebraic pattern $\mathfrak{P}$ for (generalized) $\mO$-operads.
	
\item An $\mO$-monoidal $\infty$-category is a $\mathfrak{P}$-monoidal $\infty$-category for the algebraic pattern $\mathfrak{P}$ for $\mO$-operads.

% \item An $\mO$-Segal $\infty$-category is a $\mathfrak{P}$-monoidal $\infty$-category for the algebraic pattern $\mathfrak{P}$ for generalized $\mO$-operads.
		
\end{itemize}
	
\end{example}

\subsubsection{Weakly enriched $\infty$-categories}\label{weten}
\begin{definition}
The algebraic pattern for weak enrichment on $\Delta^\op$
is the algebraic pattern, where $\mE$ is the wide subcategory $\max \subset \Delta^\op$ of inert morphisms that preserve the maximum, and
$(\Delta^\op)^\circ$ is the full subcategory spanned by $[0]$.

% $\mK$ is the collection of functors $[1] \to \Delta^\op$ classifying the morphism
% $[0] \simeq \{\n\} \subset [\n]$ for some $\n \geq 0.$

% \vspace{1mm}

% The categorical pattern for weak right enrichment on $\Delta^\op$
% is the categorical pattern, where $\mE$ is the full subcategory $\min \subset \Fun([1],\Delta^\op)$ of inert morphisms that preserve the minimum, and
% $\mK$ is the collection of functors $[1] \to \Delta^\op$ classifying the morphism
% $[0] \simeq \{0\} \subset [\n]$ for some $\n \geq 0.$

\end{definition}

\begin{definition}
	
Let $\mV$ be a monoidal $\infty$-category.
The algebraic pattern for weak enrichment on $\mV$ is the categorical algebraic induced by $\mV^\ot \to \Delta^\op$ from the algebraic pattern for weak enrichment.
\end{definition}

\begin{example}\label{wla}
Let $\mV$ be a monoidal $\infty$-category.
An $\infty$-category weakly enriched in $\mV$ is a $\mathfrak{P}$-operad for the algebraic pattern $\mathfrak{P}$ for weak enrichment on $\mV$.

\end{example}

Next we define monads on higher algebraic structures.

\begin{definition}
Let $\rS$ be an $\infty$-category and $\mathfrak{P}$ an algebraic pattern on $\rS$.

\begin{itemize}
\item A $\mathfrak{P}$-operadic (co)monad is a (co)monad in the $(\infty,2)$-category
$\infty\Cat_{/ \rS}^{\mathfrak{P}}.$

\item A $\mathfrak{P}$-monoidal (co)monad
is a (co)monad in $\infty\Cat_{/ \rS}^{\mathfrak{P}, \mon}.$

\end{itemize}
\end{definition}

\begin{definition}Let $\mV$ be a monoidal $\infty$-category.
A $\mV$-enriched (co)monad is a (co)monad in $\mV \mathrm{-}\omega\Enr.$

\end{definition}

\begin{example}
Let $\mV$ be a monoidal $\infty$-category.
A $\mV$-enriched (co)monad is precisely a $\mathfrak{P}$-operadic (co)monad for the algebraic pattern for weak enrichment on $\mV.$   

\end{example}

% \begin{definition}Let $\mV$ be a monoidal $\infty$-category.
	
% \begin{itemize}
% \item Let $\mM$ be a left $\mV$-tensored $\infty$-category.
% A $\mV$-linear (co)monad on $\mM$
% is a (co)monad on $\mM$ in $\mV\mathrm{-}\LMod.$
		
% %\item A $\mV$-enriched (co)monad on $\mM^\circledast \to \mV^\ot$is a (co)monad on $\mM^\circledast \to \mV^\ot$ in $\mV \mathrm{-}\omega\Enr.$
		
% %\item An op$\mV$-enriched (co)monad on $\mM^\circledast \to \mV^\ot$is a (co)monad on $\mM^\circledast \to \mV^\ot$ in $\LMod^\oplax_\mV.$.

% % \item 
% % A $\mV$-enriched (co)monad 
% % is a (co)monad in $\mV \mathrm{-}\omega\Enr,$
% % which is a $\mathfrak{P}$-operadic (co)monad for the algebraic pattern for weak enrichment on $\mV.$

% \item Let $\mM $ be a $\mV$-enriched $\infty$-category.
% A $\mV$-enriched (co)monad on $\mM  $
% is a $\mV$-enriched (co)monad on $\mM$ in $\mV \mathrm{-}\Cat.$
% \end{itemize}
% %\item A $\mV$-enriched (co)monad on $\mM^\circledast \to \mV^\ot $is a $\mV$-enriched (co)monad on $\mM^\circledast \to \mV^\ot $such that $\mM^\circledast \to \mV^\ot$ is a $\mV$-enriched $\infty$-category\end{itemize}
% \end{definition}

\begin{notation}

Let $\rS$ be an $\infty$-category and $\mathfrak{P}$ an algebraic pattern on $\rS$. Let $$\infty\Cat_{/ \rS}^{\mathfrak{P}, \mon, \oplax}$$
be the pullback of the $(\infty,2)$-category $(\infty\Cat_{/ \rS}^{\mathfrak{P}, \mon, \lax})^\co$
along the map of spaces $$ (\infty\Cat_{/ \rS}^{\mathfrak{P}, \mon})^\simeq \simeq (\infty\Cat_{/ \rS}^{\mathfrak{P}, \mon})^\simeq$$
that forms the opposite $\mathfrak{P}$-monoidal $\infty$-category.

%be the source of the equivalence of $(\infty,2)$-categories $$\infty\Cat_{/ \rS}^{\mathfrak{P}, \mon, \oplax} \simeq (\infty\Cat_{/ \rS}^{\mathfrak{P}, \mon, \lax})^\mathrm{co}$$ lifting the equivalence of spaces $(\infty\Cat_{/ \rS}^{\mathfrak{P}, \mon})^\simeq\simeq (\infty\Cat_{/ \rS}^{\mathfrak{P}, \mon})^\simeq$ taking the opposite cocartesian fibration. 

\end{notation}

\begin{remark}
By definition $\infty\Cat_{/ \rS}^{\mathfrak{P}, \mon, \oplax}$
is an $(\infty,2)$-category canonically equivalent to the $(\infty,2)$-category $(\infty\Cat_{/ \rS}^{\mathfrak{P}, \mon, \lax})^\co.$
The difference is only that we identify 
$\infty\Cat_{/ \rS}^{\mathfrak{P}, \mon, \oplax}$ and $(\infty\Cat_{/ \rS}^{\mathfrak{P}, \mon, \lax})^\co$ not by the identity but by
viewing $\mathfrak{P}$-monoidal $\infty$-categories in
$\infty\Cat_{/ \rS}^{\mathfrak{P}, \mon, \oplax}$ as their opposite
$\mathfrak{P}$-monoidal $\infty$-categories in
$\infty\Cat_{/ \rS}^{\mathfrak{P}, \mon, \lax}.$

\end{remark}

\begin{definition}\label{monpat}

Let $\rS$ be an $\infty$-category, $\mathfrak{P}$ an algebraic pattern on $\rS$
and $\mC \to \rS $ a $\mathfrak{P}$-monoidal $\infty$-category.

\begin{itemize}
\item A lax $\mathfrak{P}$-monoidal (co)monad on $\mC \to \rS$
is a (co)monad on $\mC \to \rS$ in $\infty\Cat_{/ \rS}^{\mathfrak{P}, \mon, \lax}.$

\item An oplax $\mathfrak{P}$-monoidal (co)monad on $\mC \to \rS$
is a (co)monad on $\mC \to \rS$ in $\infty\Cat_{/ \rS}^{\mathfrak{P}, \mon, \oplax}.$
\end{itemize}

\end{definition}

\begin{remark}\label{Remmmm}
	
An oplax $\mathfrak{P}$-monoidal comonad on $\mC \to \rS$
is by definition a lax $\mathfrak{P}$-monoidal monad on the fiberwise opposite cocartesian fibration of $\mC \to \rS.$
An oplax $\mathfrak{P}$-monoidal monad on $\mC \to \rS$
is by definition a lax $\mathfrak{P}$-monoidal comonad on the opposite cocartesian fibration of $\mC \to \rS.$

\end{remark}

%\begin{definition}Let $\n \geq 1$ and $\mV^\ot \to \bE_\n$ a $\n$-monoidal $\infty$-category.

%\begin{itemize}\item A $\n$-monoidal (co)monad on $\mV^\ot \to \bE_\n$is a (co)monad on $\mV^\ot \to \bE_\n$ in $\Op_\infty^{\bE_\n, \mon}.$
		
%\item A lax $\n$-monoidal (co)monad on $\mV^\ot \to \bE_\n$is a (co)monad on $\mV^\ot \to \bE_\n$ in $\Op_\infty^{\bE_\n, \mon, \lax}.$
		
%\item An oplax $\n$-monoidal (co)monad on $\mV^\ot \to \bE_\n$is a (co)monad on $\mV^\ot \to \bE_\n$ in $\Op_\infty^{\bE, \mon, \oplax}.$\end{itemize}\end{definition}

\begin{definition}
	
Let $\psi: \mO \to \Fin_*$ be a symmetric $\infty$-operad.
A (lax, oplax) $\mO$-monoidal (co)monad is a (lax, oplax) $\mathfrak{P}$-monoidal (co)monad for the algebraic pattern for $\mO$-operads.
	
\end{definition}

\begin{lemma}

Let $\F: \infty\Cat \to \mC$ be a left adjoint monoidal functor
that exhibits $\mC$ as enriched in $\infty\Cat$.
For every associative algebra $\X $ in $\mC$ there is a canonical monoidal functor
$\Mor_\mC(\tu,\X) \to \Mor_\mC(\X,\X)$
that sends $\alpha: \tu \to \X$ to $\X \simeq \tu \ot \X \xrightarrow{\alpha \ot \X} \X \ot \X \to \X,$
where the monoidal structure on the left hand side is induced by the algebra structure and the monoidal structure on the right hand side is the endomorphism monoidal structure.

\end{lemma}

\begin{proof}
%Since $\F \simeq (-)\ot \tu: \infty\Cat \to \mC$,the $\infty\Cat$-linear 

The functor $\F$ is right adjoint to the functor $\Mor_\mC(\tu,-),$ 
which gets canonically lax monoidal.
Hence the monoidal adjunction $\F:  \infty\Cat \rightleftarrows \mC:\Mor_\mC(\tu,-)$ induces an adjunction on associative algebras
so that the counit $\F(\Mor_\mC(\tu,\X)) \to \X$ lifts to a map of associative algebras in $\mC.$
Restricting the left action of $\X$ on itself along the counit we obtain a left action of $\F(\Mor_\mC(\tu,\X))$ on $\X$ in $\mC.$
Since $\F$ exhibits $\mC$ as enriched in $\infty\Cat$,
this left action corresponds to a monoidal functor $\Mor_\mC(\tu,\X) \to \Mor_\mC(\X,\X).$

\end{proof}

\begin{corollary}\label{polk}
Let $\rS$ be an $\infty$-category, $\mathfrak{P}$ an algebraic pattern on $\rS$
and $\phi: \mC \to \rS$ a monoid in $\infty\Cat^{\mathfrak{P}}_{/\rS}$.
There is a canonical monoidal functor 
$\Alg_{\mathfrak{P}}(\mC) \to \Alg_{\mC/\mathfrak{P}}(\mC),$
where the monoidal structure on the left hand side is induced by the monoid structure and the monoidal structure on the right hand side is the endomorphism monoidal structure.
So for every associative algebra $\A$ in $\Alg_{\mathfrak{P}}(\mC)$
the functor $\A \ot (-): \mC \to \mC$ over $\rS$ refines to a $\mathfrak{P}$-operadic monad.

\end{corollary}

%\begin{corollary}\label{Coopr}Let $\mO$ be an algebraic pattern and $ \mC \to \mO$ a monoid in $\Op_\infty^{\mO,\mon}$.There is a canonical monoidal functor $\Alg_{\mO}(\mC) \to \Alg_{\mC/\mO}(\mC),$where the monoidal structure on the left hand side is induced by the monoid structure and the monoidal structure on the right is the endomorphism monoidal structure.So for every associative algebra $\A$ in $\Alg_{\mO}(\mC)$the lax $\mO$-monoidal functor $\A \ot (-): \mC \to \mC$ is a lax $\mO$-monoidal monad.

%Note that $\Mon(\Mon_\mO) \simeq \Mon_{{\Delta^\op} \otimes \mO}$and $\Alg(\Fun^{\ot,\lax}_{\mO}(\mO,\mC)) \simeq \Alg_{{\Delta^\op} \ot \mO}(\mC).$\end{corollary}

\begin{definition}Let $\mV$ be a monoidal $\infty$-category and 
$\mM$ a left $\mV$-tensored $\infty$-category.
A $\mV$-linear (co)monad on $\mM$ is a (co)monad on $\mM$ in $\mV\mathrm{-}\LMod.$
\end{definition}

\begin{remark}\label{Fama}Let $\mV, \mW$ be monoidal $\infty$-categories.
We apply Example \ref{exck} to the $(\infty,2)$-category $\mW\mathrm{-}\R\Mod$ of \cref{enr2cat}.

\begin{enumerate}
% \item A weakly bienriched $\infty$-category $ \mM^\circledast \to \mV^\ot \times \mW^\ot $ that exhibits $\mM$ as left tensored over $\mV$
% classifies and is classified by a left $\mV$-module in $\mW^\rev\mathrm{-}\omega\Enr$
% \cite[Proposition 3.66.]{HEINE2023108941}.
% % with respect to the left action of $\infty\Cat$ on $\omega\mW\mathrm{-}\RMod$.
% We obtain a %The latter is the pullback of the endomorphism left action of $\Enr\Fun_{\mW}(\mM,\mM)$ on $\mM_{[0]} \to \mW^\ot$ along a canonical 
% monoidal functor $$ \mV \to\mW^\rev\mathrm{-}\Fun(\mM,\mM)$$ that sends an associative algebra $\A$ in $\mV$ to a monad $\A \ot (-)$ on the underlying weakly $\mW^\rev$-enriched $\infty$-category of $\mM$. Consequently, there is a canonical $\mW^\rev$-enriched equivalence $$\LMod_\A(\mM) \simeq \LMod_{\A \ot (-)}(\mM).$$

\item A $\mV, \mW$-bitensored $\infty$-category $ \mM$ corresponds to a left $\mV$-module in $\mW\mathrm{-}\R\Mod$
\cite[Proposition 3.66.]{HEINE2023108941}.
% with respect to the left action of $\infty\Cat$ on $\omega\mW\mathrm{-}\RMod$.
We obtain a %The latter is the pullback of the endomorphism left action of $\Enr\Fun_{\mW}(\mM,\mM)$ on $\mM_{[0]} \to \mW^\ot$ along a canonical 
monoidal functor $$ \mV \to\mW^\rev\mathrm{-}\LinFun(\mM,\mM)$$ that sends an (co)associative algebra $\A$ in $\mV$ to a $\mW$-linear (co)monad $\A \ot (-)$ on the underlying right $\mW$-tensored $\infty$-category of $\mM$. 

%Consequently, there is a canonical $\mW^\rev$-linear equivalence $$\LMod_\A(\mM) \simeq \LMod_{\A \ot (-)}(\mM).$$

% \item If $ \mM$ is an $\infty$-category bitensored over $\mV,\mW$, the monoidal functor $ \mV \to\mW\mathrm{-}\Fun(\mM,\mM)$
% lands in $\mW\mathrm{-}\LinFun(\mM,\mM)$ and there is a canonical equivalence $$\LMod_\A(\mM) \simeq \LMod_{\A \ot (-)}(\mM)$$ in $\mW\mathrm{-}\RMod.$

\item The monoidal functor $\mV \to \mV \mathrm{-} \LinFun(\mV,\mV)$ associated to the biaction of $\mV$ on itself, is an equivalence by \cite[Corollary 4.2.4.7.]{lurie.higheralgebra}.

\end{enumerate}

\end{remark}

Next we consider Eilenberg-Moore objects of higher algebraic structures.
The next corollary follows immediately from \cref{patty}.
\begin{corollary}\label{tgvwxlkm0}

Let $\rS$ be an $\infty$-category and $\mathfrak{P}$ an algebraic pattern on $\rS$.

\begin{enumerate}
\item For every $\mathfrak{P}$-operad $ \mC \to \rS$ and 
$\mathfrak{P}$-operadic monad $\T$ on $\mC \to \rS$ the forgetful functor $$ \LMod^{\rS}_\T(\mC) \to \mC $$
is an Eilenberg-Moore object for $\T $ in ${\infty\Cat}^{\mathfrak{P}}_{/\rS}.$

\item For every $\mathfrak{P}$-operad $ \mC \to \rS$ and $\mathfrak{P}$-operadic comonad $\R $ on $\mC \to \rS$ the forgetful functor $$\coLMod^{\rS}_\R(\mC) \to \mC $$
is a coEilenberg-Moore object for $\R $ in ${\infty\Cat}^{\mathfrak{P}}_{/\rS}.$

\item For every $\mathfrak{P}$-monoidal $\infty$-category $ \mC \to \rS$ and 
oplax $\mathfrak{P}$-monoidal monad $\T$ on $\mC \to \rS$ the forgetful functor $$ \widetilde{\LMod}^\rS_\T(\mC):= \coLMod^\rS_\T(\mC^\rev)^\rev\to \mC $$ is an Eilenberg-Moore object for $\T $ in ${\infty\Cat}^{\mathfrak{P},\mon,\oplax}_{/\rS}.$

\item For every $\mathfrak{P}$-monoidal $\infty$-category $ \mC \to \rS$ and lax $\mathfrak{P}$-monoidal comonad $\R $ on $\mC \to \rS$ the forgetful functor $$\coLMod^{\rS}_\R(\mC) \to \mC $$
is a coEilenberg-Moore object for $\R $ in ${\infty\Cat}^{\mathfrak{P},\mon,\lax}_{/\rS}.$

\end{enumerate}

\end{corollary}

\cref{tgvwxlkm0} specializes to the following two corollaries:

\begin{corollary}\label{exkq} Let $\mO$ be a (generalized) non-symmetric  $\infty$-operad and $\mC$ an $\mO$-monoidal $\infty$-category. 
\begin{enumerate}
\item Let $\T$ be a monad on $\mC $ in $\infty\Op^{\mO}$.
The functor $$\LMod_\T(\mC):= \LMod_\T^{\mO}(\mC)\to \mO$$ is an $\mO$-monoidal $\infty$-category. 
		
\vspace{1mm}
\item Let $\T$ be an oplax $\mO$-monoidal monad on $\mC. $
%(i.e. a comonad on $(\mC^\ot)^\rev \to \mO^\ot$ in $\infty\Op^{\mO}$).
The functor $$\widetilde{\LMod}_\T(\mC):= \LMod_\T^{\mO}(\mC^\rev)^\rev\to \mO$$ is an $\mO$-monoidal $\infty$-category. 
		
\end{enumerate}
\end{corollary}

\begin{corollary}\label{exk}\label{eNr} Let $\mV$ be a monoidal $\infty$-category.

\begin{enumerate}
\item Let $\mM$ be a weakly left $\mV$-enriched $\infty$-category and $\T$ a $\mV$-enriched monad on $\mM $.
The functor $$\LMod_\T(\mM)^\circledast:= \LMod_\T^{\mV^\ot}(\mM^\circledast)\to \mV^\ot$$ is a weakly $\mV$-enriched $\infty$-category.

\vspace{1mm}
\item Let $\mM$ be a left $\mV$-tensored $\infty$-category
and $\T$ a $\mV$-linear monad on $\mM$.
The functor $$\LMod_\T(\mM)^\circledast:= \LMod_\T^{\mV^\ot}(\mM^\circledast)\to \mV^\ot$$ is a left $\mV$-tensored $\infty$-category.

\item Let $\mM$ be a pseudo-$\mV$-enriched $\infty$-category. The functor $$\LMod_\T(\mM)^\circledast:= \LMod_\T^{\mV^\ot}(\mM^\circledast)\to \mV^\ot$$ is a pseudo-$\mV$-enriched $\infty$-category.

\vspace{1mm}
\item Let $\mM$ be a $\mV$-enriched $\infty$-category. 
If $\mV$ admits totalizations, the functor $$\LMod_\T(\mM)^\circledast:= \LMod_\T^{\mV^\ot}(\mM^\circledast)\to \mV^\ot$$ is a $\mV$-enriched $\infty$-category.
\end{enumerate}
\end{corollary}

\begin{proof}

(1) and (2) follow from \cref{tgvwxlkm0}.
We prove (3) and (4).
By \cite[Corollary 3.61.]{heine2024bienriched} the respective full subcategories of the 2-category $\mV\mathrm{-}\omega\Enr$
of weakly $\mV$-enriched $\infty$-categories spanned by the $\mV$-enriched
$\infty$-categories and spanned by the pseudo-$\mV$-enriched $\infty$-categories are both accessible 2-localizations, and so closed under Eilenberg-Moore objects by \cref{rembras0}.
Hence (3) and (4) follow from (1).
\end{proof}

We also have a statement about Kleisli objects of higher algebraic structures:

% Next we study Kleisli objects in the $(\infty,2)$-category $\infty\Cat_{/ \rS}$ (\cref{tgvwxlkm}) and more generally in the $(\infty,2)$-category $\infty\Cat_{/ \rS}^{\mathfrak{P}}$ of $\mathfrak{P}$-operads for every discrete algebraic pattern $\mathfrak{P}$ on $\rS$.

\begin{corollary}\label{coKlei0}

Let $\rS$ be an $\infty$-category, $\mathfrak{P}$ a discrete algebraic pattern on $\rS$
and $ \mC \to \rS$ a $\mathfrak{P}$-operad.

\begin{enumerate}
\item For every $\mathfrak{P}$-operadic monad $\T$ on $\mC \to \rS$ let $ \LMod^{\rS}_\T(\mC)' \subset \LMod^{\rS}_\T(\mC)$ be the essential image of the free functor.
The functor $\mC \to \LMod^{\rS}_\T(\mC)'$ is a Kleisli-object in $\infty\Cat^{\mathfrak{P}}_{/\rS}.$

\vspace{1mm}

\item For every $\mathfrak{P}$-operadic comonad $\R $ on $\mC \to \rS$ 
let $ \coLMod^{\rS}_\R(\mC)' \subset \coLMod^{\rS}_\R(\mC)$ be the essential image of the cofree functor. The functor $\mC \to \coLMod^{\rS}_\R(\mC)'$ is a coKleisli-object in $\infty\Cat^{\mathfrak{P}}_{/\rS}.$
\end{enumerate}

\end{corollary}

\begin{proof}

This follows immediately from \cref{coKlei} since for every discrete algebraic pattern on $\rS$
the essential image of any morphism of fibrous objects is again fibrous.

\end{proof}

\subsection{Monadicity of higher algebraic structures}\label{monpara}

Next we prove a monadicity theorem for higher algebraic structures (Theorem \ref{strmon}).

\begin{lemma}\label{splil}
Let $\mC \to \rS$ be a cartesian fibration and $\s \in \rS$.
The functor $\mC_\s \to \mC$ preserves weakly contractible colimits.	

\end{lemma}

\begin{proof}

Let $\rH: \K^{\triangleright} \to \mC_\s$ be a colimit diagram and $\K$ weakly contractible.
By \cite[Proposition 4.3.1.12.]{lurie.HTT} and \cite[Proposition 4.3.1.5. (2)]{lurie.HTT} the constant functor $\K^{\triangleright} \to * \to \rS$ is a colimit diagram since $\K$ is weakly contractible.
Hence the functor $\K^{\triangleright} \to \mC_\s \to \mC$ lifts a colimit diagram
and so by \cite[Corollary 4.3.1.16.]{lurie.HTT} and \cite[Proposition 4.3.1.5. (2)]{lurie.HTT} is a colimit diagram.

%This is equivalent to say that for every $\Z \in \mC$ the canonical map$$ \alpha: \mC(\rH(-\infty),\Z) \to \lim (\mC(-,\Z) \circ \rH_{\mid\K}) $$ is an equivalence.If $\Z $ lies over $\rt \in \rS$, the latter map is a map over $ \rS(\s,\rt)$using that $\K$ is weakly contractible.
%where $\K^\sim \simeq \colim_\K *$ is the space arising from $\K$ by formally inverting all morphisms.The map $\alpha$ induces on the fiber over any morphism $\varphi: \s \to \rt$the map $$ \mC_\s(\rH(-\infty),\varphi^*(\Z)) \to \lim (\mC_\s(-,\varphi^*(\Z)) \circ \rH_{\mid\K}),$$which is an equivalence by assumption.

\end{proof}

\begin{lemma}\label{splitt}

Let $\rS$ be an $\infty$-category, $\mC \to \rS$ a functor and $\T$ a monad on $\mC \to \rS$ in $\infty\Cat_{/ \rS}$. Let $\nu: \LMod^\rS_\T(\mC) \to \mC$ be the forgetful functor.
For every $\s \in \rS$ every $\nu_\s$-split simplicial object of $\LMod^\rS_\T(\mC)_\s$ admits a colimit that is preserved by $\nu_\s$ and the functor $\LMod^\rS_\T(\mC)_\s \to \LMod^\rS_\T(\mC).$

\end{lemma}

\begin{proof}

Assume first that $\mC \to \rS$ is a cartesian fibration whose fibers admit small colimits and $\T: \mC \to \mC$ preserves fiberwise small colimits.
Then for every $\s \in \rS$ by \cite[Corollary 4.2.3.5.]{lurie.higheralgebra} the $\infty$-category $\LMod_\T^\rS(\mC)_\s \simeq \LMod_{\T_\s}(\mC_\s)$ admits small colimits and the forgetful functor
$\nu_\s: \LMod_{\T_\s}(\mC_\s) \to \mC_\s$ preserves small colimits.
Since $\mC \to \rS$ is a cartesian fibration, by \cref{zhgbfdc} (1) 
the functor $\LMod_\T^\rS(\mC) \to \rS$ is a cartesian fibration.
So by Lemma \ref{splil} the functor $\LMod_\T^\rS(\mC)_\s \to \LMod_\T^\rS(\mC)$ preserves small colimits. So the claim follows.

For the general case we construct a cartesian fibration $\mD \to \rS$
whose fibers admit small colimits and a monad $\P$ on $\mD \to \rS$ that preserves fiberwise small colimits and a pullback square over $\rS$:
\begin{equation*} 
\begin{xy}
\xymatrix{
\LMod^\rS_\T(\mC) \ar[d] 
\ar[r]^{  } 
&\LMod^\rS_{\P}(\mD) \ar[d] 
\\
\mC  \ar[r]^{ }  & \mD.
}
\end{xy} 
\end{equation*} 
This will imply the claim.	
Let $\Env^\rS(\mC) \to \rS$ be the functor $$ \Fun([1], \rS) \times_{\Fun(\{1\}, \rS)} \mC \to \Fun([1], \rS) \to \Fun(\{0\}, \rS),$$ which is a cartesian fibration.
The diagonal embedding $\rS \to \Fun([1], \rS)$ yields an embedding $\mC \to \Env^\rS(\mC)$ over $\rS$.
By \cite[Proposition 2.2.4.9.]{lurie.higheralgebra} for every cartesian fibration $\mB \to \rS$ restriction along $\mC \subset \Env^\rS(\mC)$ induces an equivalence
$ \Fun_\rS^\cart(\Env^\rS(\mC), \mB) \to \Fun_\rS(\mC,\mB).$

Let $\infty\Cat^{\cart, \rc\rc}_{/\rS} \subset \infty\widehat{\Cat}^\cart_{/\rS}$
be the subcategory of cartesian fibrations whose fibers admit small colimits and whose fiber transports preserve small colimits and maps of cartesian fibrations over $\rS$ that preserve fiberwise small colimits.
The equivalence $ \infty\widehat{\Cat}^\cart_{/\rS} \simeq  \Fun(\rS^\op,\infty\widehat{\Cat})$ restricts to an equivalence
$$\infty\Cat^{\cart, \rc\rc}_{/\rS} \simeq \Fun(\rS^\op, \infty\Cat^{\rc\rc}).$$

By \cite[Corollary 5.3.6.10.]{lurie.HTT} the inclusion $\infty\Cat^{\rc\rc} \subset \infty\widehat{\Cat}$ admits a left adjoint, which by \cite[Theorem 5.1.5.6]{lurie.HTT} sends a small $\infty$-category $\mC$ to $\mP(\mC),$
where the unit $\mC \to \mP(\mC)$ is the Yoneda-embedding.
So the inclusion $\infty\Cat^{\cart, \rc\rc}_{/\rS} \simeq \Fun(\rS^\op, \infty\Cat^{\rc\rc}) \subset \Fun(\rS^\op,\infty\widehat{\Cat}) \simeq \infty\widehat{\Cat}^\cart_{/\rS}$ admits a left adjoint,
which we denote by $\mP^\rS(-).$
By adjointness for every functor $\mB \to \rS$ and cartesian fibration $\mD \to \rS$ whose fibers admit small colimits and whose fiber transports preserve small colimits the induced functor 
$ \Fun^{\cart,\rc\rc}_\rS(\mP^\rS(\mB), \mD) \to \Fun^\cart_\rS(\mB,\mD)$ is an equivalence. 
So the functor $\mB \to \Env^\rS(\mB) \to \mP^\rS(\Env^\rS(\mB))$ induces an equivalence $ \Fun^{\cart,\rc\rc}_\rS(\mP^\rS(\Env^\rS(\mB)), \mD) \to \Fun_\rS(\mC,\mD).$
So the inclusion $\infty\Cat^{\cart, \rc\rc}_{/\rS}\subset \infty\widehat{\Cat}^\cart_{/\rS} \subset \infty\widehat{\Cat}_{/\rS}$ of $(\infty,2)$-categories admits a $\infty\Cat$-enriched left adjoint.
The endomorphism left action of $\Fun_\rS(\mC,\mC)$ on $\mC \to \rS$
gives rise to a left action of $\Fun_\rS(\mC,\mC)$ on $\mP^\rS(\Env^\rS(\mC)) \to \rS$
that is the pullback along a canonical monoidal functor
$\Fun_\rS(\mC,\mC) \to \Fun_\rS^{\cart, \rc\rc}(\mP^\rS(\Env^\rS(\mC)), \mP^\rS(\Env^\rS(\mC)))$ lifting the canonical embedding $$\Fun_\rS(\mC,\mC) \subset \Fun_\rS(\mC,\mP^\rS(\Env^\rS(\mC))) \simeq \Fun_\rS^{\cart, \rc\rc}(\mP^\rS(\Env^\rS(\mC)), \mP^\rS(\Env^\rS(\mC))).$$
The unit $\mC \to \mP^\rS(\Env^\rS(\mC))$, which is a $\Fun_\rS(\mC,\mC)$-linear embedding, gives rise to a pullback square:
\begin{equation*}\label{hhgbhhjvcf} 
\begin{xy}
\xymatrix{
\LMod^\rS_\T(\mC) \ar[d] 
\ar[r]^{  } 
&\LMod^\rS_{\mP^\rS(\Env^\rS(\T))}(\mP^\rS(\Env^\rS(\mC))) \ar[d] 
\\
\mC  \ar[r]^{ }  & \mP^\rS(\Env^\rS(\mC)).
}
\end{xy} 
\end{equation*} 

\end{proof}

\begin{lemma}\label{lemu2}

Let $\F,\G:\mJ \to \infty\Cat$ be functors and $\alpha: \F \to \G$ a natural transformation such that for every $\Z \in \mJ$
the functor $\alpha_\Z: \F(\Z) \to \G(\Z)$ admits a left adjoint $\beta^\Z$
and for every map $\kappa: \Y \to \Z$ the canonical functor
$\beta^\Z \circ \G(\kappa) \to \F(\kappa) \circ \beta_\Y$ is an equivalence.
If the functor $\alpha_\Z: \F(\Z) \to \G(\Z)$ is monadic for every $\Z \in \mJ$, the induced functor $\lim(\alpha): \lim(\F) \to \lim(\G)$ is monadic.

\end{lemma}

\begin{proof}We use the equivalence
$\Fun(\mJ,\infty\Cat) \simeq \infty\Cat^\cart_{/ \mJ^\op}$
of $(\infty,2)$-categories provided by the Grothendieck construction.
Let $\mF \to \mJ^\op, \mG \to \mJ^\op$ be the cartesian fibrations classifying
$\F,\G,$ respectively.
Let $\rho: \mF \to \mG$ be the map of cartesian fibrations over $\mJ^\op$ classifying $\alpha: \F \to \G.$
The map $\rho: \mF \to \mG$ of cartesian fibrations over $\mJ^\op$ is fiberwise monadic and so monadic in $\infty\Cat^\cart_{/ \mJ^\op}$ by \cref{lemu} and \cref{tgvwxlkm}, \cref{rightadj2}, \cref{rembras}.
Hence $\alpha$ is monadic. Applying the right adjoint 2-functor
$\lim: \Fun(\mJ,\infty\Cat) \to \infty\Cat$
we find that $\lim(\alpha)$ is monadic (\cref{rembras}).

% Let $\phi$ be the left adjoint of $\alpha,$ which exists by (1),and $\T \simeq \alpha \circ \phi$ the monad associated to $\alpha$.Applying the 2-functor$\lim: (\infty\Cat^\circledast)^\mJ \to \infty\Cat^\circledast$we see that $\lim(\T) \simeq \lim(\alpha)\circ \lim(\phi)$ is the monad associated to $\lim(\alpha).$The induced functor $ \lim(\F) \to \LMod_{\lim(\T)}(\lim(\G))$ identies with the functor$ \lim(\F) \to \lim(\LMod_{\T}(\G))$. 

\end{proof}

\begin{lemma}\label{monlemmas}

Let $\rS$ be an $\infty$-category and $\G: \mD \to \mC$ a functor
over $\rS$ that admits a left adjoint relative to $\rS$.
The following conditions are equivalent:

\begin{enumerate}
\item The functor $\G: \mD \to \mC$ is monadic in $\infty\Cat_{/ \rS}$.	

\vspace{1mm}
\item For every functor $[1] \to \rS$ the induced functor $[1] \times_\rS \mD \to [1] \times_\rS \mC $ is monadic.

\vspace{1mm}

\item For every functor $[1] \to \rS$ the induced functor $\Fun_{\rS}([1], \mD) \to \Fun_{\rS}([1], \mC)$ is monadic.

\end{enumerate}

\end{lemma}

\begin{proof}

(1) implies (2) by \cref{monaT}. (2) implies (3) by definition of monadicity and
because the induced functor $\Fun_{\rS}([1], \mD) \to \Fun_{\rS}([1], \mC)$
identifies with the induced functor $$\Fun_{[1]}([1], [1] \times_\rS \mD) \to \Fun_{[1]}([1], [1] \times_\rS \mC).$$

We prove that (3) implies (1).
By \cref{monaT} the functor $\G: \mD \to \mC$ is monadic in $\infty\Cat_{/ \rS}$ if for every functor $\mB \to \rS$ the induced functor
$\rho: \Fun_{\rS}(\mB,\mD) \to \Fun_{\rS}(\mB,\mC)$ is monadic.
Let $\mW \subset \infty\Cat$ be the full subcategory of all $\mB$ such that
for every functor $\mB \to \rS$ the functor $\rho$ is an equivalence.
\cref{lemu2} guarantees that $\mW$ is closed under small colimits in $\infty\Cat.$ Since $\infty\Cat$ is genrerated under small colimits by $[0],[1]$, we can reduce to the case $\mB=[0],[1].$
The case of $\mB=[1]$ is (3). So it remains to see the case $\mB=[0],$
i.e. that for every $\s \in \rS$ the induced functor $\mD_\s \to \mC_\s$ is monadic. Let $[1] \to \rS$ be the constant functor with value $\s$.
The induced functor $\Fun_{\rS}([1], \mD) \to \Fun_{\rS}([1], \mC)$ 
is equivalent to the functor $\beta: \Fun([1], \mD_\s) \to \Fun([1], \mC_\s),$
which is monadic if (3) holds.
The functor $\mD_\s \to \mC_\s$ is a retract of the monadic functor $\beta$ and so monadic, too:
indeed, since $\beta$ is conservative, the functor $\mD_\s \to \mC_\s$ is conservative. This implies that the functor $\mD_\s \to \mC_\s$ is the pullback of $\beta$ along the embedding $\mC_\s \to \Fun([1], \mC_\s)$ and so is monadic.

\end{proof}

\begin{proposition}\label{relbar}

Let $\rS$ be an $\infty$-category and $\G: \mD \to \mC$ a functor
over $\rS$ that admits a left adjoint relative to $\rS$.
The following conditions are equivalent:

\begin{enumerate}
\item The functor $\G: \mD \to \mC$ is monadic in $\infty\Cat_{/ \rS}$.	

\vspace{1mm}
\item For every $\s \in \rS$ the induced functor $\G_\s: \mD_\s \to \mC_\s$ is conservative and every $\G_\s$-split simplicial object of $\mD_\s$ admits a colimit that is preserved by $\G_\s$ and the functor $\mD_\s \to \mD.$

% \vspace{1mm}
% \item For every morphism $\alpha: \s \to \rt$ in $\rS$ the induced functor $\G_\s: \mD_\s \to \mC_\s$ is conservative and every $\G_\s$-split simplicial object of $\mD_\s$ admits a colimit that is preserved by $\G_\s$ and the embedding $$\mD_\s= \{0\}\times_\rS \mD \to [1]\times_\rS \mD,$$ where the pullback is along the functor $[1]\to \rS$ taking $ \alpha.$

% \vspace{1mm}
% \item For every obejct $ \s $ of $\rS$ the induced functor $\G_\s: \mD_\s \to \mC_\s$ is conservative and every $\G_\s$-split simplicial object of $\mD_\s$ admits a colimit that is preserved by $\G_\s$.
% For every $\s \in \rS$, every
% $\G_\s$-split simplicial object $\alpha: (\Delta^\op)^\triangleright \to \mD_\s$ and every $\X \in \mD$ lying over some $\rt \in \rS$
% the functor
% $$ \Delta^\triangleleft \xrightarrow{\alpha^\op} \mD_\s^\op \to \mD^\op \xrightarrow{\Map_\mD(-,\Y)} \mS $$
% is a limit diagram.

% \vspace{1mm}

% \item For every functor $[1] \to \rS$ the induced functor $\Fun_{\rS}([1], \mD) \to \Fun_{\rS}([1], \mC)$ is monadic.

\end{enumerate}

\end{proposition}

\begin{proof}
The description of Eilenberg-Moore objects in $\infty\Cat_{/ \rS}$ of \cref{tgvwxlkm} combined with \cref{splitt} guarantees that (1) implies (2). 
%(2) trivially implies (3).
We prove that (2) implies (1). 

By \cref{monlemmas} we can reduce to the case that $\rS=[1].$
By (2) for every $\bi=0,1$ the induced functor $\mD_\bi \to \mC_\bi$ is monadic and thus also the functor $\Fun([1], \mD_\bi) \to \Fun([1], \mC_\bi)$,
which identifies with the functor
$\theta:\Fun_{[1]}(\rS, \mD) \to \Fun_{[1]}(\rS, \mC)$ for $\rS \to [1]$
the constant functor $[1] \to [1]$ with value $\bi.$
Therefore by \cref{monlemmas} it suffices to prove that the induced functor
$\theta:\Fun_{[1]}([1], \mD) \to \Fun_{[1]}([1], \mC)$ is monadic.

The functor $\theta$ is conservative because the induced functor $\mD_\bi \to \mC_\bi$ and the functor $\Fun_{[1]}([1], \mD) \to \mD_0 \times \mD_1$ are conservative.
So the functor $\theta$ is monadic if every $\theta$-split simplicial object of $\Fun_{[1]}([1], \mD)$ admits a colimit that is preserved by $\theta$. 
Let $\theta': \Fun([1],\mD) \to \Fun([1],\mC)$ be the induced functor.
The functor $\theta$ is the pullback of $\theta'$.
Moreover the canonical functor $\Fun_{[1]}([1], \mD) \to \Fun([1],\mD)$ is fully faithful since the functor $* \to \Fun([1],[1])$ taking the identity if fully faithful. 
Therefore it is enough to see that every $\theta'$-split simplicial object $\X$ of $\Fun([1], \mD)$ admits a colimit that is preserved by $\theta'$
if the image $\X^\bi$ under the functor $\Fun([1], \mD) \to \Fun(\{\bi\},\mD) $ lands in $\mD_\bi$ for $\bi=0,1.$
Since $\X^\bi$ is a split $\G_\bi$-simplicial object, by (2) the simplicial object $\X^\bi$ admits a colimit in $\mD_\bi$ that is preserved by $\G_\bi$
and the embedding $\mD_\bi \subset \mD.$
So $\X^\bi$ admits a colimit in $\mD$ that is preserved by $\G: \mD \to \mC$ for every $\bi=0,1.$
Consequently, $\X$ admits a colimit in $ \Fun([1],\mD)$ that is preserved by $\theta'$.

% It remains to see that (4) implies (1). By \cref{monaT} the functor $\G: \mD \to \mC$ is monadic in $\infty\Cat_{/ \rS}$ if for every functor $\mB \to \rS$ the induced functor
% $\rho: \Fun_{\rS}(\mB,\mD) \to \Fun_{\rS}(\mB,\mC)$ is monadic.
% Let $\mW \subset \infty\Cat$ be the full subcategory of all $\mB$ such that
% for every functor $\mB \to \rS$ the functor $\rho$ is an equivalence.
% \cref{lemu2} guarantees that $\mW$ is closed under small colimits in $\infty\Cat.$ Since $\infty\Cat$ is genrerated under small colimits by $[0],[1]$, we can reduce to the case $\mB=[0],[1].$
% The case of $\mB=[1]$ is (4). So it remains to see the case $\mB=[0],$
% i.e. that for any $\s \in \rS$ the induced functor $\mD_\s \to \mC_\s$ is monadic. Let $[1] \to \rS$ be the constant functor with value $\s$.
% The induced functor $\Fun_{\rS}([1], \mD) \to \Fun_{\rS}([1], \mC)$ 
% is equivalent to the functor $\beta: \Fun([1], \mD_\s) \to \Fun([1], \mC_\s).$
% If (4) holds, $\beta$ is monadic.
% The functor $\mD_\s \to \mC_\s$ is a retract of the monadic functor $\beta$ and so monadic, too:
% indeed, since $\beta$ is conservative, the functor $\mD_\s \to \mC_\s$ is conservative. This implies that the functor $\mD_\s \to \mC_\s$ is the pullback of
% $\beta$ along the embedding $\mC_\s \to \Fun([1], \mC_\s)$ and so is monadic.
% So (4) implies (1). 
 
\end{proof}

Next we define a notion of colimit in $\mathfrak{P}$-operads.

\begin{definition}
Let $\mathfrak{P}:= (\rS, \mE, \rS^\circ)$ be an algebraic pattern,
$\mC \to \rS$ a $\mathfrak{P}$-operad, $\K$ an $\infty$-category
and $\s \in \rS.$
A functor $\alpha: \K^\triangleright \to \mC_\s$ is a $\mathfrak{P}$-colimit
if for every $\Y \in \mC^\circ$ the following functor is a limit diagram: $$(\K^\op)^\triangleleft \simeq (\K^\triangleright)^\op \xrightarrow{\alpha^\op}\mC_\s^\op \to (\mC^\act)^\op \xrightarrow{\Map_{\mC^\act}(-,\Y)} \mS.$$

\end{definition}

\begin{example}Let $\mO^\ot \to \Delta^\op$ be a non-symmetric $\infty$-operad, $\mC^\ot \to \mO^\ot$ an $\mO$-operad, $\K$ an $\infty$-category and $\X \in \mO^\ot$ lying over $[\n].$
Let $\mathfrak{P}$ be the algebraic pattern for $\mO$-operads.

A functor $\alpha: \K^\triangleright \to \mC^\ot_\X$ is a $\mathfrak{P}$-colimit
if and only if for every $\Y \in \mC^\circ$ the following functor is a limit diagram: $$(\K^\op)^\triangleleft \simeq (\K^\triangleright)^\op \xrightarrow{\alpha^\op}(\mC^\ot_\X)^\op \simeq \mC^\op_{\X_1} \times ... \times \mC^\op_{\X_\n} \xrightarrow{\Mul_{\mC}(-,\Y)} \mS.$$

A $\mathfrak{P}$-colimit is an $\mO$-operadic colimit in the sense of \cite[Definition 3.1.1.2.]{lurie.higheralgebra}.

\end{example}

\begin{example}Let $\mV^\ot \to \Delta^\op$ be monoidal $\infty$-category,
$\mM^\circledast \to \mV^\ot$ a weakly $\mV$-enriched $\infty$-category, 
$\K$ an $\infty$-category
and $\V \in \mV^\circledast$ lying over $[\n].$
Let $\mathfrak{P}$ be the algebraic pattern for weak $\mV$-enrichment.

A functor $\alpha: \K^\triangleright \to \mM^\circledast_\X$ is a $\mathfrak{P}$-colimit
if and only if for every $\Y \in \mM$ the following functor is a limit diagram: $$(\K^\op)^\triangleleft \simeq (\K^\triangleright)^\op \xrightarrow{\alpha^\op}(\mM^\circledast_\V)^\op \simeq \mM^\op \xrightarrow{\Mul_{\mM}(\V_1,...,\V_\n,-,\Y)} \mS.$$

A $\mathfrak{P}$-colimit a $\mV$-enriched conical colimit in the sense of
\cite[Definition 3.25.]{heine2024higher}.

\end{example}

\begin{theorem}\label{strmon} Let $\mathfrak{P}:= (\rS, \mE, \rS^\circ)$ be an algebraic pattern.

A map $\G: \mD \to \mC$ in $\infty\Cat^{\mathfrak{P}}_{/\rS}$ is monadic if and only if the following hold:

\begin{enumerate}
\item The functor $\G: \mD \to \mC$ admits a left adjoint relative to $\rS.$

%\item The functor $\G: \mD \to \mC$ is conservative.

\item For every $\s \in \rS^\circ$ the induced functor $\G_\s: \mD_\s \to \mC_\s$ is conservative.
	
\item For every $\s \in \rS^\circ$ every $\G_\s$-split simplicial object of $\mD_\s$ admits a colimit that is preserved by $\G_\s$.

\item For every $\s \in \rS$ the colimit of every $\G_\s$-split simplicial object of $\mD_\s$ is a $\mathfrak{P}$-colimit.
	
\end{enumerate}

\end{theorem}

\begin{proof}

Let $\phi: \mD \to \rS, \psi: \mC \to \rS$ be the canonical functors.
By \cref{relbar} the functor $\G: \mD \to \mC$ is monadic if and only if
$\G: \mD \to \mC$ admits a left adjoint relative to $\rS$, 
for every object $\s \in \rS$ the induced functor $\G_\s: \mD_\s \to \mC_\s$ is conservative and every $\G_\s$-split simplicial object of $\mD_\s$ admits a colimit that is preserved by $\G_\s$ and the functor $\mD_\s \to \mD.$
By the first axiom of a $\mathfrak{P}$-operad
the induced functor $\G_\s: \mD_\s \to \mC_\s$ is conservative for every
$\s \in \rS$ if and only if the induced functor $\G_\s: \mD_\s \to \mC_\s$ is conservative for every $\s \in \rS^\circ$.
Again by the first axiom of a $\mathfrak{P}$-operad 
for every $\s \in \rS$ every $\G_\s$-split simplicial object of $\mD_\s$ admits a colimit that is preserved by $\G_\s$ if and only if this holds only for
every $\s \in \rS^\circ.$

Consequently, it remains to prove that for every $\s \in \rS$ the colimit of 
every $\G_\s$-split simplicial object of $\mD_\s$ is preserved by the functor $\mD_\s \to \mD$ if and only if for every $\s \in \rS$ the colimit of every $\G_\s$-split simplicial object of $\mD_\s$ is a $\mathfrak{P}$-colimit.

Let $\s \in \rS$ and $\beta : \Delta^\op \to \mD_\s$ a
$\G_\s$-split simplicial object.
The colimit $\colim(\beta)$ of $\beta$ is a preserved by the functor $\mD_\s \to \mD$ if and only if for every $\X \in \mD$ the induced map
$$ \rho: \Map_\mD(\colim(\beta),\X) \to \lim(\Map_\mD(-,\X) \circ \beta) $$ 
is an equivalence.
The map $\rho$ is a map over $ \Map_\rS(\s,\phi(\X)).$
So $\rho$ is an equivalence if and only if it induces an equivalence on the fiber over every object $\alpha$ of $ \Map_\rS(\s,\phi(\X)).$
Taking $\X \in \mD^\circ$ and $\alpha$ an active morphism,
we find the colimit of $\beta$ is a $\mathfrak{P}$-colimit if it is preserved by the functor $\mD_\s \to \mD$.
Hence every monadic morphism of $\infty\Cat^{\mathfrak{P}}_{/\rS}$ satisfies conditions (1)- (4).
We prove the converse. We assume that for every $\s \in \rS$ the colimit of every $\G_\s$-split simplicial object of $\mD_\s$ is a $\mathfrak{P}$-colimit. 
It remains to see that for every $\s \in \rS$ the colimit of every $\G_\s$-split simplicial object of $\mD_\s$ is preserved by the functor $\mD_\s \to \mD$.

Let $\s \in \rS$ and $\beta : \Delta^\op \to \mD_\s$ a
$\G_\s$-split simplicial object. We have to see that $\rho$ is an equivalence
for every $\X \in \mD.$
We prove first that $\rho$ is an equivalence for every $\X \in \mD$
if $\rho$ is an equivalence for every $\X \in \mD^\circ.$

By (3) of \cref{fibrous} the functor $$(\rS^\circ \times_{\rS}\mE_{\s /})^{\triangleleft} \to \mC, \ (\alpha: \s \to \rt) \mapsto \alpha_!(\X) $$ 
induces a pullback square
\begin{equation*} 
\begin{xy}
\xymatrix{
\Map_{\mD}(\colim(\beta),\X)
\ar[d] 
\ar[r]  & \lim_{\alpha \in \rS^\circ \times_{\rS}\mE_{\s /}}\Map_{\mD}(\colim(\beta),\alpha_!(\X))
\ar[d]
\\
\Map_{\rS}(\phi(\colim(\beta)),\phi(\X)) \ar[r] & \lim_{\alpha \in \rS^\circ \times_{\rS}\mE_{\s /}}\Map_{\rS}(\phi(\colim(\beta)),\phi(\alpha_!(\X))).
}
\end{xy}
\end{equation*}

We paste this square with the following commutative square, which is a pullback square if we assume the case of $\X \in \mD^\circ:$

\begin{equation*} 
\begin{xy}
\xymatrix{
\lim_{\alpha \in \rS^\circ \times_{\rS}\mE_{\s /}}\Map_{\mD}(\colim(\beta),\alpha_!(\X))
\ar[d] 
\ar[r]  & \lim_{\alpha \in \rS^\circ \times_{\rS}\mE_{\s /}}\lim(\Map_\mD(-,\alpha_!(\X)) \circ \beta)
\ar[d]
\\
\lim_{\alpha \in \rS^\circ \times_{\rS}\mE_{\s /}}\Map_{\rS}(\phi(\colim(\beta)),\phi(\alpha_!(\X))) \ar[r] & \lim_{\alpha \in \rS^\circ \times_{\rS}\mE_{\s /}}\lim(\Map_\rS(-,\phi(\alpha_!(\X))) \circ \phi \circ \beta).
}
\end{xy}
\end{equation*}

The pasting of both squares is therefore a pullback square which is the pasting of the following two squares:
\begin{equation}\label{ohol}
\begin{xy}
\xymatrix{
\Map_{\mD}(\colim(\beta),\X)
\ar[d] 
\ar[r]  & \lim(\Map_\mD(-,\X) \circ \beta)
\ar[d]
\\
\Map_{\rS}(\phi(\colim(\beta)),\phi(\X)) \ar[r]^\simeq & \lim(\Map_\rS(-,\phi(\X)) \circ\phi \circ \beta)
}
\end{xy}
\end{equation}
and
\begin{equation*} 
\begin{xy}
\xymatrix{
\lim(\Map_\mD(-,\X) \circ \beta)
\ar[d] 
\ar[r]  & \lim(\lim_{\alpha \in \rS^\circ \times_{\rS}\mE_{\s /}}\Map_\mD(-,\alpha_!(\X)) \circ \beta)
\ar[d]
\\
\lim(\Map_\rS(-,\phi(\X)) \circ\phi \circ \beta) \ar[r] & \lim(\lim_{\alpha \in \rS^\circ \times_{\rS}\mE_{\s /}}\Map_\rS(-,\phi(\alpha_!(\X))) \circ\phi \circ \beta)
}
\end{xy}
\end{equation*}

The last commutative square is a pullback square as a consequence of (3) of \cref{fibrous}.
Thus the pasting law implies that \ref{ohol} is a pullback square.

So $\rho$ is an equivalence for every $\X \in \mD$
if $\rho$ is an equivalence for every $\X \in \mD^\circ.$

We finish the proof by showing that $\rho$ is an equivalence for every
$\X \in \mD^\circ$. Let $\X \in \mD^\circ$.
Every morphism $\alpha: \s \to \phi(\X)$ of $\rS$ factors an inert morphism $\gamma: \s \to \rt$  followed by an active morphism $\alpha': \rt \to \phi(\X)$.

We consider the commutative diagram
\begin{equation*}
\begin{xy}
\xymatrix{
\Map_{\mD}(\gamma_!(\colim(\beta)),\X)
\ar[d] 
\ar[r]  &\Map_{\mD}(\colim(\beta),\X)
\ar[d] 
\ar[r]  & \lim(\Map_\mD(-,\X) \circ \beta)
\ar[d]
\\
\Map_{\rS}(\rt,\phi(\X)) \ar[r] & \Map_{\rS}(\phi(\colim(\beta)),\phi(\X)) \ar[r]^\simeq & \lim(\Map_\rS(-,\phi(\X)) \circ\phi \circ \beta)
}
\end{xy}
\end{equation*}

The left hand square is a pullback square by the definiton of cocartesian morphism. So by the pasting law it suffices to see that the outer square is a pullback square. The outer square identifies with the outer square in the following commutative diagram
\begin{equation*}
\begin{xy}
\xymatrix{
\Map_{\mD}(\gamma_!(\colim(\beta)),\X)
\ar[d] 
\ar[rr]  &&\lim(\Map_\mD(-,\X) \circ \gamma_! \circ \beta)
\ar[d] 
\ar[r]  & \lim(\Map_\mD(-,\X) \circ \beta)
\ar[d]
\\
\Map_{\rS}(\rt,\phi(\X)) \ar[rr]^\simeq && \lim(\Map_{\rS}(-,\phi(\X)) \circ \phi \circ \gamma_! \circ \beta) \ar[r] & \lim(\Map_\rS(-,\phi(\X)) \circ\phi \circ \beta)
}
\end{xy}
\end{equation*}

The right hand square is a pullback square by the definiton of cocartesian morphism. The left hand square is a pullback square if $\gamma_!: \mD_\s \to \mD_\rt$ sends the colimit of $\beta$ to a $\mathfrak{P}$-colimit.
We prove that this holds.
Since $\gamma$ is inert, the functor $\gamma_!: \mD_\s \to \mD_\rt$ covers the functor $\gamma_!: \mC_\s \to \mC_\rt$ and so sends the colimit diagram of the $\G_\s$-split simplicial object $\beta: \Delta^\op \to \mD_\s $ to a $\G_\rt$-split simplicial diagram $\gamma_! \circ \beta$. The latter is a $\mathfrak{P}$-colimit diagram since $\gamma_! \circ \beta: \Delta^\op \to \mD_\rt$ admits a $\mathfrak{P}$-colimit that is preserved by $\G_\rt$ and $\G_\rt$ is conservative.

\end{proof}

We obtain an operadic version of the monadicity theorem.

% \begin{definition}
% Let $\mC^\ot \to \Fin_*$ be a symmetric $\infty$-operad. A functor $\alpha: \K^\triangleright \to \mC$ is an operadic colimit diagram if for every $\Z_1,...,\Z_\n, \Z \in \mC$ for $\n \geq 0$ the following functor is a limit diagram: $$(\K^\op)^\triangleleft \simeq (\K^\triangleright)^\op \xrightarrow{\alpha^\op}\mC^\op \xrightarrow{\Mul_\mC(\Z_1,...,\Z_\n,-;\Z)} \mS.$$
% \end{definition}
% Note that an operadic colimit diagram is a colimit diagram in $\mC$ taking $\n=0.$

% \begin{lemma}\label{splittOp}
	
% Let $\mO^\ot \to \Fin_*$ be an $\infty$-operad, $\mC^\ot \to \mO^\ot$ be an $\mO$-operad and $\T$ a monad on $\mC^\ot \to \mO^\ot$ in $\infty\Op^\mO$. 
% Let $\nu: \LMod^\mO_\T(\mC)^\ot \to \mC^\ot$ be the forgetful functor.
% For every $\X \in \mO$ every $\nu_\X$-split simplicial object of $\LMod^\mO_\T(\mC)_\X$ admits an $\mO$-operadic colimit that is preserved by $\nu_\X$.

% \end{lemma}

\begin{corollary}\label{relbart}

Let $\mO$ be a non-symmetric $\infty$-operad and $\mC, \mD$ be $\mO$-operads.

A map $\G: \mD \to \mC$ of $\mO$-operads is monadic in $\infty\Op^\mO $ if and only if the following hold: 
%relative to $\mO$.The following conditions are equivalent:
\begin{enumerate}
\item The map $\G$ admits a left adjoint in $\infty\Op^\mO$.

\item For every color $\X$ of $\mO$ the induced functor $\G_{\X}: \mD_\X \to \mC_\X$ is conservative.

\item For every color $\X$ of $\mO$ every $\G_\X$-split simplicial object of $\mD_\X$ admits an $\mO$-operadic colimit preserved by $\G_\X$.
% and that is an $\mO$-operadic colimit.

%the functor $\mD_\X \to \mD.$
\end{enumerate}
\end{corollary}

% We obtain a monoidal version of the monadicity theorem.

% \begin{corollary}\label{relbart3}

% Let $\mO$ be a non-symmetric $\infty$-operad and $\mC, \mD$ be $\mO$-monoidal $\infty$-categories.

% A lax $\mO$-monoidal functor $\G: \mD \to \mC$ is monadic in $\infty\Op^\mO $ if and only if the following hold: 
% \begin{enumerate}
% \item The map $\G$ admits a left adjoint in $\infty\Op^\mO$.

% \item For every color $\X$ of $\mO$ the induced functor $\G_{\X}: \mD_\X \to \mC_\X$ is conservative.

% \item For every color $\X$ of $\mO$ every $\G_\X$-split simplicial object of $\mD_\X$ admits a colimit that is preserved by $\G_\X$ and the tensor product in each variable.

% \end{enumerate}
% \end{corollary}

We obtain an enriched version of the monadicity theorem.

% \begin{definition}
% Let $\mV$ be a monoidal $\infty$-category and $\mM$ a weakly $\mV$-enriched $\infty$-category. A functor $\alpha: \K^\triangleright \to \mM$ is a $\mV$-enriched colimit diagram if for every $\V_1,...,\V_\n \in \mV$ for $\n \geq 0$ and $\Y \in \mM$
% the following functor is a limit diagram: $$(\K^\op)^\triangleleft \simeq (\K^\triangleright)^\op \xrightarrow{\alpha^\op}\mM^\op \xrightarrow{\Mul\Mor_\mM(\V_1,...,\V_\n;-,\Y)} \mS.$$
% \end{definition}
% Note that a $\mV$-enriched colimit diagram is a colimit diagram in $\mM$ taking $\n=0.$

% \begin{remark}\label{remboq}

% A functor $\alpha: \K^\triangleright \to \mM$ is a $\mV$-enriched colimit diagram if and only if for every $\V_1,...,\V_\n \in \mV$ for $\n \geq 0$ and $\Y \in \mM^\circledast$ the functor $$(\K^\op)^\triangleleft \simeq (\K^\triangleright)^\op \xrightarrow{\alpha^\op}\mM^\op \xrightarrow{\mM^\circledast((\V_1,...,\V_\n,-),\Y)} \mS $$ is a limit diagram.
% This is equivalent to say that for every $\V \in \mV^\ot$
% the canonical functor
% $\mM \simeq \mM^\circledast_{\V} \to \mM^\circledast$ sends
% $\alpha: \K^\triangleright \to \mM$ to a colimit diagram.
% Note that for $\K=\Delta^\op$ a colimit $\alpha: \K^\triangleright \to \mM$ 
% that extends to a split simplicial object in $\mM$, is preserved by any functor
% and so is a $\mV$-enriched colimit diagram.

% \end{remark}

\begin{corollary}\label{EnrBaro}\label{EnrBar}
Let $\mV$ be a monoidal $\infty$-category and $\mM, \mN$ (weakly) $\mV$-enriched $\infty$-categories.
A $\mV$-enriched functor $\G: \mN \to \mM$ is monadic if and only if the following hold:

\begin{enumerate}
\item The $\mV$-enriched functor $\G: \mN \to \mM$ admits a $\mV$-enriched left adjoint.

\item The underlying functor of $\G$ is conservative.
	
\item Every $\G$-split simplicial object of $\mN$ admits a $\mV$-enriched colimit that is preserved by $\G$.
	
\end{enumerate}

\end{corollary}

We obtain a double $\infty$-categorical version of the monadicity theorem.

\begin{corollary}\label{relbart3}

Let $\mC, \mD$ be double $\infty$-categories.

A lax functor $\G: \mD \to \mC$ of double $\infty$-categories
is monadic in the $(\infty,2)$-category of double $\infty$-categories and lax functors if and only if the following hold: 
\begin{enumerate}
\item The map $\G$ admits a left adjoint.

\item The map $\G$ induces a monadic functor on $\infty$-categories of objects
and on $\infty$-categories of morphisms.

%\item For every objects $\A,\B$ of $\mD$ the induced functor $\G_{\A,\B}: \Mor_\mD(\A,\B) \to \Mor_\mC(\G(\A),\G(\B))$ is conservative.

% \item The composition functor $$\mD_{[1]} \times_{\mD_{[0]}} \mD_{[1]} \to \mD_{[1]}$$ preserves the colimit of every $\G_{[1]} \times_{\G_{[0]}} \G_{[1]}$-split simplicial object.

\item For every $\n \geq 0$ the composition functor $$\mD_{[1]} \times_{\mD_{[0]}} \mD_{[1]} \times_{\mD_{[0]}} ... \times_{\mD_{[0]}} \mD_{[1]} \to \mD_{[1]}$$ preserves the colimit of every $\G_{[1]} \times_{\G_{[0]}} \G_{[1]}\times_{\G_{[0]}} ... \times_{\G_{[0]}} \G_{[1]}$-split simplicial object.

% \item For every objects $\A,\B$ of $\mD$ every $\G_{\A,\B}$-split simplicial object of $\Mor_\mD(\A,\B)$ admits a colimit that is preserved by $\G_{\A,\B}$ and the functor $\Mor_\mD(\A,\B) \to \Mor_\mD(\A',\B') $ composing along any morphisms
% $\A' \to \A, \B'\to \B$ in $\mD.$

\end{enumerate}
\end{corollary}

\begin{proof}

In view of \cref{strmon} it suffices to prove that condition (3) is equivalent to the condition that for every $\n \geq 0$ and $\X \in \mD^\circ$ the colimit of every $\G_{[\n]}$-split simplicial object is preserved by the functor
\begin{equation}\label{mapiz}
\mD_{[\n]} \to \mD^\act \xrightarrow{\Map_{\mD^\act_{[\n]}}(-,\X)} \mS^\op.\end{equation} 

If $\X \in \mD^\circ$, then $\X$ lies over $[0]$ or $[1].$
If $\X$ lies over $[0]$, there is nothing to show since there is no active morphism $[0] \to [\n] $ in $\Delta$ unless $\n=0$. But if $\n=0$, there is nothing to show since the functor $\mD_{[0]} \subset \mD$ is fully faithful.
So we consider the case that $\X$ lies over $[1]$.
There is a unique active morphism $\alpha^\n: [1] \to [\n] $ in $\Delta$.
In this case the functor (\ref{mapiz}) identifies with the functor
$$ \mD_{[\n]} \xrightarrow{\alpha^\n_!} \mD_{[1]} \xrightarrow{\Map_{\mD_{[1]}}(-,\X)} \mS^\op.$$
The latter functor preserves the colimit of every $\G_{[\n]}$-split simplicial object for every $\X \in \mD_{[1]}$ if and only if condition (3) holds.

% We prove by induction on $\n \geq 0$ that the latter condition is equivalent to
% condition (3), the case of $\n=2$.
% If $\n=1$, there is nothing to show.

% The unique active morphism $\alpha^\n: [1] \to [\n] $ in $\Delta$
% factors as the active morphism $\alpha^{\n-1}: [1] \to [\n-1] $ followed by the active injection $\beta: [\n-1] \to [\n]$ omitting $\n-1.$
% Hence the functor $\alpha^\n_!: \mD_{[\n]} \to \mD_{[1]}$ factors as
% $\beta_!: \mD_{[\n]} \to \mD_{[\n-1]}$
% followed by 
% $\alpha^{\n-1}_!: \mD_{[\n-1]} \to \mD_{[1]}$.
% By induction hypothesis the functor $\alpha^{\n-1}_!$ preserves the 
% colimit of every $\G_{[\n-1]}$-split simplicial object.
% The functor $\beta_!: \mD_{[\n]} \simeq \mD_{[\n-2]} \times_{\mD_{[0]}} \mD_{[2]}\to \mD_{[\n-1]} \simeq \mD_{[\n-2]} \times_{\mD_{[0]}} \mD_{[1]}$
% identifies with the functor
% $\mD_{[\n-2]} \times_{\mD_{[0]}}\alpha^2_!: \mD_{[\n-2]} \times_{\mD_{[0]}} \mD_{[2]}\to \mD_{[\n-2]} \times_{\mD_{[0]}} \mD_{[1]}$
% and so preserves the 
% colimit of every $\G_{[\n]}$-split simplicial object.

\end{proof}

% \begin{corollary}\label{EnrBar}
% Let $\mV$ be a monoidal $\infty$-category and $\mM$ a $\mV$-enriched $\infty$-category.
% A $\mV$-enriched functor $\G: \mN \to \mM$ is monadic 
% in $\mV \mathrm{-}\Cat$ if and only if the following hold:

% \begin{enumerate}
% \item The $\mV$-enriched functor $\G: \mN \to \mM$ admits a $\mV$-enriched left adjoint.

% \item The underlying functor of $\G: \mN \to \mM$ is conservative.

% \item Every $\G$-split simplicial object of $\mN$ admits a $\mV$-enriched colimit that is preserved by $\G$.

% \end{enumerate}

% \end{corollary}	

\section{Monadic higher algebra}

\subsection{The relative tensor product of monadic algebras}\label{rEl}

In the following we construct a relative tensor product for algebras of a structured monad that extends the given tensor product along the free functor. 
To state our theorem we use the following terminology:

\begin{definition}Let $\K$ be an $\infty$-category.
A functor $\mC \to \rS$ is compatible with $\K$-indexed colimits if
the diagonal functor $\mC \to \mC^\K$ over $\rS$ admits a left adjoint relative to $\rS$.	

\end{definition}

\begin{example}
Let $\K$ be an $\infty$-category.
A locally cocartesian fibration $\mC \to \rS$ is compatible with $\K$-indexed colimits if and only if every fiber admits $\K$-indexed colimits, which are preserved by the fiber transports. This follows immediately from \cite[Proposition 7.3.2.11.]{lurie.higheralgebra}.	

\end{example}

\begin{theorem}\label{tgfdcccccc} Let $\rS$ be an $\infty$-category, $\mathfrak{P}$ an algebraic pattern on $\rS$ and $\mC \to \rS$ a $\mathfrak{P}$-monoidal $\infty$-category compatible with geometric realizations
(small colimits) and $\T$ a $\mathfrak{P}$-operadic monad on $\mC \to \rS$ that preserves fiberwise geometric realizations.
Then $\LMod^\rS_\T(\mC) \to \rS $ is a $\mathfrak{P}$-monoidal $\infty$-category compatible with geometric realizations (small colimits).

\end{theorem}

\cref{tgfdcccccc} immediately follows from \cref{Rell} and \cref{tgvwxlkm} (1).

%\begin{theorem}\label{tgfdcccccc}Let $\mE \subset \rS$ be a subcategoryand $\mC \to \rS$ a cocartesian fibration relative to $\mE$ and locally cocartesian fibration whose fibers admit geometric realizations and whose fiber transports preserve geometric realizations.Let $\T$ be a monad on $\mC \to \rS$ in $\infty\Cat_{/ \rS}$ such that for every $\s\in \rS$ the functor $\T_\s: \mC_\s \to \mC_\s$ preserves geometric realizations.	The functor $\LMod^{\rS}_\T(\mC) \to \rS$ is a cocartesian fibration relative to $\mE$ and locally cocartesian fibration whose fiber transports preserve geometric realizations. If the fiber transports of $\mC \to \rS$ preserve small colimits, thefiber transports of $\LMod^{\rS}_\T(\mC) \to \rS$ preserve small colimits.\end{theorem}

Before proving Theorem \ref{Rell} we consider important corollaries. % of Theorem \ref{tgfdcccccc}.
%\begin{corollary}Let $\mO^\ot \to \Fin_*$ be a symmetric $\infty$-operad and $\mC^\ot \to \mO^\ot$ an $\mO$-monoidal $\infty$-categorycompatible with geometric realizations.Let $\T$ be a lax $\mO$-monoidal monad on $\mC$ %in $\infty\Op^\mO$ that preserves geometric realizations.The $\mO$-operad $\LMod_\T(\mC)^\ot \to \mO^\ot $ is an$\mO$-monoidal $\infty$-category compatible with geometric realizations.\end{corollary}and Theorem \ref{tgfdcccccc} give the following theorem:\begin{proof}
%By Proposition \ref{corso} the functor $\LMod^\rS_\T(\mC) \to \rS $ is a $\mathfrak{P}$-operad. By Theorem \ref{tgfdcccccc} for $\mE=\rS$ the functor $\LMod^\rS_\T(\mC) \to \rS $ is a cocartesian fibration.	\end{proof}

\begin{corollary}\label{atono}
Let $\mO$ be an $\bE_{\bk}$-operad for $1 \leq \bk \leq \infty$ and $\mC$ an $\mO$-monoidal $\infty$-category compatible with geometric realizations. Let $\T$ be a monad on $\mC$ in $\infty\Op^{\mO}$ that preserves fiberwise geometric realizations.
Then $\LMod^\mO_\T(\mC) $ is an $\mO$-monoidal $\infty$-category
and the free functor $\mC \to \LMod^\mO_\T(\mC)$ is an $\mO$-monoidal functor.
	
\end{corollary}

\begin{corollary}Let $1 \leq \bk \leq \infty$ and $\mV$ % \simeq ({\Delta^\op})^{\times \bk} \times {\Delta^\op}$ 
an $\bE_{\bk+1}$-monoidal $\infty$-category compatible with geometric realizations
%(seen as $({\Delta^\op})^{\times \bk}$-family of monoidal $\infty$-categories)
%corresponding to a cocartesian $\bE_\bk$-family of monoidal $\infty$-categories
and $\A$ % \in \Fun_{({\Delta^\op})^{\times \bk} }(({\Delta^\op})^{\times \bk},\Alg^{({\Delta^\op})^{\times \bk}}(\mV^\ot)) \hookrightarrow  \Fun_{({\Delta^\op})^{\times \bk+1} }(({\Delta^\op})^{\times \bk+1},\mV^\ot)  $ 
an $\bE_{\bk+1}$-algebra in $\mV$.
Then $\LMod_\A(\mV)$ carries a canonical $\bE_{\bk}$-monoidal structure compatible with geometric realizations and the free functor
$\mV \to \LMod_\A(\mV)$ is $\bE_{\bk}$-monoidal.
	
\end{corollary}

\begin{proof}
An $\bE_{\bk+1}$-monoidal $\infty$-category $\mV$
corresponds to an associative algebra structure for the cartesian structure on
$\infty\Op^{\bE_{\bk}, \mon}$ on the underlying $\bE_{\bk}$-monoidal $\infty$-category $\mV'$ of $\mV$ 
and there is a canonical equivalence $\Alg_{\bE_{\bk+1}}(\mV) \simeq \Alg(\Alg_{\bE_{\bk}}(\mV'))$. By Corollary \ref{polk} the monad $\A \ot (-):\mV \to \mV$ is an $\bE_{\bk}$-monoidal monad on $\mV$, i.e. a monad on $\mV'$ in $\infty\Op^{\bE_{\bk}, \mon}$.
Thus by Theorem \ref{Rell} the functor $ \LMod^{\bE_\bk}_{\A \ot(-)}(\mV'^\ot) \to \bE_\bk $ is an $\bE_{\bk}$-monoidal $\infty$-category compatible with geometric realizations.
	
\end{proof}

%\begin{corollary}Let $\mV^\ot \to \Delta^\op, \mW^\ot \to \Delta^\op$ be monoidal $\infty$-categories compatible with geometric realizations, $\mM^\circledast \to \mV^\ot \times \mW^\ot$ an $\infty$-category bitensored over $\mV, \mW$ compatible with geometric realizations and $\A$ an associative algebra in $\mW$ 
% in $\mV \mathrm{-}\omega\Enr.$such that the functor $(-)\ot \A: \mM \to \mM$ preserves geometric realizations. The weakly left enriched $\infty$-category $\RMod_\A(\mM)^\circledast \to \mV^\ot $ is left tensored over $\mV$ compatible with geometric realizations.end{corollary}\begin{proof}The bitensored $\infty$-category $\mM^\circledast \to \mV^\ot \times \mW^\ot$corresponds to a left $\tau(\mW)$-action on $\mM_{[0]}^\circledast \to \mV^\ot$ with respect to the left action of $\infty\Cat$ on $\mV\mathrm{-}\LMod$.The latter gives a monoidal functor $\tau(\mW) \to \LinFun_\mV(\mM,\mM)$that sends an associative algebra $\A$ in $\mW$ to a $\mV$-linear monad $(-)\ot \A$ on $\mM$ such that $\RMod_\A(\mM)^\circledast \simeq \LMod_{\A \ot (-)}(\mM)^\circledast.$\end{proof}

%\subsection{The construction of the relative tensor product}In the following we use the theory of the latter subsection to construct the relative tensor product (Theorem \ref{tgfdccc}) and prove that it is closed (Proposition \ref{closed}). 

We construct the relative tensor product by extending it from free algebras. That this is possible is the content of Proposition \ref{ztgvfdwq}.

\begin{proposition}\label{ztgvfdwq}
	
Let 
\begin{equation*}
\begin{xy}
\xymatrix{
\Y 
\ar[rd]_\varphi
\ar[rr]^\xi && \X
\ar[ld]^\phi
\\
& \rS
}
\end{xy}
\end{equation*}
be a commutative triangle, where $\varphi: \Y \to \rS$ is a cocartesian fibration relative to a subcategory $\mE$ of $\rS$ and $\xi$ sends $\varphi$-cocartesian lifts of morphisms of $\mE$ to $\phi$-cocartesian morphisms. Let $\K$ be an $\infty$-category. The functor $\phi: \X \to \rS$ is a cocartesian fibration relative to $\mE$ if for every object $\s$ of $\rS$ the fiber $\X_\s$ is generated 
under $\K$-indexed colimits by the essential image of $\xi_\s :\Y_\s \to \X_\s$
and the diagonal functor 
$ \X \to \X^{  \K }  $ over $\rS$ admits a left adjoint relative to $\rS.$
 	
\end{proposition}

Proposition \ref{ztgvfdwq} follows from the following three lemmas.

\begin{notation}
	
Let $ \phi: \X \to \rS$ be a functor and $\mE \subset \rS$ a subcategory.
Let $\tilde{\X} \subset \X$ be the full subcategory spanned by the objects $\A$ lying over some object $\s $ of $\rS$ such that for every morphism 
$\F: \s \to \rt $ of $\mE $ there exists a $\phi$-cocartesian lift $\A \to \B$ of $\F.$
Let $\phi': \tilde{\X} \subset \X \xrightarrow{\phi} \rS$ be the restriction.
\end{notation}

\begin{lemma}\label{ghbnnhj}Let $ \phi: \X \to \rS$ be a functor and $\mE \subset \rS$ a subcategory.
The functor $\phi': \tilde{\X} \to \rS$ is a cocartesian fibration relative to $\mE$. The embedding $\tilde{\X} \subset \X$
sends $\phi'$-cocartesian morphisms to $\phi$-cocartesian morphisms.

\end{lemma}

\begin{proof}

Let $\F: \s \to \rt$ be a morphism of $\mE$ and $\A \in \tilde{\X}_\s$. 
Then there is a $\phi$-cocartesian lift $\A \to \F_\ast(\A)$ of $\F. $ We show that $\F_\ast(\A)$ belongs to $\tilde{\X}.$ Let $\G: \rt \to \br$ be a morphism of $\mE.$ 
Since $\A \in \tilde{\X}_\s,$ there is a $\phi$-cocartesian lift $\A \to (\G \circ \F)_\ast(\A)$ of $\G \circ \F: \s \to \rt \to \br. $ 
As $\A \to \F_\ast(\A)$ is $\phi$-cocartesian, the morphism $\A \to (\G \circ \F)_\ast(\A)$ factors as the morphism 
$\A \to \F_\ast(\A)$ followed by a lift $\F_\ast(\A) \to  (\G \circ \F)_\ast(\A) $ of $\G:\rt \to \br.$ 
Because $\A \to \F_\ast(\A)$ and  $\A \to (\G \circ \F)_\ast(\A)$ are $\phi$-cocartesian, the morphism $\F_\ast(\A) \to  (\G \circ \F)_\ast(\A) $ is $\phi$-cocartesian, too. 
Thus $\F_\ast(\A)$ belongs to $\tilde{\X}.$ 

\end{proof}

\begin{lemma}\label{levcoc}

Let $ \phi: \X \to \rS$ be a functor and $\K$ an $\infty$-category.
%$\phi$ induces a functor $\Fun(\K, \phi): \Fun(\K, \X) \to \Fun(\K,\rS).$
A morphism of $\Fun(\K, \X) $ is $\Fun(\K, \phi)$-cocartesian if it is objectwise $\phi$-cocartesian. %, i.e. for every $\bk \in \K$ the component $\tau(\bk)$ is $\phi$-cocartesian.

%Especially we have the following:Let $\psi: \X^\K \simeq \rS \times_{\Fun(\K, \rS) } \Fun(\K,\X) $ the cotensor of the $\infty$-category $\K$ with the $\infty$-category $\X $ over $\rS.$ Every levelwise $\phi$-cocartesian morphism of $\X^\K $ is $\Fun(\K, \phi)$-cocartesian and thus especially $\psi$-cocartesian.

\end{lemma}

\begin{proof}

Let $\mW \subset {\infty\Cat}$ be the full subcategory of $\infty$-categories $\K$ for which the statement holds.
We show that $\mW= {\infty\Cat}. $
Since $[0]\in\mW$ and ${\infty\Cat} $ is generated under small colimits by the $\infty$-categories $[0],[1] $, it remains to see that $\mW$ contains $[1] $ and is closed in ${\infty\Cat} $ under small colimits.

The functor $\Fun(-,\X): \infty\Cat^\op \to  {\infty\Cat} $ sends small colimits to limits. So closedness under small colimits follows from the fact that for every functor $\F: \mJ \to \infty\Cat$
a morphism in $\lim \Fun(\F(-),\X) $ is $ \lim \Fun(\F(-),\phi)$-cocartesian if for every $\bj \in \mJ$ its image in $\Fun(\F(\bj),\X)$ is $\Fun(\F(\bj),\phi)$-cocartesian and for every morphism $\bi \to \bj$ in $\mJ$
the induced functor $\Fun(\F(\bj),\X) \to \Fun(\F(\bi),\X)$ preserves cocartesian morphisms, using that mapping spaces in limits are limits of mappings spaces.

It remains to show that $[1]$ belongs to $\mW.$ 
For that we like to see that every objectwise $\phi$-cocartesian morphism of $ \Fun([1],\X)  $ corresponding to a commutative square
\begin{equation*}\label{trwds} 
\begin{xy}
\xymatrix{
\A 
\ar[d]^\g 
\ar[r]  & \B
\ar[d]^\h
\\
\C \ar[r] & \D
}
\end{xy}
\end{equation*}
in $\X,$ whose horizontal morphisms are  $\phi$-cocartesian, is $\Fun([1],\phi)$-cocartesian. 

For every morphism $\bk: \E \to \F $ of $\X$ the commutative square
\begin{equation*}\label{trwds} 
\begin{xy}
\xymatrix{
\Fun([1],\X)( \h, \bk) 
\ar[d] 
\ar[r]  & \Fun([1],\rS)(\phi(\h), \phi(\bk)) 
\ar[d]
\\
\Fun([1],\X)(\g, \bk)  \ar[r] & 	\Fun([1],\rS)(\phi(\g), \phi(\bk)) 
}
\end{xy}
\end{equation*}
is equivalent to the commutative square
\begin{equation*}\label{trwds} 
\begin{xy}
\xymatrix{
\X(\D, \F) \times_{ \X(\B, \F)} \X(\B, \E)
\ar[d] 
\ar[r]  & \rS(\phi(\D), \phi(\F)) \times_{ 	\rS(\phi(\B), \phi(\F))} 	\rS(\phi(\B), \phi(\E)) 
\ar[d]
\\
\X(\C, \F) \times_{ \X(\A, \F)} \X(\A, \E) \ar[r] &  	\rS(\phi(\C), \phi(\F)) \times_{ 	\rS(\phi(\A), \phi(\F))} 	\rS(\phi(\A), \phi(\E)) 
}
\end{xy}
\end{equation*}
This square is a pullback square since the morphisms $\A \to \B$ and $\C \to \D$ of $\X$ are $\phi$-cocartesian.
% and pullbacks commute with pullbacks.

\end{proof}	

\begin{lemma}\label{hrwdc}

Let $ \phi: \X \to \rS$ be a functor, $\K$ an $\infty$-category and $\mE \subset \rS$ a subcategory. If the diagonal functor 
$ \X \to \X^{  \K }  $ over $\rS$ admits a left adjoint relative to $\rS,$
for every $\s \in \rS$ the fiber $\tilde{\X}_\s$ is closed in $\X_\s$ under $\K$-indexed colimits.
\end{lemma}

\begin{proof}

Let $\K^\vartriangleright \to \X_\s $ be a colimit diagram whose restriction $\rH: \K \subset \K^\vartriangleright \xrightarrow{} \X_\s $ factors through $\tilde{\X}_\s.$
We want to see that $ \colim(\rH)$ belongs to $ \tilde{\X}_\s. $ 

Let $\F: \s \to \rt $ be a morphism of $\mE.$ We have to find a $\phi$-cocartesian lift $ \colim(\rH) \to \Z $ of $\F.$ 
Let $\phi': \tilde{\X} \subset \X \xrightarrow{\phi} \rS$ be the restriction and $\psi: \tilde{\X}^\K \simeq \rS \times_{\Fun(\K, \rS) } \Fun(\K,\tilde{\X}) \to \rS$ the cotensor.
The functor $\phi': \tilde{\X} \to \rS$ is a cocartesian fibration relative to $\mE$.
Thus the functor $\tilde{\X}^\K \to \rS$ is a cocartesian fibration relative to $\mE$, whose cocartesian morphisms are objectwise $\phi'$-cocartesian.
So we obtain a $\psi$-cocartesian morphism $\alpha: \rH \to \F_\ast(\rH)$ lying over $\F. $ By assumption the diagonal functor 
$  \X \to  \X^{  \K }  $ over $ \rS $ admits a left adjoint $\chi:  \X^{  \K }  \to \X$ relative to $ \rS. $ 
The functor $\chi$ sends $\alpha$ to a morphism
$\beta: \colim(\rH) \to \colim( \F_\ast(\rH)) $ of $ \X$ lying over $\F.$ 
The morphism $\beta$ is $\phi$-cocartesian since the composition $  \tilde{\X}^\K \subset  \X^\K \xrightarrow{\chi}   \X$ sends $\psi$-cocartesian to $\phi$-cocartesian morphisms:
as a relative left adjoint $\chi: \X^\K \to \X$ sends morphisms that are cocartesian with respect to the functor $ \X^\K \to \rS $ to $\phi$-cocartesian morphisms.
The embedding $\tilde{\X} \subset \X$ sends $\phi'$-cocartesian morphisms to $\phi$-cocartesian morphisms. So the embedding $\tilde{\X}^\K \subset \X^\K$ sends $\psi$-cocartesian morphisms to objectwise  $\phi$-cocartesian morphisms, which are cocartesian with respect to the functor $ \X^\K \to \rS $ by Lemma \ref{levcoc}.

\end{proof}

\begin{proof}[Proof of Proposition \ref{ztgvfdwq}]	
Since $ \varphi: \Y \to \rS$ is a cocartesian fibration relative to $\mE$ and $\xi$ sends $\varphi$-cocartesian lifts of morphisms of $\mE$ to $\phi$-cocartesian morphisms, for every object $\s $ of $\rS$ the fiber $\tilde{\X}_\s$ contains the essential image of $\xi_\s :\Y_\s \to \X_\s.$
If the diagonal functor 
$ \X \to \X^{  \K }  $ over $\rS$ has a left adjoint relative to $\rS,$
by Lemma \ref{hrwdc} the full subcategory $\tilde{\X}_\s$ is closed in $\X_\s$ under $\K$-indexed colimits.
So by assumption we have $\tilde{\X}_\s = \X_\s $ and so $\tilde{\X} = \X.$
Thus the claim follows from Lemma \ref{ghbnnhj}.

\end{proof}

%Next we apply Proposition \ref{ztgvfdwq} to prove Theorem \ref{Rell}.

%\begin{theorem}\label{tgfdccc}Let $\mC \to \rS$ be a functor compatible with geometric realizationsand $\T$ a monad on $\mC \to \rS$ that preserves fiberwise geometric realizations.Let $\mE \subset \rS$ be a subcategory and $\mT \subset \Fun([2],\rS)$ a full subcategory such that the functor $\mC \to \rS$ is a cocartesian fibration relative to $\mE, \mT$.Then $\LMod^{\rS}_\T(\mC) \to \rS$ is a cocartesian fibration relative to $\mE, \mT$ compatible with geometric realizations. 
%$\mC \to \rS$ is a cocartesian fibration compatible with small colimits, then $\LMod^{\rS}_\T(\mC) \to \rS$ is a cocartesian fibration compatible with small colimits.\end{theorem}

\begin{theorem}\label{Rell} Let $\rS$ be an $\infty$-category and $\mE \subset \Fun([1],\rS)$ a full subcategory.
Let $\mC \to \rS$ be a cocartesian fibration relative to $\mE$
compatible with geometric realizations and $\T$ a monad on $\mC \to \rS$ in
$\infty\Cat_{/\rS}$ that preserves fiberwise geometric realizations.
The functor $\LMod^\rS_\T(\mC) \to \rS $ is a cocartesian fibration relative to $\mE$ compatible with geometric realizations.

Let $\alpha: \s \to \rt$ be a morphism of $\mE$ such that
$\alpha_!: \mC_\s \to \mC_\rt$ preserves small colimits. Then $$\alpha_!: \LMod^{\rS}_\T(\mC)_\s \to \LMod^{\rS}_\T(\mC)_\rt$$ preserves small colimits.

\end{theorem}

\begin{proof}
	
%In view of Theorem \ref{corso} we can assume that $\mK$ is empty.	
% By definition a functor $\mD \to \rS$ is a cocartesian fibration relative to $\mE, \mT$ if and only if for every functor $\sigma: [2]\to \rS$ that belongs to $\mT$ the pullback $[2]\times_\rS \mD \to [2]$ along $\sigma$ is a cocartesian fibration relative to $\sigma^{-1}(\mE)$. So we can reduce to the case that $\mT=\Fun([2],\rS).$	

% The functor $\mC \to \rS$ is a locally cocartesian fibration relative to $\mE$ if and only if for every functor $\sigma: [1]\to \rS$ that belongs to $\mE$ the pullback $\mC':= [1]\times_\rS \mC \to [1]$ along $\sigma$ is a cocartesian fibration. Moreover the pullback of the functor
% $\LMod^{\rS}_\T(\mC) \to \mC$ along $\sigma$ is $\LMod^{[1]}_{\T'}(\mC') \to [1],$ where $\T':= [1]\times_\rS \T.$
% Therefore the statement about locally cocartesian fibrations relative to $\mE$ follows from the statement about cocartesian fibrations over $[1]$
% and so cocartesian fibrations relative to $\mE.$ So it suffices to prove the statement about cocartesian fibrations relative to $\mE.$ 

We apply Proposition \ref{ztgvfdwq} to the free functor $\mC \to \LMod^{\rS}_\T(\mC)$ over $\rS, $ which preserves cocartesian morphisms as a left adjoint relative to $\rS.$
By assumption the diagonal functor $ \mC \to \mC^{\Delta^\op} $ over $\rS$ admits a left adjoint relative to $\rS.$
Thus for any functor $\mD \to \rS$ the diagonal functor
$\Fun_\rS(\mD, \mC) \to \Fun_\rS(\mD, \mC^{\Delta^\op}) \simeq \Fun_\rS(\mD, \mC)^{\Delta^\op}$ admits a left adjoint.
Let $\Fun'_\rS(\mC,\mC) \subset \Fun_\rS(\mC,\mC)$ be the full subcategory of functors preserving fiberwise geometric realizations. 
The diagonal functor $\Fun_\rS(\mD, \mC) \to \Fun_\rS(\mD, \mC)^{\Delta^\op}$ is $\Fun_\rS(\mC,\mC)$-linear and so $\Fun'_\rS(\mC,\mC)$-linear
and the left adjoint is $\Fun'_\rS(\mC,\mC)$-linear.
Hence we obtain an adjunction $$\LMod_\T(\Fun_\rS(\mD, \mC)^{\Delta^\op}) \simeq \LMod_\T(\Fun_\rS(\mD, \mC))^{\Delta^\op} \rightleftarrows \LMod_\T(\Fun_\rS(\mD, \mC)),$$
which identifies with an adjunction
$$\Fun_\rS(\mD,  \LMod^{\rS}_\T(\mC))^{\Delta^\op} \simeq \Fun_\rS(\mD,  \LMod^{\rS}_\T(\mC)^{\Delta^\op}) \rightleftarrows \Fun_\rS(\mD, \LMod^{\rS}_\T(\mC)).$$
This implies by Proposition \ref{prosta} that the functor $ \LMod^{\rS}_\T(\mC)^{\Delta^\op} \to \LMod^{\rS}_\T(\mC) $ over $\rS$ admits a left adjoint relative to $\rS.$
By \cite[Example 4.7.2.7.]{lurie.higheralgebra} for every $\s\in\rS$ the fiber $ \LMod^{\rS}_\T(\mC)_\s \simeq \LMod_{\T(\s)}(\mC_\s)$ is generated under geometric realizations by the free $\T(\s)$-algebras.
Thus we can apply Proposition \ref{ztgvfdwq} to deduce that the functor $\LMod^{\rS}_\T(\mC) \to \rS$ is a cocartesian fibration relative to $\mE.$

%Let $\mC \to \rS$ be a cocartesian fibration compatible with small colimits. Then $\LMod^{\rS}_\T(\mC) \to \rS$ is a cocartesian fibration compatible with geometric realizations. % whose fibers admit small colimits.
To see that $\alpha_!$ preserves small colimits, it is enough to see that
$\alpha_! $ preserves arbitrary coproducts.
For every morphism $\alpha: \s \to \rt$ the $\infty$-category $\LMod^{\rS}_\T(\mC)_\s \simeq \LMod_{\T(\s)}(\mC_\s)$ admits arbitrary coproducts:
%and the induced functor $\alpha_! : \LMod_{\T(\s)}(\mC_\s) \to \LMod_{\T(\rt)}(\mC_\rt)$ preserves arbitrary coproducts.
every $\A \in \LMod_{\T(\s)}(\mC_\s)$ is a geometric realization of a simplicial object $\A_\bullet$ of free modules and the coproduct of a family
$(\A^\bj)_{\bj \in \mJ}$ in $\LMod_{\T(\s)}(\mC_\s) $ is $\mid \coprod_{\bj \in \mJ} \A^\bj_\bullet \mid, $
where we use that the family $ (\A^\bj_\bullet)_{\bj \in \mJ} $ admits a coproduct in $\Fun(\Delta^\op, \LMod_{\T(\s)}(\mC_\s))$ (computed objectwise) since families of free modules admit a coproduct.
The induced functor $\alpha_! $ preserves arbitrary coproducts since it preserves geometric realizations. 
\end{proof}

Corollary \ref{exk} and \cref{Rell} give the following corollary:

\begin{corollary}\label{enrtensor}

Let $\mV$ be a monoidal $\infty$-category compatible with geometric realizations, $\mM$ an $\infty$-category left tensored over $\mV$ compatible with geometric realizations and $\T$ a $\mV$-enriched monad on $\mM$
that preserves geometric realizations. The weakly $\mV$-enriched $\infty$-category $\LMod_\T(\mM) $ is left tensored over $\mV$
and for every $\V \in \mV$ the functor
\begin{equation}\label{eqos}
\V \ot (-):\LMod_\T(\mM) \to \LMod_\T(\mM)
\end{equation} preserves geometric realizations.

If moreover for every $\V \in \mV$ the functor
$ \V \ot (-):\mM \to \mM$ preserves small colimits, for every $\V \in \mV$ 
the functor (\ref{eqos}) preserves small colimits.

\end{corollary}

\begin{proof}
	
We apply \cref{Rell} for $\mE \subset \mV^\ot$ the subcategory of cocartesian lifts of morphisms of $\Delta^\op$, where we use that the functor 
$\mM^\circledast \to \mV^\ot$ is a locally cocartesian fibration and cocartesian fibration relative to $\mE$ whose fibers are all $\mM$ and whose fiber transports preserve geometric realizations.
\end{proof}

\subsection{A tensor product for algebras over Hopf $\infty$-operads}\label{Hhopf}

%\subsection{Hopf monads}\label{Hhopf}define enriched $\infty$-operads following \cite[\S 2]{heine2024restricted} and 

In the following we prove that the underlying monad of every Hopf $\infty$-operad is a Hopf monad (Theorem \ref{dfghjbbvl}, Corollary \ref{hopfo}).
Via Corollary \ref{exkq} we conclude that the $\infty$-category of algebras over every Hopf $\infty$-operad inherits a symmetric monoidal structure (\cref{Hoopf}).

% We start with defining the composition product on symmetric sequences in $\mC.$
\begin{notation}
Let $\Sigma \simeq  \coprod_{\n \geq 0} B(\Sigma_\n)$ be the groupoid of finite sets and bijections, which is the free symmetric monoidal $\infty$-category generated by the point. 

\end{notation}

For the next notation we use \cref{operat}:

\begin{notation}

For every $\infty$-category $\mC$ let 
$$\sSeq(\mC):= \sSeq_*(\mC)= \Fun(\Sigma, \mC) \simeq \prod_{\n \geq 0} \Fun(B(\Sigma_\n), \mC)$$ be the monoidal $\infty$-category of symmetric sequences in $\mC$ endowed with the composition product.
\end{notation}

%and $\triv \in \mC^\Sigma$ for the symmetric sequence concentrated in degree 1with value the tensor unit of $\mC.$

% The following \cite[Corollary 4.2.9.]{articles} defines the composition product
% and enriched $\infty$-operads:

% \begin{proposition}\label{compo}
	
% Let $\mC$ be a symmetric monoidal $\infty$-category compatible with small colimits.	
% The $\infty$-category $\sSeq(\mC)$ carries a monoidal structure 
% such that the tensor product $\circ$ on $\sSeq(\mC)$ admits the following description:
% for every $\X, \Y \in \mC$ and $\n \geq 0$ there is a canonical equivalence
% $$ (\X \circ \Y ) _\n \simeq \underset{ \bk \geq 0 }{\coprod} ( \underset{\n_1 +...+ \n_\bk = \n }{\coprod} \Sigma_\n \times_{(\Sigma_{\n_1} \times ... \times \Sigma_{\n_\bk})} (\X_\bk \ot (\bigotimes_{1 \leq \bj \leq \bk } \Y_{\n_\bj })))_{\Sigma_\bk } .$$	
	
% \end{proposition}

% By work of Haugseng \cite[Corollary 4.2.9.]{Rune} the $\infty$-category $\Alg(\sSeq(\mC))$ is
% equivalent to the $\infty$-category of $\mC$-enriched $\infty$-operads in the sense of Chu-Haugseng \cite{chu2020enriched}. This motivates the following definition:

\begin{definition}Let $\mC$ be a presentably symmetric monoidal $\infty$-category.
A $\mC$-enriched $\infty$-operad is an associative algebra in $\sSeq(\mC)$
with respect to the composition product. 
We set $$\Op(\mC):=\Alg(\sSeq(\mC)).$$

\end{definition}

\begin{remark}

Let $\mC$ be a presentably symmetric monoidal $\infty$-category.
The left action of $\sSeq(\mC)$ on itself restricts to a left action of
$\sSeq(\mC)$ on $\mC$, which sends $(\X, \Y)$ to $\underset{\bk \geq 0 }{\coprod} \X_\bk \ot_{\Sigma_\bk} \Y^{\ot \bk}.$
The latter gives rise to a monoidal functor $$ \sSeq(\mC) \to \Fun(\mC,\mC), \mO \mapsto \mO \circ (-)$$
by the universal property of endomorphism algebra,
which sends any $\mC$-enriched $\infty$-operad $\mO$ to a monad $ \mO \circ (-) $ on $\mC$.

\end{remark}

\begin{definition}
Let $\mC$ be a presentably symmetric monoidal $\infty$-category.
For every $\mC$-enriched $\infty$-operad $\mO$ we define the $\infty$-category of $\mO$-algebras in $\mC$ by $$\Alg_\mO(\mC):=\LMod_\mO(\mC),$$
where we use the left action of $\sSeq(\mC)$ on $\mC$.
\end{definition}

\begin{remark}\label{remol}
By definition for every $\mC$-enriched $\infty$-operad $\mO$ there is a canonical equivalence $$ \Alg_\mO(\mC) =\LMod_\mO(\mC) \simeq \LMod_{\mO \circ (-)}(\mC). $$
	
\end{remark}

\begin{notation}
For every symmetric monoidal $\infty$-category $\mC$ let $$\Coalg_{\bE_\infty}(\mC):=\Alg_{\Fin_*}(\mC^\op)^\op$$ be the $\infty$-category of $\bE_\infty$-coalgebras in $\mC$.
\end{notation}
Since the functor $\Coalg_{\bE_\infty}: \infty\Op^{\Fin_*,\mon} \to \infty\Cat$ preserves finite products and $ \infty\Op^{\Fin_*,\mon}$ is preadditive by \cite[Proposition 3.2.4.7.]{lurie.higheralgebra}, we find that $\Coalg_{\bE_\infty}(\mC)$ inherits a symmetric monoidal structure from $\mC$ such that the forgetful functor
$\Coalg_{\bE_\infty}(\mC)\to \mC$ is symmetric monoidal.
By \cite[Corollary 3.2.2.5.]{lurie.higheralgebra} the symmetric monoidal structure on $\Coalg_{\bE_\infty}(\mC)$ is compatible with small colimits if $\mC$ is compatible with small colimits.
So we can make the following definition:

\begin{definition}Let $\mC$ be a symmetric monoidal $\infty$-category compatible with small colimits.
	
A $\mC$-enriched Hopf $\infty$-operad is an $\infty$-operad enriched in $\Coalg_{\bE_\infty}(\mC)$. Let
$$\Hopf(\mC):=\Op(\Coalg_{\bE_\infty}(\mC)).$$
	
\end{definition}

Now we are able to state the main proposition:

\begin{theorem}\label{dfghjbbvl}
Let $\mC$ be a symmetric monoidal $\infty$-category compatible with small colimits. There is a monoidal functor
$$ \sSeq(\Coalg_{\bE_\infty}(\mC)) \to \Fun^{\ot, \oplax}(\mC, \mC) $$ that fits into a
commutative square of monoidal $\infty$-categories:
\begin{equation*} 
\begin{xy}
\xymatrix{
\sSeq(\Coalg_{\bE_\infty}(\mC))  \ar[d] 
\ar[r]^{  } 
& \Fun^{\ot, \oplax}(\mC, \mC)  \ar[d] 
\\
\sSeq(\mC)  \ar[r]^{ }  &  \Fun(\mC, \mC).
}
\end{xy} 
\end{equation*} 

\end{theorem}

Passing to associative algebras we obtain the following corollary:

%The latter commutative square of monoidal $\infty$-categories yields 
\begin{corollary}\label{hopfo}
Let $\mC$ be a symmetric monoidal $\infty$-category compatible with small colimits. There is a commutative square
\begin{equation*}\label{hhgbhhjvcf} 
\begin{xy}
\xymatrix{
\Hopf(\mC) \ar[d] 
\ar[r]^{  } 
& \Alg(\Fun^{\ot, \oplax}(\mC, \mC))  \ar[d] 
\\
\Op(\mC)  \ar[r]^{ }  &  \Alg(\Fun(\mC, \mC)).
}
\end{xy} 
\end{equation*} 
\end{corollary}

In other words the monad associated to a Hopf $\infty$-operad in $\mC$ is an oplax symmetric monoidal monad.
Corollary \ref{patty} implies %that the $\infty$-category $\Alg_\mH(\mC) \simeq \Alg_{\T_\mH}(\mC) $ of $\mH$-algebras in $\mC$ carries a canonical symmetric monoidal structure such that the forgetful functor$\Alg_\mH(\mC) \to \mC$ is symmetric monoidal.%We define operads in $\mC$ as associative algebras in $(\mC^\Sigma)^\ot $
the following theorem:

\begin{corollary}\label{Hoopf}
	
Let $\mC$ be a symmetric monoidal $\infty$-category compatible with small colimits and $\mO$ a $\mC$-enriched Hopf $\infty$-operad.
The $\infty$-category $\Alg_\mO(\mC) $ carries a canonical symmetric monoidal structure such that the forgetful functor $\Alg_\mO(\mC) \to \mC$ is symmetric monoidal.

\end{corollary}

%\begin{proof}By the functoriality of constr. \ref{laxx} the symmetric monoidal Yoneda-embedding $\mC\subset \mC':= \mP(\mC)$ yields a canonical equivalence $\Alg_\mH(\mC) \simeq \mC \times_{\mC'} \Alg_\mH(\mC') $ over $\mC$ that makes the forgetful functor$\Alg_\mH(\mC) \to \mC$ symmetric monoidal as $\mC':= \mP(\mC)$ is compatible with small colimits.\end{proof}\vspace{4mm}

In the following we prove Proposition \ref{dfghjbbvl}.
%As next we prepare the proof of Proposition \ref{dfghjbbvl}.

\begin{notation}
For every functor $\mD \to \T$ let $\sSeq^\T(\mD) \to \T$ be the cotensor
$$\mD^\Sigma:= \T \times_{\Fun(\Sigma, \T)} \Fun(\Sigma, \mD) \to \T.$$
\end{notation}
We have the following functoriality (see \cite[3.2.]{articles}):

\begin{proposition}
		
Let $\T$ be an $\infty$-category and $\mC \to \T$, $\mD \to \T $
commutative monoids in $\infty\Cat_{/ \T}^\cart$ 
%cartesian $\T$-families of symmetric monoidal $\infty$-categories 
fiberwise compatible with small colimits.

\begin{enumerate}
\item The functor $\sSeq^\T(\mD) \to \T$ carries the structure of a 
%is a cartesian $\T$-family $\sSeq^\T(\mC)^\ot \to \T \times {\Delta^\op} $ of monoidal $\infty$-categories with the following properties:
monoid in $\infty\Cat_{/ \T}^\cart$.
For $\T$ contractible the monoid structure on $\sSeq^\T(\mD)\to \T$ is the composition monoidal structure.
For every functor $\rS \to \T$ the canonical equivalence $$\sSeq^\rS(\rS \times_\T \mC) \simeq \rS \times_\T \sSeq^\T(\mC) $$ over $\rS$ refines to an equivalence of monoids over $\rS.$

%\begin{itemize}
%\item the underlying functor $\sSeq^\T(\mC) \to \T$ is the cotensor $\mC^\Sigma \to \T$.
%\item for any functor $\rS \to \T$ there is a canonical equivalence over $\rS \times {\Delta^\op}:$
%$$\sSeq^\rS(\rS \times_\T \mC)^\ot \simeq \rS \times_\T \sSeq^\T(\mC)^\ot. $$\item for $\T=*$ we have $\sSeq^\T(\mC)^\ot \simeq \sSeq(\mC)^\ot.$\end{itemize}

\item Let $\mC \to \mD$ be a map of commutative monoids in $\infty\Cat_{/ \T},$ 
preserving fiberwise small colimits.
The canonical functor $\sSeq^\T(\mC) \to \sSeq^\T(\mD)$ over $\T$ is a map of monoids in $\infty\Cat_{/ \T}$ and 
%of monoidal $\infty$-categories such that the underlying functor $\sSeq^\T(\mC) \to \sSeq^\T(\mD)$ is the induced functor $\mC^\Sigma \to \mD^\Sigma$
for every functor $\rS \to \T$ the equivalence between the pullback of $\sSeq^\T(\mC) \to \sSeq^\T(\mD)$ along $\rS \to \T$ and $\sSeq^\rS(\rS \times_\T \mC) \to \sSeq^\rS(\rS \times_\T \mD)$ is an equivalence of monoids
in $\infty\Cat_{/ \T}$.

\end{enumerate}
	
\end{proposition}

%We fix the following notation:

\begin{notation}
Let $\T \to \rS$ be a functor and $\mE \subset \T^{[1]}$ a full subcategory.
For any cartesian fibrations $\mB \to \T, \mD \to \T $ let
$$\Fun^{\rS, \mE}_{\T}(\mB, \mD) \subset \Fun^{\rS}_{\T}(\mB, \mD) $$ 
be the full subcategory spanned by the maps
$ \mB_\s \to \mD_\s$ of cartesian fibrations relative to $\mE_\s $
for some $\s \in \rS$.
For $\mE =\T^{[1]}$ we write $\Fun^{\rS, \cart}_{\T}(\mB, \mD) $ for $\Fun^{\rS, \mE}_{\T}(\mB, \mD).$
For $\rS$ contractible we drop $\rS$.
%write $ \Fun^{\mE}_{\T}(\mB, \mD)$ for $\Fun^{\rS, \mE}_{\T}(\mB, \mD)$ and $ \Fun^{\cart}_{\T}(\mB, \mD)$ for $\Fun^{\rS, \cart}_{\T}(\mB, \mD).$

\end{notation}

\begin{remark}\label{remrem}
	
If the cartesian fibration $\mD \to \T$ admits fiberwise small colimits, by \cite[Proposition 5.4.7.11.]{lurie.higheralgebra} the $\infty$-category $\Fun^{\rS}_{\T}(\mB, \mD)$ admits fiberwise small colimits which are formed objectwise. Moreover if the fiber transports of $\mD \to \T$ along morphisms of $\mE$ lying over equivalences in $\rS$ preserve small colimits, the full subcategory $$\Fun^{\rS, \mE}_{\T}(\mB, \mD)\subset \Fun^{\rS}_{\T}(\mB, \mD)$$ is fiberwise closed under small colimits.

\end{remark}

\begin{construction}
	
Let $\T \to \rS$ be a functor and $\mE \subset \T^{[1]}$ a full subcategory.
Let $\mC \to \T $ be a commutative monoid in $\infty\Cat^\cart_{/ \T} $ fiberwise compatible with small colimits.
%If for every $\rt \in \T$ the induced symmetric monoidal $\infty$-category $\mC_{\rt}$ is compatible with small colimits, 
%The cartesian $\rS$-family of symmetric monoidal $\infty$-categories $\Fun^{\rS}_{\T}(\T, \mC)^\ot :=\Fun^{\rS \times \Fin_*}_{\T \times \Fin_*}(\T \times \Fin_*, \mC^\ot \to \rS \times \Fin_* $ is fiberwise compatible with small colimits and restricts to $ \Fun^{\rS, \mE}_{\T}(\T, \mC)$, where $ \Fun^{\rS, \mE}_{\T}(\T, \mC)$ is fiberwise closed in$ \Fun^{\rS}_{\T}(\T, \mC)$ under small colimits.
The functor $\Fun^{\rS}_{\T}(\T, -):\infty\Cat^\cart_{/ \T} \to \infty\Cat^\cart_{/ \rS} $ preserves finite products and thus 
$\Fun^{\rS}_{\T}(\T, \mC) \to \rS$ inherits the structure of a commutative monoid in $\infty\Cat^\cart_{/\rS} $, which by Remark \ref{remrem}	is fiberwise compatible with small colimits.
The commutative monoid structure on $\Fun^{\rS}_{\T}(\T, \mC) \to \rS$ restricts to $\Fun^{\rS, \mE}_{\T}(\T, \mC) \to \rS$.
Hence $\Fun^{\rS, \mE}_{\T}(\T, \mC) \to \rS$ inherits the structure of a commutative monoid in $\infty\Cat^\cart_{/\rS} $ that is fiberwise compatible with small colimits since $\Fun^{\rS, \mE}_{\T}(\T, \mC) \to \rS$ is fiberwise closed under small colimits in $\Fun^{\rS}_{\T}(\T, \mC) \to \rS$ by Remark \ref{remrem}. 

%Thus $\Fun^{\rS, \mE}_{\T}(\T, \mC) \to \rS$ is fiberwise compatible with small colimits.

%the commutative monoid $\Fun^{\rS}_{\T}(\T, \mC) $ in $ \Cat^\cart_{\infty/ \rS}$gives rise to a monoid $$\sSeq^\rS(\Fun^{\rS}_{\T}(\T, \mC)) \to \rS $$ in $ \Cat_{\infty/  \rS}$ and the composition product on $ \sSeq^\rS(\Fun^{\rS}_{\T}(\T, \mC)) \to \rS $restricts to $ \sSeq^\rS(\Fun^{\rS, \mE}_{\T}(\T, \mC)).$

\end{construction}

Now we are able to state the main proposition, from which we deduce Corollary \ref{gghfghjnh}.

\begin{proposition}\label{dfghjklk}
Let $\T \to \rS$ be a functor and $\mC \to \T$ a commutative monoid in $\infty\Cat^\cart_{/ \T} $ fiberwise compatible with small colimits.

\begin{enumerate}
\item There is a map of monoids
$ \rho: \sSeq^\rS(\Fun^{\rS}_{\T}(\T, \mC)) \to \Fun^{\rS}_{\T}(\mC, \mC)$
in $\infty\Cat_{/ \rS} $
that sends a symmetric sequence $\A $ in $\Fun_{\T_\s}(\T_\s,\mC_\s)$
for $\s \in \rS$ and $ \X \in \mC_{\rt} $ for $\rt \in \T_\s$
to $ \A(\rt) \circ \X \in \mC_{\rt}.$  

\vspace{1mm}
\item Let $\mE \subset \T^{[1]}$ be a full subcategory. 
If for every map $\F: \rt \to \rt'$ in $\T$ that belongs to $\mE $ the induced functor
$\mC_{\rt'} \to \mC_{\rt} $ preserves small colimits, $\rho$ restricts to a map of monoids in $\infty\Cat_{/ \rS}:$
$$ \sSeq^\rS(\Fun^{\rS, \mE}_{\T}(\T, \mC)) \to \Fun^{\rS, \mE}_{\T}(\mC, \mC).$$
\end{enumerate}	
\end{proposition}

\begin{proof}
(1): The counit transformation $ \Fun^{\rS}_{\T}(\T, -) \times_\rS \T \to \id $ of finite products preserving functors
$\infty\Cat_{/ \T} \to \infty\Cat_{/ \T}$
yields a map $\alpha: \Fun^{\rS}_{\T}(\T, \mC) \times_\rS \T \to \mC  $ of commutative monoids in $\infty\Cat_{/ \T}$ between commutative monoids in $ \infty\Cat^\cart_{/ \T}$,
which induces on the fiber over every $\rt \in \T$ lying over $\s \in \rS$ the small colimits preserving functor $ \Fun_{\T_\s}(\T_\s, \mC_\s) \to \Fun_{\T_\s}( \{ \rt \}, \mC_\s) \simeq \mC_{\rt}.$
The map $\alpha$ yields a map  of monoids $$\phi: \sSeq^\rS(\Fun^{\rS}_{\T}(\T, \mC)) \times_\rS \T \to \sSeq^\T(\mC)$$ in $\infty\Cat_{/ \T}$ between
monoids in $ \infty\Cat^\cart_{/ \T}$.

The evaluation map $ \sSeq^\T(\mC) \to \mC^{\{0\}} $ of cartesian fibrations over $\T$ induces on the fiber over any $\rt \in \T$ 
the evaluation functor $ \sSeq(\mC_{\rt}) \to \mC_{\rt}^{\{0\}} $ right adjoint
to the fully faithful functor that views an object of $\mC_{\rt} $ as a symmetric sequence concentrated in degree 0. 
As a map of cartesian fibrations over $\T$ the functor
$ \sSeq^\T(\mC) \to \mC^{\{0\}} $ has a fully faithful left adjoint $ \mC \to \sSeq^\T(\mC) $ relative to $\T$.
The monoid structure on $ \sSeq^\T(\mC) \to \T $ in 
$\infty\Cat_{/ \T}$ endows $ \sSeq^\T(\mC) \to \T $ with a left action over itself in $ \infty\Cat_{/\T}$ that restricts to a left action on $ \mC \to \T $ over $ \sSeq^\T(\mC) \to \T $ as it restricts fiberwise.
%for every $\s \in \rS,\rt \in \T$ the induced left action on $ \mC_{\s, \rt}^\Sigma $ over itself restricts to a left action on $\mC_{\s, \rt}$ over $\mC_{\s, \rt}^\Sigma. $
Pulling back along $\phi$ we get a left action on $\mC \to \T  $ over $ \sSeq^\rS(\Fun^{\rS}_{\T}(\T, \mC))  $
with respect to the $\infty\Cat_{/ \rS} $-left action on $ \infty\Cat_{/ \T}$ corresponding to a map of monoids $ \sSeq^\rS(\Fun^{\rS}_{\T}(\T, \mC)) \to \Fun^{\rS}_{\T}(\mC, \mC)$ in $\infty\Cat_{/ \rS}.$

(2): To show (2) we can assume that $\rS$ is contractible.
For every functor $\T' \to \T$ there is a commutative square of monoidal $\infty$-categories:
\begin{equation*}\label{hhgbvcf} 
\begin{xy}
\xymatrix{
\sSeq(\Fun_\T(\T,\mC))   \ar[d] 
\ar[r]^{  } 
& \Fun_{\T}(\mC, \mC)  \ar[d] 
\\
\sSeq(\Fun_{\T'}(\T', \T' \times_\T \mC)) \ar[r]^{ }  &  \Fun_{\T'}(\T' \times_\T \mC, \T' \times_\T \mC)
}
\end{xy} 
\end{equation*} 

So to prove (2) we can reduce to the case that $\mE = \Fun([1], \T)$.
If for every morphism $\F: \rt \to \rt'$ in $\T$ the induced functor
$\mC_{\rt'} \to \mC_{\rt}$ preserves small colimits, the monoid
$ \sSeq^\T(\mC) \to \T $ in $\infty\Cat_{/ \T}$ is a monoid in $ \infty\Cat^\cart_{/ \T}$.
Moreover the embedding $\mC \subset \sSeq^\T(\mC)$ is a map of cartesian fibrations over $\T $ so that $\mC \to \T $ is a left module over 
$ \sSeq^\T(\mC) \to \T $ in $ \infty\Cat^\cart_{/ \T} .$
The functor $\alpha: \Fun_{\T }(\T, \mC) \times \T \to \mC $
over $\T $ restricts to a map
$ \Fun^{\cart}_{ \T }( \T ,\mC) \times \T \to \mC$
of cartesian fibrations over $\T.$ 
Thus the restriction $$\sSeq(\Fun^{\cart}_{\T}( \T, \mC)) \times \T \subset \sSeq(\Fun^{}_{ \T}(\T, \mC)) \times \T \xrightarrow{\phi} \sSeq^\T(\mC)$$ is a map of cartesian fibrations over $\T.$
The composition product on $ \sSeq(\Fun_{\T}(\T, \mC)) $ restricts to the full subcategory $ \sSeq(\Fun^{\cart}_{ \T}(\T, \mC)).$
Hence the left $ \sSeq(\Fun_{\T}(\T, \mC)) $-action on $\mC \to \T $ with respect to the canonical left $\infty\Cat $-action on $\infty\Cat_{/ \T}$ restricts to a left $ \sSeq(\Fun^{\cart}_{\T}(\T, \mC)) $-action on $\mC \to \T $ with respect to the canonical left $\infty\Cat $-action on $ \infty\Cat^\cart_{/ \T}$. So the monoidal functor
$$ \sSeq(\Fun_{\T}( \T, \mC)) \times \T  \to \Fun_{\T}(\mC, \mC)$$
restricts to a monoidal functor
$ \sSeq(\Fun^{\cart}_{ \T}( \T, \mC)) \times \T \to \Fun^{\cart}_{\T}(\mC, \mC).$

\end{proof}

\begin{notation} Let $ \mE \subset \Fun([1],\rS \times\Fin_*)$ be the full subcategory of pairs of equivalences and inert morphisms.
For cocartesian $\rS$-families $\mC^\ot \to \rS \times\Fin_*, \mD^\ot \to \rS \times\Fin_* $ of symmetric monoidal $\infty$-categories we set $$ \Fun^{\rS^\op, \ot, \oplax}_{ }(\mC, \mD) := \Fun^{\rS^\op, \mE^\op}_{\rS^\op \times\Fin_*^\op}( ((\mC^\ot)^\rev)^\op, ((\mD^\ot)^\rev)^\op).$$

For $ \mC^\ot \to  \rS \times\Fin_* $ the identity we set 
$$\Coalg_{\bE_\infty}^{\rS^\op}(\mD) :=  \Fun^{\rS^\op, \ot, \oplax}_{ }(\rS, \mD).$$

\end{notation}

\begin{corollary}\label{gghfghjnh}
	
Let $\rS$ be an $\infty$-category and $\mC \to \rS$ a commutative monoid in 
$ \Cat^\cocart_{\infty/\rS} $ compatible with small colimits.
There is a map of associative monoids in $ \infty\Cat_{/\rS^\op} $:
$$ \sSeq^{\rS^\op}(\Coalg_{\bE_\infty}^{\rS^\op}(\mC)) \to \Fun^{\rS^\op, \ot, \oplax}(\mC, \mC).$$ 

\end{corollary}

\begin{proof}
	 
The forgetful functor $ \Cmon(\Cmon(\Cat^\cocart_{\infty/\rS})) \to \Cmon(\Cat^\cocart_{\infty/\rS} )$ is an equivalence by \cite[Corollary 2.4.3.10., Proposition 3.2.4.10.]{lurie.higheralgebra}.
Thus the commutative monoid $\mC \to \rS$ in 
$ \Cat^\cocart_{\infty/ \rS}$ fiberwise compatible with small colimits 
uniquely lifts to a commutative monoid $\bar{\mC} \to \rS \times \Fin_*$ in 
$$ \Cmon(\Cat^\cocart_{\infty/\rS}) \subset \Fun(\Fin_*, \Cat^\cocart_{\infty/\rS}) \simeq \Cat^\cocart_{\infty/ \rS \times\Fin_* }$$ fiberwise compatible with small colimits. 
By Proposition \ref{dfghjklk} 1. there is a map 
$$ \sSeq^{\rS^\op}(\Fun^{\rS^\op}_{\rS^\op \times\Fin_*^\op }(\rS^\op \times\Fin_*^\op,(\bar{\mC}^\rev)^\op)) \to \Fun^{\rS^\op}_{\rS^\op \times\Fin_*^\op }((\bar{\mC}^\rev)^\op, (\bar{\mC}^\rev)^\op) $$
of monoids in $ \infty\Cat_{/\rS^\op} $, which is canonically equivalent to a map of monoids in $ \infty\Cat_{/\rS^\op}:$
$$ \sSeq^{\rS^\op}(\Fun^{\rS}_{\rS \times\Fin_*}(\rS \times\Fin_*,\bar{\mC}^\rev)^\op) \to \Fun^{\rS}_{\rS \times\Fin_* }(\bar{\mC}^\rev, \bar{\mC}^\rev)^\op.$$
By Proposition \ref{dfghjklk} 2. this map restricts to the desired map of monoids in $ \infty\Cat_{/\rS^\op} $
because for every $\s \in \rS$ and inert map $\langle \n \rangle \to \langle \m \rangle $ the induced functor
$\bar{\mC}_{\s, \langle \n \rangle  } \to \bar{\mC}_{ \s, \langle \m \rangle  }  $
preserves small colimits.

\end{proof}

\subsection{Enrichment of the $\infty$-category of algebras}

%\subsection{Enriched monads}

In the following we prove that for every presentable $\infty$-category $\mV$ and presentably symmetric monoidal $\mV$-enriched $\infty$-category, where $\mV$ carries the cartesian structure,
and every $\infty$-operad $\mO$ in $\mC$ the $\infty$-category of $\mO$-algebras in $\mC$ underlies a $\mV$-enriched $\infty$-category (\cref{enrist}).

By \cite[Corollary 4.64.]{heine2024higher} for every presentably symmetric monoidal $\infty$-category $\mV$ the $\infty$-category
$\mV \mathrm{-}\Cat$ inherits a presentably symmetric monoidal structure.

\begin{definition}
Let $\mV$ be a presentably symmetric monoidal $\infty$-category.
A symmetric monoidal $\mV$-enriched $\infty$-category is a commutative algebra in $\mV \mathrm{-}\Cat.$
A symmetric monoidal $\mV$-enriched functor is a morphism of commutative algebras in $\mV \mathrm{-}\Cat.$
	
\end{definition}

% \begin{definition}
% Let $\mV$ be a presentably symmetric monoidal $\infty$-category.
% A presentably symmetric monoidal $\mV$-enriched $\infty$-category is a
% presentably symmetric monoidal $\infty$-category equipped with a left adjoint symmetric monoidal functor from $\mV.$

% \end{definition}

\begin{remark}\label{remip}

By \cite[Theorem 7.21.]{HEINE2023108941} every presentably symmetric monoidal $\mV$-linear $\infty$-category has an underlying large symmetric monoidal $\mV$-enriched $\infty$-category, which uniquely determines it.
Moreover every left adjoint symmetric monoidal $\mV$-linear functor between 
presentably symmetric monoidal $\mV$-linear $\infty$-categories has an underlying $\mV$-enriched symmetric monoidal functor, which uniquely determines it.
\end{remark}

\begin{notation}Let $\mV$ be a presentably symmetric monoidal $\infty$-category and 
$\mA, \mB, \mC$ presentably symmetric monoidal $\mV$-linear $\infty$-categories and
$\mC \to \mA, \mC \to \mB$ left adjoint symmetric monoidal  $\mV$-linear functors.
	
Let $$\mV\mathrm{-}\Fun_\mC^{\ot,\L}(\mA,\mB):=* \times_{\mV\mathrm{-}\Fun_\mC^{\ot,\L}(\mC,\mB)}\mV\mathrm{-}\Fun_\mC^{\ot,\L}(\mA,\mB)  $$ be the $\infty$-category of left adjoint symmetric monoidal $\mV$-linear functors under $\mC.$
\end{notation}

\cite[Proposition 2.16.]{heine2024topologicalmodelcellularmotivic} and Remark \ref{remip} 
imply the following:

\begin{lemma}\label{lemom} Let $\mV$ be a presentably symmetric monoidal $\infty$-category and
$\mC$ a presentably symmetric monoidal $\mV$-linear $\infty$-category
and $\mD$ a presentably symmetric monoidal $\mC$-linear $\infty$-category.

The Day-convolution monoidal structure on $\sSeq(\mC)$ refines to a
presentably symmetric monoidal $\mV$-linear $\infty$-category and evaluation at the point induces an equivalence 
$$ \mV\mathrm{-}\Fun_{\mC}^{\ot,\L}(\sSeq(\mC),\mD) \simeq  \mD. $$

%where $\sSeq(\mC)$ on the left hand side carries the $\mV$-enriched Day-convolution monoidal structure.
	
\end{lemma}

\begin{theorem}\label{enrlift}
Let $\mV$ be a cartesian closed presentable $\infty$-category
and $\mC$ a presentably symmetric monoidal $\mV$-linear $\infty$-category.
%\begin{enumerate}
% \item The monoidal $\infty$-category $\sSeq(\mC)$ underlies a monoidal $\mV$-enriched $\infty$-category
% and the left action of $\sSeq(\mC)$ on $\mC$ underlies a $\mV$-enriched left action.
% \vspace{1mm}\item 
The monoidal functor $\sSeq(\mC) \to \Fun(\mC,\mC)$ lifts to a $\mV$-enriched monoidal functor
$$\sSeq(\mC) \to\mV\mathrm{-}\Fun(\mC,\mC).$$

%\end{enumerate}
	
\end{theorem}

%\begin{remark}\label{opol}
%	
%For every presentably symmetric monoidal $\infty$-categories $\mA,\mB$ under $\mC$ let $\Fun^{\ot,\L}(\mA,\mB) $ be the $\infty$-category of small colimits preserving symmetric monoidal functors under $\mC.$
%
%There is a canonical monoidal equivalence $\sSeq(\mC) \simeq \Fun^{\ot,\L}_\mC(\sSeq(\mC),\sSeq(\mC)),$ where $\sSeq(\mC)$ on the left hand side carries the composition product monoidal structure and $\sSeq(\mC)$ on the right hand side carries the Day convolution monoidal structure.
%	
%\end{remark}

\begin{proof}
%By Remark \ref{opol} there is a canonical monoidal equivalence $\sSeq(\mC) \simeq \Fun^{\ot,\L}_\mC(\sSeq(\mC),\sSeq(\mC)),$ where $\sSeq(\mC)$ on the left hand side carries the composition product monoidal structure and $\sSeq(\mC)$ on the right hand side carries the Day convolution monoidal structure.

%The 2-functor  $  \Calg(\Pr^L)  $  induces a monoidal equivalence $$ \Enr\Fun_{\mV, \ot,\L}_\mC(\sSeq(\mC),\sSeq(\mC)) \to \Fun^{\ot,\L}_\mC(\sSeq(\mC),\sSeq(\mC)).$$

By Lemma \ref{lemom} there is a canonical equivalence
$$\mV\mathrm{-}\Fun_\mC^{\ot,\L}(\sSeq(\mC),\sSeq(\mC)) \simeq \sSeq(\mC).$$

Let $\mB$ be a presentably symmetric monoidal $\infty$-category.
By \cite[Theorem 4.6.]{Gepner_2015} the $\infty$-category $\Cmon(\mB)$ of commutative monoid objects in $\mB$ inherits a presentably symmetric monoidal structure and the free functor $\mB \to \Cmon(\mB)$
refines to a symmetric monoidal functor. 

Let $\mW$ be a presentably symmetric monoidal $\infty$-category.
By \cite{haugseng2023tensorproductenrichedinftycategories} the $\infty$-category $\mW\mathrm{-}\Cat$ inherits a presentably symmetric monoidal structure.
Hence by \cite[Theorem 4.6.]{Gepner_2015} the $\infty$-category $\Cmon(\mW\mathrm{-}\Cat)$ inherits a presentably symmetric monoidal structure and the free functor $$\mW\mathrm{-}\Cat \to \Cmon(\mW\mathrm{-}\Cat)$$ refines to a symmetric monoidal functor. 
So $\Cmon(\mW \mathrm{-}\Cat)$ is a presentably $\mW \mathrm{-}\Cat$-enriched
$\infty$-category and the free and forgetful adjunction $\mW \mathrm{-}\Cat \rightleftarrows \Cmon(\mW \mathrm{-}\Cat)$ is a $\mW \mathrm{-}\Cat$-enriched adjunction.

By \cite[A.4.]{heine2021realktheorywaldhauseninfinity}
this implies that $\Cmon(\mW \mathrm{-}\Cat)_{\mC/ }$ is enriched in $\mW \mathrm{-}\Cat$ for every $\mC \in \Cmon(\mW \mathrm{-}\Cat)$
the forgetful functor $\Cmon(\mW \mathrm{-}\Cat)_{\mC/ } \to\Cmon(\mW \mathrm{-}\Cat) $ is a $\mW$-enriched functor.

% In particular, for every $\mA, \mB \in \Cmon(\mW \mathrm{-}\Cat)_{\mC/ }$
% there is a canonical $\mW$-enriched functor
% $$ \mV\mathrm{-}\Fun_\mC^{\ot}(\mA,\mB)\to \mV\mathrm{-}\Fun^{\ot}(\mA,\mB).$$

We apply this when $\mW$ is a cartesian closed presentable $\infty$-category $\mV$ viewed as a presentably symmetric monoidal $\infty$-category.
In this case the symmetric monoidal structure on $\mV \mathrm{-}\Cat$ is cartesian and $\Cmon(\mV \mathrm{-}\Cat)$ is the $\infty$-category of symmetric monoidal $\mV$-enriched $\infty$-categories and symmetric monoidal $\mV$-enriched functors.

Since $\Cmon(\mV \mathrm{-}\Cat)_{\mC/ }$ is enriched in $\mV \mathrm{-}\Cat$ for every $\mC \in \Cmon(\mV \mathrm{-}\Cat),$ for every symmetric monoidal $\mV$-enriched $\infty$-categories $\mA,\mB$ under $\mC$ the
$\infty$-category $\mV\mathrm{-}\Fun_\mC^{\ot}(\mA,\mB)$ underlies a $\mV$-enriched $\infty$-category.
So also the full subcategory $\mV\mathrm{-}\Fun_\mC^{\ot,\L}(\mA,\mB)$ underlies a $\mV$-enriched $\infty$-category.
By the universal property of the endomorphism algebra the monoidal $\infty$-category $\mV\mathrm{-}\Fun_\mC^{\ot}(\mA,\mA)$ underlies a $\mV$-enriched monoidal $\infty$-category.
So also the full monoidal subcategory $\mV\mathrm{-}\Fun_\mC^{\ot,\L}(\mA,\mA)$ underlies a $\mV$-enriched monoidal $\infty$-category.

Let $ \mD$ be presentably symmetric monoidal $\mC$-linear $\infty$-category.
The canonical 
$\mV$-enriched functor
$$ \mV\mathrm{-}\Fun_{\mC}^{\ot,\L}(\sSeq(\mC),\mD) \to \mV\mathrm{-}\Fun(\sSeq(\mC),\mD) \to \mD $$
evaluating at the point is an equivalence.

Restricting the endomorphism left action of the $\mV$-enriched monoidal $\infty$-category $\mV\mathrm{-}\Fun(\sSeq(\mC),\sSeq(\mC))$ on $ \sSeq(\mC) $ along the $\mV$-enriched monoidal functor 
$$\sSeq(\mC) \simeq \mV\mathrm{-}\Fun_\mC^{\ot,\L}(\sSeq(\mC),\sSeq(\mC))\to \mV\mathrm{-}\Fun(\sSeq(\mC),\sSeq(\mC))$$ 
gives a left action of $\sSeq(\mC)$ on itself in $\mV\mathrm{-}\widehat{\Cat}$, which restricts
to a left action of $\sSeq(\mC)$ on $\mC$ in $\mV\mathrm{-}\widehat{\Cat}$.
By the  universal property of the endomorphism algebra
the latter left action induces a monoidal $\mV$-enriched functor
$\sSeq(\mC) \to \mV\mathrm{-}\Fun(\mC,\mC).$

\end{proof}

We obtain the following:

\begin{corollary}\label{enrist}
Let $\mV$ be a cartesian closed presentable $\infty$-category and $\mC$ a presentably symmetric monoidal $\mV$-linear $\infty$-category.
For every $\mC$-enriched $\infty$-operad $\mO$ the $\infty$-category $\Alg_\mO(\mC)$ of $\mO$-algebras in $\mC$ underlies a $\mV$-enriched $\infty$-category that is tensored and cotensored over $\mV$.

\end{corollary}

\begin{proof}
	
By Proposition \ref{enrlift} the monoidal functor $\sSeq(\mC) \to \Fun(\mC,\mC)$ lifts to a monoidal functor
$\sSeq(\mC) \to\mV\mathrm{-}\Fun(\mC,\mC),$
which sends a $\mC$-enriched $\infty$-operad $\mO$ to a $\mV$-enriched monad $\T:=\mO \circ (-)$ on $\mC,$ which by description of the composition product preserves small sifted colimits and so geometric realizations.
By Proposition \ref{eNr} the $\infty$-category $ \LMod_{\T}(\mC)$
underlies a $\mV$-enriched $\infty$-category, which by \cref{enrtensor}
is a left $\mV$-tensored $\infty$-category and for every $\V \in \mV$
the functor
$\V \ot (-):\LMod_\T(\mC) \to \LMod_\T(\mC)$ preserves small colimits.
By \cite[Proposition 4.2.3.4.]{lurie.higheralgebra} 
the $\infty$-category $\LMod_\T(\mC)$ is presentable since $\T: \mC \to \mC$ is accessible and $\mM$ is presentable.
By presentability of $\LMod_\T(\mC)$ the functor $\V \ot (-):\LMod_\T(\mC) \to \LMod_\T(\mC)$ admits a right adjoint.
Thus $\LMod_\T(\mC)$ is also cotensored over $\mV.$
By \cref{remol} there is a canonical equivalence $$\Alg_\mO(\mC) \simeq \LMod_{\T}(\mC).$$
	
\end{proof}

\appendix

\vspace{1mm}

\section{Locally cocartesian fibrations}

\begin{proposition}\label{swi}

Let $\rS,\T$ be $\infty$-categories, $\mE \subset \Fun([1], \rS)$ a full  subcategory and $\phi:\mC \to \rS \times \T$ a map of locally cocartesian fibrations over $\T.$

\begin{enumerate}
\item The functor $\mC \to \rS$ is a cocartesian fibration relative to
$\mE$ and $\phi$ is a map of locally cocartesian fibrations relative to $\mE$ if for every $\rt \in \T$ the induced functor $\mC_\rt \to \rS$
is a cocartesian fibration relative to $\mE$ and for every morphism $\rt \to \rt'$ of $\T$ the induced functor $\mC_\rt \to \mC_{\rt'}$ over $\rS$ is a map of locally cocartesian fibrations relative to $\mE$.

\vspace{1mm}

\item The functor $\mC \to \rS$ is a cartesian fibration relative to
$\mE$ and $\phi$ is a map of locally cartesian fibrations relative to $\mE$ if for every $\rt \in \T$ the induced functor $\mC_\rt \to \rS$
is a cartesian fibration relative to $\mE$.

\end{enumerate}

\end{proposition}

\begin{proof}
(1): Let $\alpha: \s \to \s'$ be a morphism of $\rS$ that belongs to $\mE.$	
For any $\X \in \mC$ lying over $\s \in \rS$ and some $\rt \in \T$ there is a
$\phi_\rt$-cocartesian lift $\X \to \Y$ of $\alpha$ in $\mC_\rt.$
For any $\Z \in \mC$ lying over $\s'' \in \rS$ and $\rt' \in \T$
the commutative square
\begin{equation*} 
\begin{xy}
\xymatrix{
\mC(\Y,\Z)
\ar[d] 
\ar[r] & \mC(\X,\Z)
\ar[d]
\\
\rS(\s',\s'') \times \T(\rt,\rt') \ar[r] & \rS(\s,\s'') \times \T(\rt,\rt')
}
\end{xy}
\end{equation*}
over $ \T(\rt,\rt')$
induces on the fiber over every map
$\beta:\rt \to \rt'$ of $\T$
the square
\begin{equation*} 
\begin{xy}
\xymatrix{
\mC_{\rt'}(\beta_*(\Y),\Z)
\ar[d] 
\ar[r] & \mC_{\rt'}(\beta_*(\X),\Z)
\ar[d]
\\
\rS(\s',\s'') \ar[r] & \rS(\s,\s''),
}
\end{xy}
\end{equation*}
which is a pullback square by assumption.

(2): Let $\alpha: \s' \to \s$ be a morphism of $\rS$ in $\mE.$	
For every $\X \in \mC$ lying over $\s \in \rS$ and some $\rt \in \T$ there is a
$\phi_\rt$-cartesian lift $\Y \to \X$ of $\alpha$ in $\mC_\rt.$
For every $\Z \in \mC$ lying over $\s'' \in \rS$ and $\rt' \in \T$
the commutative square
\begin{equation*} 
\begin{xy}
\xymatrix{
\mC(\Z,\Y)
\ar[d] 
\ar[r] & \mC(\Z,\X)
\ar[d]
\\
\rS(\s'',\s') \times \T(\rt',\rt) \ar[r] & \rS(\s'',\s) \times \T(\rt',\rt)
}
\end{xy}
\end{equation*}
over $ \T(\rt',\rt)$
yields on the fiber over any map
$\beta:\rt' \to \rt$ of $\T$ the pullback
\begin{equation*} 
\begin{xy}
\xymatrix{
\mC_\rt(\beta_*(\Z),\X)
\ar[d] 
\ar[r] & \mC_{\rt}(\beta_*(\Z),\Y)
\ar[d]
\\
\rS(\s'',\s') \ar[r] & \rS(\s'',\s).
}
\end{xy}
\end{equation*}

\end{proof}

\begin{lemma}\label{heopo}

Let $\mV, \mW$ be monoidal $\infty$-categories, $\F,\G: \mV \to \mW$ lax monoidal functors and $\alpha: \F \to \G$ a lax monoidal natural transformation. Let $\mM$ be a left $\mW$-tensored $\infty$-category.
There is a canonical $\mV$-enriched functor
$\G^*(\mM) \to \F^*(\mM) $ that fits into a commutative triangle of enriched $\infty$-categories
\begin{equation*}
\begin{xy}
\xymatrix{
\mM \ar[rd]
\ar[rr]
&&\G^*(\mM) \ar[ld]
\\
& \F^*(\mM),}
\end{xy} 
\end{equation*}  
where both other enriched functors are the cartesian lifts.

For every associative algebra $\A$ in $\mV$ the induced functor
$$ \LMod_{\G(\A)}(\mM) \simeq \LMod_\A(\G^*(\mM)) \to \LMod_\A(\F^*(\mM)) \simeq \LMod_{\F(\A)}(\mM)$$ 
is the functor that restricts along the map of associative algebras
$\F(\A) \to \G(\A)$ induced by $\alpha.$

\end{lemma}

\begin{proof}

The lax monoidal natural transformation $\alpha: \F \to \G$ is classified by a functor $\theta: [1]\times \mV^\ot \to \mW^\ot$ over $\Delta^\op.$
The pullback of the functor $\mM^\circledast \to \mW^\ot $ along $\theta$
is a functor $\mN^\circledast \to [1]\times \mV^\ot $ whose fibers over 0 and 1 are $\F^*(\mM)$ and $ \G^*(\mM), $ respectively.
The functor $\mN^\circledast \to [1]\times \mV^\ot $ is a locally cocartesian fibration since $\mM^\circledast \to \mW^\ot $ is a locally cocartesian fibration. This implies that the functor $\mN^\circledast \to [1]\times \mV^\ot $ is a map of locally cocartesian fibrations over $\mV^\ot$. For every $\V \in \mV^\ot$ the fiber $\mN_\mV \to [1]$ identifies with the projection $[1]\times \mM \to [1],$ which is a cartesian fibration. \cref{swi} implies that the functor 
$\mN^\circledast \to [1]\times \mV^\ot $ is a map of cartesian fibrations over $[1]$ and so classifies a functor $\G^*(\mM)^\circledast \to \F*(\mM)^\circledast$ over $\mV^\ot$ that encodes a $\mV$-enriched functor
$\G^*(\mM) \to \F*(\mM).$
By construcion there is the claimed commutative triangle and the description of the functor on modules.

\end{proof}

\section{Adjunctions in $(\infty,2)$-categories}

\begin{lemma}\label{lemu}

Let $\F,\G:\mJ \to \infty\Cat$ be functors and
$\alpha: \F \to \G$ a natural transformation such that for every $\Z \in \mJ$
the functor $\alpha_\Z: \F(\Z) \to \G(\Z)$ admits a left adjoint $\beta^\Z$
and for every map $\kappa: \Y \to \Z$ the canonical functor
$\beta^\Z \circ \G(\kappa) \to \F(\kappa) \circ \beta_\Y$ is an equivalence.

\begin{enumerate}

\item The map $\alpha$ admits a left adjoint in $\Fun(\mJ,\infty\Cat).$

\item The induced functor
$\lim(\alpha): \lim(\F) \to \lim(\G)$ admits a left adjoint.

\end{enumerate}

\end{lemma}

\begin{proof}

(1): We use the canonical equivalence
$\Fun(\mJ,\infty\Cat) \simeq  \infty\Cat^\cart_{/ \mJ^\op}$
of $(\infty,2)$-categories provided by the Grothendieck construction.
Let $\mF \to \mJ^\op, \mG \to \mJ^\op$ be the cartesian fibrations classifying
$\F,\G,$ respectively.
Let $\rho: \mF \to \mG$ be the map of cartesian fibrations over $\mJ^\op$ classifying $\alpha: \F \to \G.$
By assumption the functor $\rho$ over $\rS$ admits fiberwise a left adjoint.
This implies by \cite[Proposition 7.3.2.6.]{lurie.higheralgebra} that $\rho$ admits a left adjoint $\tau$ relative to $\mJ^\op$, which is a left adjoint in the $(\infty,2)$-category $ \infty\Cat_{/ \mJ^\op}$.
Since for every map $\kappa: \Y \to \Z$ the canonical functor
$\beta^\Z \circ \G(\kappa) \to \F(\kappa) \circ \beta_\Y$ is an equivalence, $\tau$ is a map of cartesian fibrations over $\mJ^\op$. (2) follows from applying the 2-functor $\lim: \Fun(\mJ,\infty\Cat) \to \infty\Cat$ right adjoint to the diagonal 2-functor. 

\end{proof}

\begin{lemma}\label{prostos}

Let $\mC$ be an $(\infty,2)$-category and $\F: \X \to \Y, \G: \Y \to \X$ morphisms in $\mC$ and $\eta: \id \to \G \circ \F$ 
a 2-morphism.
The following conditions are equivalent:
\begin{enumerate}

\item The 2-morphism $\eta$ exhibits $\F$ as left adjoint to $\G.$

\item For every $\Z \in \mC$ the transformation $\Mor_{\mC}(\Z, \eta)$ 
exhibits $\Mor_{\mC}(\Z, \F)$ as left adjoint to $\Mor_{\mC}(\Z, \G).$
\end{enumerate}

\end{lemma}

\begin{proof}
Condition (1) trivially implies (2). So let (2) be satisfied.
For every $\Z \in \mC$ let $\varepsilon^\Z$ be the counit of the adjunction
$$\Mor_{\mC}(\Z, \F):  \Mor_{\mC}(\Z, \X) \rightleftarrows \Mor_{\mC}(\Z, \Y) : \Mor_{\mC}(\Z, \G)$$ and $\varepsilon:= \varepsilon^\Y_\id.$

For every morphism $\tau: \Z \to \Z'$ in $\mC$ the natural transformation
$$ \varepsilon^{\Z} \circ \Mor_{\mC}(\tau, \Y) : \Mor_{\mC}(\Z, \F) \circ \Mor_{\mC}(\Z, \G) \circ \Mor_{\mC}(\tau, \Y) \to \Mor_{\mC}(\tau, \Y)$$
factors as the canonical equivalence
$$ \Mor_{\mC}(\Z, \F) \circ \Mor_{\mC}(\Z, \G) \circ \Mor_{\mC}(\tau, \Y) 
\simeq \Mor_{\mC}(\tau, \Y)\circ \Mor_{\mC}(\Z', \F) \circ \Mor_{\mC}(\Z', \G) $$
followed by $$ \Mor_{\mC}(\tau, \Y) \circ \varepsilon^{\Z'}:  \Mor_{\mC}(\tau, \Y)\circ\Mor_{\mC}(\Z', \F) \circ \Mor_{\mC}(\Z', \G) \to \Mor_{\mC}(\tau, \Y).$$
For $\tau=\F:\X \to \Y$ and the identity of $\Y$ we find that
$ \varepsilon^{\X}_\F: \F \circ \G \circ \F \to \F$ is canonically equivalent to
$ \Mor_{\mC}(\F, \Y)(\varepsilon^{\Y}_\id)= \varepsilon \circ \F.$

By the triangle identities the transformation
$\Mor_{\mC}(\Z, \eta) \circ \Mor_{\mC}(\Z, \G) $
is a section of $\Mor_{\mC}(\Z, \G) \circ \varepsilon^\Z$
and $\Mor_{\mC}(\Z, \F) \circ \Mor_{\mC}(\Z, \eta)$
is a section of $\varepsilon^\Z \circ \Mor_{\mC}(\Z, \F)$.
For the identity of $\Z=\Y$ we find that $\eta \circ \G $ is a section of $\G \circ \varepsilon$, for the identity of $\Z=\X$ we find that $\F \circ \eta$ is a section of $\varepsilon^\X_\F= \varepsilon \circ \F$.

\end{proof}

\begin{proposition}\label{prosta}

Let $\mC$ be an $(\infty,2)$-category and $\G: \Y \to \X$ a morphism in $\mC.$ 
The following conditions are equivalent:
\begin{enumerate}

\item The morphism $\G$ admits a left adjoint.

\item For every $\Z \in \mC$
the induced functor $\Mor_{\mC}(\Z, \G)$ admits a left adjoint $\F^\Z$ and for every morphism $\tau: \Z \to \Z'$ in $\mC$
the following canonical 2-morphism is an equivalence:
$$\lambda: \F^{\Z} \circ \Mor_{\mC}(\tau, \X) \to \Mor_{\mC}(\tau, \Y) \circ \F^{\Z'}.$$
\end{enumerate}

\end{proposition}

\begin{proof}

Condition (1) trivially implies (2). So let (2) be satisfied.

Embedding $\mC$ into $\FUN(\mC^\op,\infty\Cat)$ we can assume that $\mC$ is tensored and cotensored over $\infty\Cat.$
For every $\Z \in \mC$ let $\eta^\Z$ be the unit and $\varepsilon^\Z$ the counit of the adjunction
$$\F^\Z: \Mor_{\mC}(\Z, \X) \rightleftarrows  \Mor_{\mC}(\Z, \Y) : \Mor_{\mC}(\Z, \G).$$
Let $\F:= \F^\X(\id):\X \to \Y$ and $\eta:= \eta^\X_\id: \id \to \G \circ \F,$
where $\eta$ corresponds to a morphism $\rho: \X \to \X^{[1]}$ in $\mC.$
Let $\G'$ be the image of $\G$ under the 2-functor 
$\mC \hookrightarrow \FUN(\mC^\op,\infty\Cat) \to \Fun(\mC^\op,\infty\Cat)$, where the latter 2-functor forgets the enrichment.
By Lemma \ref{lemu} the map $\G'$ admits a left adjoint
$\F'$. Let $\eta': \id \to \G' \circ \F'$ be the unit of this adjunction.
Then for every $\Z \in \mC$ we have $\F'_\Z \simeq \F^\Z$ and the natural transformation $$\eta'_\Z : \id \to \Mor_\mC(\Z,\G) \circ \F'_\Z \simeq \Mor_\mC(\Z,\G) \circ \F^\Z$$ exhibits $\F^\Z$ as left adjoint to $\Mor_\mC(\Z,\G).$

We prove that for every $\Z \in \mC$ there is an equivalence $\Mor_\mC(\Z,\F) \simeq \F^\Z$ and the transformation $$\eta'_\Z : \id \to \Mor_\mC(\Z,\G) \circ \F^\Z \simeq \Mor_\mC(\Z,\G) \circ \Mor_\mC(\Z,\F)$$
is equivalent to $\Mor_\mC(\Z,\eta)$. Then we apply Lemma \ref{prostos}.
The transformation $\eta'_\Z$ corresponds to $$\rho'_\Z:  \Mor_\mC(\Z,\X) \to \Mor_\mC(\Z,\X^{[1]}) \simeq \Fun([1],\Mor_\mC(\Z,\X))$$
and we need to construct an equivalence $\rho'_\Z \simeq \Mor_\mC(\Z,\rho)$
that induces the identity under the functor $\Mor_\mC(\Z,\ev_0)$
and the equivalence $$\Mor_\mC(\Z,\G) \circ \F^\Z \simeq \Mor_\mC(\Z,\G) \circ \Mor_\mC(\Z,\F)$$ under $\Mor_\mC(\Z,\ev_1)$.
By the Yoneda-lemma the map $\Mor_\mC(\Z,\F)^\simeq$ identifies with the map
$(\F'_\Z)^\simeq = (\F^\Z)^\simeq.$
Since $ \Fun(\K, \Mor_\mC(\Z,\G)) \simeq \Fun(\K \ot \Z,\G),$
there is an equivalence
$ \Fun(\K, \F^\Z)  \simeq \F^{\K \ot \Z}$ and so an equivalence
$$ \Fun(\K, \F^\Z)^\simeq  \simeq (\F^{\K \ot \Z})^\simeq \simeq \Mor_\mC(\K \otimes \Z,\F)^\simeq \simeq \Fun(\K,\Mor_\mC(\Z,\F))^\simeq $$
representing an equivalence $\F^\Z \simeq \Mor_\mC(\Z,\F).$

Similarly, by the Yoneda-lemma there is an equivalence of maps $(\rho'_\Z)^\simeq \simeq \Mor_\mC(\Z,\rho)^\simeq$ 
inducing the identity under the map $\Mor_\mC(\Z,\ev_0)^\simeq$
and under the map $\Mor_\mC(\Z,\ev_1)^\simeq$ the equivalence $$\Mor_\mC(\Z,\G)^\simeq \circ (\F^\Z)^\simeq \simeq \Mor_\mC(\Z,\G)^\simeq \circ \Mor_\mC(\Z,\F)^\simeq.$$

Using the equivalence $ \Fun(\K, \Mor_\mC(\Z,\G))  \simeq \Fun(\K \ot \Z,\G)$
condition (2) implies that there is a canonical equivalence
$ \Fun(\K, \rho'_\Z) \simeq \rho'_{\K \ot \Z}:$$$ \Fun(\K, \Mor_\mC(\Z,\X)) \simeq \Fun(\K \ot \Z,\X)\to \Fun([1],\Fun(\K, \Mor_\mC(\Z,\X))) \simeq  \Fun([1],\Fun(\K \ot \Z,\X))$$ 
that induces the equivalence 
$ \Fun(\K, \Mor_\mC(\Z,\X)) \simeq \Fun(\K \ot \Z,\X)$ under $ \Fun(\K,\Mor_\mC(\Z,\ev_0))$
and the equivalence $ \Fun(\K, \Mor_\mC(\Z,\G)) \circ \Fun(\K, \F^\Z)  \simeq \Fun(\K \ot \Z, \G) \circ \F^{\K \ot \Z} $ under $\Mor_\mC(\Z,\ev_1)$.
We obtain an equivalence representing the desired equivalence
$$\Fun(\K, \rho'_\Z)^\simeq \simeq (\rho'_{\K \ot \Z})^\simeq \simeq \Mor_\mC(\K \ot \Z,\rho)^\simeq \simeq \Fun(\K, \Mor_\mC(\Z,\rho))^\simeq. $$ 

\end{proof}

\begin{corollary}\label{propos}

Let $\mC, \mD$ be $(\infty,2)$-categories and $\alpha: \F \to \G$ a morphism in $\FUN(\mC,\mD).$ The following conditions are equivalent:
\begin{enumerate}

\item The morphism $\alpha$ admits a left adjoint.

\item For every $\Y \in \mC$ the morphism $\alpha(\Y):\F(\Y) \to \G(\Y)$
admits a left adjoint $\beta^\Y$ and for every morphism $\rho: \Y \to \Y'$ in $\mC$ the following canonical 2-morphism is an equivalence:
$$ \beta^{\Y'} \circ \G(\rho) \to \F(\rho) \circ \beta^\Y.$$

\end{enumerate}

\end{corollary}

\begin{proof}

(1) clearly implies (2). So let (2) be satisfied.
We first reduce to the case that $\mD=\infty\Cat.$ The 2-Yoneda-embedding
$\mD \hookrightarrow \FUN(\mD^\op,\infty\Cat) $
induces an embedding
$$\kappa: \FUN(\mC,\mD) \hookrightarrow \FUN(\mC,\FUN(\mD^\op,\infty\Cat)) \simeq \FUN(\mC \times \mD^\op,\infty\Cat).$$
If (2) holds, for every $\Y \in \mC, \Z \in \mD$ the functor $\Mor_\mD(\Z,\alpha(\Y)) : \Mor_\mD(\Z,\F(\Y)) \to \Mor_\mD(\Z,\G(\Y))$
admits a left adjoint and and for every morphisms $\rho: \Y \to \Y'$ in $\mC$
and $\tau: \Z \to \Z'$ in $\mD$ the canonical natural transformation
$$ \Mor_\mD(\Z,\beta^{\Y'}) \circ \Mor_\mD(\Z,\G(\rho)) \to \Mor_\mD(\Z,\F(\rho)) \circ \Mor_\mD(\Z,\beta^\Y) $$
is an equivalence.
Hence by assumption $\kappa(\alpha)$ admits a left adjoint so that $\alpha$ does.

So we can assume that $\mD=\infty\Cat.$
By Proposition \ref{prosta} the morphism $\alpha$ admits a left adjoint if
for every $\Z \in \FUN(\mC,\infty\Cat)$
the induced functor $\Mor_{\FUN(\mC,\infty\Cat)}(\Z, \alpha)$ admits a left adjoint $\Gamma^\Z$ and for every morphism $\tau: \Z \to \Z'$ in $\FUN(\mC,\infty\Cat)$
the canonical natural transformation
$$\lambda: \Gamma^{\Z} \circ \Mor_{\FUN(\mC,\infty\Cat)}(\tau, \G) \to \Mor_{\FUN(\mC,\infty\Cat)}(\tau, \F) \circ \Gamma^{\Z'} $$
is an equivalence.

If $\Z$ is in the image of the 2-Yoneda-embedding $\mC^\op \hookrightarrow \FUN(\mC,\infty\Cat)$, the functor $\Mor_{\FUN(\mC,\infty\Cat)}(\Z, \alpha)$ admits a left adjoint by the 2-Yoneda-lemma (\cref{enryol}).
Let $\mW \subset \FUN(\mC,\infty\Cat)$
be the full subcategory of $\Z$, for which the induced functor $\Mor_{\FUN(\mC,\infty\Cat)}(\Z, \alpha)$ admits a left adjoint.
The full subcategory $\mW$ is closed in $\FUN(\mC,\infty\Cat)$ under the left $\infty\Cat$-action, and condition (2) and Lemma \ref{lemu} imply that $\mW$ is closed under small colimits in $\FUN(\mC,\infty\Cat)$.
Since $\FUN(\mC,\infty\Cat)$ is generated by the essential image of the 2-Yoneda-embedding under small colimits and the left $\infty\Cat$-action \cite[Remark 5.6.]{HEINE2023108941}, for every $\Z \in \FUN(\mC,\infty\Cat)$
the induced functor $\Mor_{\FUN(\mC,\infty\Cat)}(\Z, \alpha)$ admits a left adjoint.

Let $\Z$ belong to the essential image of the 2-Yoneda-embedding $\mC^\op \hookrightarrow \FUN(\mC,\infty\Cat)$
and let $\mW_\Z \subset \FUN(\mC,\infty\Cat)$ be the full subcategory spanned by those $\Z'$ such that for every morphism $\tau: \Z \to \Z'$ in $\FUN(\mC,\infty\Cat)$
the natural transformation $\lambda$ is an equivalence. 
By condition (2) the full subcategory $\mW_\Z$ contains the essential image of the 2-Yoneda-embedding. Moreover condition (2) implies that $\mW_\Z$ is closed in  $\FUN(\mC,\infty\Cat)$ under small colimits and the left $\infty\Cat$-action since 
$\Mor_{\FUN(\mC,\infty\Cat) }(\Z,-)$
preserves small colimits and the left $\infty\Cat$-action by the 2-Yoneda-embedding.
Hence $\mW_\Z = \FUN(\mC,\infty\Cat)$.
So the full subcategory $\mV \subset \FUN(\mC,\infty\Cat)$
spanned by those $\Z$ such that for every morphism $\tau: \Z \to \Z'$ in $\FUN(\mC,\infty\Cat)$
the natural transformation $\lambda$ is an equivalence, contains the essential image of the 2-Yoneda-embedding.
Since $\mV$ is closed in $\FUN(\mC,\infty\Cat)$
under small colimits and the left $\infty\Cat$-action, we find that $\mV=\FUN(\mC,\infty\Cat).$ 

\end{proof}

\bibliographystyle{plain}

\bibliography{ma}

\end{document}